\documentclass[12pt,oneside,a4paper,dvipsnames]{book} 
\usepackage{afterpage}
\newcommand\blankpage{%
    \null
    \thispagestyle{empty}%
    \addtocounter{page}{-1}%
    \newpage}

\usepackage{sectsty}
\usepackage{pgfplots}
\pgfplotsset{compat=1.18}
\usepackage{silence}
\PassOptionsToPackage{dvipsnames}{xcolor}
\usepackage{bm}
\usepackage{dsfont}
\usepackage{lettrine}
\usepackage[russian, english]{babel}
\usepackage[nohints,tight]{minitoc}
\usepackage{pdfpages}
\usepackage{mathtools,amssymb}
\usepackage{amsmath}
\usepackage{amssymb, amsthm}
\usepackage{multirow}
\usepackage{booktabs}
\usepackage{varwidth}
\usepackage{listings}
\usepackage{relsize}
\usepackage{float} 
\usepackage{subfig}
\usepackage{multirow}
\usepackage[left=1in, right=1in, top=1in, bottom=1in]{geometry}
\usepackage []{csquotes}
\usepackage{enumitem}
\MakeOuterQuote{"}
\usepackage[referable]{threeparttablex}
\usepackage[ruled,vlined]{algorithm2e}  
\usepackage{changepage}   
\usepackage{epsfig}

\usepackage[dvipsnames]{xcolor}
\usepackage{color}
\usepackage{latexsym}
\usepackage{mathrsfs}
\usepackage[titletoc]{appendix}
\usepackage{alltt}
\usepackage{longtable}
\usepackage{graphicx}
\usepackage{tikz}
\usepackage{forest}
\usepackage{stanli}
\usetikzlibrary{positioning,shapes,shadows,arrows ,automata,arrows.meta,decorations, decorations.text,calligraphy}
\usepackage{pgfplots}
\pgfplotsset{compat=1.18}
\usepackage{array}
\usepackage{arydshln}
\usepackage{acronym}
\usepackage{csquotes}
\usepackage[colorlinks]{hyperref}
\definecolor{colorlink}{RGB}{110,42,175}
\hypersetup{
    allcolors  = {colorlink},
}

\usepackage{xpatch}

\usepackage[subfigure]{tocloft}
\usepackage{titletoc}
\usepackage{calligra}

\usepackage{times}

\usepackage{titlesec}

\usepackage{textcase}

\sectionfont{\normalfont}
\titleformat{\section}
{\normalfont\fontsize{14pt}{16pt}\selectfont}{\thesection}{1em}{}
\titleformat{\subsection}
{\normalfont\fontsize{12pt}{14pt}\selectfont}{\thesubsection}{1em}{}
\titleformat{\subsubsection}
{\normalfont\fontsize{12pt}{14pt}\selectfont}{\thesubsubsection}{1em}{}
\usepackage{caption}
\newtheorem{definition}{Definition}
\newcommand{\brdefinition}{\begin{definition}}
	\newcommand{\erdefinition}{\end{definition}}
\DeclareTextSymbolDefault{\ae}{T1}
\newtheorem{corollary}{Corollary}
\newcommand{\bcorollary}{\begin{corollary}}
	\newcommand{\ecorollary}{\end{corollary}}
\newtheorem{example}{Example}
\newcommand{\bexample}{\begin{example}}
	\newcommand{\eexample}{\end{example}}
\newtheorem{remark}{Remark}
\newcommand{\bremark}{\begin{remark}}
	\newcommand{\eremark}{\end{remark}}

\newtheorem{theorem}{Theorem}
\newcommand{\btheorem}{\begin{theorem}}
	\newcommand{\etheorem}{\end{theorem}}
\newtheorem{lemma}{Lemma}
\newcommand{\blemma}{\begin{lemma}}
	\newcommand{\elemma}{\end{lemma}}

\newcommand{\brdef}{\begin{defi}}
	\newcommand{\erdef}{\end{defi}}
\newcommand{\bcor}{\begin{cor}}
	\newcommand{\ecor}{\end{cor}}
\newcommand{\bth}{\begin{thm}}
	\newcommand{\ble}{\begin{lem}}
		\newcommand{\ele}{\end{lem}}
	\newcommand{\bcha}{\end{cha}}

\renewcommand{\thechapter}{\arabic{chapter}}
\renewcommand{\thesection}{\thechapter.\arabic{section}}
\newcommand{\inp}{\text{in}}
\newcommand{\outp}{\text{out}}

\renewcommand{\&}{and}

\theoremstyle{definition}

\usepackage{sectsty}
\chapterfont{\centering }
\titleformat{\chapter}[display]
{\bf\centering}
{\chaptertitlename\ \thechapter}{18pt}{\large}
\sectionfont{\normalfont}

\titleformat{\section}
{\normalfont\fontsize{14}{18}\bfseries}{\thesection}{1em}{}

\cftsetindents{chapter}{0em}{3em}      
\cftsetindents{section}{0em}{3em}

\usepackage[sort&compress,numbers]{natbib}
\usepackage[acronym,toc=false,section=section, nogroupskip]{glossaries}
\makeglossaries

\newcolumntype{L}[1]{>{\raggedright\let\newline\\\arraybackslash\hspace{0pt}}m{#1}}
\newcolumntype{C}[1]{>{\centering\let\newline\\\arraybackslash\hspace{0pt}}m{#1}}
\newcolumntype{R}[1]{>{\raggedleft\let\newline\\\arraybackslash\hspace{0pt}}m{#1}}

\usepackage{etoolbox}

\makeatletter
\patchcmd{\ttlh@hang}{\parindent\z@}{\parindent\z@\leavevmode}{}{}
\patchcmd{\ttlh@hang}{\noindent}{}{}{}
\makeatother

\newcommand{\eqnref}[1]{Eq.~(\ref{#1})}
\newcommand{\figref}[1]{Fig.~\ref{#1}}
\newcommand{\tabref}[1]{Tab.~\ref{#1}}
\newcommand{\algref}[1]{Alg.~\ref{#1}}

\def\abs#1{\left\lvert#1\right\rvert}

\definecolor{curcolor}{RGB}{181,126,220} 
\definecolor{curcolordark}{RGB}{142,57,201} 

\definecolor{cclcolor}{RGB}{165,220,126} 
\definecolor{cclcolordark}{RGB}{124,204,67} 

\dominitoc

\newcommand{\objectif}[1]{%
\vspace{.1cm}
\begin{singlespace}
\tikzstyle{titlebox}=[rectangle,inner sep=10pt,inner ysep=10pt,draw=curcolor,draw]%
\tikzstyle{title}=[fill=white]%
\bigskip\noindent\begin{tikzpicture}
\node[titlebox] (box){%
    \begin{minipage}{0.9\textwidth}
#1
    \end{minipage}
};
\node[title] at (box.north west) {\color{curcolordark}  Objectives};
\end{tikzpicture}\bigskip%
\end{singlespace}
}

\newcommand{\recap}[1]{%
\vspace{.1cm}
\begin{singlespace}
\tikzstyle{titlebox}=[rectangle,inner sep=10pt,inner ysep=10pt,draw=cclcolor,draw]%
\tikzstyle{title}=[fill=white]%
\bigskip\noindent\begin{tikzpicture}
\node[titlebox] (box){%
    \begin{minipage}{0.9\textwidth}
#1
    \end{minipage}
};
\node[title] at (box.north west) {\color{cclcolordark}  Summary };
\end{tikzpicture}\bigskip%
\end{singlespace}
}

\newcommand{\tstar}[5]{
\pgfmathsetmacro{\starangle}{360/#3}
\draw[#5] (#4:#1)
\foreach \x in {1,...,#3}
{ -- (#4+\x*\starangle-\starangle/2:#2) -- (#4+\x*\starangle:#1)
}
-- cycle;
}

\newcommand{\meta}{\mathrm{met}}
\newcommand{\cat}{\mathrm{cat}}

\newcommand{\quant}{\mathrm{qnt}}

\newcommand{\neutral}{\mathrm{neu}}
\newcommand{\decreed}{\mathrm{dec}}
\newcommand{\acting}{\mathrm{inc}}

\newcommand{\y}{{\bm{y}}}

\newcommand{\hyi}{\hat{y}_i(\bm{x})}

\newcommand{\As}{{\mathcal{A}}}
\newcommand{\A}{\bm{A}}

\DeclareCaptionLabelFormat{customa}
{%
      #1 \textbf{(a)}
}
\DeclareCaptionLabelFormat{customb}
{%
      #1 \textbf{(b)}
}

\newlength{\chaptertopskip}
\newlength{\chapterbottomskip}
\makeatletter
\def\@makechapterhead#1{%
  \vspace*{\chaptertopskip}%
  {\parindent \z@ \raggedright \normalfont
    \ifnum \c@secnumdepth >\m@ne
      \if@mainmatter
         \@chapapp \\
         \tiny \thechapter\ 
      \fi
    \fi
    \interlinepenalty\@M
    \Huge \bfseries #1\par\nobreak
    \vskip \chapterbottomskip
  }}
\makeatother
\chapterfont{\centering }
\titleformat{\chapter}[display]
{\normalfont%
	\centering 
     \Large
	\bfseries}
 { 
  \vspace*{-2.5cm}
 \centering \chaptertitlename\  \thechapter}{14pt}{%
	 \large 
}
\newacronym{GP}{GP}{Gaussian Process}
\newacronym{BO}{BO}{Bayesian Optimization}
\newacronym{MDO}{MDO}{Multidisciplinary Design Optimization}
\newacronym{MDA}{MDA}{Multidisciplinary Design Analysis}
\newacronym{EGO}{EGO}{Efficient Global Optimization}
\newacronym{PLS}{PLS}{ Partial Least Squares}
\newacronym{MBSE}{MBSE}{Model-Based Systems Engineering}
\newacronym{SEGO}{SEGO}{Super Efficient Global Optimization}
\newacronym{SEGOMOE}{SEGOMOE}{Super Efficient Global Optimization with Mixture Of Experts}
\newacronym{UTB}{UTB}{Upper Trust Bound}
\newacronym{SMT}{SMT}{Surrogate Modeling Toolbox}
\newacronym{LHS}{LHS}{Latin Hypercube Sampling}
\newacronym{KPLS}{KPLS}{Kriging with Partial Least Squares}
\newacronym{FAST-OAD}{FAST-OAD}{Future Aircraft Sizing Tool with Overall Aircraft Design}
\newacronym{EU}{EU}{European Union}
\newacronym{EGORSE}{EGORSE}{Efficient Global Optimization with Random and Supervised Embeddings}
\newacronym{WB2S}{WB2s}{Watson and Barnes $2^{\text{nd}}$ criterion Scaled}
\newacronym{NSGA2}{NSGA-II}{Non-dominated Sorting Genetic Algorithm number II}
\newacronym{AIAA}{AIAA}{American Institute of Aeronautics and Astronautics}
\newacronym{DOE}{DoE}{Design of Experiments}
\newacronym{EI}{EI}{Expected Improvement}
\newacronym{MDF}{MDF}{MultiDisciplinary Feasible}
\newacronym{IDF}{IDF}{Individual Disciplinary Feasible}
\newacronym{XDSM}{XDSM}{eXtended Design Structure Matrix}

\newacronym{DRAGON}{"\texttt{DRAGON}"}{Distributed fans Research Aircraft with electric Generators by ONera}
\newacronym{KPLSK}{KPLSK}{Kriging with Partial Least Squares + Kriging}
\newacronym{TLAR}{TLAR}{Top Level Aircraft Requirements}

\usepackage{tlsflyleaf}

\title{High-dimensional multidisciplinary design optimization for aircraft eco-design}
\author{Paul SAVES}
\defencedate{19/01/2024}

\lab{ISAE-ONERA MOIS - MOdélisation et Ingénierie des Systèmes}

\nboss{2}
\makesomeone{boss}{1}{Nathalie BARTOLI}{}{} 
\makesomeone{boss}{2}{Youssef DIOUANE}{}{}  
\nreferee{2}
\makesomeone{referee}{1}{Sébastien DA VEIGA}{}{}
\makesomeone{referee}{2}{Emmanuel VAZQUEZ}{}{}
\njudge{9}
\makesomeone{judge}{1}{Olivier ROUSTANT}{Professeur, INSA}{Président du Jury}

\makesomeone{judge}{2}{Sébastien DA VEIGA}{Professeur Associé, ENSAI}{Rapporteur}
\makesomeone{judge}{3}{Emmanuel VAZQUEZ}{Professeur, CentraleSupélec}{Rapporteur}

\makesomeone{judge}{5}{Delphine SINOQUET} {Ingénieure de Recherche, IFPEN}{Examinatrice}
\makesomeone{judge}{4}{Julien PELAMATTI} {Ingénieur de Recherche, EDF R$ \hspace{-0.07cm} \And \hspace{-0.09cm}$D}{Examinateur}

\makesomeone{judge}{7}{Nathalie BARTOLI} {Directrice de recherche, ONERA}{Directrice}
\makesomeone{judge}{8}{Youssef DIOUANE}{Professeur Associé, Polytechnique Montréal}{Co-Directeur}

\makesomeone{judge}{6}{Joseph MORLIER}{Professeur, ISAE-SUPAERO}{Membre du Jury}

\makesomeone{judge}{9}{Thierry LEFEBVRE}{Ingénieur de Recherche, ONERA}{Invité}

\begin{document}
\clearpage\thispagestyle{empty}
\makeflyleaf
\blankpage

\thispagestyle{plain}
\pagenumbering{roman}
\newpage

\begin{spacing}{1.5}
	\tableofcontents	
\end{spacing}

\newpage
\thispagestyle{plain}
\addcontentsline{toc}{section}{\bf Résumé}
{\centerline { {\textbf{Résumé}}}}
~\\
\vspace{-0.5cm}
\lettrine[lines=2, lhang=0.33, loversize=0.25, findent=1.5em]{D}{e nos jours}, un intérêt significatif et croissant pour améliorer les processus de conception de véhicules s'observe dans le domaine de l'optimisation multidisciplinaire grâce au développement de nouveaux outils et de nouvelles techniques. Concrètement,  en conception aérostructure, les variables aérodynamiques et structurelles s'influencent mutuellement et ont un effet conjoint sur des quantités d'intérêt telles que le poids ou la consommation de carburant. L'optimisation multidisciplinaire se présente alors comme un outil puissant pouvant effectuer des compromis inter-disciplinaires.
Dans le cadre de la conception aéronautique, le processus multidisciplinaire implique généralement des variables de conception mixtes, continues et catégorielles. Par exemple, la taille des pièces structurelles d'un avion peut être décrite à l'aide de variables continues, le nombre de panneaux est associé à un entier et la liste des sections transverses ou le choix des matériaux correspondent à des choix catégoriels.
L'objectif de cette thèse est de proposer une approche efficace pour optimiser un modèle multidisciplinaire boîte noire lorsque le problème d'optimisation est contraint et implique un grand nombre de variables de conception mixtes (typiquement 100 variables). 
L'approche d'optimisation bayésienne utilisée consiste en un enrichissement séquentiel adaptatif d'un métamodèle pour approcher l'optimum de la fonction objectif tout en respectant les contraintes. Les modèles de substitution par processus gaussiens sont parmi les plus utilisés dans les problèmes d'ingénierie pour remplacer des modèles haute fidélité coûteux en temps de calcul. L'optimisation globale efficace est une méthode heuristique d'optimisation bayésienne conçue pour la résolution globale de problèmes d'optimisation coûteux à évaluer permettant d'obtenir des résultats de bonne qualité rapidement. Cependant, comme toute autre méthode d'optimisation globale, elle souffre du fléau de la dimension, ce qui signifie que ses performances sont satisfaisantes pour les problèmes de faible dimension, mais se détériorent rapidement à mesure que la dimension de l'espace de recherche augmente. Ceci est d'autant plus vrai que les problèmes de conception de systèmes complexes intègrent à la fois des variables continues et catégorielles, augmentant encore la taille de l'espace de recherche. Dans cette thèse, nous proposons des méthodes pour réduire de manière significative le nombre de variables de conception comme, par exemple, des techniques d'apprentissage actif telles que la régression par moindres carrés partiels. Ainsi, ce travail adapte l'optimisation bayésienne aux variables discrètes et à la grande dimension pour réduire le nombre d'évaluations lors de l'optimisation de concepts d'avions innovants moins polluants comme la configuration hybride électrique "\texttt{DRAGON}".
\vskip 0.5cm
\noindent
\textbf{Mots Clefs :} {\it Processus gaussien}, {\it Optimisation boîte noire}, {\it Inférence bayésienne}, {\it Variables hiérarchiques et catégorielles}, {\it Conception d'avions décarbonés}.
\newpage
\thispagestyle{plain}
\addcontentsline{toc}{section}{\bf Abstract}
{\centerline { {\textbf{Abstract}}}}
~\\
\vspace{-0.5cm}
\lettrine[lines=2, lhang=0.33, loversize=0.25, findent=1.5em]{N}{owadays}, there is a significant and growing interest in improving the efficiency of vehicle design processes through the development of tools and techniques in the field of \gls{MDO}. Specifically, in aerostructure design, aerodynamic and structural variables influence each other and have a joint effect on quantities of interest like weight or fuel consumption and, as such, \gls{MDO} arises as a powerful tool for automatically making interdisciplinary trade-offs. In the aircraft design context, the process generally involves mixed continuous and categorical design variables. For instance, the size of an aircraft's structural parts can be described using continuous variables, while discrete variables may include either integer variables, like the number of panels, or categorical variables, like cross-sections or material choices. 
 The objective of this \textit{Philosophiae Doctor} (Ph.D) thesis is to propose an efficient approach for optimizing a multidisciplinary black-box model when the optimization problem is constrained and involves a large number of mixed integer design variables (typically 100 variables).
The targeted optimization approach, called \gls{EGO}, is based on a sequential enrichment of an adaptive surrogate model and, in this context, \gls{GP} surrogate models are one of the most widely used in engineering problems to approximate time-consuming high fidelity models. \gls{EGO} is a heuristic \gls{BO} method that performs well in terms of solution quality. However, like any other global optimization method, \gls{EGO} suffers from the curse of dimensionality, meaning that its performance is satisfactory on lower dimensional problems, but deteriorates as the dimensionality of the optimization search space increases. For realistic aircraft design problems, the typical size of the design variables can even exceed 100 and, thus, trying to solve directly the problems using \gls{EGO} is ruled out. The latter is especially true when the problems involve both continuous and categorical variables increasing even more the size of the search space. In this Ph.D thesis, effective parameterization tools are investigated, including techniques like partial least squares regression, to significantly reduce the number of design variables.
Additionally, Bayesian optimization is adapted to handle discrete variables and high-dimensional spaces in order to reduce the number of evaluations when optimizing innovative aircraft concepts such as the "\texttt{DRAGON}" hybrid airplane to reduce their climate impact.
\vskip 0.5cm
\noindent
\textbf{Keywords:} {{\it Gaussian process}, {\it Black-box optimization}, {\it Bayesian inference},  \it Multidisciplinary design optimization},   {\it Mixed hierarchical and categorical inputs},  {\it Eco-friendly aircraft design}.  

\newpage

\thispagestyle{plain}
\addcontentsline{toc}{section}{\bf Remerciements}
{\centerline { {\textbf{Remerciements }}}}

\vspace*{1cm}

\lettrine[lines=2, lhang=0.33, loversize=0.25, findent=1.5em]{C}{e} parcours  doctoral a été une aventure extraordinaire, et mes remerciements vont à ceux qui y ont joué un rôle essentiel, tant sur le plan académique que personnel. Cette page est trop courte pour vous tous, mais je vais essayer ici de vous exprimer ma gratitude même si je pense que chacun d'entre vous mériterait plusieurs lignes pour vous rendre honneur comme vous le mériterait.

Nathalie Bartoli et Youssef Diouane, mes directeurs de thèse, ont toute ma gratitude pour leur engagement inébranlable depuis notre première collaboration en 2020. Leurs conseils, leur dévouement et leur soutien ont été des fondations présentes tout au long de cette expérience, même à travers les défis de la pandémie. Joseph Morlier et Thierry Lefebvre, mes encadrants, ont ajouté une dimension significative à ma recherche. Votre expertise et vos conseils éclairés à tous et toutes ont façonné cette thèse d'une manière que je n'aurais jamais imaginée.

La réussite de cette soutenance a été possible grâce aux relectures minutieuses et aux conseils de Sébastien Da Veiga et Emmanuel Vazquez, qui ont accepté d'être rapporteurs pour le manuscrit.  Un immense merci à Delphine Sinoquet et Julien Pelamatti qui ont fait le déplacement jusqu'à Toulouse pour examiner mes travaux. Enfin, un merci tout particulier à  Olivier Roustant, que j'ai eu l'occasion de rencontrer lors de mes études d'ingénieur, pour avoir présidé ma soutenance. Merci également pour leurs recommandations à mes professeurs et enseignants à l'ENAC et à l'INSA, c'est aussi grâce à vous si j'en suis là ! 

A tous mes amis de l'INSA Toulouse pour avoir été là, depuis notre rencontre en 2015 et à mes amis plus anciens, je suis fier de vous connaître et d'avoir partagé tant de moments en votre compagnie ces 15 dernières années. 

À mes collègues de l'ISAE-SUPAERO, de Polytechnique Montréal et de l'ONERA, merci pour les échanges enrichissants. Chacun d'entre vous a contribué à créer un environnement de travail stimulant et collaboratif et cela a toujours été une joie que d'aller au labo le matin. 
A Rémy Charayron, Andrés, et Luiz pour avoir partagé mon bureau et aussi à tous les membres du club de tarot à l'ONERA, merci, vous avez fait de ce lab  un espace chaleureux et convivial. Vos amitiés et votre compagnie ont rendu ces années de recherche particuliérement inoubliables et agréables.

À ma famille, pilier de ma vie, qui mérite une reconnaissance spéciale, je vous aime tous et toute et je vous embrasse bien fort. À mes parents et à mon frère, dont le soutien inconditionnel a été ma boussole et à ma tante, la première docteure de la famille pour m'avoir donné envie de perséverer dans la recherche.

Enfin, à mon aimée Dorielle, ta présence et tes encouragements m'ont aidé à garder la motivation jusqu'au bout et au-delà !

Cette réussite est le résultat de nombreuses collaborations et de soutiens inestimables, et je suis profondément reconnaissant envers chacun d'entre vous. Je pense à vous tous et je vous remercie profondément !

\thispagestyle{plain}

\medskip 

Cette thèse est dédiée à feu Prof. Feodor Aleksandrovich Murzin (1952-2021), chef du laboratoire de modélisation et de systèmes complexes, membre honoraire en sciences physiques et mathématiques de l'institut A.P. Ershov en systèmes informatique, directeur scientifique adjoint aux sciences, membre de la branche sibérienne de l'académie des sciences de la fédération russe, qui nous a récemment quitté.
C'est grâce à lui que j'ai découvert la science et la recherche et il m'a toujours soutenu et encouragé à poursuivre et à réaliser cette thèse de doctorat. 

\vskip 0.8cm
\noindent
$\quad$ Fait à Toulouse,
\vskip 0.3cm
\noindent
$\quad$ Le 19/01/2024 

\quad \quad \quad \quad \quad \quad 
\quad \quad \quad \quad \quad \quad 
\quad \quad \quad \quad \quad \quad 
\quad \quad \quad \quad \quad \quad \includegraphics[scale=0.085]{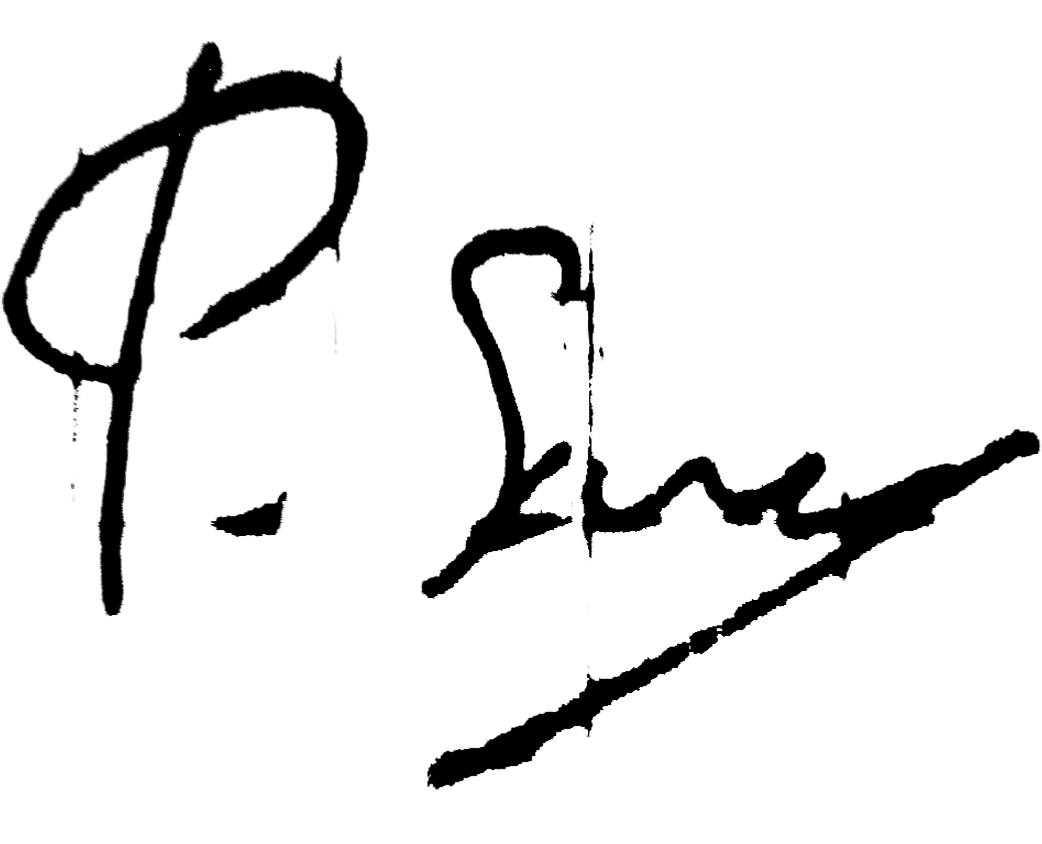}
\newpage

\thispagestyle{plain}
\addcontentsline{toc}{section}{\bf Declaration}
{\centerline { {\textbf{Declaration}}}}

\vspace*{1cm}

\lettrine[lines=2, lhang=0.33, loversize=0.25, findent=1.5em]{I}{ hereby} declare that this thesis represents my own work which has been done after registration for the degree of \textit{Philosophiae Doctor} (Ph.D) at ISAE-SUPAERO. The contents of this report have not been submitted and will not be submitted either in part or in full, for the award of any other degree or diploma in this institute or any other institute or university. 
In the following, I acknowledge all the people who collaborated significantly to this work.

Firstly, the main contributors are my supervisors: N. Bartoli, Y. Diouane, T. Lefebvre and J. Morlier.  

Then, comes the other co-authors: 
J. H. Bussemaker,
C. David,
S. Defoort,
P. Della Vecchia,
G. Donelli,
S. Dubreuil,
J. M. Gomes de Mello,
R. Grapin,
J. T. Hwang,
R. Lafage,
B. Nagel,
E. Nguyen Van, 
M. Mandorino,
J. R. R. A. Martins and
R. Priem. 
I also want to thanks those with whom our worked are not finished yet, namely, E. Hall\'e-Hannan, S. Le Digabel and C. Audet. 

Next, I want to acknowledge the people who contributed indirectly to the results presented in this work,
O. Altinault,
A. Averland,
G. Berthelin, 
M. A. Bouhlel, 
I. Cardoso,
R. Carreira Rufato,
G. Cataldo,
A. Chan, 
R. Charayron, 
R. Conde Arenzana,
C. Döll, 
D. Donjat,
V. Drouet,
J. Hermetz,
A. Kenens,
A. F. López-Lopera,
M. Meliani,
M. Menz,
N. Moëllo,
M. Mongeau,
B. Paluch,
M. Ridel,
E. Roux,
C. Tardieu,
A. Thouvenot,
H. Valayer and
F. Vergnes.

We acknowledge the partners institutions: ONERA, NASA Glenn, ISAE-SUPAERO, Institut Clement Ader (ICA), the  University of Michigan, Polytechnique Montr\'eal, the Research Group in Decision Analysis of Montréal (GERAD), the University of California San Diego, the German Aerospace Center (DLR), Embraer S.A. and the University of Naples Federico II (UNINA). 
From 2019 to 2023, the research presented in this paper was performed in the framework of the AGILE 4.0 project (Towards Cyber-physical Collaborative Aircraft Development) and received funding from the European Union Horizon 2020 Programme under grant agreement n$^{\circ}$ 815122. 
Since 2023, the research presented in this paper was performed in the framework of the COLOSSUS project (Collaborative System of Systems Exploration of Aviation Products, Services and Business Models) and has received funding from the European Union Horizon Programme under grant agreement n$^{\circ}$ 101097120. 
This work is part of the activities of ONERA - ISAE - ENAC joint research group.

I certify that I have no known competing financial interests or personal relationships that could have appeared to influence the work reported in this document.
\newpage


\let\origaddvspace\addvspace
\renewcommand{\addvspace}[1]{}
\begin{spacing}{1.5}
	\addcontentsline{toc}{section}{\bf List of figures}
	\setcounter{lofdepth}{1} \listoffigures
\end{spacing}
\newpage
\begin{spacing}{1.5}
	\addcontentsline{toc}{section}{\bf List of tables}
	\setcounter{lotdepth}{1} \listoftables 	
\end{spacing}

\thispagestyle{plain}
\clearpage

\vspace*{0.5cm}
\begin{spacing}{1.5}
\addcontentsline{toc}
{section}{\bf List of terms and abbreviations}

\thispagestyle{plain}
\printglossary[type=\acronymtype, title={\bf ~~~~~~~~~~~~~~~~~~List of terms and abbreviations}]
\end{spacing}
\renewcommand{\addvspace}[1]{\origaddvspace{#1}}

\thispagestyle{plain}
\mainmatter

\pagenumbering{arabic}

\chapter{Introduction} \label{c1}

\setlength{\fboxrule}{0pt}
\hspace{7cm} \noindent\fbox{%
    \parbox{0.55\textwidth}{%
      \selectlanguage{russian}
  \hspace*{4.2cm} Эшелон за эшелоном, \\
  \hspace*{4.2cm} Эшелон за эшелоном, \\
  \hspace*{4.2cm} Путь-дорога широка...
     \\
     \hrule \vspace{0.2cm}
     \hspace*{\fill} Эшелонная, Осип Яковлевич Колычев }%
}

\objectif{ 
\lettrine[lines=2, lhang=0.33, loversize=0.25, findent=1.5em]{T}{his} chapter introduces the context and objectives of this \textit{Philosophiae Doctor} (Ph.D) thesis and explains how and why this work was carried out, what the goals were and how we addressed them.  
}



\setcounter{section}{-1}
\section{Synthèse du chapitre en français}
Le chapitre 
"Introduction"
introduit le travail présenté dans ce manuscrit. Il revient sur le contexte dans lequel ce travail a eu lieu et sur ses principales contributions.
Cette thèse reprend des articles pour le contenu de ses chapitres dont les thématiques et les motivations sont ici exposées. 

Ce travail de recherche a été réalisé dans le cadre de l'Ecole Doctorale Mathématiques, Informatique et Télécommunications de Toulouse (EDMITT) en collaboration avec l'ONERA et l'ISAE-SUPAERO. L'objectif principal était d'investiguer et de développer des méthodologies pour optimiser la conception d'aéronefs plus respectueux de l'environnement. Les projets AGILE 4.0 (2019-2023) et COLOSSUS (2023-2026), financés par l'Union européenne ont été initiés pour relever ces défis écologiques en étendant les efforts de recherche sur l'ensemble du cycle de vie des avions. Lors de la conception de systèmes complexes comme les aéronefs, l'utilisation de méthodes de conception optimale multidisciplinaire (MDO pour Multidisciplinary Design Optimization) et de technologies d'ingénierie basée sur les modèles (MBSE pour Model-Based System Engineering) est essentielle pour réduire l'impact environnemental lié à la consommation de carburant, à la génération de déchets et aux émissions tout au long des activités réalisées par l'avion. La complexité de la conception moderne des avions, couplée à la multitude de disciplines impliquées, présente des défis importants pour l'optimisation,  qui devient coûteuse, de grande dimension, en entiers mixtes, multi-objectif, sous contraintes,... 
Différentes approches existent pour réaliser une MDO, dont, en particulier, les méthodes MDF (pour MultiDisciplinary Feasible) et IDF (pour Individual Disciplinary Feasible), chacune avec leurs propres avantages. La MDF résout l'analyse multidisciplinaire (MDA pour Multidisciplinary Design Analysis) à chaque étape d'optimisation, permettant d'avoir des réponses exactes ayant un sens physique. Au contraire, l'IDF sépare les solveurs disciplinaires mais résout l'optimisation et la MDA simultanément, ce qui la rend plus efficace mais cela se fait au risque d'avoir des solutions non réalisables si les itérations sont interrompues avant la convergence. Néanmoins, lorsque les dérivées des disciplines ne sont pas disponibles, comme pour des systèmes complexes ou méconnus, les méthodes sans gradient se basant sur la MDF sont cruciales. Ces méthodes sont générales, adaptables, et utilisables pour de nombreux problèmes d'ingénierie. 
Dans ce contexte, et en particulier pour optimiser des boîtes noires coûteuses, les méthodes d'optimisation basées sur des modèles de substitution, telles que l'optimisation bayésienne, sont utilisées pour explorer efficacement l'espace de conception et identifier des solutions prometteuses tout en minimisant le nombre d'évaluations coûteuses. Cette recherche vise donc à fournir des informations, des méthodologies et des outils pour l'industrie aérospatiale afin de concevoir des avions plus respectueux de l'environnement et plus efficaces selon un ensemble de critères. Ainsi, cette thèse fait suite aux thèses de doctorat de M.-A. Bouhlel (2012-2016) et R. Priem (2017-2020), qui ont développé une stratégie adaptative d'optimisation globale sous contraintes ainsi que des modèles réduits permettant de traiter des systèmes complexes de grande dimension.

Tout d'abord, le logiciel Super Efficient Global Optimization with Mixture Of Experts (SEGOMOE), développé par l'ONERA et l'ISAE-SUPAERO, est un optimiseur bayésien capable de prendre en compte efficacement les problèmes contraints en utilisant des données parcimonieuses. Il repose sur la méthode Super Efficient Global Optimization (SEGO) qui construit un modèle de substitution basé sur des processus gaussiens (GP pour Gaussian Process) à la fois pour les fonctions objectifs et les contraintes. Le processus d'optimisation itératif se base sur le modèle de substitution de la fonction objectif et permet de reformuler le problème d'optimisation originel comme la maximisation d'un critère d'enrichissement qui respecte les prédictions des modèles des contraintes. 

Ensuite, le logiciel open-source Surrogate Modeling Toolbox (SMT) regroupe plusieurs techniques de modélisation et notamment certaines qui sont améliorées par l'ajout de différentes informations comme les dérivées locales. En particulier, les modèles GP (ou krigeage ou Kriging en anglais) sont particulièrement intéressants et leurs dérivées analytiques peuvent être facilement utilisées dans un contexte de systèmes couplés, par exemple avec des méthodes adjointes. SMT implémente également des modèles de haute dimension tels que le krigeage avec moindres carrés partiels (KPLS pour Kriging with Partial Least Squares) qui combine le krigeage avec des projections en espace réduit pour construire un modèle dans un sous-espace de petite dimension, réduisant ainsi le temps de calcul nécéssaire à la construction du modèle et facilitant l'optimisation de ses paramètres internes. Suivant cette idée, la méthode KPLSK construit d'abord un modèle KPLS initial comme configuration de départ, puis construit un modèle de krigeage complet en repartant du modèle réduit pour accélérer le processus de construction.

Cette thèse étend et améliore les capacités de modélisation par processus gaussiens d'une part, et d'optimisation bayésienne avec et sans contraintes d'autre part, sur de nombreux aspects, dont les principaux sont énumérés ci-dessous.
\begin{itemize}
    \item 
    Optimisation bayésienne très haute dimension~: le développement d'une méthode d'optimisation bayésienne appelée EGORSE par R. Priem permet de résoudre des problèmes avec plusieurs centaines de variables de manière plus efficace que les méthodes existantes. A titre d'information, la contribution apportée par ce travail s'est restreinte à la prise en main de l'optimiseur SEGOMOE afin de poursuivre les expériences numériques initiées par R. Priem. 

    \item 
    Optimisation bayésienne sous contraintes avec objectifs multiples~: une généralisation du critère d'enrichissement WB2s a été développée avec R. Grapin pour l'optimisation bayésienne multi-objectif. Cette extension permet d'obtenir des fronts de Pareto avec 20 à 50 fois moins d'évaluations que les méthodes existantes. 
    Plus précisément, la contribution de ce travail à l'optimisation multi-objectif s'est concentrée sur l'encadrement du stage de recherche de R. Grapin (02/2021-02/2022) et sur le travail conjoint d'implémentation algorithmique effectué dans l'optimiseur SEGOMOE.

     \item 
   Modèles de processus gaussiens avec variables mixtes et hiérarchiques~: des modèles GP adaptés aux variables catégorielles et discrètes ont été développés, permettant d'étendre les méthodes d'optimisation bayésienne aux variables continues, discrètes ou catégorielles. De plus, des modèles hiérarchiques ont été développés pour prendre en compte des structures hiérarchiques entre les variables. Tous ces modèles sont implémentés dans le logiciel open-source SMT. 
   Ce travail a été initié durant une mobilité internationale à l'école Polytechnique de  Montréal en collaboration avec le doctorant E. Hallé-Hannan (05/2022-09/2022). Il a été poursuivi à Toulouse en collaboration avec le doctorant J. Bussemaker (DLR) lors de sa mobilité à l'ONERA (03/2023-06/2023).

      \item 
  Applications à la conception d'aéronefs plus respectueux de l'environnement~: les développements réalisés ont été appliqués à plusieurs problèmes de conception d'aéronefs. 
  Tout d'abord, une configuration de référence basée sur un A320 a été optimisée avec un ou plusieurs objectifs grâce à l'outil de conception multidisciplinaire \gls{FAST-OAD}.
  Ensuite et plus concrètement, un avion innovant hybride électrique long-courrier avec propulsion distribuée (\texttt{DRAGON}) a été optimisé en collaboration avec l'équipe M2CI de l'ONERA et en particulier grâce au concours de E. Nguyen Van. 
  Enfin, et avec la participation des partenaires européens et de l'ONERA, des optimisations ont permis de réduire la consommation de carburant, les émissions de gaz à effet de serre, le bruit émis, les coûts de développement et les risques de dysfonctionnement, tout en améliorant la qualité de fabrication de différentes configurations aéronautique afin de les rendre moins polluantes. 
  \end{itemize}

En résumé, cette partie présente les contributions et les développements réalisés dans le domaine de l'optimisation bayésienne sous contraintes en haute dimension, avec des variables mixtes et hiérarchiques. Les méthodes développées ont été appliquées à la conception d'aéronefs plus respectueux de l'environnement contribuant, par la même occasion, à l'amélioration des processus de conception aéronautique.
\section{General aircraft design context}

%
This work was carried out within the Toulouse Mathematics, Computer Science and Telecommunications doctoral school (EDMITT) and was co-funded by the Office National d'Etudes et de Recherches Aérospatiales (ONERA) and the Institut Supérieur de l'Aéronautique et de l'Espace - SUPAERO (ISAE-SUPAERO). This Ph.D was directed by Nathalie BARTOLI, Senior Research Director in the Multidisciplinary Methods and Integrated Concepts (M2CI) team of the Department of Information Processing and Systems (DTIS) at ONERA. It was co-directed by Youssef DIOUANE, Associate Professor in the Department of Mathematics and Industrial Engineering (MAGI) at Polytechnique Montréal. Furthermore, this work was supervised at ONERA by Thierry LEFEBVRE, Research Engineer in the M2CI team and, at ISAE-SUPAERO, it was supervised by Joseph MORLIER, Professor in the Department of Materials and Stuctural Mechanics (DMSM). This research work was carried out in both the Applied Mathematics (MA) unit of the Complex Systems Engineering Department (DISC) at ISAE-SUPAERO and the M2CI team at ONERA. Additionally, four months were spent in the MAGI department at Polytechnique Montréal as part of an international collaboration.

%
%

Being aerospace research laboratories, one of the main objectives pursued by ONERA and ISAE-SUPAERO is to investigate and develop future aircraft configurations. In the past few years, there has been an increasing emphasis on the development of more sustainable and eco-friendly solutions in various industries, and the aerospace sector is no exception~\cite{bravo2022unconventional}. With the ever-growing demand for air travel, it has become crucial to address the environmental impact of aircraft and to seek for more efficient and greener designs. \gls{MDA} plays a significant role in achieving this objective by integrating various disciplines into a unique framework based on coupled systems. For optimization purposes, the aeronautical industry has primarily focused on incremental improvements through aircraft design optimization. However, to meet future expectations in terms of environmental impact, noise reduction, and cost-effectiveness, substantial advancements are required. The AGILE 4.0\footnote{\url{https://www.agile4.eu/}} (2019-2023) and COLOSSUS\footnote{\url{https://colossus-sos-project.eu/}} (2023-2026) \gls{EU} funded H2020 projects address these challenges by extending the research efforts over the whole lifecycle of an aircraft, including production, operation, and end-of-life disposal~\cite{roussel2023assembly}. To achieve this ambitious goal, \gls{MDO} methods and \gls{MBSE} technologies are being employed~\cite{cramer1994problem,ramos}. The objective is to reduce the environmental impact associated with fuel consumption, waste generation, and emissions throughout the aircraft system's activities and operations. This requires collaborative efforts involving not only the aircraft design domain but also other industrial domains such as manufacturing, supply chain management, maintenance, and certification~\cite{chan2022aircraft,bartoli2023Agile,Effectiveness}. Additionally, the evaluation of new architecture systems requires consideration of categorical variables, such as on-board system architecture (conventional, hybrid, electric) and material properties (\textit{e.g.}, aluminum, titanium), which lack a defined order. When design problems include both discrete variables and continuous variables, they are said to have mixed variables. In consequence, the objective of this Ph.D thesis is to investigate and propose methodologies for high-dimensional \gls{MDO} including mixed integer inputs. The complexity of modern aircraft design, coupled with the multitude of disciplines involved, presents significant challenges for optimization. The first one is the computational time of the \gls{MDA} meaning that evaluating {the performances and constraint functions at a given operating point} is expensive. The second one is the difficulty to capture how the overall models behave. In terms of optimization, it leads to expensive-to-evaluate black-box problems in which no additional information such as the derivatives is available. This work builds upon existing \gls{MDO} frameworks and optimization techniques and adapts them to handle the challenges posed by high-dimensional mixed integer eco-design optimization problems. 
To perform \gls{MDO}, several formulations can be used~\cite{Lambe2012}, among which the \gls{MDF} and \gls{IDF} approaches are of high interest~\cite{cramer1994problem}. Each approach offers its unique set of advantages and drawbacks.
The \gls{MDF} approach is characterized by its agnostic approach for optimization where the \gls{MDA} is solved at every iteration of the optimizer. \gls{MDF} uses a non-intrusive coupling between the disciplinary solvers offering modality, adaptability, and, as a result, a physically relevant solution is obtained at every optimization stage.
Conversely, \gls{IDF} takes a different route, decoupling disciplinary solvers but solving together the optimization and the \gls{MDA} resolution. While often more computationally efficient, it may yield non-feasible solutions when the iterations are stopped before convergence. Moreover, coupling both \gls{IDF} and gradient-based optimization, when the derivative information is available, is known to have better performance and reduced computation time~\cite{Mader_ad,Martins2021}.
Notwithstanding, when we have no access to disciplinary derivatives or even to disciplines {as} a whole, gradient-free methods are of high interest and \gls{MDF} appears to be the most {appropriate} approach. The latter is particularly true for complex systems and novel configurations where little is known, as for eco-design optimization. The interest of gradient-free \gls{MDF} method is that, thanks to its generality and because it has very few requirements, it can be applied directly to many engineering problems without adaptation. This is why this method is often applied in a \gls{MDO} context for aircraft design~\cite{brevault2020multi,dufour2015trajectory,serna2009advanced,dubreuil2021development}. 
In particular, special attention is given to the development of surrogate-based optimization methods, such as \gls{BO}, to efficiently explore the design space and identify promising solutions while minimizing the number of computationally expensive evaluations.
Consequently, to optimize, without derivatives, expensive-to-evaluate black-box problems, whether coupled systems or more generally complex systems, \gls{BO} appears to be the more fitting method. 
As a matter of fact, \gls{BO} relies on surrogate models to obtain both a prediction of the unknown model and a quantification of the uncertainty in our knowledge. Mathematically, we need a surrogate model providing both a mean prediction and a variance prediction, and, in this setting, \gls{GP} based metamodels are of particular interest~\cite{KIM,DACE}. In fact, \gls{BO} leverages \gls{GP} models to predict the behavior of the objective function, enabling the identification of promising regions in the search space. \gls{BO} uses a Bayesian infill criterion that balances the exploration of unexplored regions and the exploitation of known promising regions, allowing for \gls{EGO} even with mixed, non-convex or high-dimensional functions. The versatility and effectiveness of \gls{BO} have been demonstrated in various scientific domains, ranging from materials science to environmental monitoring. Its applications include optimizing experimental conditions, designing materials with desired properties, and fine-tuning the parameters of machine learning models. By efficiently exploring the search space and accommodating limited evaluations, Bayesian optimization has proven to be a valuable tool for addressing the challenges faced by the optimization of complex systems such as aircraft designs. 

The outcomes of this research will contribute to advance the field of high-dimensional \gls{MDO} for eco-design aircraft, providing valuable insights, methodologies, and tools to the aerospace industry. The developed optimization algorithms and frameworks will enable engineers and designers to explore a vast design space, identify eco-efficient aircraft configurations, and make informed decisions that balance performance, safety, and environmental considerations, in presence of mixed variables. This work also adapts the \gls{MDO} framework to the so-called hierarchical variables. For example, consider different aircraft propulsion architectures: a conventional gas turbine would not require a variable to represent a choice in the electrical power source, while hybrid or pure electric propulsion would require such a variable. Consequently, the architectural choice influences the total number of sizing variables and therefore induces a hierarchy between the inputs variables.

In conclusion, this Ph.D thesis will address the pressing need for more eco-friendly aircraft designs by tackling the challenges of high-dimensional hierarchical and mixed \gls{MDO} under constraints. By leveraging advanced optimization techniques and integrating sustainability aspects, this research aims to improve the aircraft design process, ultimately leading to greener and more efficient aircraft that meet the incentives for a sustainable future.

\section{Previous works}
\label{sec:sota}
This work follows from the Ph.D theses of M. A. Bouhlel (2012-2016) and R. Priem (2017-2020), in which an adaptive strategy for global optimization under constraints (enrichment-based algorithm) was developed together with reduced order models adapted for high dimension.

First, the \gls{SEGOMOE} software~\cite{bartoli:hal-02149236} developed by ONERA and ISAE-SUPAERO~\cite{bartoli:hal-02149236} is a Bayesian optimizer able to take efficiently into account constrained black-box problems in a sparse data context. It is based on the \gls{SEGO}~\cite{sasena2002exploration} method that builds a \gls{GP} surrogate for the objective {function} as well as for the constraint {functions}. These constraints are taken into account during the infill criterion optimization where, for example, the optimization problem can be posed as maximizing a criterion based on the objective surrogate such that the constraints, known by surrogate predictions, are respected. Also, criteria like the \gls{UTB} criterion~\cite{SEGO-UTB} have been developed to take into account the GP models of the constraints less restrictively for multimodal constraints handling {purposes}. 

Second, the \gls{SMT} open-source software~\cite{SMT2019} regroups different modeling techniques and notably, some that are enhanced thanks to several derivatives information. In particular, the \gls{GP} models (or Kriging) are of high interest and their derivatives can be used in the context of coupled systems, with adjoint methods, for example. \gls{SMT} also implements high-dimensional models like \gls{KPLS} that combines Kriging with projective information to build the model in a small subspace, leading to both a reduced computational time and an easier optimization of the internal parameters of the aforementioned model. Furthermore, the \gls{KPLSK} method builds a first \gls{KPLS} model as a starting configuration and then builds a full Kriging model on top of it to speed up the building process. 

At the start of this Ph.D thesis, several technological locks were identified over both the modeling capabilities and the Bayesian optimization under constraints that relies on them.

\begin{itemize}
    \item During the Bayesian optimization process, the infill criteria are still to be optimized in the original high dimension to predict the best infill point to evaluate. This optimization is costly, scales poorly and {can quickly become intractable in time}. Moreover, there is a need to extend smooth infill criteria from mono to multi-objective optimization under constraints to account for opposite search interests.
    \item \gls{GP} models are well-suited for modeling continuous inputs but need to account for categorical or discrete variables that can only take a finite number of values. Moreover, these models have to be adapted to the high-dimension, as for the continuous models. 
    \item There is a need for \gls{GP} that can take into account hierarchical structure between inputs variables: for example, if we add a motor, we have one more propeller to optimize. In other terms, the bigger the number of motors, the bigger the number of sizing variables. These problems are also called "variable-size" problems and generally involve mixed categorical variables.  
\end{itemize}

The aim of this thesis is to show the recent developments done in these directions. Various test cases available in the team served as aeronautical and engineering benchmarks. The extension to mixed variables also extended the scope of applications (choice of materials, choice of electrical architectures, number of electric motors) for the \gls{EU} funded H2020 projects AGILE 4.0 (2019-2023) and COLOSSUS (2023-2026). 

\section{Contributions and developments}
\label{sec:work}

This work led to the development of high-dimensional mixed and hierarchical Bayesian optimization based on Gaussian processes. Several publications have been produced for each aspect of the thesis and the aspects over which I contributed are listed in the following.
First, the contributions to the theory of Gaussian process and Bayesian optimization can be distinguished as such: 
\begin{enumerate}
    \item  Gaussian processes with mixed variables 
    \item  High-dimensional Gaussian processes for mixed variables
    \item  Gaussian processes with hierarchical variables
    \item  Constrained Bayesian optimization with multiple objectives
    \item  Constrained Bayesian optimization in high dimension
\end{enumerate}
Then, the application to engineering test cases and aircraft design can be enumerated below.
\begin{enumerate}
    \setcounter{enumi}{5}
    \item  Software developments in \gls{SMT} for surrogate models practical application of the aforementioned GP with mixed and hierarchical variables. 
    \item  Software developments in \gls{SEGOMOE} for efficient global optimization practical application of the aforementioned Bayesian optimization methods with constraints, multi-objective and a high number of variables.
    \item  Application to the design of more environmentally-friendly aircraft. 
    \item Applications within AGILE 4.0 project for manufacturing many aircraft.
\end{enumerate}
Concerning the extensions of Bayesian optimization, in~\cite{EGORSE}, R. Priem \textit{et al.} developed a method allowing for Bayesian optimization on problems with up to 600 variables, relying on both supervised and random embeddings. This work showed that \gls{EGORSE} is more efficient and scale to high dimension better than {similar high-dimensional Bayesian optimization methods like} RREMBO~\cite{RREMBO}, TurBO~\cite{Turbo} and HESBO~\cite{Hesbo}. Through computer experiments, it has been shown that the best method is to combine \gls{PLS} regressions as supervised embeddings with Gaussian random matrices as random embeddings.
It must be noted that the personal contribution brought within the scope of this thesis was restricted to running numerical experiments to assist R. Priem. 

In~\cite{grapin_constrained_2022}, R. Grapin \textit{et al.} developed a method that generalizes \gls{WB2S} to smooth infill criteria in the context of multi-objective Bayesian  optimization. It has been shown that the developed criteria are easier to optimize than the previous ones and, with them, multi-objective Bayesian optimization was efficiently applied to aircraft design problems. Notably, Pareto fronts of the objectives have been obtained with between 20 to 50 times less evaluations than the \gls{NSGA2}~\cite{nsga2}. The contribution of this work to multi-objective optimization concerns algorithms implementation within SEGOMOE and supervision of R. Grapin internship (02/2021-02/2022).

Concerning the extensions of \gls{GP} to mixed variables, in~\cite{Mixed_Paul}, we developed a new Gaussian process model for categorical and discrete inputs to extend the Bayesian optimization methods aforementioned to  mixed integer variables. In particular, we proposed a new exponential correlation kernel that unifies both distance-based approaches and matrix-based approaches through a unified formulation. This led to a unique framework for various preexisting methods. Next, we extended this unified model to high-dimension in~\cite{Mixed_Paul_PLS} (under review) using \gls{PLS} regression. Namely, we develop a method that approximates the correlation matrix between the levels of a categorical variable with a small number of hyperparameters. This formulation gives a reduced order model to build a \gls{GP} in high-dimension for an affordable time cost.
Another milestone brought to \gls{GP} models was to take into account the hierarchy, and therefore develop new hierarchical models.
The hierarchical works have been initiated during a secondment at Polytechnique Montréal, school of Engineering, in collaboration with E. Hallé-Hannan (05/2022-09/2022). It has been followed by another collaboration with J. Bussemaker (DLR) at the time of {his} international mobility at ONERA (03/2023-06/2023). 

In~\cite{saves2023smt},
we presented a new hierarchical model and we developed a new correlation kernel peculiar to this type of variables. Moreover, the mixed integer \gls{GP}  can be coupled with the hierarchy and all these models are free-to-use and implemented in the \gls{SMT} open-source software for reproducibility. 

The new \gls{GP} models and the new Bayesian optimization capabilities have been combined and applied to various aircraft design problems. 
%
%
First, in~\cite{saves2021constrained, grapin_constrained_2022}, we optimized a reference configuration based on a A320 named CeRAS (Central Reference Aircraft System)~\cite{risse2016central}. 
Then, in~\cite{SciTech_cat}, we optimized an hybrid electric aircraft for long range missions~\cite{schmollgruber2} and we developed an adaptive criterion to select automatically the number of effective dimensions when building the models during the optimization process. This work received the \gls{AIAA} \gls{MDO} Best Paper Award 2022 in the aerospace design and structures group. In the context of AGILE 4.0, the works presented in this manuscript have been applied to several aircraft problems. Namely, in~\cite{bartoli2023Agile}, we present applications of  Bayesian optimization to multi-objective, mixed discrete and hierarchical problems. Aircraft design being multidisciplinary at its core, a lot of actors worked together for the realization of these works. At ONERA, the MDA developments have been realized by E. Nguyen Van, T. Lefebvre, N. Bartoli, C. David, S. Defoort, R. Lafage and many others.
Extensive uses have been made of the open-source softwares WhatsOpt~\cite{lafage2019whatsopt} and \gls{FAST-OAD}~\cite{David_2021}.
In fact, WhatsOpt is a web application allowing the ONERA experts to define collaboratively aircraft multidisciplinary analyses in terms of disciplines and data exchanges. 
Concerning \gls{FAST-OAD}, it is a software program developed by ONERA and ISAE-SUPAERO for aircraft sizing analysis and optimization that aims to ease the user experience and code modularity. Its aircraft sizing code is based on multidisciplinary design optimization techniques and on a point mass approach to estimate the required fuel and energy consumption for a given set of \gls{TLAR}.

\medbreak

This Ph.D thesis consists in four main chapters. Chapter~\ref{c2} corresponds to the article~\cite{Mixed_Paul} for mixed categorical \gls{GP} that is the point 1 in the list of contributions.
Chapter~\ref{c3} corresponds to the article~\cite{Mixed_Paul_PLS}, under review, for the extension of these mixed categorical \gls{GP} to high dimension with \gls{PLS} that is the point 2 in the list of contributions. 
Chapter~\ref{c4} presents the open-source software in which these models are implemented and the development of new hierarchical \gls{GP} models as in the article~\cite{saves2023smt}
, corresponding to the points 3 and 6 in the list of contributions. 
{Lastly}, Chapter~\ref{c5} regroups all the optimized application cases tackled during the last three years and correspond to our articles~\cite{EGORSE,bartoli2023Agile,SciTech_cat,Aviation_model,grapin_constrained_2022, Mixed_Paul_PLS, saves2021constrained} for really high dimension, mixed hierarchical and multi-objective optimization and regroups the points 4, 5, 7, 8 and 9 in the list of contributions.  
This has been made possible thanks to the several co-workers of this thesis that applied the developed algorithms to industrial test case of high interest, notably for green aircraft developments. 
In these works, a rover trajectory, an A320 based aircraft, an hybrid aircraft with distributed propulsion
, a coupled system for manufacturing, supply chain and design of aircraft, a regional aircraft and even a family of similar business jets were optimized. In practice, these multi-objective optimizations resulted in decreased fuel burned, cumulative emission index, emitted noise, costly expenses, manufacturing time and risk, while simultaneously increasing the manufacturing quality of aircraft, in order to "decarbonate" aeronautics.

\chapter{A mixed-categorical correlation kernel for Gaussian process} 
\chaptermark{A\MakeLowercase{ mixed-categorical correlation kernel for} G\MakeLowercase{aussian process} }
 \label{c2}

\setlength{\fboxrule}{0pt}
\hspace{3cm} \noindent\fbox{%
     \parbox{0.8\textwidth}{%
        A proposition is in itself neither probable nor improbable. An event occurs or does not occur, there is no middle course. Probability is a generalization. It involves a general description of a propositional form. Only in default of certainty do we need probability. \\
        \hrule \vspace{0.15cm}
     \hspace*{\fill} Tractatus logico-philosophicus, Ludwig Wittgenstein }%
} 



\objectif{ 
\lettrine[lines=2, lhang=0.33, loversize=0.25, findent=1.5em]{R}{ecently}, there has been a growing interest for mixed-categorical metamodels based on Gaussian process (GP) surrogates.
The main objective of this chapter is to develop a correlation kernel that can handle mixed-categorical variables to use within GP. These chapter goals are listed below.
\begin{itemize}
    \item To bridge the gap between the distance based kernels used for continuous inputs and the matrix based kernels used for categorical inputs by developing a unified and more general kernel. 
    \item To prove that the developed kernel is Symmetric Definite Positive (SPD).
    \item To use the unified kernel to build a mixed GP that generalizes continuous GP.
    \item To establish a framework based on the unified kernel that can be declined in several other state-of-the-art kernels to give new insights on these kernels.
    \item To implement the developed kernel in the open-source software SMT.
    \item To test and compare the state-of-the-art models to show their behaviours 
    and compare their performances.    
\end{itemize}

}

\minitoc
\setcounter{section}{-1}
\section{Synthèse du chapitre en français}

Le chapitre "A mixed-categorical correlation kernel for Gaussian process" propose une approche pour modéliser les variables d'entrée mixtes-catégorielles dans le cadre des processus gaussiens (GP pour Gaussian Process), également appelé modèles de krigeage. Les simulations boîtes noires coûteuses à évaluer sont couramment rencontrées dans de nombreuses applications industrielles et, dans de tels cas, les modèles de substitution sont couramment employés afin d'approcher la fonction coûteuse avec peu de données. Or, dans ce cadre, les problèmes à modéliser peuvent impliquer des variables d'entrée discrètes et catégorielles, et c'est pourquoi, l'objectif de ce chapitre  est d'apprendre un modèle de substitution peu coûteux à partir d'une fonction boîte noire donnée, tout en prenant en compte des variables continues et discrètes. En l'occurrence, nous utilisons une approche basée sur un GP qui modélise le modèle coûteux en tout point de l'espace de conception à partir d'une valeur moyenne et d'un écart-type.

La principale contribution de ce travail est le développement d'un noyau de corrélation mixte-catégoriel qui permet de modéliser avec précision les corrélations entre les variables mixtes d'entrée en combinant deux types d'approches, à savoir les approches basées distance et les approches basées matrice. 
Ce noyau de corrélation peut ensuite être utilisé pour estimer la moyenne et l'écart-type du GP.
En particulier, notre nouveau noyau de corrélation permet d'expliciter une formulation "homogène" qui unifie plusieurs approches préexistantes, qu'elles soient basées sur des distances ou sur des matrices.
Les principaux résultats obtenus montrent l'efficacité du modèle proposé sur divers cas d'étude analytiques et d'ingénierie. Ces résultats illustrent également les avantages et les inconvénients des noyaux matriciels par rapport aux noyaux basés sur des distances pour prendre en compte des variables de conception catégorielles. 

En conclusion, ce chapitre présente une approche novatrice pour modéliser les variables d'entrée mixtes-catégorielles dans le cadre des GP. Le modèle proposé offre des résultats prometteurs et ouvre la voie à de futures recherches dans ce domaine en unifiant plusieurs approches à travers une formulation unique pour exprimer les corrélations entre les variables catégorielles.

\newpage

\section{Introduction}
\renewcommand*\footnoterule{}

\label{sec:intro}

Expensive-to-evaluate black-box simulations play a key role for many engineering and industrial applications. In this context, surrogate models have shown great interest for a wide range of applications, \textit{e.g.}, aircraft design~\cite{SciTech_cat}, deep neural networks~\cite{snoek2015scalable}, coastal flooding prediction~\cite{lopez}, agriculture forecasting~\cite{MLP}, turtle retinas modeling~\cite{retina} or seismic imaging~\cite{YDiouane_SGratton_XVasseur_LNVicente_HCalandra_2016}.
These black-box simulations are generally complex and may involve mixed-categorical input variables. Typically, an aircraft design tool has to take into account variables such as the number of panels, the list of cross sectional areas or the material choices.

In this work, we target to learn an inexpensive surrogate model $ \hat{f}$ from a mixed-categorical black-box function given by
\begin{equation}
f :  \Omega \times S \times \mathbb{F}^l \to \mathbb{R}.
  \label{eq:opt_prob}
\end{equation}
This function $f$ is typically an expensive-to-evaluate simulation with no exploitable derivative information.
 $\Omega \subset \mathbb{R}^n$ represents the bounded continuous design set for the $n$ continuous variables.  $S \subset \mathbb{Z}^m$ represents the bounded integer set where $L_1, \ldots, L_m$ are the numbers of levels of the $m$ quantitative integer variables on which we can define an order relation and $ \mathbb{F}^l = \{1, \ldots, L_1\} \times \{1, \ldots, L_2\} \times  \ldots \times \{1, \ldots, L_l\}$ is the design space for the $l$ categorical qualitative variables with their respective  $L_1, \ldots, L_l$ levels.
Typical examples of $f$ can be found in different  engineering contexts. Mechanical performance of hybrid discontinuous composite materials~\cite{RaulAIAA} is an example where the mixed-categorical function $f$ represents the stiffness value which depends on a set of input variables $w=(x,c)\in \Omega \times \mathbb{F}^2$. The continuous part $x$ has two components, the length of the fibers $x_1$ and the proportion of carbon fibers $x_2$ (\textit{i.e.}, $\Omega =
[515,12000] \times [0,1]$). The categorical choices $c$ represent the types of carbon fibers $c_1$ and glass ones  $c_2$ (\textit{i.e.}, $\mathbb{F}^2 = \{\text{XN-90},\text{T800H}\} \times \{  \text{GF}, \text{T300},\text{C100},\text{C320} \} $).

For that purpose, \textit{Gaussian process} (GP)~\cite{williams2006gaussian}, also called Kriging model~\cite{krige1951statistical}, is known to be a good modeling strategy to learn a response surface model from a given dataset. Namely,
we will consider that our unknown black-box function $f$ follows a Gaussian process of mean $\mu^{f}$ and of standard deviation $\sigma^f$, $\textit{i.e.}$,
\begin{equation}
f \sim \hat{f}=
\mathcal{GP} \left(\mu^{f}, [\sigma^f]^2\right). \label{eq:GP:f}\end{equation}
For a general problem involving categorical or integer variables, several modeling strategies to build a mixed-categorical GP have been proposed~\cite{Pelamatti, Zhou, Deng, Roustant,GMHL,Gower,cuesta2021comparison,SciTech_cat}. Compared to a continuous GP, the major changes are in the estimation of the correlation matrix, the latter being essential to build estimates of $\mu^{f}$ and $\sigma^f$. Similarly to the process of constructing a GP with continuous inputs, relaxation techniques~\cite{GMHL,SciTech_cat}, continuous latent variables~\cite{cuesta2021comparison} and Gower distance based models~\cite{Gower} use a kernel-based approach to estimate the correlation matrix.
Other recent approaches try to estimate the correlation matrix independently of a kernel choice by modeling  directly the possible correlation entries of the correlation matrix~\cite{Pelamatti, Zhou, Deng, Roustant}. 

Using GP surrogates is not the only possible approach whatsoever. Random forests are often used instead of GP as they also can model both mean and variance~\cite{SMAC} and tree-structured Parzen estimators have been shown to be well-adapted for such problems~\cite{TPE}. Other surrogate models for black-box include ReLU functions~\cite{Relu-surr}, piecewise linear neural network~\cite{nn-surr} or categorical regression splines~\cite{splines}. Models other than GP could also be based on a mixed integer kernel as for support vector regression~\cite{herrera} or on a mixed integer distance as for radial basis functions~\cite{RBF_geo}. Another classical modeling strategy is to consider a different continuous model for every possible categorical choice and to build another model peculiar to the categorical variables besides the continuous models. This categorical model can be, for instance, a probability law~\cite{CAT-EGO}, a multi-arm bandit~\cite{Bandit-BO} or an integer model~\cite{AMIEGO}. Also, in case of prior information, latent variables approaches~\cite{cuesta2021comparison} and user-defined neighbourhood~\cite{Mixed_Abramson} based models are of great interest.

In this chapter, we target to extend the classical paradigm for continuous inputs (where a kernel is used to build the GP) to cover the mixed-categorical case. Namely, we will present a kernel-based approach that will lead to a unified model for existing approximation strategies~\cite{Pelamatti, GMHL,Gower}. Namely, this work unifies both distance based kernels and matrix based kernels into a unique homogeneous formulation.
This work generalizes existing methods that were already proven to be efficient over deep learning models~\cite{GMHL} and analytical test cases~\cite{Gower}. 
A similar kernel for the estimation of the correlation matrix could be applied to continuous, integer and categorical inputs. The good potential of the proposed approach is shown and analyzed over analytical and industrial test cases. 

Another main benefit behind the use of specific kernels~\cite{Pelamatti, Zhou, Deng, Roustant} to handle mixed-categorical inputs is to model accurately correlations between the  variables; which is required to get accurate GP models. It might be possible to use continuous kernels to model categorical data but in this case one needs to define a distance function. Such function is not trivial to define on categorical data; only simple distances are possible (\textit{e.g.}, Gower distance~\cite{Gower}) which in general leads to poor GP models. This chapter shows in particular the utility of mixed-categorical kernel over continuous based ones on both analytical and industrial test cases.

The GP models and the Bayesian Optimization (BO) that could be performed with them are implemented in the Surrogate Modeling Toolbox (SMT) v2.0\footnote{\url{https://smt.readthedocs.io/en/latest/}}~\cite{SMT2019}.
Our modeling software is free and open-source and has been used regularly in the aircraft industry, for example with a deep learning model~\cite{DL1, DL2, DL3, DL4} or with a deep gaussian process~\cite{DGP1,DGP2}. 

The remainder of this chapter is as follows.
In Section~\ref{sec:GP}, a detailed review of the GP model for continuous and for categorical inputs is given.
The extended kernel-based approach for constructing the correlation matrix is presented in Section~\ref{sec:model_uni}.
Section~\ref{sec:Results} presents academical tests as well as the obtained results.
Conclusions and perspectives are finally drawn in Section~\ref{sec:conclu}.

\section{GP for mixed-categorical inputs}
\renewcommand{\footnoterule}{ \hrule width5cm \vspace*{0.1cm} }    
\label{sec:GP}

In this section, we will present the mathematical background associated with GP for mixed-categorical variables. This part also introduces the notations that will be used throughout the chapter. In this section, we are considering the general case involving mixed integer variables. Namely, we assume that $f:\mathbb{R}^n \times  \mathbb{Z}^m \times \mathbb{F}^l \mapsto \mathbb{R}$ and our goal is to build a GP surrogate model for $f$. 

Given a set of data points, called a \gls{DOE}~\cite{forrester}, Bayesian inference learns the GP model that explains the best this data set. A GP model consists of a mean response hypersurface $\mu^{f}$, as well as an estimation of its variance $[\sigma^f]^2$. In the following, $n_t$ denotes the size of the given DoE data set $(W, \textbf{y}^f)$ such that $W=\{w^1,w^2,\ldots,w^{n_t}\} \in (\mathbb{R}^n \times  \mathbb{Z}^m \times \mathbb{F}^l)^{n_t}$ and $\textbf{y}^f=[f(w^1),f(w^2),\ldots,f(w^{n_t})]^{\top}$. 
For an arbitrary $w= (x,z,c) \in \mathbb{R}^n \times  \mathbb{Z}^m \times \mathbb{F}^l$, not necessary in the DoE, the GP model prediction at $w$ writes as $\hat f(w) = \mu ({w})+\epsilon(w) \in \mathbb{R} $, with $\epsilon$ being the error between $f$ and the model approximation $\mu$~\cite{GP14}. The considered error terms are random variables of variance $\sigma^2$.  Using the DoE, the expression of  $\mu^{f}$ and the estimation of its variance $[\sigma^f]^2$ are given as follows:

\begin{equation} \label{eq:mean:GP}
\mu^f(w)= \hat{\mu}^f+r(w)^\top  [R(\Theta)]^{-1}(\textbf{y}^f-\mathds{1} \hat{\mu}^f), 
\end{equation}
and
\begin{equation}
\label{eq:std:GP}
[\sigma^f(w)]^2=[\hat{\sigma}^f]^2\left[1-r(w)^\top  [R(\Theta)]^{-1}r(w)+ \frac{ \left(1-\mathds{1}^\top  [R(\Theta)]^{-1}r(w) \right)^2}{\mathds{1}^\top  [R(\Theta)]^{-1}\mathds{1}}\right], \end{equation}
where $\hat{\mu}^f$ and $\hat{\sigma}^f$, respectively, are the maximum likelihood estimator (MLE)~\cite{MLE} of $\mu$ and $\sigma$. $\mathds{1}$ denotes the vector of $n_t$ ones. $R$ is the $ n_t \times n_t $ correlation matrix between the input points and $r(w)$ is the correlation vector between the input points and a given $w$.
The correlation matrix $R$ is defined, for a given couple  $(r,s) \in (\{1,\ldots,n_t\})^2$, by \begin{equation}
\label{eq:R}
 [R(\Theta)]_{r,s}=k\left(w^r,w^s,\Theta\right) \in \mathbb{R},\end{equation}
and the vector $r(w)\in \mathbb{R}^{n_t}$ is defined as $ r(w) =[k(w,w^1), \ldots , k(w,w^{n_t})]^{\top}$,
where $k$ is a given correlation kernel that relies on a set of hyperparameters $\Theta$~\cite{de2015optimizing,bernal2022criteria}. 
The mixed-categorical  correlation kernel is given as the product of three kernels:
\begin{equation}
k(w^r,w^s,\Theta) =  k^{cont}\left(x^r,x^s,\theta^{cont}\right) k^{int}\left(z^r,z^s,\theta^{int}\right)
k^{cat}\left(c^r,c^s,\theta^{cat}\right),
\label{eq:decomp_mix}
\end{equation}
where $k^{cont}$  and $\theta^{cont}$ are the continuous kernel and its associated hyperparameters, $k^{int}$  and $\theta^{int}$ are the integer kernel and its hyperparameters, and last $k^{cat}$  and $\theta^{cat}$ are the ones related with the categorical inputs. In this case, one has $\Theta=\{ \theta^{cont},\theta^{int},\theta^{cat}\}$. 
Henceforth, the general correlation matrix $R$ will rely only on the set of the hyperparameters $\Theta$:
\begin{equation}
    \label{eq:corel:mat}
    [R(\Theta)]_{r,s} = [R^{cont}(\theta^{cont})]_{r,s}  [R^{int}(\theta^{int})]_{r,s}
    [R^{cat}(\theta^{cat})]_{r,s},
\end{equation}
where $[R^{cont}(\theta^{cont})]_{r,s} =k^{cont}(x^r,x^s,\theta^{cont}) $, $[R^{int}(\theta^{int})]_{r,s} =k^{int}(z^r,z^s,\theta^{int}) $ and  $[R^{cat}(\theta^{cat})]_{r,s}=k^{cat}(c^r,c^s,\theta^{cat})$. 
The set of hyperparameters $\Theta$ could be estimated using the DoE data set $({W},\textbf{y}^f)$ through the MLE approach on the following way
\begin{equation}
\Theta^*= \arg\max_{\Theta} \mathcal{L}(\Theta):=\left( - \frac{1}{2} {\textbf{y}^f}^\top [R(\Theta)]^{-1} {\textbf{y}^f}   - \frac{1}{2} \log 	\abs{  [R(\Theta)]} - \frac{n_t}{2} \log 2 \pi    \right),
\label{eq:likelihood}
\end{equation}
where $R(\Theta)$ is computed using~\eqnref{eq:corel:mat}. 
To construct the correlation matrix, several choices for the correlation kernel are possible. Usual families of kernels include exponential kernels or Matern kernels~\cite{Lee2011}. 
In the rest of this section, we will focus mainly on the exponential kernels and describe in details the construction of the continuous $R^{cont}(\theta^{cont})$, the integer $R^{int}(\theta^{int})$ and the categorical $R^{cat}(\theta^{cat})$ correlation matrices. 

\subsection{Correlation matrices for continuous and integer inputs}
\label{sec:GP_cont}

The construction of the correlation matrix $R^{cont}(\theta^{cont})$ for continuous inputs, based on an exponential kernel, can be described as follows. 
For a couple of continuous inputs $x^r\in \mathbb{R}^n$ and $x^s\in \mathbb{R}^n$, one sets 
\begin{equation}
[R^{cont}(\theta^{cont})]_{r,s}= {\displaystyle \prod_{j=1}^{n}{  \exp \left(-   \theta^{cont}_j  \abs{x^r_j - x^s_j}^p \right) }}. 
\label{eq:E_cont}
\end{equation}
Different values for $p$ can be used. Typically, when $p=1$, one gets the absolute exponential kernel (Ornstein-Uhlenbeck process~\cite{Lee2011}) and, when $p=2$, the squared exponential kernel (or Gaussian kernel~\cite{williams2006gaussian}) is obtained. Clearly, in the continuous case, constructing $R^{cont}(\theta^{cont})$ would require the estimation of $n$ non-negative hyperparameters, $\textit{i.e.}$, $\theta^{cont} \in \mathbb{R}^n_+$.

Thanks to a continuous relaxation technique that transforms integer inputs into continuous ones, the integer inputs can be naturally handled with continuous kernels. On this base, in what comes next, there will be no distinction between continuous and integer inputs; the two of them will be handled in the same way. In fact, for integer variables, the distance defined in the continuous case is still valid. Thus, for an integer couple $z^r\in \mathbb{Z}^m$ and $z^s\in \mathbb{Z}^m$, a natural extension of the exponential kernel that handles integer variables can be given as follows:
\begin{equation}
[R^{int}(\theta^{int})]_{r,s} = {\displaystyle \prod_{{j}=1}^{m}{  \exp \left(-   \theta^{int}_{j}  \abs{z^r_{j} - z^s_{j}}^p \right) }}.
\end{equation}
In a similar fashion, constructing $R^{int}(\theta^{int})$ would require the estimation of $m$ non-negative hyperparameters, $\textit{i.e.}$, $\theta^{int} \in \mathbb{R}^m_+$.

\subsection{Correlation matrices for categorical inputs}
\label{subsec:mi_kriging}

For categorical inputs, different choices can be made to build the correlation matrix $R^{cat}(\theta^{cat})$. Some choices are sophisticated and can therefore lead to better GP models, but are known to be computationally expensive (particularly as the number of categorical inputs increases)~\cite{Roustant,Pelamatti}. On the contrary, simple extensions of the well-known continuous kernels based on the Gower distance~\cite{Gower} or on the continuous relaxation techniques~\cite{one-hot} would be less expensive.
In the rest of this section, we will describe three known techniques to build correlation matrices for categorical inputs that are based on kernels. 

\subsubsection{Gower distance based kernel}

The Gower distance based kernel dedicates one hyperparameter per categorical input variable~\cite{Gower,RaulAIAA}. 
Namely, for two given inputs $c^r \in \mathbb{F}^l$ and $c^s \in \mathbb{F}^l$, the Hamming distance, or score, $s$ between the $i^{th}$ component of $c^r$ and $c^s$ is defined as: $s(c_i^r, c_i^s)=0$ if $c_i^r = c_i^s$, otherwise $s(c_i^r, c_i^s)=1$. Thanks to the Hamming distance, one can straightforwardly uses a continuous kernel  to define $R^{cat}(\theta^{cat})$. For instance, in the case of an exponential kernel, the Gower distance based correlation matrix will be given by 
$$
[R^{cat}(\theta^{cat})]_{r,s} =k^{cat}(c^r,c^s,\theta^{cat}) = {\displaystyle \prod_{i=1}^{l}{   \exp \left(- \theta_i^{cat} s(c_i^r,c_i^s)^p  \right) }}.  
$$
Similarly to the continuous and integer correlation matrices, the construction of the categorical correlation matrix based on the Gower distance kernel requires the estimation of $l$ hyperparameters ($\theta^{cat} \in \mathbb{R}^l_+$). Note that, as the Hamming distance can only take the values $0$ and $1$, all the exponential kernels lead to the same result independently of the value of $p$.

\subsubsection{Continuous relaxation based kernel}

To handle categorical variables through continuous relaxation, the design space $ \mathbb{F}^l $ is relaxed to a continuous space $ \Omega^l$  constructed in the following way. 
For a given $i\in \{1,\ldots,l\}$, let $c_i$ be the $i^{th}$ categorical variable with $L_i$ levels, and, for a given input point $c^r$, let $\ell_r^i$ be the index of the level taken by $c^r$ on the variable $i$.
Denote $e_{c^r_i}$ the one-hot encoding~\cite{one-hot} of $c^r_i$ that takes value $0$ everywhere but on the dimension $\ell_r^i$:  $e_{c^r_i} \in \mathbb{R}^{L_i} $ such that $\left( e_{c^r_i} \right)_{\ell^i_r}=1$ and $\left( e_{c^r_i} \right)_k=0$ for $k \neq \ell^i_r$. 
For example, if the $i^{th}$ component is the color with $L_i=3$ levels being $ \{ \mbox{red}, \mbox{blue}, \mbox{green} \}$ and if the $r^{th}$ variable takes value $\mbox{blue}$ ($c_i^r = \mbox{blue}$), then the corresponding index is $\ell_r^i =2$ and the corresponding one-hot encoding is $e_{c^r_i} = (0,1,0)$.
The continuous relaxation idea is as follows. At the beginning we set $ \Omega^l$ to be empty, then, for each $i  \in \{1, \ldots,l\}$, a relaxed one-hot encoding is used for $c_i$. The latter increases the dimension of the relaxed continuous space $ \Omega^l$ by $L_i$ and, at the end of the relaxation, we get the final continuous design space $\Omega^l\subseteq \{0,1\}^{n^l}$, where $n^l=\sum_{i=1}^l L_i >l$. 
Like the Gower distance based kernel, the continuous relaxation based kernel adapts continuous kernels to handle categorical variables, $\textit{i.e.}$, for a couple of categorical inputs $c^r$ and $c^s$,
\begin{equation*}
[R^{cat}(\theta^{cat})]_{r,s} =k^{cat}(c^r,c^s,\theta^{cat}) ={\displaystyle \prod_{i=1}^{l}{ k^{cat}( c_i^r, c_i^s,\theta^{cat})}}= {\displaystyle \prod_{i=1}^{l}{ {\displaystyle \prod_{j=1}^{L_i} k^{cont}( [e_{c_i^r}]_j, [e_{c_i^s}]_{j},\theta^{cat})}}}.    
\end{equation*}
Typically, for an exponential continuous kernel, one has
\begin{equation}
\begin{split}
[R^{cat}(\theta^{cat})]_{r,s} &= {\displaystyle \mathlarger{\prod}_{i=1}^{l} {\displaystyle \mathlarger{\prod}_{j=1}^{L_i} \exp\left(-\theta_{ \sum_{i'=1}^{i-1} {L_{i'}+j}}^{cat}\abs{[e_{c^r_i}]_j - [e_{c^s_i}]_j} ^p\right) }},
\end{split}
\label{eq:kernel}
\end{equation}
and, by using the one-hot encoding structure of $e_{c^r_i}$ and $e_{c^s_i}$, it leads to
\begin{equation*}
\begin{split}
[R^{cat}(\theta^{cat})]_{r,s} &= {\displaystyle \mathlarger{\prod}_{i=1}^{l}  \exp\left(- \theta_{ \sum_{i'=1}^{i-1} {L_{i'}+ \ell^i_r}}^{cat} - \theta_{ \sum_{i'=1}^{i-1} {L_{i'}+ \ell^i_s}}^{cat} \right). } 
\end{split}
\label{eq:kernel:2}
\end{equation*}
Hence, this kernel relies on $n^l=\sum_{i=1}^l L_i$ hyperparameters ($\theta^{cat} \in \mathbb{R}^{n^l}_+$) which can be much more higher than the number of hyperparameters required to build the Gower distance based kernel. Due to one-hot encoding strategy, the value of $p$ is also irrelevant for the construction of the continuous relaxation based kernel.

\newpage

\subsubsection{Homoscedastic hypersphere kernel}
\label{subsec:ho_hs}
The idea of the homoscedastic hypershere kernel~\cite{Roustant,Pelamatti} is to directly model the correlation matrix instead of looking for a kernel function.
The use of a kernel function guarantees the related correlation matrix $R^{cat}$ to be symmetric positive definite (SPD). However, with the homoscedastic hypershere kernel, one will directly construct an SPD matrix with the desired properties. Namely, for a given $i \in \{1, \ldots, l\}$, let $c^r_{i} $  and $c^s_{i} $  be a couple of categorical variables taking respectively the $\ell^i_r$ and the $\ell^i_s$ level on the categorical variable $c_i$, $[R^{cat}(\theta^{cat})]_{r,s}$ can be formulated in a  level-wise form~\cite{Pelamatti} as: 
\begin{equation}
[R^{cat}(\theta^{cat})]_{r,s}=k^{cat}(c^r,c^s,\theta^{cat}) 
= {\displaystyle \prod_{i=1}^{l}  [R_i(\Theta_i)]_{\ell^i_r,\ell^i_s} }= {\displaystyle \prod_{i=1}^{l}  [C(\Theta_i)C(\Theta_i)^{\top}]_{\ell^i_r,\ell^i_s} }.
\label{eq:homo_HS}
\end{equation}
For all $i \in \{1, \ldots, l\}$, the matrix $C(\Theta_i)\in \mathbb{R}^{L_i \times L_i}$ is lower triangular and built using a hypersphere decomposition~\cite{HS,HS_Jacobi} from a symmetric matrix $\Theta_i \in \mathbb{R}^{L_i \times L_i}$ of hyperparameters. For any $k, k' \in \{1,\ldots, L_i \}$, the  matrix $C(\Theta_i)$ is given by: 
\begin{equation}
 \begin{cases}
 [C(\Theta_i)]_{1,1} = 1, \\
 [C(\Theta_i)]_{k,1} =  \cos\left([\Theta_i]_{k,1}\right) ~\mbox{for any}~  2 \leq k \leq L_i \\
 [C(\Theta_i)]_{k,k'} =  \cos\left([\Theta_i]_{k,k'}\right)  \prod_{j=1}^{k'-1} \sin\left([\Theta_i]_{k,j}\right),~\mbox{for any}~ 2 \le k'< k \leq L_i \\
 [C(\Theta_i)]_{k,k} =    \prod_{j=1}^{k-1} \sin\left([\Theta_i]_{k,j}\right),~\mbox{for any}~ 2 \le k \leq L_i,  \\
\end{cases} 
\label{eq:hy_decomp}
\end{equation}
where the hyperparameters are set such that  $[\Theta_i]_{k,k'} \in [0,\pi]$ for all $1 \le k'< k \leq L_i $. 
For this kernel, the hyperparameters $\theta^{cat}$ can be seen as a concatenation of the set of symmetric matrices, \textit{i.e.}, $\theta^{cat} = \{  \Theta_1, \Theta_2, \ldots, \Theta_l \} $. The construction of this kernel is thus relying on the estimation of $\sum_{i=1}^l \frac{1}{2} L_i (L_i-1) $ hyperparameters. Unlike the previous kernels where the elements of the correlation matrix are non-negative, the correlation values for the homoscedastic hypersphere kernel can be negative, \textit{i.e.}, $[R^{cat}(\theta^{cat})]_{r,s} \in [-1,1]$.

\section{An exponential kernel-based model for categorical inputs}
\label{sec:model_uni}

In this section, we propose an  extension of the classical exponential kernels (used for continuous inputs) to handle categorical variables.
Thanks to the one-hot encoding, we can replace the distance-based approach by an hyperparameter-based approach. 
This extension will naturally lead to a generalization of both continuous relaxation and Gower distance based kernels.  

Distance based approaches (like Gower distance or continuous relaxation) can not model every possible correlation between the various categorical choices. Therefore, these methods do not lead to an exhaustive GP model but to an imprecise approximation. In what follows, we propose to introduce a new formulation that includes a correlation matrix so that we could reach a higher accuracy for the resulting distance-based GP model.
To begin with, the continuous relaxation kernel described in~\eqnref{eq:kernel} can be reformulated as:
\begin{equation}
\begin{split}
[R^{cat}(\theta^{cat})]_{r,s}
&= {\displaystyle \mathlarger{\prod}_{i=1}^{l} {\displaystyle \mathlarger{\prod}_{j=1}^{L_i} \exp\left(- \abs{[e_{c^r_i} - e_{c^s_i}]_j} ^{p/2} [{\Theta}_{i}]_{j,j} \abs{[e_{c^r_i} - e_{c^s_i}]_j} ^{p/2}\right) }},
\end{split}
\label{eq:CR2}
\end{equation}
where, for all $i=1, \ldots,l$, the matrix $\Theta_i \in \mathbb{R}^{L_i\times L_i}$ is diagonal  such that $[\Theta_i]_{j,j}=\theta_{ \sum_{i'=1}^{i-1} {L_{i'}+j}}^{cat} \in \mathbb{R}_+$, and $\theta^{cat}$ is defined as the list of hyperparameter matrices $\theta^{cat} = \{\Theta_1, \ldots, \Theta_l\}$.
The idea of the new kernel is the following: we start from the reformulation of the continuous relaxation kernel
of~\eqnref{eq:CR2}. Then, as for the kernel of~\eqnref{eq:homo_HS}, we consider, for every categorical variable $i=1, \ldots,l$, a SPD matrix $\Phi( \Theta_i) \in  \mathbb{R}^{L_i\times L_i}$ used to build a kernel associated with the correlation matrix $R_i(\Phi(\Theta_i))$.
Let $c^r_{i} $  and $c^s_{i} $  be a couple of categorical variables taking respectively the $\ell^i_r$ and the $\ell^i_s$ level of the variable $c_i$, we set
\begin{equation}
[R^{cat}(\theta^{cat})]_{r,s} = {\displaystyle \prod_{i=1}^{l}  [R_i(\Phi(\Theta_i))]_{\ell^i_r,\ell^i_s} },
\label{eq:nat_ext_gaussian_ker:0}
\end{equation}
and, for all $i=1, \ldots,l$, one has
\begin{equation}
[R_i(\Phi(\Theta_i))]_{\ell^i_r,\ell^i_s}  
=   {\displaystyle \prod_{j=1}^{L_i}  {\displaystyle \prod_{j'=1}^{L_i}  
\exp \left(-   \abs{ [e_{c^r_i} -  e_{c^s_i}]_j}^{p/2} [\Phi(\Theta_i)]_{j,j'}  \abs{ [e_{c^r_i}   -  e_{c^s_i}]_{j'}}^{p/2} \right) }},
\label{eq:nat_ext_gaussian_ker}
\end{equation}
where $\ell_r^i$ and $\ell_s^i$ are the indices of the levels taken by the variables $c^r$ and $c^s$, respectively, on the $i^{th}$ categorical variable and the coefficient $[\Phi(\Theta_i)]_{{\ell^i_r},{\ell^i_s}}$ is characterizing the correlation between these two levels. 

\begin{remark}
One can easily see that~\eqnref{eq:nat_ext_gaussian_ker} generalizes the continuous relaxation approach. In fact, by setting  $\Phi( \Theta_i)=\Theta_i$ to be a diagonal matrix, we recover~\eqnref{eq:CR2}.
\end{remark}

Now, by using the one-hot encoding nature of the vectors $e_{c_i^r}$ and $e_{c_i^s}$, we get naturally what follows. Namely, if $c^r_i =c^s_i$, one deduces that $ [R_i(\Phi(\Theta_i))]_{\ell^i_r,\ell^i_s}   = \exp (0) = 1$. Otherwise, if $ c^r_i \neq c^{s}_i$, we get
\begin{equation}
\begin{split}
[R_i(\Phi(\Theta_i))]_{\ell^i_r,\ell^i_s}   &= \exp  \left( - {\displaystyle \sum_{j=1}^{L_i} {\displaystyle \sum_{j'=1}^{L_i}{    \abs{  [e_{c^r_i} -  e_{c^s_i}]_j }^{p/2} [\Phi(\Theta_i)]_{j,j'}  \abs{[e_{c^r_i}   -  e_{c^s_i}]_{j'}}^{p/2}  }}} \right) \\
&= \exp \left( - \left([\Phi(\Theta_i)]_{{\ell^i_r},{\ell^i_r}}+[\Phi(\Theta_i)]_{{\ell^i_s},{\ell^i_s}} +[\Phi(\Theta_i)]_{{\ell^i_r},{\ell^i_s}} +[\Phi(\Theta_i)]_{{\ell^i_s},{\ell^i_r}} \right) \right)\\ 
&= \exp \left( - [\Phi(\Theta_i)]_{{\ell^i_r},{\ell^i_r}}-[\Phi(\Theta_i)]_{{\ell^i_s},{\ell^i_s}} - 2[\Phi(\Theta_i)]_{{\ell^i_r},{\ell^i_s}} \right). 
\end{split}
\label{eq:mat_ker_i}
\end{equation}

\begin{remark}
Note that the resulting correlation matrix $R_i(\Phi(\Theta_i))$ does not depend on the chosen parameter $p$ (used within the definition of the exponential kernels). Therefore, in our case, when dealing with categorical variables kernels, there will be no distinction between squared or absolute exponential  kernels. 
\end{remark}

In addition, as far as the matrices $\Theta_i$ respect a specific parameterization, we will show that our approach guarantees that the correlation matrix $R$ is SPD with a unit diagonal and off-diagonal terms values in $[0,1]$~\cite{PDUDE}. In general, the latter properties are required to be satisfied by the correlation matrices. Otherwise, one may get numerical issues to build the GP model, see~\eqnref{eq:mean:GP} and~\eqnref{eq:std:GP}. For that purpose,
for a given $i \in \{1,\ldots, l\}$, we propose to use the following parameterization for the hyperparameter matrix $\Phi(\Theta_i)$:
\begin{equation}
\begin{split}
[\Phi(\Theta_i)]_{j,j} &:= [\Theta_i]_{j,j} ~  \geq 0  \\
[\Phi(\Theta_i)]_{j,j'} &:= \frac{\log \epsilon }{2} ([C(\Theta_i) C(\Theta_i)^\top]_{j,j'} -1)  ~~~\mbox{if}~{j\neq j'},   \\
\end{split}
\label{eq:hs_exp}
\end{equation}
where the parameter $\epsilon$ is chosen as a small positive tolerance ($ 0<\epsilon \ll 1$) and the matrix
$C(\Theta_i)$ is a Cholesky lower triangular matrix that relies on the symmetric matrix $\Theta_i$ of $L_i(L_i-1)/2$ elements in $[0,\frac{\pi}{2}]$. The elements of  $C(\Theta_i)$ represent the coordinates of a point on the surface of a unit radius sphere as in~\cite{ Roustant,Pelamatti}. They are described in~\eqnref{eq:hy_decomp}.  Note that, by taking into consideration the symmetry of the matrix $\Theta_i$, the total number of hyperparameters for the categorical variable $i$ is $\frac{L_i(L_i+1)}{2}$.

In the next theorem, we will show that the parameterization given by~\eqnref{eq:hs_exp} guarantees the desirable properties for the correlation matrices $R_i(\Theta_i)$ and therefore for the matrix $R^{cat}$. In particular, we will show  that the matrix $R^{cat}(\theta^{cat})$ is SPD with elements in $[0,1]$, $\textit{i.e.}$, for all $s,r \in \{1,\ldots, n_t\}$,~ $[R^{cat}(\theta^{cat})]_{r,s} \in [0,1]$.

\begin{theorem}
\label{th:definiteness}
 Assume that, for all $i \in \{1, \ldots, l\}$, $ \Phi(\Theta_i)$ satisfies the parameterization of~\eqnref{eq:hs_exp}.
Then the matrix $R^{cat}(\theta^{cat})$, given by~\eqnref{eq:nat_ext_gaussian_ker:0}, is SPD with  elements in $[0,1]$.
\end{theorem}
\begin{proof}Indeed, for all $i \in \{1, \ldots, l\}$, by using~\eqnref{eq:mat_ker_i} and~\eqnref{eq:hs_exp}, one has
\begin{equation*}
\begin{split}    
[R_i(\Phi(\Theta_i))]_{\ell^i_r,\ell^i_s}&=[W_i]_{\ell^i_r,\ell^i_s} [T_i]_{\ell^i_r,\ell^i_s}, ~~~~\mbox{if}~ \ell^i_r \neq \ell^i_s \\
[R_i(\Phi(\Theta_i))]_{\ell^i_r,\ell^i_r}&= 1,
\end{split}
\end{equation*}

where $[W_i]_{\ell^i_r,\ell^i_s}=\exp \left( - [\Theta_i]_{{\ell^i_r},{\ell^i_r}}-[\Theta_i]_{{\ell^i_s},{\ell^i_s}}\right)$ and $[T_i]_{\ell^i_r,\ell^i_s}=\exp \left( - 2 [\Phi(\Theta_i)]_{{\ell^i_r},{\ell^i_s}} \right)$.
The matrix $R_i(\Phi(\Theta_i))$ is thus defined as a Hadamard product  ($\textit{i.e.}$, element-wise product of matrices)~\cite{hadamard}. Hence, by application of the Schur product theorem~\cite[Lemma~3.7.1]{schurmatapp}, it suffices to show that the matrices $W_i$ and $T_i$ are SPD to prove that $R_i$ is also SPD~\cite{Schur}. 
Taking into account that, for all $s,r \in \{1,\ldots, n_t\}$, $e_{c^r_i}$ and $e_{c^s_i}$ are one-hot encoding elements of $\mathbb{R}^{L_i}$,  the matrix $W_i$ corresponds to the correlation matrix associated with the exponential kernel in the continuous space, \textit{i.e.},
$$[W_i]_{\ell^i_r,\ell^i_s}=\exp \left( - [\Theta_i]_{{\ell^i_r},{\ell^i_r}}-[\Theta_i]_{{\ell^i_s},{\ell^i_s}}\right)={\displaystyle \prod_{j=1}^{L_i} {  \exp \left(-   [\Theta_i]_{j,j}  \abs{ [e_{c^r_i}   -  e_{c^s_i}]_{j} }^p \right) }}.$$ 
Hence, since  the diagonal elements of $\Theta_i$ are positive, the matrix $W_i$ is SPD. In fact, the kernel function $\phi(x) = \exp [- \theta |x|^p]$ is positive definite for a given positive $\theta$ if $0<p \leq 2$~\cite[Corollary~3]{Schoenberg}. 
Regarding the matrix $T_i$, by using the fact that $\Theta_i$ satisfies~\eqnref{eq:hs_exp}, one has
\begin{equation}
\label{eq:EHH}
[T_i]_{\ell^i_r,\ell^i_s}= \epsilon  \exp \left( -(\log \epsilon) [C(\Theta_i) C(\Theta_i)^{\top}]_{\ell^i_r,\ell^i_s} \right).
\end{equation}
For an $\epsilon \in (0,1)$,  the matrix $-(\log \epsilon) [C(\Theta_i) C(\Theta_i)^{T}]$ is SPD as a Cholesky like-decomposition matrix. Thus, $T_i$ is also SPD as the Hadamard exponential of an SPD matrix~\cite[Theorem 7.5.9]{horn2012matrix}.

For the second part of the proof, the matrix $C(\Theta_i)$ is constructed by hypersphere decomposition such that the values of $C(\Theta_i) C(\Theta_i)^\top$ belong to $[0,1]$~\cite{hypersphere}. Hence,
\begin{equation*}
\begin{split}
 &[T_i]_{\ell^i_r,\ell^i_s}=\epsilon  \exp \left(  -(\log \epsilon)  [C(\Theta_i) C(\Theta_i)^\top]_{\ell^i_r,\ell^i_s} \right) \geq \epsilon  \exp 0 =  \epsilon,  \\ 
&[T_i]_{\ell^i_r,\ell^i_s}=\epsilon  \exp \left(  -(\log \epsilon)  [C(\Theta_i) C(\Theta_i)^\top]_{\ell^i_r,\ell^i_s} \right) \leq \frac{\epsilon}{\epsilon} =1.
\end{split}
\end{equation*}
Also, the elements of the matrix $W_i$ are in $[0,1]$ since the diagonal elements of $\Theta_i$ are chosen to be positive. Consequently, the extra-diagonal elements of $R_i$ are in $[0,1]$.
Finally, the Hadamard product being conservative for those two latter properties, one concludes that the correlation matrix $R^{cat}(\theta^{cat})$ is SPD and all its elements are in $[0,1]$. 
\end{proof}

\begin{remark}
\label{rmk:2}
For a given small $\epsilon>0$, the transformation
\begin{equation}
\label{eq:bijection_alpha}
\alpha \to \epsilon \exp[-\log(\epsilon) (\alpha-1)]    
\end{equation}
is a bijection over $[\epsilon,1]$, thus one can deduce that (when $[T_i]_{j,j'}> \epsilon$ for all $j,j'$) there exists a unique  matrix $\hat \Theta_i$ such that 
$$
T_i = C(\hat \Theta_i)C(\hat \Theta_i)^{\top}.
$$
This, in particular, shows that if we set $W_i$ to identity in our parameterization, then, as far as the correlations are larger than $\epsilon$, the homoscedastic hypersphere parameterization of Zhou \textit{et al.}~\cite{Zhou, pinheiro_unconstrained_1996} is equivalent to our proposed one.
\end{remark}

In the next theorem, using the hypersphere decomposition properties~\cite{HS_Jacobi}, we will show that the correlation matrix $R_i$, as given by ~\eqnref{eq:mat_ker_i}, can be built in an equivalent way without the diagonal elements of the matrix  $\Phi(\Theta_i)$. Such result is of high interest as it reduces the number of hyperparameters from $\frac{L_i(L_i+1)}{2}$ to  $\frac{L_i(L_i-1)}{2}$ per categorical variable $i$ without any loss in the accuracy in the final model.

\newpage

\begin{theorem}
\label{lemma_eq_ho_hs}
The correlation matrix $R_i$, as given by ~\eqnref{eq:mat_ker_i}, can be rewritten as follows
\begin{equation}
\label{new:R}
\begin{split}
[R_i(\Phi(\bar \Theta_i))]_{\ell^i_r,\ell^i_s}=& \exp \left( - 2 [\Phi(\bar \Theta_i)]_{{\ell^i_r},{\ell^i_s}} \right), ~~~~\mbox{if}~ \ell^i_r \neq \ell^i_s \\
[R_i(\Phi(\bar \Theta_i))]_{\ell^i_r,\ell^i_r}=& 1,
\end{split}
\end{equation}
where $[\Phi(\bar \Theta_i)]_{{\ell^i_r},{\ell^i_s}}= \frac{\log \epsilon }{2} ([C(\bar\Theta_i) C(\bar\Theta_i)^\top]_{{\ell^i_r},{\ell^i_s}} -1)$ and $\bar\Theta_i$ is a symmetric matrix whose diagonal elements are set to zero (\textit{i.e.}, $[\bar\Theta_i]_{j,j} =0$ for all $j=1, \ldots, L_i$).
\end{theorem}
\begin{proof}
Indeed, by using the hypersphere decomposition~\cite{HS_Jacobi}, any SPD matrix $T_i(\Theta_i)$ with unitary diagonal and values in $[\epsilon,1]$ can be modeled as $T_i(\Theta_i) = [C({\hat \Theta}_i) C({\hat\Theta}_i)^\top]$ from a certain symmetric matrix ${\hat \Theta_i}$ without using additional diagonal elements (\textit{i.e.}, $[\hat \Theta_i]_{j,j}=0$ for all $j=1, \ldots, L_i$).
Thus, using the fact that $R_i(\Theta_i)$ is written as the image of this SPD matrix $T_i(\Theta_i)$ by the element-wise transformation of ~\eqnref{eq:bijection_alpha} bijective over $[\epsilon,1]$, one can deduce that there must exist a symmetric matrix ${\bar \Theta_i}$ whose diagonal elements are set to zero (\textit{i.e.}, $[\bar \Theta_i]_{j,j}=0$ for all $j=1, \ldots, L_i$) such that 
\begin{equation*}
\begin{split}
[R_i(\Phi(\bar \Theta_i))]_{\ell^i_r,\ell^i_s}=& \exp \left( - 2 [\Phi(\bar \Theta_i)]_{{\ell^i_r},{\ell^i_s}} \right), ~~~~\mbox{if}~ \ell^i_r \neq \ell^i_s \\
[R_i(\Phi(\bar \Theta_i))]_{\ell^i_r,\ell^i_r}=& 1,
\end{split}
\end{equation*}
where $[\Phi(\bar \Theta_i)]_{{\ell^i_r},{\ell^i_s}}= \frac{\log \epsilon }{2} ([C(\bar\Theta_i) C(\bar\Theta_i)^\top]_{{\ell^i_r},{\ell^i_s}} -1)$. 
\end{proof}
 
 In what comes next, we will refer to our kernel when it uses the parameterization of~\eqnref{new:R} as the Exponential Homoscedastic Hypersphere (EHH) kernel (see Remark~\ref{rmk:2}). We will call the original parameterization of the correlation matrix (as given by~\eqnref{eq:mat_ker_i}) as the Fully Exponential (FE) kernel. Note that, as explained in~\ref{apendix:EHH2CR}, whenever the matrix $\Theta_i$ is chosen to be diagonal, the matrix $\Phi(\Theta_i)$ used within FE kernel will also be diagonal. Thus, we are able to recover the continuous relaxation kernel~\cite{GMHL} and  this parameterization will be called the Continuous Relaxation (CR) kernel. Similarly, if we choose $\Theta_i $ to be of the form $\theta_i \times I_{L_i}$ where $\theta_i \in \mathbb{R}^+$ and $I_{L_i}$ is the identity matrix of size $L_i$, we are able to recover the Gower distance based kernel~\cite{Gower}: this parameterization will be called the Gower Distance (GD) kernel. 
 
 As mentioned earlier, the EHH kernel is similar to the FE kernel. Therefore, we can deduce that the EHH kernel generalizes the CR kernel and also that the CR kernel generalizes the GD kernel. Table~\ref{tab:General scheme} gives all the details associated with the four categorical kernels described above, \textit{i.e.}, GD, CR, EHH and FE.

\begin{table}[H]
   \caption{Description of the four categorical kernels (GD, CR, EHH and FE) using our proposed exponential parameterization.}
   \vspace*{-0.2cm}
   \begin{center}
   \resizebox{1\columnwidth}{!}{%
      \begin{tabular}{cccc}
       \hline
  \textbf{Kernel} &    $\Theta_i=$ & 

 $[R_i(\Theta_i)]_{{\ell^i_r},{\ell^i_s}}=$ 
          
       &  \# of Hyperparam. \\
       \hline 
 & & &  \\
 \textbf{GD} 
 & 
 \scalebox{1.2}{$ \displaystyle \frac{\theta_i}{2} \  \times $  } $  
\begin{bmatrix}
1 & \textcolor{white}{9} & \hspace{2em} { \textbf{\textit{ Sym.}}}  \textcolor{white}{9} & \\
0  &1 & \textcolor{white}{9} \\
\vdots & \ddots &  \ddots & \textcolor{white}{9}  \\
0 &  \ldots & 0 &1 \\
\end{bmatrix}$  &
$ \exp \left( - [\Phi(\Theta_i)]_{{\ell^i_r},{\ell^i_r}}-[\Phi(\Theta_i)]_{{\ell^i_s},{\ell^i_s}}  \right) $
&1\\
 & & &  \\
 & & &  \\
 \textbf{CR }
 & 
  $   
\begin{bmatrix}
[\Theta_i]_{1,1}  & \textcolor{white}{9} & \textcolor{white}{9} & \hspace{-2em} { \textbf{\textit{ Sym.}}} \\
0 &[\Theta_i]_{2,2}  & \textcolor{white}{9} \\
\vdots & \ddots &  \ddots & \textcolor{white}{9}  \\
0 &  \ldots & 0 & [\Theta_i]_{L_i,L_i} \\
\end{bmatrix}$ & 
  $ \exp \left( - [\Phi(\Theta_i)]_{{\ell^i_r},{\ell^i_r}}-[\Phi(\Theta_i)]_{{\ell^i_s},{\ell^i_s}}  \right) $
& $ L_i$ \\
 & & &  \\
 & & &  \\
 \textbf{EHH}  & 
$ 
\begin{bmatrix}
0 & \textcolor{white}{9} & \hspace{2em} { \textbf{\textit{ Sym.}}}  \textcolor{white}{9} & \\
[\Theta_i]_{1,2}  &0 & \textcolor{white}{9} \\
\vdots& \ddots &  \ddots & \textcolor{white}{9}  \\
[\Theta_i]_{1,L_i} &  \ldots & [\Theta_i]_{L_i-1,L_i} &0 \\
\end{bmatrix}$ &  
$\exp \left(  - 2  [\Phi(\Theta_i)]_{{\ell^i_r},{\ell^i_s}} \right)$
&  $\frac{1}{2} L_i(L_i-1) $ \\
 & & &  \\
      & & &  \\
\textbf{FE} & 

$\begin{bmatrix}
[\Theta_i]_{1,1} & \textcolor{white}{9} & \hspace{2em} { \textbf{\textit{ Sym.}}}  \textcolor{white}{9} & \\
[\Theta_i]_{1,2}  & [\Theta_i]_{2,2} & \textcolor{white}{9} \\
\vdots &\ddots & \ddots & \textcolor{white}{9}  \\
[\Theta_i]_{1,L_i} &  \ldots & [\Theta_i]_{L_i-1,L_i} &[\Theta_i]_{L_i,L_i} \\ 
\end{bmatrix} $  & 
$\exp \ ( - [\Phi(\Theta_i)]_{{\ell^i_r},{\ell^i_r}}-[\Phi(\Theta_i)]_{{\ell^i_s},{\ell^i_s}}- 2[\Phi(\Theta_i)]_{{\ell^i_r},{\ell^i_s}} )$ &

  $\frac{1}{2} L_i(L_i+1) $ \\
   & & &  \\
        \hline  
      \end{tabular}
      }
   \end{center}
   \label{tab:General scheme}
\end{table}

{\color{black}
To sum up, we have seen that the HH kernel can be more general than the EHH  one (as it can deal with negative correlations) and that EHH generalizes both CR and GD. All these categorical models can be unified in a single formulation as follows. 
For each $i \in \{1, \ldots, l\}$,  a hyperparameter matrix $\Theta_i$ is associated with each variable $c_i$, \textit{i.e.},
$$\Theta_i= \begin{bmatrix}
[\Theta_i]_{1,1} & \textcolor{white}{9} & \hspace{2em} { \textbf{\textit{ Sym.}}}  \textcolor{white}{9} & \\
[\Theta_i]_{1,2}  & [\Theta_i]_{2,2} & \textcolor{white}{9} \\
\vdots &\ddots & \ddots & \textcolor{white}{9}  \\
[\Theta_i]_{1,L_i} &  \ldots & [\Theta_i]_{L_i-1,L_i} &[\Theta_i]_{L_i,L_i} \\ 
\end{bmatrix}.$$
The correlation term $[R_i(\Theta_i)]_{\ell^i_r,\ell^i_s} $ associated with $c_i$ can be formulated in the following level-wise form: 
\begin{eqnarray*}
[R_i(\Theta_i)]_{\ell^i_r,\ell^i_s} &= &\  \kappa ( 2 [ \Phi(\Theta_i) ]_{{ \ell_i^r},{\ell_i^s}} ) \  \kappa ( [ \Phi(\Theta_i) ]_{{ \ell_i^r},{\ell_i^r}} ) \  \kappa ( [ \Phi(\Theta_i) ]_{{ \ell_i^s},{\ell_i^s}} ), 
\label{eq:homogeneous}
\end{eqnarray*}
where $\kappa$ can be set either to any positive definite kernel (to get  GD, CR, EHH or FE kernels) or to identity (to get HH kernel). The transformation 
function $\Phi(.)$ is selected such that, for any SPD matrix $\Theta_i$, the output matrix $\Phi(\Theta_i)$ is also SPD. \tabref{tab:kernels} gives a list of possible choices for $\Phi$ when the is function $\kappa$ is set to exponential or identity. 
For all categorical variables $i \in \{1, \ldots, l\}$, the matrix $C(\Theta_i)\in \mathbb{R}^{L_i \times L_i}$ (lower triangular) is  built using a hypersphere decomposition. 

\begin{table}[H]
\caption{Kernels using different choices for the function $\Phi$.}
\begin{tabular}{ccl}
\hline
\textbf{Kernel} & $\kappa(\phi)$  &  \hspace{3cm} ${\centering \Phi(\Theta_i)}$  \\
\hline
\textbf{GD}   &  $\exp(-\phi) $ & ${ \displaystyle [\Phi(\Theta_i)]_{j,j} := \frac{1}{2}  \theta_{i} \quad  ~;~ [\Phi(\Theta_i)]_{j \neq j'} := 0 }$    \\
\textbf{CR}  & $\exp(-\phi) $ &  $  { \displaystyle [\Phi(\Theta_i)]_{j,j} := [\Theta_i]_{j,j} ~;~ [\Phi(\Theta_i)]_{j \neq j'} := 0 } $   \\
\textbf{EHH}  & $\exp(-\phi)$ & 
$  { \displaystyle [\Phi(\Theta_i)]_{j,j} := 0 \quad \quad ~;~ [\Phi(\Theta_i)]_{j \neq j'} := \frac{\log \epsilon }{2} ([C(\Theta_i) C(\Theta_i) ^\top]_{j,j'} -1)  }$ \\
\textbf{HH} &  $\phi$ &    $  { \displaystyle [\Phi(\Theta_i)]_{j,j} := 1 \quad \quad ~;~ [\Phi(\Theta_i)]_{j \neq j'} := \frac{1}{2} [C(\Theta_i) C(\Theta_i)^\top]_{j,j'} }$  \\
\hline
\end{tabular}
\label{tab:kernels}
\end{table}
}

In the next section, we will see how these kernels perform on different test cases. In particular, we study numerically the trade-off between the kernel efficiencies and their respective computational efforts (related directly to the number of hyperparamters).

\section{Results and discussion}
\label{sec:Results}

In this section, we propose several illustrations and comparisons on three different test cases (from 2 to 10 continuous variables and 1 or 2 categorical variables up to 12 levels) to show the interest of our method and the equivalence with other kernels from the literature.
The likelihood value and the approximate errors are the quantities of interest considered to compare different correlation kernels.
\subsection{Implementation details}

The optimization of the likelihood as a function of the hyperparameters needs a performing gradient-free algorithm, in this work, we are using COBYLA~\cite{COBYLA} to maximize this quantity from the Python library Scipy with default termination criterion related to the trust region size. As COBYLA is a local search algorithm, a multi-start technique is used. Our models and their implementation are available in the toolbox SMT v2.0\footnote{\url{https://smt.readthedocs.io/en/latest/}}~\cite{SMT2019}.  By default, in SMT, the number of starting points for COBYLA is equal to 10 with evenly spaced starting points.

A simple noiseless Kriging with a constant prior model for the GP is used. We recall that the absolute exponential kernel and the squared exponential kernel are similar for categorical variables and  differ only for the continuous ones. The correlation values range between 2.06e-9 and 0.999999 for both continuous and categorical hyperparameters. Therefore, the constant $\epsilon$ is chosen to correspond to a correlation value of 2.06e-9.
The random DoEs are drawn by \gls{LHS}~\cite{LHS} and  the validation sets are given by some evenly spaced points. 

\subsection{ Analytic validation: a categorical cosine problem ($n= 1$, $m=0$, $l=1$ and $L_1=13$)}

In this section,  we consider the categorical cosine problem, from~\cite{Roustant}, to illustrate the behaviour of our proposed kernels. In this problem, the objective function $f$ depends on a continuous variable in $[0,1]$ and on a categorical variable with 13 levels.~\ref{subsec:cosine} provides a detailed description of this function. Let $w= (x,c)$ be a given point with  $x$ being the continuous variable and $c$ being the categorical variable, $c \in \{1, \ldots, 13\}$. There are two groups of curves corresponding to levels 1 to 9 and levels 10 to 13 with strong positive within-group correlations, and strong negative between-group correlations. 

\begin{figure}[H]
\begin{center}
	\subfloat[GD kernel]{
      \centering 
		\includegraphics[height=4.5cm, width=7.2cm]{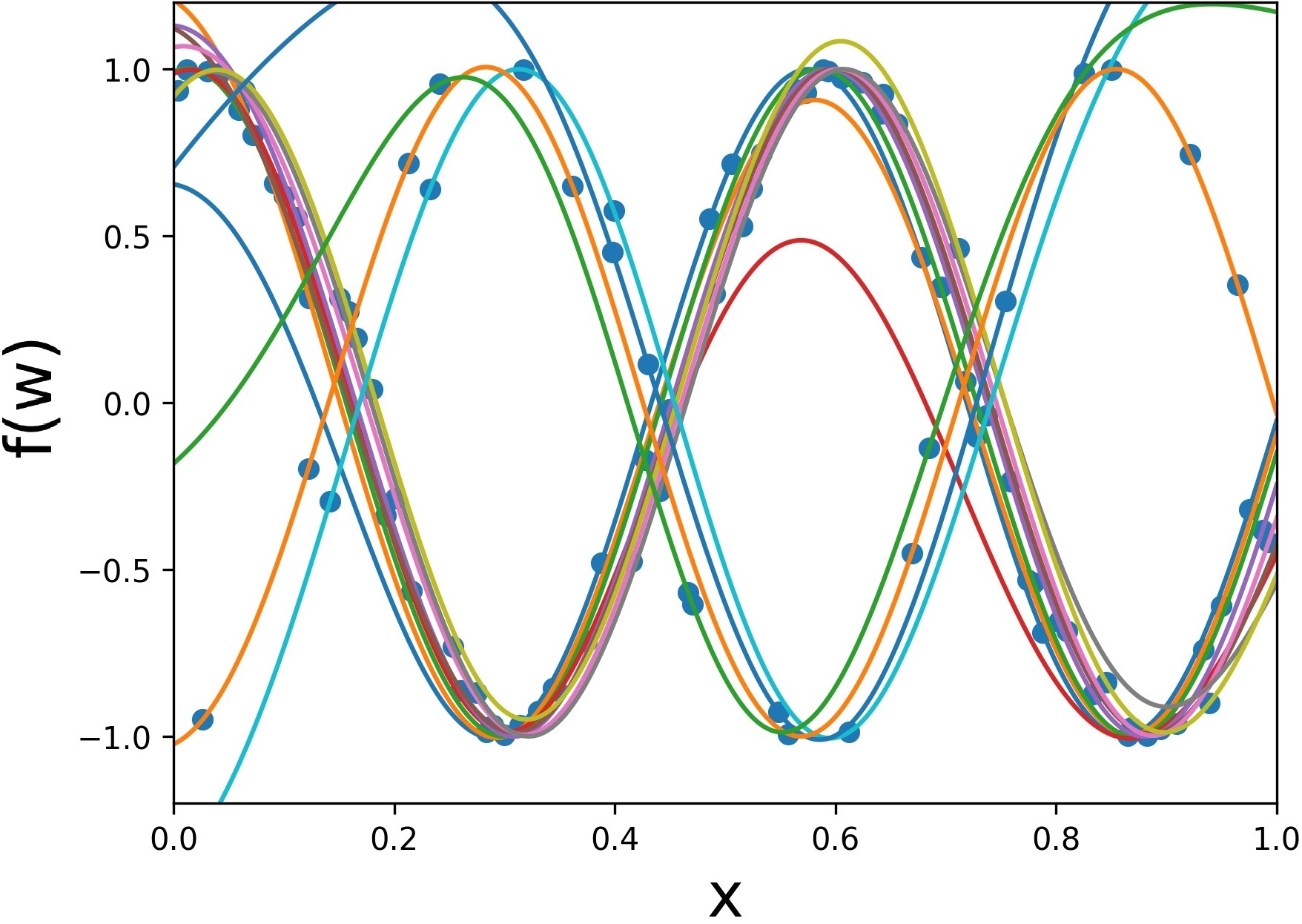}
     }   
        \subfloat[CR kernel]{
      \centering
	\includegraphics[height=4.5cm, width=8cm]{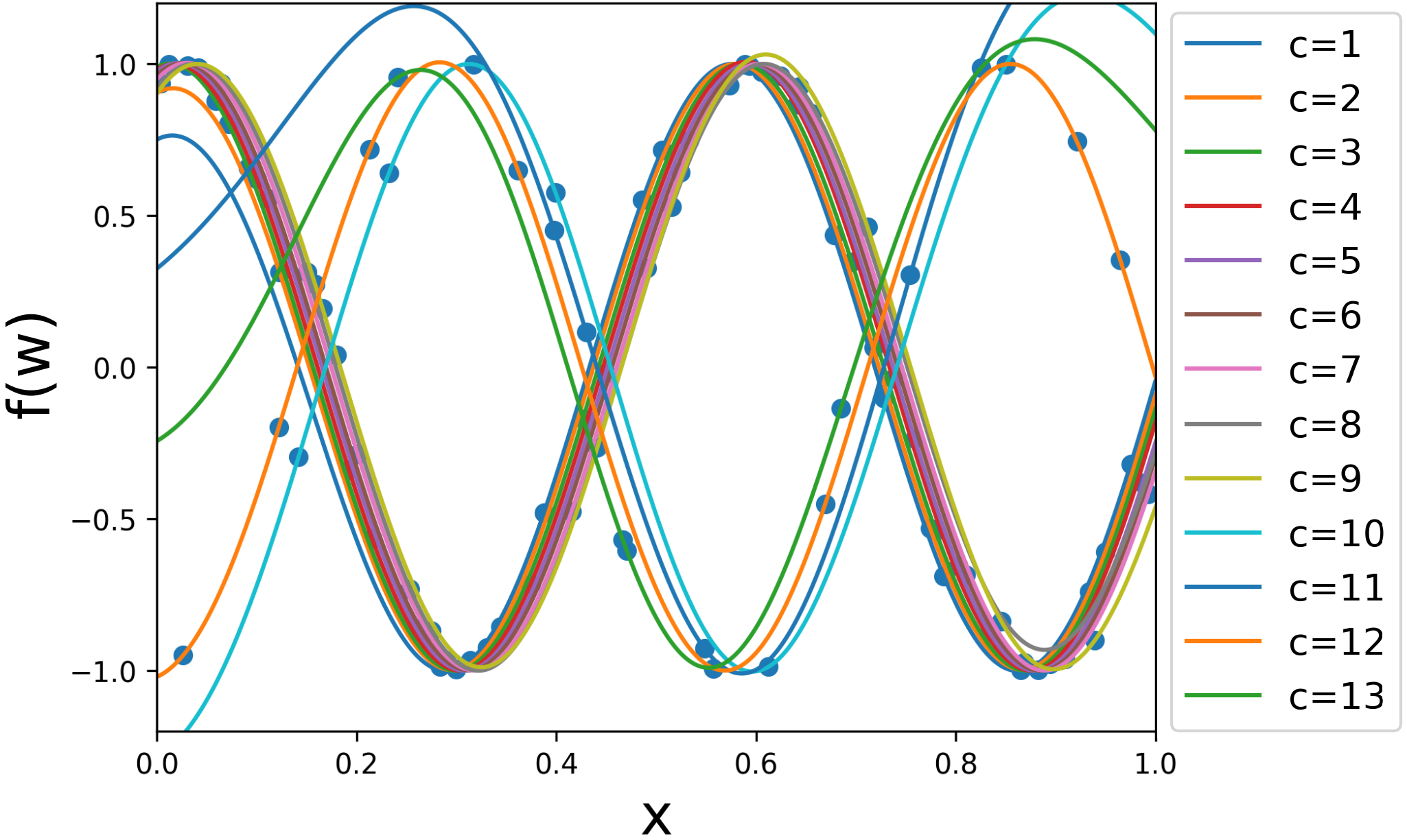}} 
 \\
        \subfloat[EHH kernel]{
      \centering
	\includegraphics[height=4.5cm, width=7.2cm]{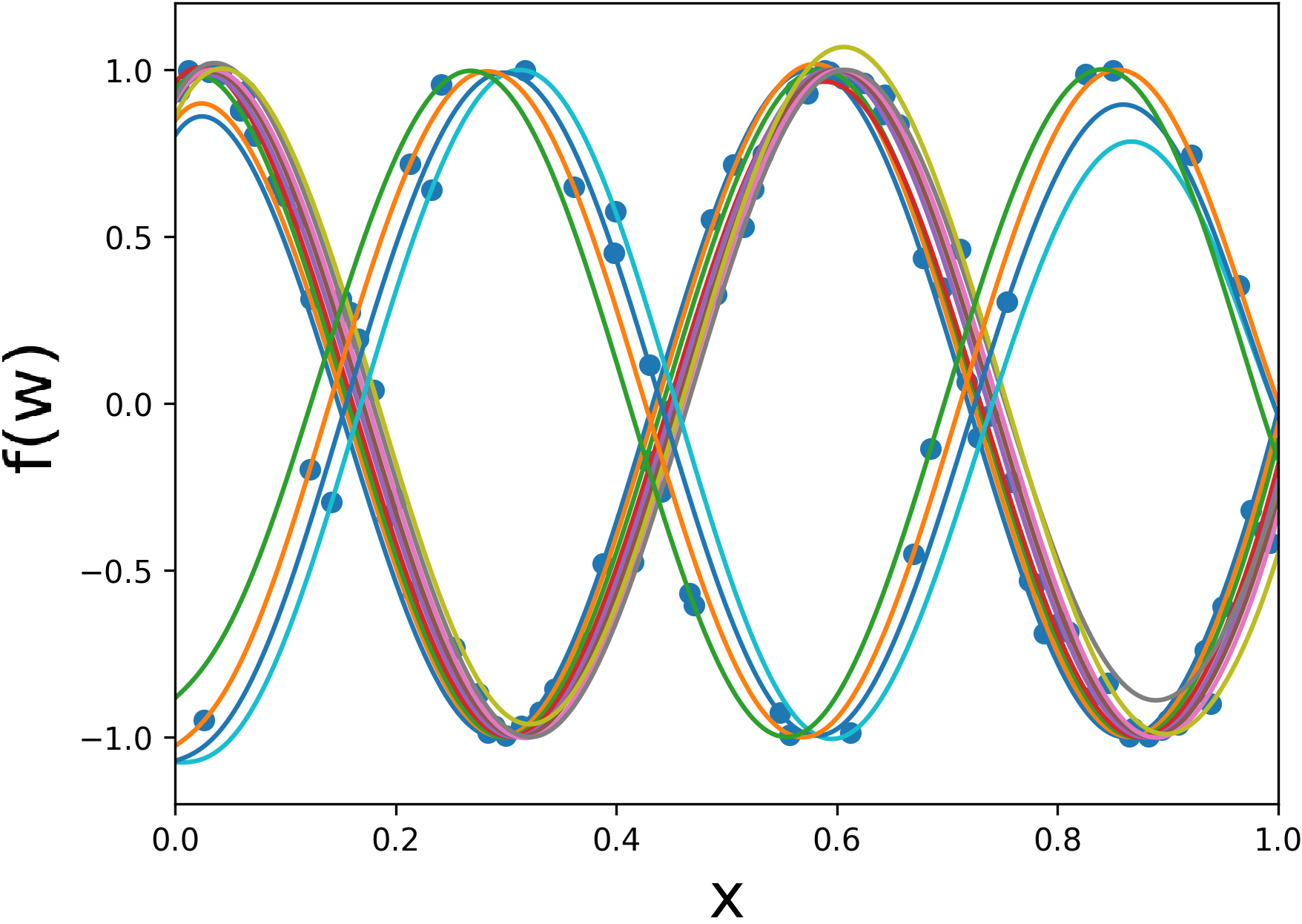}} 
 \subfloat[HH kernel]{
  \centering
	\includegraphics[ height=4.5cm, width=8cm]{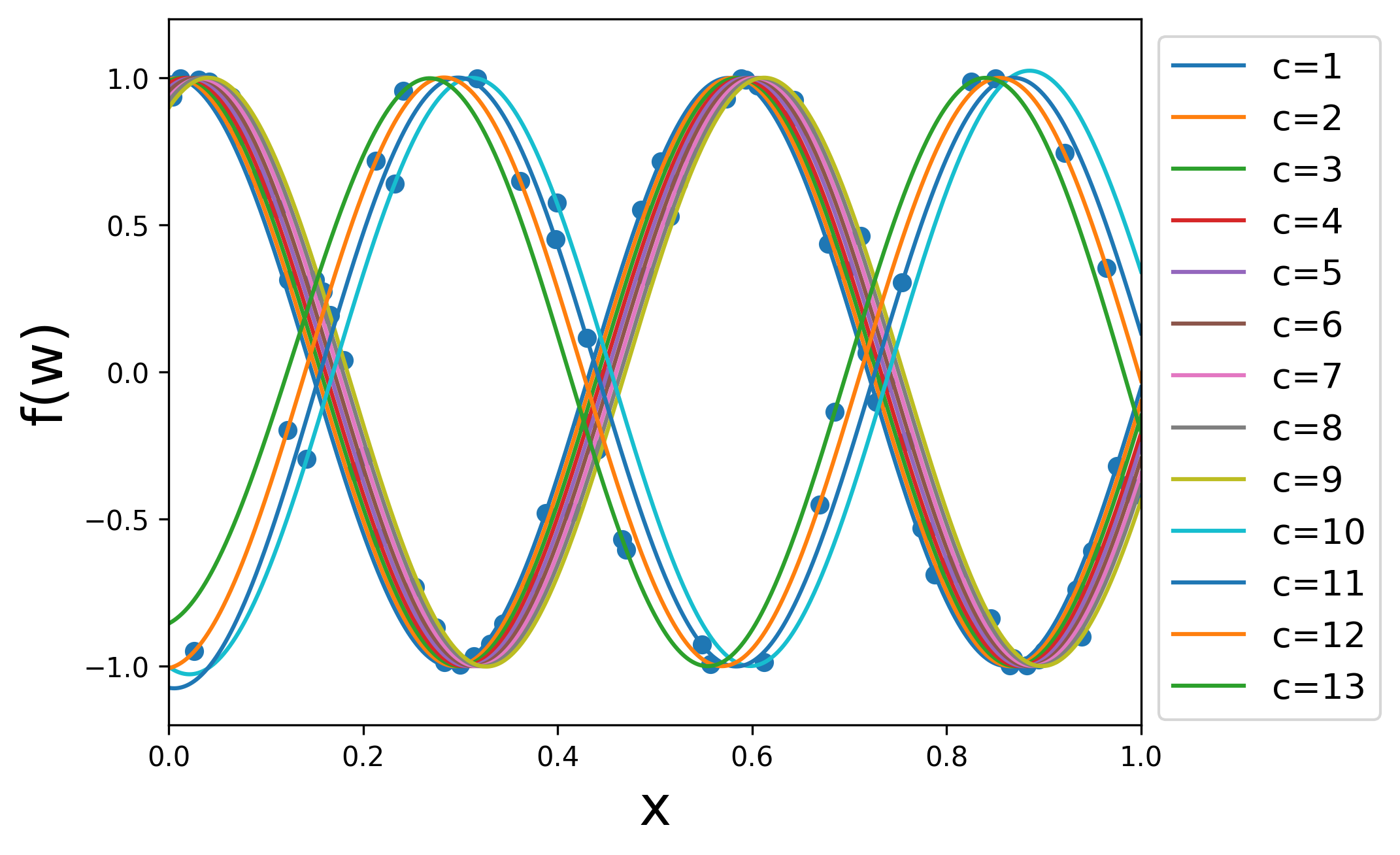} 
    \label{Roustant_hs_gsh_theo}
     }
\caption{Mean predictions for our proposed model using different kernels over the matrix $\Theta_1$ for the cosine problem with a 98 point DoE. }
\label{models_Roustant}
\end{center} 
\end{figure}

In this example, the number of relaxed dimensions for continuous relaxation is 14. A \gls{LHS} DoE with 98 points ($14\times 7$, if 7 points per dimension are considered) is chosen to built the Gaussian process models. 
The associated mean posterior models are shown on~\figref{models_Roustant} for GD, CR, EHH and HH. 
 The number of hyperparameters to optimize is therefore $2$ for GD, $14$ for CR and $79$ for EHH and HH as indicated in~\tabref{tab:resRoustant}. 
 
{\color{black}\figref{models_Roustant} shows that the predicted values remain properly within the interval $[-1,1]$ only with HH and EHH kernels. Therefore, these kernels seem to be better modeling methods. 

To better assess the accuracy of each kernel, we compute the root mean square error (RMSE) and the predictive variance adequacy (PVA)~\cite{PVA} are respectively given by 
\begin{equation*}
\mbox{RMSE}= \sqrt{\underset{i=1}{\overset{n}{\sum}} \frac{1}{n} \left( \hat{f}(w_i) - f(w_i) \right)^2} \quad \mbox{and} \quad
\mbox{PVA}= \log \left( {\underset{i=1}{\overset{n}{\sum}} \frac{1}{n} \frac{\left( \hat{f}(w_i) - f(w_i) \right)^2}{ [\sigma^f(w_i)]^2 }  } \right),
\end{equation*}
where $n$ is the size of the validation set, $\hat{f}(w_i)$ is the prediction of our GP model at point $w_i$ and $f(w_i)$ is the associated true value and the validation set consists of 13000 evenly spaced points (see~\ref{subsec:cosine}). The values, reported in~\tabref{tab:resRoustant}, show that the PVA is constant, meaning that the estimation of the variance is kept proportional to the RMSE. The RMSE decreases as the number of hyperparameters is increasing. }

\begin{figure}[H]
\begin{center}
    \subfloat[EHH kernel: 79 hyperparameters, RMSE= 1.882.]{
    \centering
	\includegraphics[  height=4.5cm, width=7.2cm]{images/HOMO_50_Roustant_curves.jpg}} 
	\subfloat[FE kernel: 92 hyperparameters, RMSE= 22.610.]{
      \centering
	\includegraphics[  height=4.5cm, width=8cm]{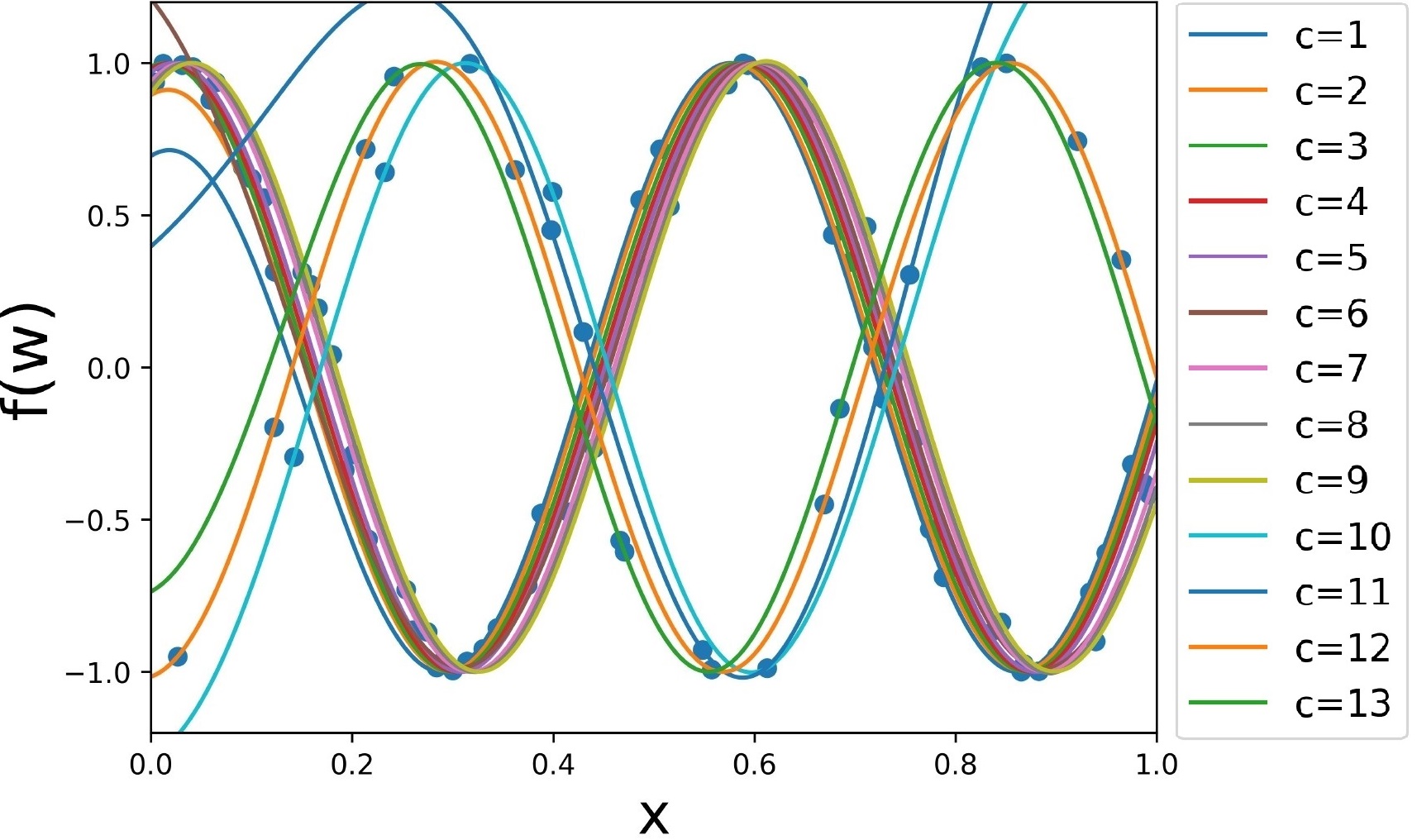}
	\label{fullroustant}} 	
\caption{Mean predictions comparison between EHH and FE kernels over the matrix $\Theta_1$ for the cosine problem with a 98 point DoE. }
\label{models_EHH}
\end{center} 
\end{figure}

\begin{table}[H]
\centering
 \caption{Kernel comparison for the cosine test case.}
\begin{tabular}{ccccc}
  \hline 
  \textbf{Kernel} &  \# of Hyperparam. & RMSE & PVA  & CPU time (s) \\
  \hline 
   \textbf{GD } &  2 & 30.079&  21.99 & 1.4 \\   
   \textbf{CR } & 14 & 22.347 & 23.04& 24.5 \\
   \textbf{EHH} & 79 & 1.882 & 23.74 &514.5 \\
   \textbf{HH}  & 79 & 1.280   & 24.31 & 514.5 \\
\hline
\end{tabular}
\label{tab:resRoustant}
\end{table}
{\color{black}
\figref{models_EHH} shows a comparison between the FE and EHH kernels. Although the two kernels are equivalent in exact precision, the EHH kernel shows  more stable and better results in term of the accuracy compared to the FE general kernel. For this reason, in what comes next, only the EHH kernel will be considered on practical use cases.}

\begin{figure}[ht]
\begin{center}

	\subfloat[GD kernel.]{
      \centering 
		\includegraphics[  height=4.5cm, width=5cm]{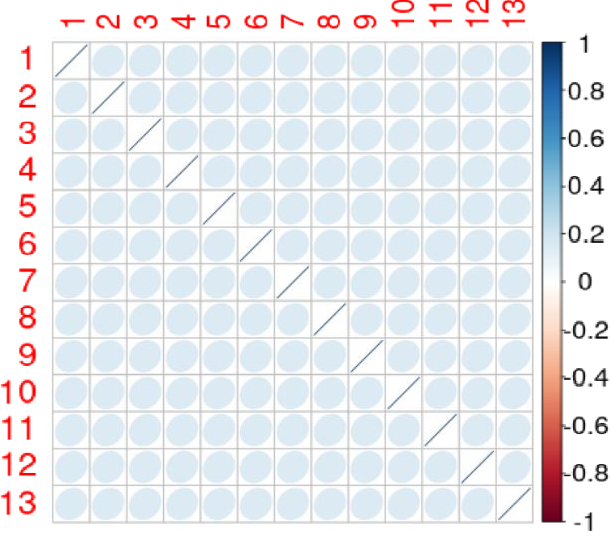} \label{corr_gower}
     }  
        \subfloat[CR kernel.]{
      \centering
		\includegraphics[  height=4.5cm, width=5cm]{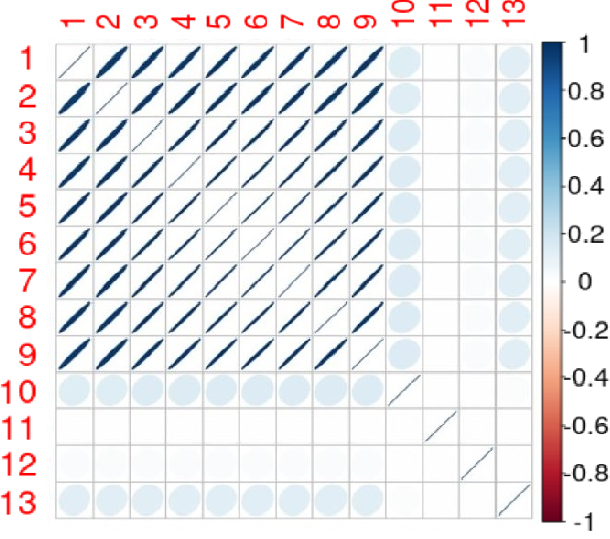} \label{corr_CR}
     }  
     
        \subfloat[EHH kernel.]{
      \centering
		\includegraphics[  height=4.5cm, width=5cm]{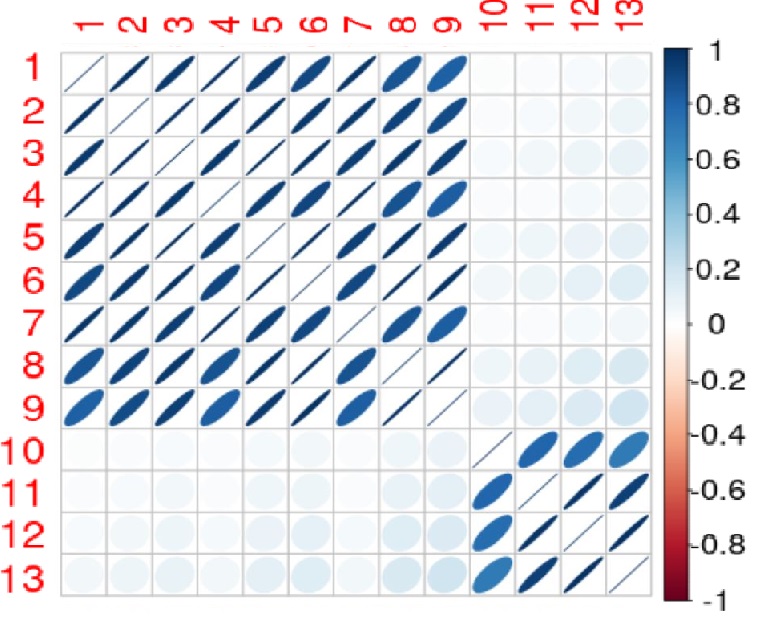}\label{corr_HS}
     }  
    \subfloat[HH kernel.]{
  \centering
	\includegraphics[  height=4.5cm, width=5cm]{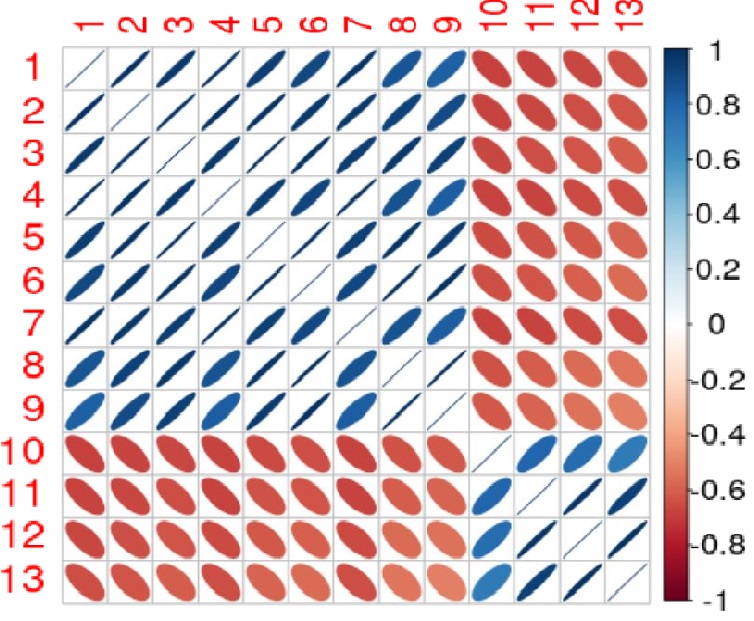}
    \label{corr_HH}
     }
\caption{Correlation matrix $R_1^{cat}$ using different choices for $\Theta_1$ on the cosine problem with a DoE of 98 points. }\label{corr_Roustant}
\end{center} 
\end{figure}

The estimated correlation matrix $R_i=R_1^{cat}$ is shown in~\figref{corr_Roustant}. For two given levels $\{\ell_r^1,\ell_s^1\}$, the correlation term $[R_1]_{\ell_r^1,\ell_s^1}$ is in blue for correlation values close to 1, in white for correlations close to 0 and in red for values close to -1; moreover the thinner the ellipse, the higher the correlation and we can see that the correlation between a level and itself is always 1. As expected, with GD kernel, there is only one estimated "mean correlation" as in~\figref{corr_gower}. For CR kernel (see~\figref{corr_CR}), the most important levels (1 to 9) are strongly correlated (in blue) with one another and the other levels (10 to 13) that should also have been correlated are badly estimated because of the kernel limitations that neglected them. 
In contrast, the EHH kernel (see~\figref{corr_HS}) gives a good approximation of the real correlations as it recovers the two groups of highly correlated levels. We recover the levels 1 to 9 as strongly similar and the levels 10 to 13 as strongly similar which is a good point but the two groups, even if less similar, are still positively correlated with one another. The latter between-group correlations should have been negative but the squared exponential kernel does not allow negative values. 
The comparison with the HH kernel, as proposed in~\cite{Zhou}, see~\figref{corr_HH}, shows that  even if the HH kernel is more general compared to EHH kernel, both kernels have an RMSE of the same order of magnitude (around 1.8 for EHH and 1.3 for HH).

\begin{figure}[H]
\centering
\subfloat[RMSE value versus DoE size.]{
\includegraphics[  height=5.2cm, width=7cm]{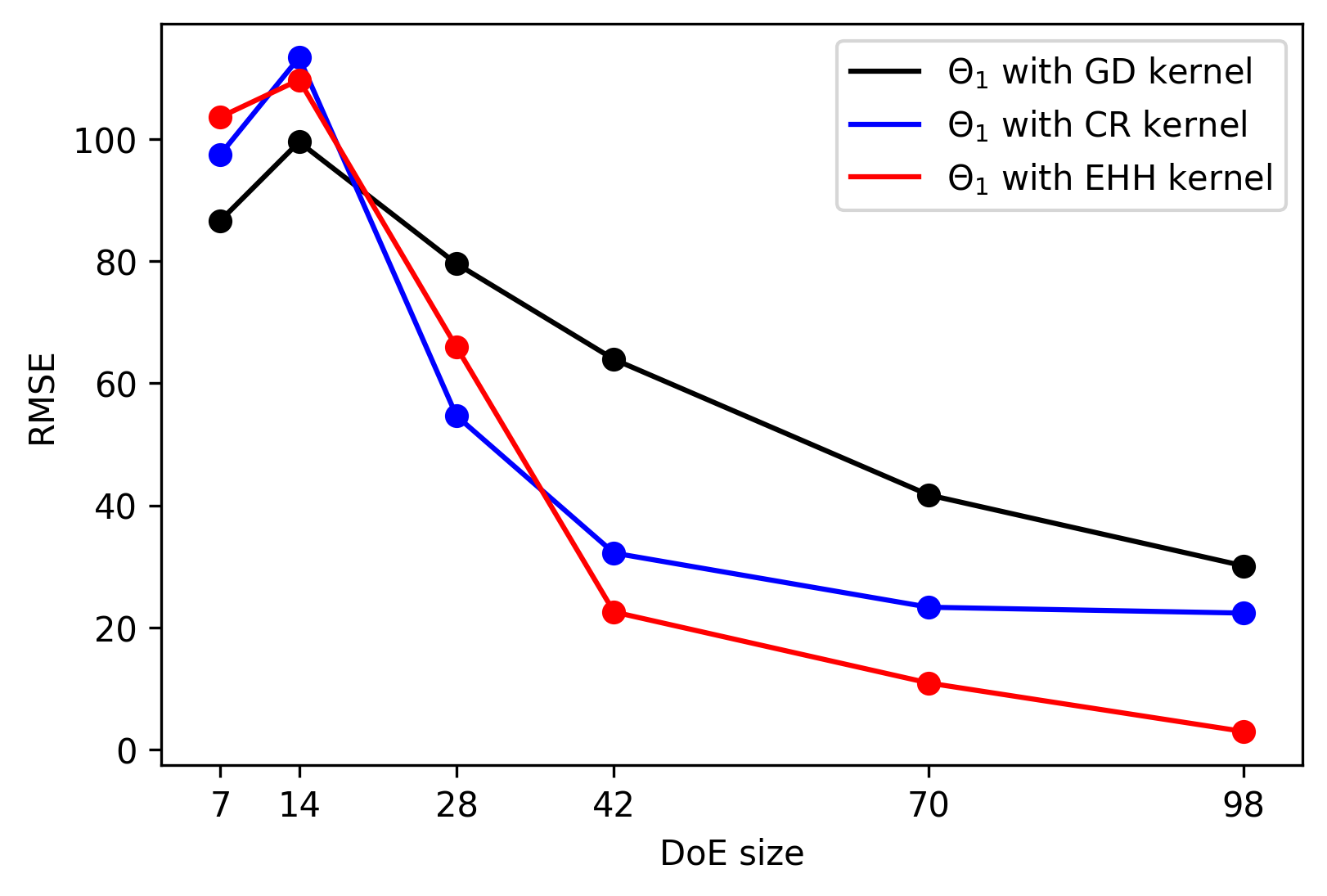}
}
\subfloat[CPU time (log scale) versus DoE size.]{
\includegraphics[  height=5.2cm, width=7cm]{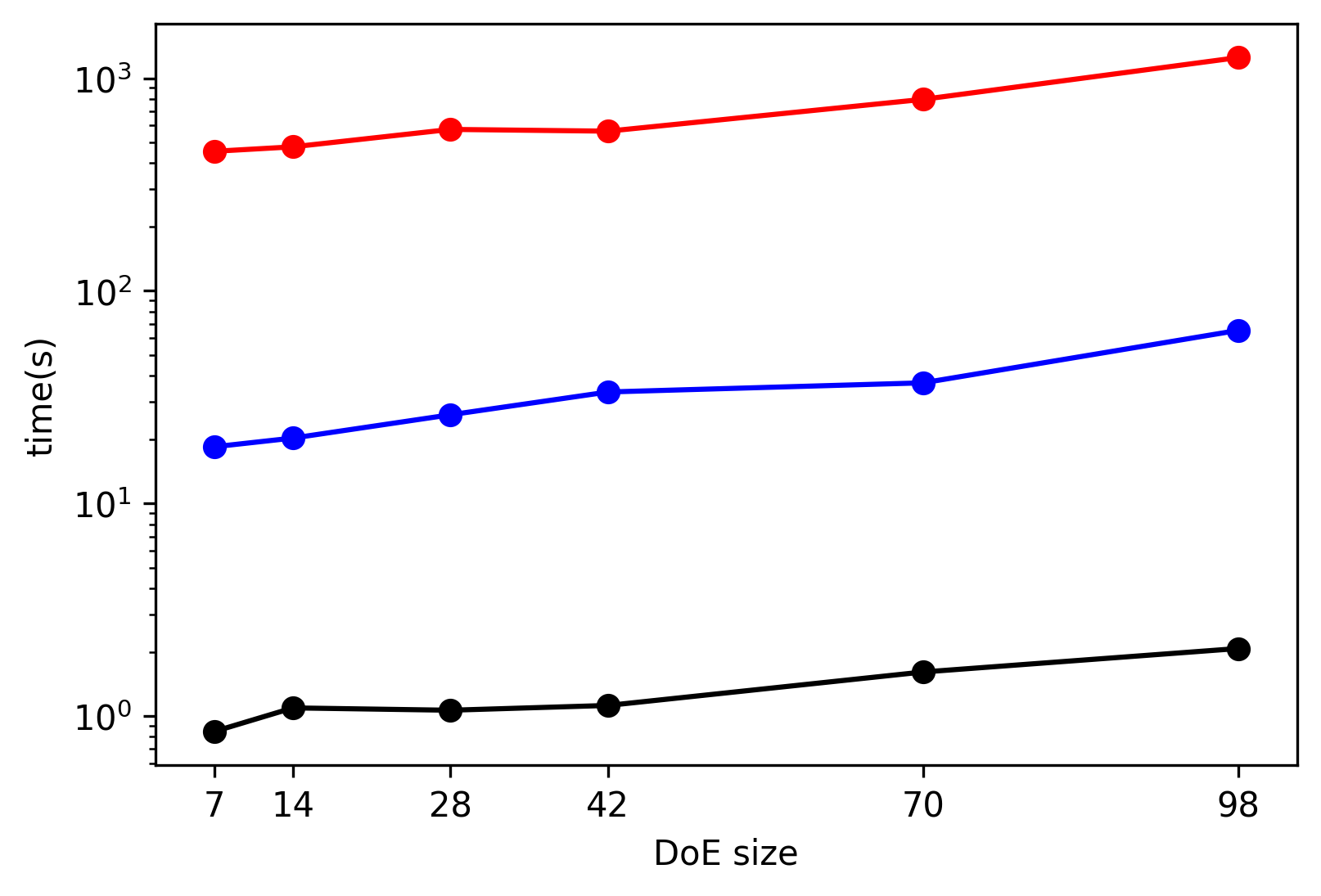}
\label{fig:time_increase}    
}
\caption{RMSE and CPU time to compute models with respect to DoE size.}
\label{fig:conv_RMSE}    
\end{figure}

On this particular test case, with a 98 point DoE, the more general the kernel, the better the performance and precision of the resulting GP. To show the DoE size impact, on~\figref{fig:conv_RMSE}, we draw 6 \gls{LHS} DoEs of different sizes and we plot the RMSE and computational time for the kernels to see how they behave for different DoE sizes.

As expected, when the size of the DoE is too small for the problem (here smaller than 15 points), the three model behaviours are similarly bad because too little information is available for the hyperparameters optimization. However, when the size of the DoE is sufficiently large, we found the same hierarchy we found with 98 points on~\figref{models_Roustant} and the more complex the model, the faster the RMSE convergence. Nevertheless, on~\figref{fig:time_increase}, we can see that the computational costs of the models scale hardly with the DoE size on a logarithmic scale.

\subsection{Application to engineering problems}

To validate and compare our method on real applications, we will consider two engineering problems of different scale to analyze the model behaviour. In Section~\ref{subsec:beam}, we present an engineering beam bending problem and in Section~\ref{subsec:aircraft}, we introduce a complex system problem from aircraft design.

\subsubsection{Cantilever beam bending problem ($n=2$, $ m=0$, $l=1$ and $L_1=12$) }
\label{subsec:beam}
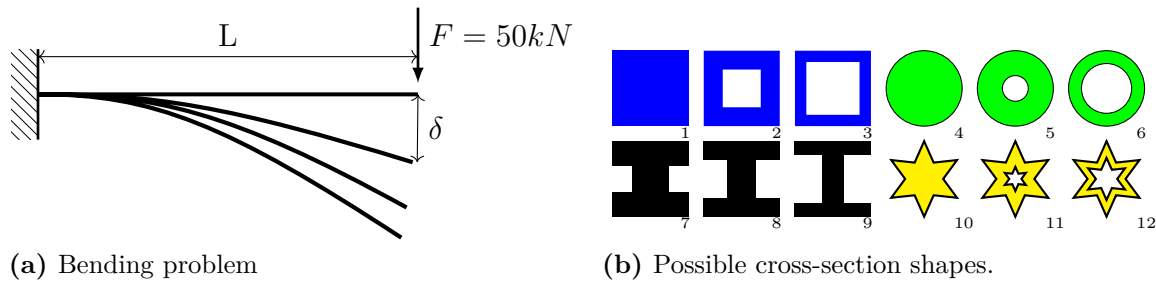
\begin{figure}[ht]
\centering

\vspace{-5pt}
\captionsetup{justification=raggedright,singlelinecheck=false}
\subfloat[Bending problem]{
\begin{tikzpicture}
    \hspace{-3pt}
    \point{origin}{-0.75}{-0.25};
    \point{begin}{0}{0};
    \point{end}{5}{0};
    \point{end_bot}{4.99}{-0.9};
    \point{end_up}{5}{0.5};
    \beam{2}{begin}{end};
    \support{3}{begin}[-90];
    \load{1}{end}[90]   ;
    \notation{1}{end_up}{$F=50kN$};

     \draw[<->] (end) -- (end_bot) node[midway, right] {$\delta$} ;
     \draw[<->] (0,0.5) -- (5,0.5) node[midway, above] {L};
     
    \draw
      [-, ultra thick] (begin) .. controls (1.5, +.01) and (2.5, -.15) .. (4.93, -0.9)
      [-, ultra thick] (begin) .. controls (1.5, +.01) and (2.5, -.2) .. (4.85, -1.5)
      [-, ultra thick] (begin) .. controls (1.5, +.01) and (2.5, -.4)   .. (4.78, -1.9);
  \end{tikzpicture}
\label{fig:beamMATHS}   
}
\subfloat[Possible cross-section shapes.]{
\centering
\begin{tikzpicture}

\hspace{-5pt}

\tstar{0.25}{0.5}{6}{0}{thick,fill=yellow,xshift=+3.6cm,yshift= -1.2cm}
\tstar{0.14}{0.28}{6}{0}{thick,fill=white,xshift=+3.6cm,yshift= -1.2cm}

\tstar{0.25}{0.5}{6}{0}{thick,fill=yellow,xshift=+2.4cm,yshift= -1.2cm}
\tstar{0.08}{0.16}{6}{0}{thick,fill=white,xshift=+2.4cm,yshift= -1.2cm}

\tstar{0.25}{0.5}{6}{0}{thick,fill=yellow,xshift=+1.2cm,yshift= -1.2cm}

\fill[green,even odd rule] (3.6,0) circle (0.5) (3.6,0) circle (0.33);
\draw (3.6,0) circle (0.5) ;
\draw (3.6,0) circle (0.33) 
; 
\fill[green,even odd rule] (2.4,0) circle (0.5)(2.4,0) circle (0.17);
\draw (2.4,0) circle (0.5) ;
\draw (2.4,0) circle (0.17) 
; 
\fill[green,even odd rule] (1.2,0) circle (0.5) ;
\draw (1.2,0) circle (0.5) ;

\def\pos{-2.4}
\fill[blue,even odd rule]  (\pos-0.5,-1.7+1.2) -- (\pos-0.5,-0.7+1.2) -- (\pos+0.5,-0.7+1.2) -- (\pos+0.5,-1.7+1.2) -- cycle ;

\def\pos{-1.2}
\fill[blue,even odd rule]  (\pos-0.5,-1.7+1.2) -- (\pos-0.5,-0.7+1.2) -- (\pos+0.5,-0.7+1.2) -- (\pos+0.5,-1.7+1.2) -- cycle   (\pos-0.25,-1.45+1.2) -- (\pos-0.25,-0.95+1.2) -- (\pos+0.25,-0.95+1.2) -- (\pos+0.25,-1.45+1.2) -- cycle ;

\def\pos{0}
\fill[blue,even odd rule]  (\pos-0.5,-1.7+1.2) -- (\pos-0.5,-0.7+1.2) -- (\pos+0.5,-0.7+1.2) -- (\pos+0.5,-1.7+1.2) -- cycle   (\pos-0.35,-1.55+1.2) -- (\pos-0.35,-0.85+1.2) -- (\pos+0.35,-0.85+1.2) -- (\pos+0.35,-1.55+1.2) -- cycle ;

\def\pos{-2.4}
\fill[black] (\pos-0.24,-0.5-1.2) -- (\pos-0.24,0.5-1.2) -- (\pos+0.24,0.5-1.2)  -- (\pos+0.24,-0.5-1.2)   -- cycle ;
\fill[black] (\pos-0.5,-0.5-1.2) -- (\pos-0.5,-0.18-1.2) -- (\pos+0.5,-0.18-1.2)  -- (\pos+0.5,-0.5-1.2)   -- cycle ; 
\fill[black] (\pos-0.5,0.18-1.2) -- (\pos-0.5,0.5-1.2) -- (\pos+0.5,0.5-1.2)  -- (\pos+0.5,0.18-1.2)   -- cycle ;  

\def\pos{-1.2}
\fill[black] (\pos-0.19,-0.5-1.2) -- (\pos-0.19,0.5-1.2) -- (\pos+0.19,0.5-1.2)  -- (\pos+0.19,-0.5-1.2)   -- cycle ;
\fill[black] (\pos-0.5,-0.5-1.2) -- (\pos-0.5,-0.25-1.2) -- (\pos+0.5,-0.25-1.2)  -- (\pos+0.5,-0.5-1.2)   -- cycle ; 
\fill[black] (\pos-0.5,0.25-1.2) -- (\pos-0.5,0.5-1.2) -- (\pos+0.5,0.5-1.2)  -- (\pos+0.5,0.25-1.2)   -- cycle ;  

\def\pos{0}
\fill[black] (\pos-0.14,-0.5-1.2) -- (\pos-0.14,0.5-1.2) -- (\pos+0.14,0.5-1.2)  -- (\pos+0.14,-0.5-1.2)   -- cycle ;
\fill[black] (\pos-0.5,-0.5-1.2) -- (\pos-0.5,-0.32-1.2) -- (\pos+0.5,-0.32-1.2)  -- (\pos+0.5,-0.5-1.2)   -- cycle ; 
\fill[black] (\pos-0.5,0.32-1.2) -- (\pos-0.5,0.5-1.2) -- (\pos+0.5,0.5-1.2)  -- (\pos+0.5,0.32-1.2)   -- cycle ;  

\point{un}{-2.15}{-0.80};
\notation{1}{un}{\tiny 1};
\point{deux}{-2.15+1.2}{-0.80};
\notation{1}{deux}{\tiny 2};
\point{trois}{-2.15+2.4}{-0.80};
\notation{1}{trois}{\tiny 3};

\point{quatre}{-2.15+3.6}{-0.80};
\notation{1}{quatre}{\tiny 4};
\point{cinq}{-2.15+4.8}{-0.80};
\notation{1}{cinq}{\tiny 5};
\point{six}{-2.15+6}{-0.80};
\notation{1}{six}{\tiny 6};

\point{sept}{-2.15}{-0.80-1.2};
\notation{1}{sept}{\tiny 7};
\point{huit}{-2.15+1.2}{-0.80-1.2};
\notation{1}{huit}{\tiny 8};
\point{neuf}{-2.15+2.4}{-0.80-1.2};
\notation{1}{neuf}{\tiny 9};

\point{dix}{-2.15+3.6}{-0.80-1.2};
\notation{1}{dix}{\tiny 10};
\point{onze}{-2.15+4.8}{-0.80-1.2};
\notation{1}{onze}{\tiny 11};
\point{douze}{-2.15+6}{-0.80-1.2};
\notation{1}{douze}{\tiny 12};

\end{tikzpicture}
\label{fig:beam_shapeMaths}    
}

\captionsetup{justification=centering,singlelinecheck=false}
\caption{Cantilever beam problem.}
\end{figure}

A first engineering problem commonly used for model validation is the beam bending problem in its linear elasticity range~\cite{Roustant, Cheng2015TrustRB}.
This problem is illustrated on~\figref{fig:beamMATHS} and consists of a cantilever beam loaded at its free extremity with a force $F$. As in Cheng $\textit{et al.}$~\cite{Cheng2015TrustRB}, we choose a constant Young modulus of $E=200$GPa and a load of $F=50$kN. Moreover, as in Roustant $\textit{et al.}$~\cite{Roustant}, we consider 12 possible cross-sections: there are 4 possible shapes, illustrated in~\figref{fig:beam_shapeMaths} that could be hollow, thick or full. For a given cross-section (shape and thickness), its size is determined by its surface $S$. Every cross section is associated with a normalized moment of inertia $\tilde{I}$ about the neutral axis. The latter is a latent variable associated to the beam shape~\cite{oune2021latent}.

Therefore, the problem to model has two continuous variables: the length $L \in [10,20]$ (in $m$) and the surface  $S \in [1,2]$ (in $m^2$) and one categorical variable $ \tilde{I}$ with 12 levels. The tip deflection, at the free end, $\delta$ is given by $$ \delta = f( \tilde{I}, L,S) = \frac{F}{3E} \frac{L^3}{S^2\tilde{I}} $$



    
To compare our models, we draw a 98 point LHS as training set and the validation set is a grid of $12\times30\times30=10800$ points. For both squared exponential and absolute exponential kernels, the RMSE, likelihood and computational time for every model are shown in~\tabref{tab:resCantilever}. We recall that squared exponential and absolute exponential kernels differ only on the continuous variables and are the same for the categorical part.
As expected, the computational time and the likelihood increase when the model is more complex. The DoE seems of sufficient size for this problem as the computed RMSE (\textit{i.e.}, the total displacement error) decreases with the model complexity.

\begin{table}[H]
\centering
 \caption{Results of the cantilever beam models.}
\resizebox{0.99\columnwidth}{!}{%
\small
\begin{tabular}{ccccc}
\hline
  \textbf{Categorical kernel} & Continuous kernel & Displacement error (cm) & $ \ $ Likelihood  &$\ $ Time (s) \\
  \hline 
  \textbf{GD} & squared exponential &1.3858 & 111.13&  8.02 \\   
  \textbf{CR} & squared exponential & 1.1604 & 162.26 & 89.1 \\
  \textbf{EHH} & squared exponential &0.1247 &
256.90 & 2769.4 \\
  \hdashline 
  \textbf{GD} & absolute exponential & 3.2403 & 74.48   & 14.71  \\
  \textbf{CR} & absolute exponential & 3.0918 & 99.00 & 260.1 \\
  \textbf{EHH} & absolute exponential & 2.0951& 102.48 & 19784\\
\hline
\end{tabular}
}
\label{tab:resCantilever}
\end{table}

\begin{figure}[H]
\begin{center}
	\subfloat[With GD kernel.]{
      \centering 
\includegraphics[  height=4.5cm, width=4.3cm]{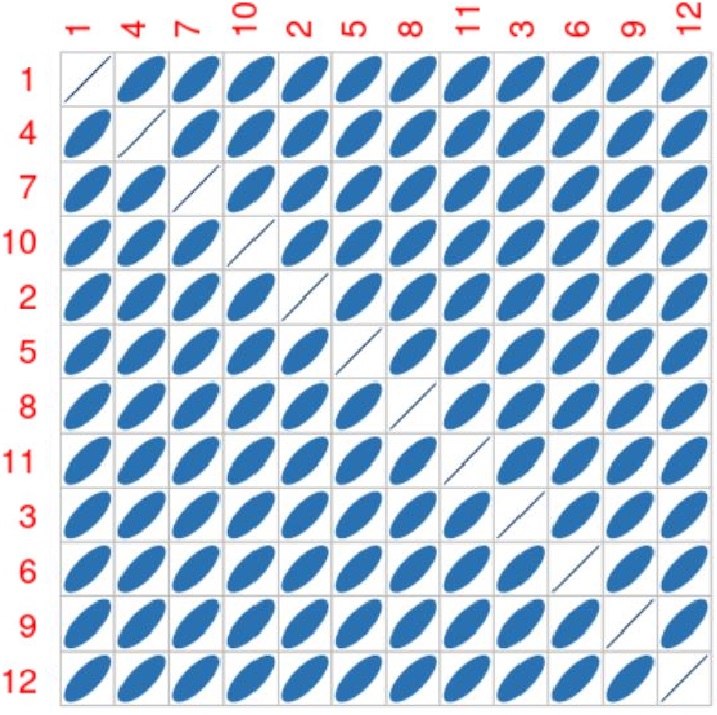} \label{corr_canti_gower}
     }  
	\subfloat[ With CR kernel.]{
      \centering
\includegraphics[  height=4.5cm, width=4.3cm]{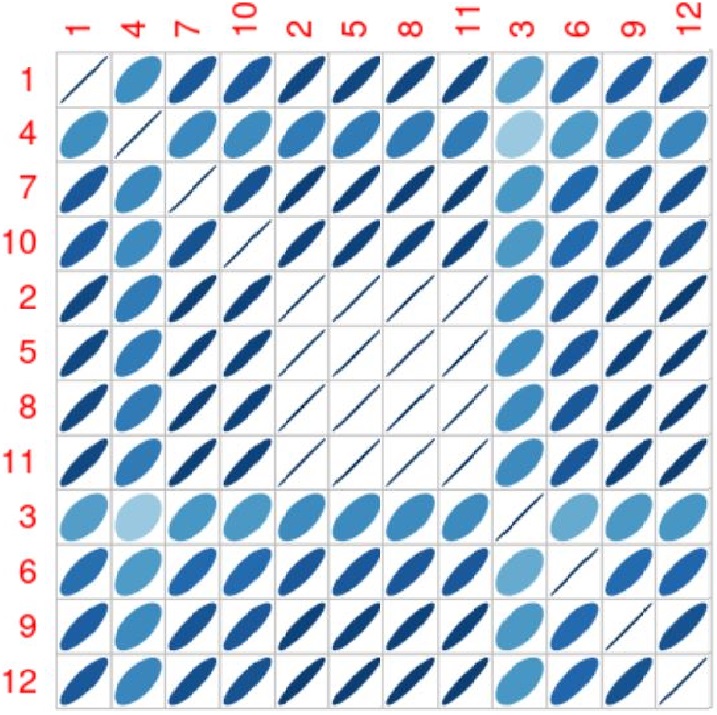} \label{corr_canti_cr}
     }  
	\subfloat[ With EHH kernel. ]{
      \centering 
\includegraphics[  height=4.5cm, width=4.3cm]{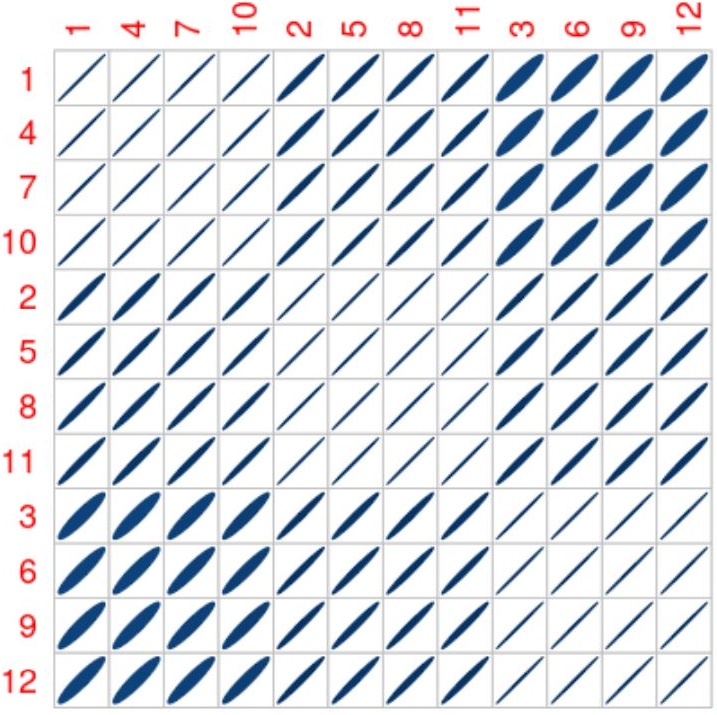} \label{corr_canti_ehh}
     }  
\includegraphics[  height=4.9cm, width=0.5cm]{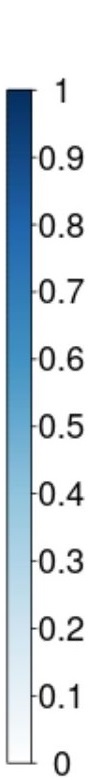}
\centering     
\caption{Correlation matrix $R_1^{cat}$  using different choices for $\Theta_1$ for the categorical variable $\tilde{I}$ from the cantilever beam problem.}
\label{corr_Cantilever}
\end{center} 
\end{figure}

In~\figref{corr_Cantilever}, we have drawn the correlation matrix found between the cross-section shape (the resulting $R_1$ correlation matrix) for the three models. On the figure below, the higher the correlation, the thinner the ellipse.

As expected, we have 3 groups of 4 shapes depending on their respective thickness (respectively, the levels \{1,4,7,10\} the levels \{2,5,8,11\} and the levels \{3,6,9,12\}). The more the thickness is similar, the higher the correlation: the thickness has more impact than the shape of the cross-section on the tip deflection. However, given the database, two points with similar $L$ and $S$ values will have similar output whatever the cross-section. The effect of the cross-section on the output is always the same (in the form of $\frac{1}{\tilde{I}}$) leading to an high correlation after maximizing the likelihood. In~\figref{corr_canti_ehh}, with the EHH kernel, we can distinguish the 3 groups of 4 shapes and, because the correlations are close to 1, the homoscedastic hyperphere model~\cite{Pelamatti} would lead to the same correlation matrix.  Also, with the CR kernel of~\figref{corr_canti_cr}, the medium thick group \{2,5,8,11\} being correlated with both the full and the hollow group, its correlation values are the higher whereas the correlation hyperparameters associated to the two other groups are smaller. 
For the GD model in~\figref{corr_canti_gower}, there is only one mean positive correlation value as before.

\subsubsection{Aircraft design application ($n=10$, $ m=0$, $l=2$ and $L_1=9$, $L_2=2$) }
\label{subsec:aircraft}

The \gls{DRAGON} aircraft concept 
has been introduced by ONERA in 2019~\cite{schmollgruber} within the scope of the European CleanSky 2 program\footnote{\href{https://www.cleansky.eu/technology-evaluator}{\color{blue}https://www.cleansky.eu/technology-evaluator}} which sets the objective of 30\% reduction of CO2 emissions by 2035 with respect to 2014 state-of-the-art.
The employment of a distributed propulsion comes at a certain cost; a turboelectric propulsive chain is necessary to power the electric fans which brings additional complexity and weight.
The turboelectric propulsive chain being an important weight penalty, it is of particular interest to optimize the chain and particularly the number and type of each component, characterized by some discrete values.  The definition of the architecture variable is given in~\tabref{tab:dragon_archi1} and the definition of the turboshaft layout is given in~\tabref{tab:dragon_archi2}. For the sake of simplicity, we restrict the optimization problem to the case of two electric cores and generators but more optimizations have been performed in~\cite{SciTech_cat}.

%
%

\begin{table}[H]
\caption{Categorical variable definition for the "\texttt{DRAGON}" test case.}
\centering
\vspace*{-0.1cm}

 \subfloat[Definition of the architecture variable and its 9 associated levels.]{
\small

\resizebox{0.8\columnwidth}{!}{%
\small

\begin{tabular}{cccc}
  \hline
  \textbf{Architecture number} $\ $ & Number of motors $\ $ & Number of cores $\ $ & Number of generators $\ $ \\
  \hline
  \textbf{1} & 8 &2 & 2\\
  \textbf{2} & 12 & 2 & 2\\
  \textbf{3} & 16 & 2 & 2\\
  \textbf{4} &20 &2 & 2\\
  \textbf{5} & 24 & 2 & 2\\
  \textbf{6} & 28 & 2 & 2\\
  \textbf{7} &32 & 2 & 2\\
  \textbf{8} & 36  & 2 & 2\\
  \textbf{9} & 40 & 2 & 2\\

\hline
\end{tabular}
}
\label{tab:dragon_archi1}
}
\centering
\vspace*{+0.1cm}

\subfloat[Definition of the turboshaft layout variable and its 2 associated levels.]{
\small

\resizebox{0.8\columnwidth}{!}{%
\small

\begin{tabular}{cccccc}
  \hline 
  \textbf{Layout} & Position & y ratio & Tail & VT aspect ratio & VT taper ratio\\
  \hline 
  \textbf{1} & under wing &0.25 & without T-tail& 1.8 & 0.3 \\
  \textbf{2} & behind & 0.34 & with T-tail& 1.2 & 0.85\\
 
\hline
\end{tabular}
}
\label{tab:dragon_archi2}
}
\end{table}
The analysis of \gls{DRAGON} is treated with Overall Aircraft Design method in \gls{FAST-OAD}~\cite{David_2021}.  We are considering the following problem described in~\tabref{tab:dragon_neucom}.
\begin{table}[H]
\centering
\vspace*{-0.3cm}

 \caption{Definition of the ``\texttt{DRAGON}'' optimization problem.}
\small

\resizebox{1.0\columnwidth}{!}{%
\small

\begin{tabular}{lllrr}
 & Function/variable & Nature & Quantity & Range\\
\hline
\hline
Model & Fuel mass & cont & 1 &\\
\hline
with respect to & \mbox{Fan operating pressure ratio} & cont & 1 & $\left[1.05, 1.3\right]$ \\  
     & \mbox{Wing aspect ratio} & cont & 1 &    $\left[8, 12\right]$ \\
    & \mbox{Angle for swept wing} & cont & 1 & $\left[15, 40\right]$  ($^\circ$) \\
     & \mbox{Wing taper ratio} & cont & 1 &    $\left[0.2, 0.5\right]$ \\
     & \mbox{HT aspect ratio} & cont & 1 &    $\left[3, 6\right]$ \\
    & \mbox{Angle for swept HT} & cont & 1 & $\left[20, 40\right]$  ($^\circ$) \\
     & \mbox{HT taper ratio} & cont & 1 &    $\left[0.3, 0.5\right]$ \\
 & \mbox{TOFL for sizing}  & cont &1 & $\left[1800, 2500\right]$ ($m$) \\
 & \mbox{Top of climb vertical speed for sizing $ \ $} & cont & 1 & $\left[300, 800\right]$ ($ft/min$) \\
 & \mbox{Start of climb slope angle} & cont & 1 & $\left[0.075, 0.15\right]$ ($rad$) \\

 & \multicolumn{2}{l}{Total  continuous variables} & 10 & \\
 \cline{2-5}
& \mbox{Architecture} & cat & 9 levels & \{1,2,3, \ldots,7,8,9\} \\
& \mbox{Turboshaft layout} & cat & 2 levels & \{1,2\} \\

 & \multicolumn{2}{l}{Total categorical variables} & 2 & \\
 \cline{2-5}

  &   \multicolumn{2}{l}{\textbf{Total relaxed variables}} & {\textbf{21}} & \\
  \hline

\end{tabular}
}
\label{tab:dragon_neucom}
\end{table}

\begin{figure}[H]
\begin{center}
	\subfloat[GD kernel.]{
      \centering 
\includegraphics[  height=4.4cm, width=4.3cm]{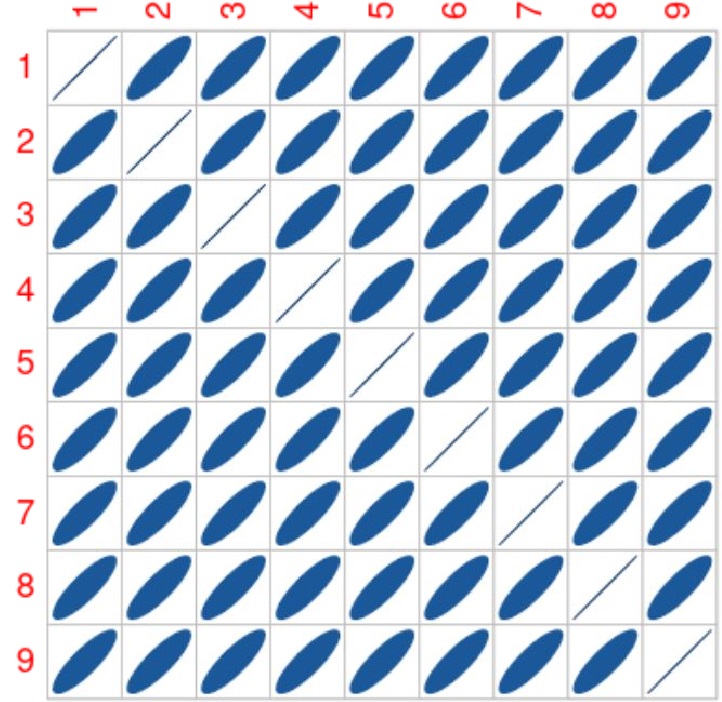} \label{corr_turboelectric_gower}
     }  
	\subfloat[CR kernel.]{
      \centering 
\includegraphics[  height=4.4cm, width=4.3cm]{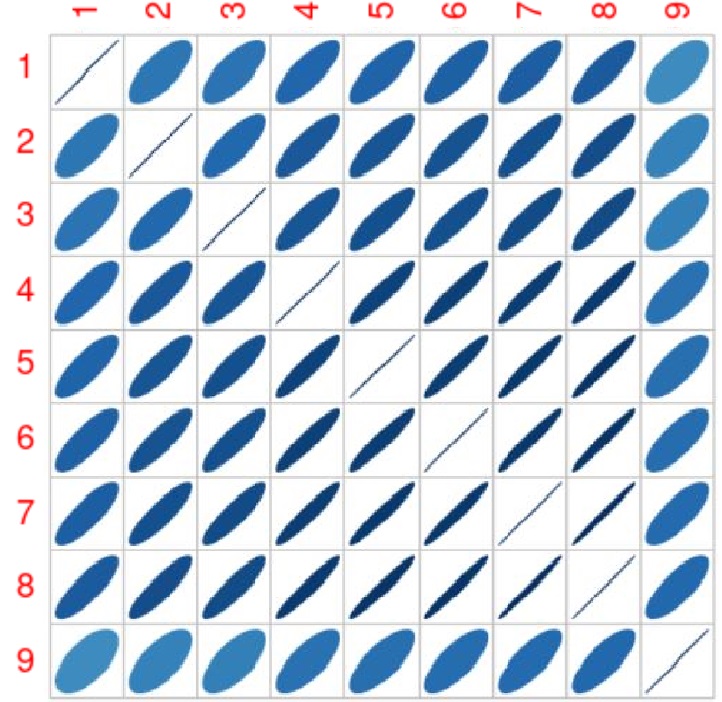} \label{corr_turboelectric_cr}
     }  
	\subfloat[EHH kernel.]{
      \centering 
\includegraphics[   height=4.4cm, width=4.3cm]{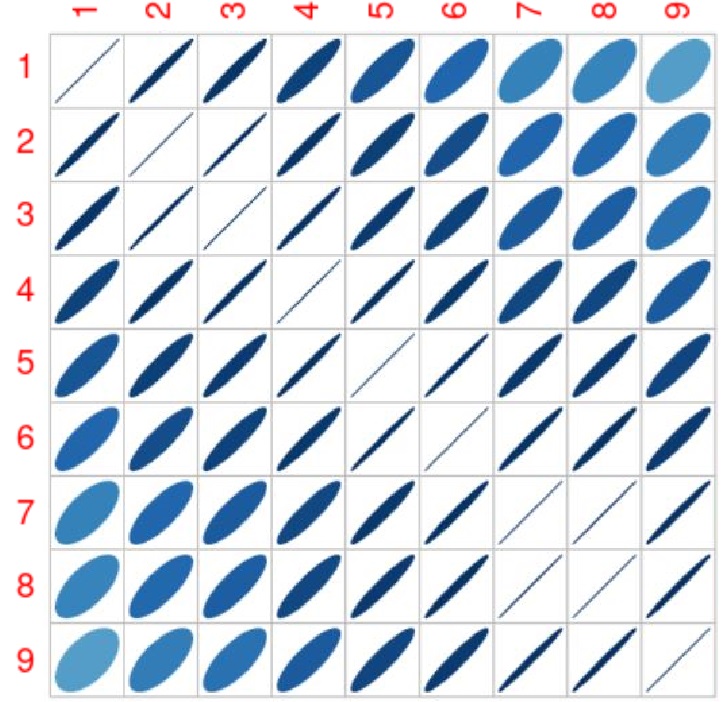}
\label{corr_turboelectric_ehh}
     }  
\includegraphics[  height=4.7cm, width=0.5cm]{images/leg.jpg}
\centering     
\caption{Correlation matrix $R_1^{cat}$  using different choices for $\Theta_1$  for the turboelectric architecture variable.}
\label{corr_turboelectric}
\end{center} 
\end{figure}

\begin{figure}[H]
\begin{center}
	\subfloat[GD kernel.]{
      \centering 
\includegraphics[  height=4.4cm, width=4.3cm]{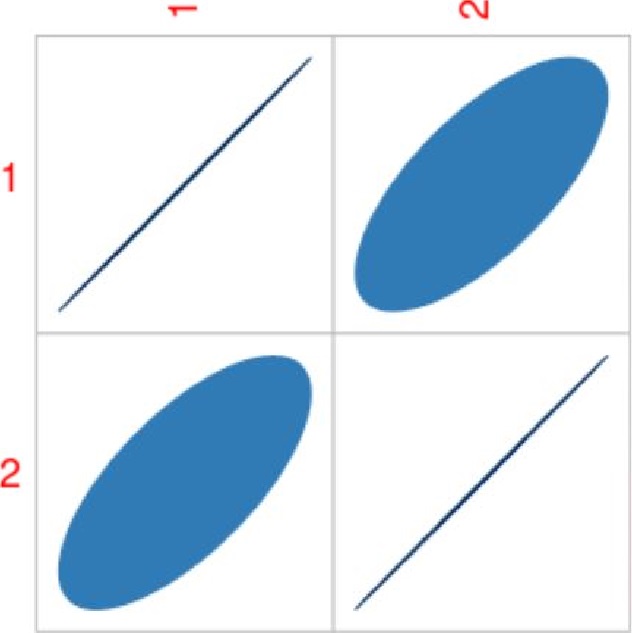} \label{corr_turboshaft_gower}
     }  
	\subfloat[CR kernel.]{
      \centering 
\includegraphics[  height=4.4cm, width=4.3cm]{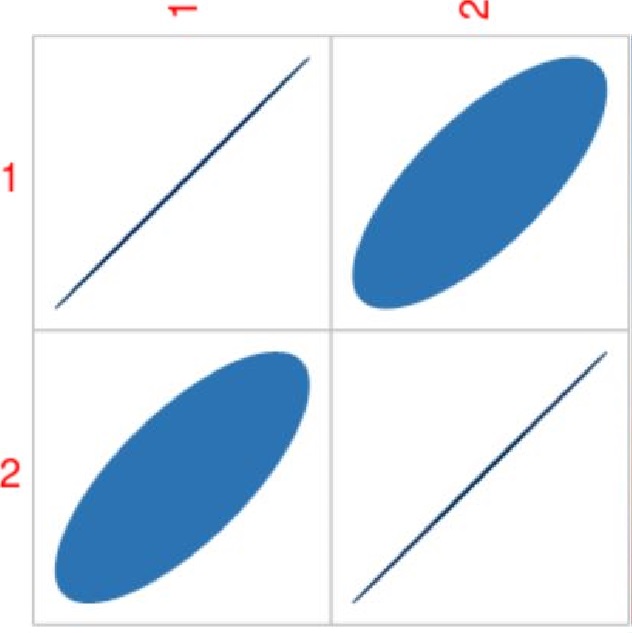} \label{corr_turboshaft_cr}
     }  
	\subfloat[EHH kernel.]{
      \centering 
\includegraphics[  height=4.4cm, width=4.3cm]{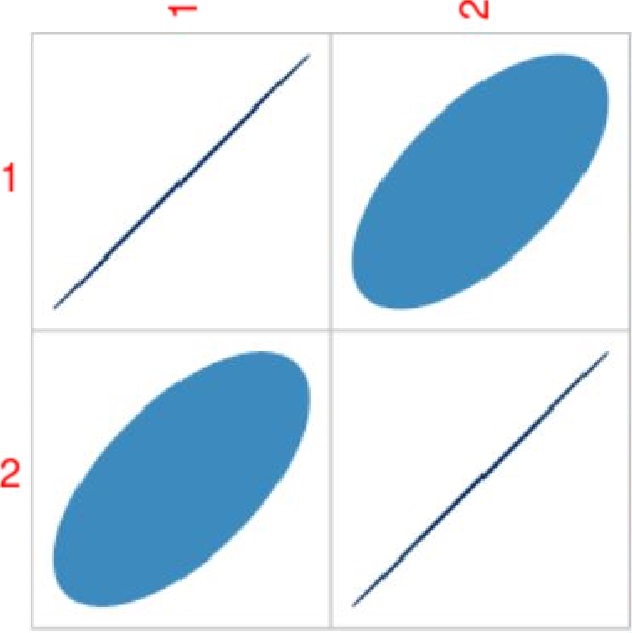} \label{corr_turboshaft_ehh}
     }  
\centering
\includegraphics[  height=4.65cm, width=0.5cm]{images/leg.jpg}
\caption{Correlation matrix $R_2^{cat}$  using different choices for $\Theta_2$  for the turboshaft layout variable.}
\label{corr_turboshaft}
\end{center} 
\end{figure}

Twice, we draw  $250$ points by LHS. Over the first DoE, that is the training set, we build the model to predict the fuel mass and over the second one, we validate our prediction and compute the RMSE reported in~\tabref{tab:resDragon}.
In this case, the number of hyperparameters is 12 for GD kernel, 21 for CR kernel and 47 for EHH kernel. Evaluating the function is costly, around 4 minutes for a single point. We observed similar performances for all models, the performance is mostly determined by the choice of the continuous kernel. 
For a problem that has that many variables, it seems useless and impractical to use a complicated model, the GD kernel being already performing well. 
On~\figref{corr_turboelectric}, we plot, for the three kernels, the approximate correlation matrices for the first categorical variable. As we can see, when considering the general EHH kernel, as in~\figref{corr_turboelectric_ehh}, the closer the levels, the higher the correlation. In fact, in this case, the only difference between two levels is the number of motors. Therefore, the more similar the number of motors, the more similar the fuel consumption. Given that, we expect, when considering CR kernel as in~\figref{corr_turboelectric_cr} that the higher correlation should appear "in the middle" \{4,5,6\} as these levels are meant to be the most correlated with the others. This is what happens to a certain extent but the levels 7 and 8 are weirdly appearing too much correlated with one another. This could be a numerical problem, the optimization being hard with that many variables and hyperparameters. As before, the GD kernel is the less precise and just give a mean correlation over the whole space as in~\figref{corr_turboelectric_gower}. In~\figref{corr_turboshaft}, we plot, for the three methods, the approximated correlation matrices for the second categorical variable. There is only two engine layouts so there is only one correlation. In this case, the correlation is positive indicating that the plane behave in the same way no matter the layout.

\begin{table}[H]    
\centering
 \caption{Results of the aircraft models based on a 250 point validation set.}
\small
\resizebox{0.9\columnwidth}{!}{%
\small

\begin{tabular}{ccccc}
\hline
  \textbf{Kernel}& $ \ $ number of hyperparameters  & $\ $  kernel   & $\ $fuel error (kg) $\ $  & time (s)   \\
  \hline 
  \textbf{GD} & 12 & squared exponential & 2115 & 65 \\
  \textbf{CR} & 21 & squared exponential & 2068  & 210  \\
  \textbf{EHH} & 47 & squared exponential & 2147 &  9450\\
   \hdashline 
  \textbf{GD} & 12 & absolute exponential &1666 & 65 \\
  \textbf{CR} & 21 & absolute exponential & 1664  & 210 \\
  \textbf{EHH} & 47 & absolute exponential & 1593 &  9295 \\
\hline
\end{tabular}
}
\label{tab:resDragon}
\end{table}

One can note that increasing the number of motors or changing a layout will not change the way an aircraft flies. For example, having more motors will only increase the fuel consumption by a given factor. The latter will always remain positive and related to the continuous variables. Hence, in this test case, we do not have opposite effects between two categorical levels.

In most industrial applications, radically opposite effects over a complex system do not occur so often. For instance, on the industrial applications that can be found on the literature, there was not a clear need for negative correlation values~\cite{Pelamatti,Roustant, cuesta2021comparison}. Therefore, in practice, the exponential model is not that limiting compared to the homoscedastic hypersphere model.

\section{Conclusion}
\label{sec:conclu}

In this work, we have proposed a class of kernels for GP models that extends the exponential continuous kernels to the mixed-categorical setting. We showed that this class of kernels generalizes Gower distance and continuous relaxation based kernels. A classification between the proposed kernels as well as a proof of the SPD nature of the resulting correlation matrices have been also proposed. Numerical illustrations on analytical toy problems showed the good potential of the proposed kernels to reduce the number of hyperparameters and thus the computational time. The implementation of our proposed method has been released in the toolbox SMT v2.0\footnote{\url{https://smt.readthedocs.io/en/latest/}}. 

When considering complex kernels, a good approach would be to use a model reduction technique such as Kriging with Partial Least Squares (KPLS)~\cite{Bouhlel18} that is derived from the construction of the correlation matrix via a kernel function. KPLS is an adaptation of the Partial Least Squares regression for exponential kernels and is used to reduce the number of hyperparameters and handle a large number of mixed inputs. Further works will consider to include such dimension reduction techniques to improve the computational efficiency of our model  and tackle higher dimensional problems.

\recap{ 
\lettrine[lines=2, lhang=0.33, loversize=0.25, findent=1.5em]{I}{n} this chapter, we proposed a new mixed-categorical metamodel based on GP  by developing a correlation kernel that can handle mixed-categorical variables to use within GP. These chapter achievements are the following. 
\begin{itemize}
    \item A new kernel combining the continuous distance-based exponential kernel with the categorical HH kernel has been developed and was named EHH.
    \item This kernel has been proven to be SPD.
    \item This kernel generalizes both continuous exponential kernels and categorical Gower distance and continuous relaxation based kernels.
    \item This chapter established a new general framework that can be declined into every other state-of-the-art kernels either for continuous or categorical kernels.
    \item  This framework leads to a classification of the kernels directly related to their respective number of hyperparameters. More precisely, we proved that the HH kernel generalizes EHH kernel which in turns generalizes the CR kernel that generalizes GD kernel. 
    \item Every kernel was implemented into the open-source software SMT.
    \item The various models have been compared and tested on analytic and engineering test cases. These test cases have shown results confirming the theory. Moreover, as we could have expected, the more complex the model, the better the performance. Notwithstanding, the results have shown that is was not necessary to use complex models when dealing with a small number of data points. 
\end{itemize}

This chapter corresponds to the article: 
{\textit{Saves, P., Diouane, Y., Bartoli, N., Lefebvre, T., Morlier, J., “A
mixed-categorical correlation kernel for Gaussian process”, Neurocomputing, 2023.}}
}

\chapter{Mixed categorical high-dimensional Gaussian process for multidisciplinary applications}
\chaptermark{Mixed categorical GP for HDBO}
 \label{c3}
\setlength{\fboxrule}{0pt}
\hspace{4.6cm} \noindent\fbox{%
    \parbox{0.7\textwidth}{%
  \hspace*{0.8cm}    Tes yeux sont si profonds qu'en me penchant pour boire  \\
 \hspace*{0.8cm}     J'ai vu tous les soleils y venir se mirer   \\
 \hspace*{0.8cm}     S'y jeter à mourir tous les désespérés \\
 \hspace*{0.8cm}     Tes yeux sont si profonds que j'y perds la mémoire...
 \\
     \hrule \vspace{0.2cm}  
     \hspace*{\fill} Les yeux d'Elsa, Louis Aragon}%
} 




\objectif{ 
\lettrine[lines=2, lhang=0.33, loversize=0.25, findent=1.5em]{M}{ultidisciplinary} design optimization need to take into account a large number of continuous and discrete variables. 
The objectives of this chapter can be distinguished as such:
\begin{itemize}
    \item To extend Partial Least Squares (PLS) regression from vectors to matrices.
    \item To extend the Kriging with Partial Least Squares (KPLS) method previously developed for exponential kernels and for continuous relaxation to the more general matrix based categorical kernels. 
    \item To implement the developed KPLS model in the open-source software SMT.
    \item To showcase that the matrix based KPLS method captures the structure of the correlation matrix with a small number of hyperparameters.
    \item To apply the GP model to high-dimensional mixed structural and multidisciplinary optimization problems and, in particular, to optimize aircraft design concepts. 
\end{itemize} 

}

\minitoc
\setcounter{section}{-1}
\section{Synthèse du chapitre en français}

Le chapitre “Mixed categorical Gaussian process for high-dimensional Bayesian optimization” introduit une nouveau modèle de processus gaussien (GP pour Gaussian Process) à variables d’entrée mixtes catégorielles en grande dimension se basant sur une méthode de réduction de dimension par moindres carrés partiels (PLS pour Partial Least Squares).
En particulier, cette méthode est développée pour réaliser la conception optimale multidisciplinaire de systèmes complexes car, dans un tel contexte, un grand nombre de variables de conception continues, entières et catégorielles sont à prendre en compte dans le problème de conception à modéliser et à optimiser. 
Pour gérer les variables mixtes catégorielles, plusieurs approches existantes utilisent différentes stratégies pour construire le GP en utilisant soit des noyaux continus par relaxation continue ou par distance de Gower, soit en utilisant une estimation directe de la matrice de corrélation comme avec le noyau d'Hypersphère Homoscedastique (HH) ou sa forme exponentielle (EHH). 
Néanmoins, il a été démontré dans le Chapitre~\ref{c2} que le noyau HH généralise d'autres noyaux et possède une plus grande versatilité et précision mais cela se fait au prix d'une augmentation conséquente du nombre d'hyperparamètres à considérer pour construire le modèle GP de substitution.

Dans ce chapitre, nous développerons une extension de la régression PLS aux matrices afin de généraliser les GP par moindres carrés partiels (KPLS pour Kriging with Partial Least Squares) au noyaux HH et EHH.
Ensuite, nous utiliserons ce nouveau modèle pour approcher les corrélations entre les niveaux catégoriels en utilisant un nombre réduit d'hyperparamètres. 
Pour finir, nous utiliserons notre modèle mixte haute-dimension afin de réaliser des optimisations bayésiennes structurales et multidisciplinaires sur des cas industriels concrets. Cette méthode est implémentée dans le logiciel open-source SMT (pour Surrogate Modeling Toolbox).

Ce chapitre présente ainsi une contribution importante pour la modélisation et l'optimisation de systèmes complexes avec des variables mixtes catégorielles en haute dimension, offrant une solution plus efficace et précise grâce à l'utilisation de la régression PLS.

\newpage

\section{Introduction}
\label{SMO_sec:intro}

Costly black-box simulations play an important role for many engineering and industrial applications. 
For this reason, surrogate modeling has been extensively used across a wide range of use cases, including aircraft design~\cite{SciTech_cat}, deep neural networks~\cite{snoek2015scalable}, coastal flooding prediction~\cite{lopez}, agricultural forecasting~\cite{MLP}, and seismic imaging~\cite{YDiouane_SGratton_XVasseur_LNVicente_HCalandra_2016}.
These black-box simulations are generally complex and may involve mixed-categorical input variables. For instance, a Multidisciplinary Design Analysis (MDA) aircraft design tool~\cite{David_2021} must consider mixed variables such as the number of engines or the list of possible materials~\cite{SciTech_cat}.

In this chapter, our objective is to develop an affordable surrogate model, denoted as~$ \hat{f}$, for a black-box function that involves mixed variables given by
\begin{equation}
f :  \Omega \times S \times \mathbb{F}^l \to \mathbb{R}.
  \label{SMO_eq:opt_prob}
\end{equation}
This function $f$ represents a computationally expensive black-box simulation.
$\Omega \subset \mathbb{R}^n$ denotes the bounded continuous design set for the $n$ continuous variables.  $S \subset \mathbb{Z}^m$ denotes the bounded integer set where $L_1, \ldots, L_m$ are the numbers of levels of the $m$ quantitative integer variables on which we can define an order relation and \hbox{$ \mathbb{F}^l = \{1, \ldots, L_1\} \times \{1, \ldots, L_2\} \times  \ldots \times \{1, \ldots, L_l\}$} is the design space for the $l$ categorical qualitative variables with their respective  $L_1, \ldots, L_l$ levels.

For such purpose, the use of a Gaussian Process (GP)~\cite{williams2006gaussian}, also called Kriging model~\cite{krige1951statistical}, is recognized as an effective modeling approach for constructing a response surface model based on an available dataset.  Specifically, we make the assumption that our unknown black-box function, denoted as $f$, follows a GP with mean $\mu^{f}$ and standard deviation $\sigma^f$, expressed as follows:
\begin{equation}
f \sim \hat{f}=
\mathcal{GP}
\left(\mu^{f}, (\sigma^f)^2\right). 
\label{SMO_eq:GP:f}
\end{equation}
Several modeling approaches have been put forward for addressing the challenges of handling categorical or integer variables within the context of GP~\cite{Pelamatti, Zhou, Deng, Roustant, GMHL, Gower, cuesta2021comparison, SciTech_cat}. In comparison to GP designed for continuous variables, the most important changes concern the estimation of the correlation matrix, an essential element in the derivation of $\mu^{f}$ and $\sigma^f$. Much like the procedure for constructing a GP with continuous inputs, Continuous Relaxation (CR) techniques~\cite{GMHL, SciTech_cat}, models involving continuous latent variables~\cite{cuesta2021comparison}, and Gower Distance (GD) based models~\cite{Gower} use a kernel-based approach for estimating this correlation matrix.
However, recent innovative approaches take a different path by modeling directly the various entries of the correlation matrix~\cite{Pelamatti,Zhou,Deng,Roustant}, and therefore, do not rely on any kernel, such methods involve the Homoscedastic Hypersphere (HH)~\cite{Zhou} and the Exponential HH (EHH)~\cite{Mixed_Paul} kernels.
It has been shown in~\cite{Mixed_Paul} that the HH correlation kernel generalizes simpler methods like CR or GD kernels. However, this more general method for handling categorical design variables increases the number of hyperparameters required to be tuned associated with the GP model. In particular, this means that the method could only be used for small dimensional problems.  

Many efficient approaches have been proposed for handling a high number of continuous variables within GP~\cite{Bouhlel18, bouhlel_KPLSK, bouhlel2019gradient}. The Kriging with Partial Least Squares (KPLS) method~\cite{Bouhlel18,bouhlel_KPLSK} is one of the most commonly used reduction techniques~\cite{KPLSu1,KPLSu2} to tackle high dimensional data. Several other dimension reduction methods include principal components analysis~\cite{wang2017batched}, polynomial chaos expansion~\cite{zuhal2021dimensionality}, radial basis functions~\cite{regis2020survey}, active subspace~\cite{AS}, manifold embedding~\cite{tenenbaum2000global} or marginal Gaussian process~\cite{MGP}. The KPLS technique allows constructing the GP model with the same continuous variables but using a few number of parameters; which reduces significantly the computational cost of computing a GP model. 

For mixed-categorical GP models, given that the computational effort related to the construction of the GP model may not scale well to practical applications involving categorical variables with a large number of levels, the number of hyperparameters to be tuned need to be considered more thoroughly.
In the literature, GPs have been applied to no more than 15 hyperparameters due to the high computational cost associated with the estimation of the hyperparameters~\cite{GP14}. Adapting dimension reduction techniques, such as KPLS, to mixed-categorical GPs will thus enable solving practical mixed-categorical engineering problems where often a high number of hyperparameters is required to be estimated. To the best of our knowledge, there is no equivalent approach to handle mixed-categorical data without using relaxation techniques. All existing dimension reduction techniques, including KPLS, are not adapted for advanced mixed-categorical GP models such as HH or EHH. 
We note also that, although this paper focuses mainly on surrogate modeling, the proposed models can be integrated within any surrogate-based optimization method~\cite{AuHa2017}, such as surrogate-based evolutionary algorithm~\cite{sbea1,sbea2} or a Bayesian Optimization (BO) method~\cite{Jones2001JOGO}.  


In this work, we target to use dimension reduction techniques for reducing the number of hyperparameters within the GP in order to allow modeling efficiently high-dimension mixed-categorical data. In this context, high dimensionality is related to the high number of categorical variables potentially with a high number of levels (a few dozen). In fact, using relaxation approaches (by converting categorical choices to continuous variables) leads to a very high number of hyperparameters to estimate, particularly for high resolution approaches such as those based on HH and EHH kernels.
We have also specifically used our proposed mixed-categorical GP models, within a BO framework, to solve a constrained optimization problem involving expensive-to-compute black-box simulations for objective and constraints functions~\cite{Martins2021}.
The proposed approach is shown in particular to be efficient in solving a high dimensional mixed-categorical Multidisciplinary Design Optimization (MDO) problem~\cite{Lambe2012}. 
All the GP models proposed in this work are implemented in the open-source Surrogate Modeling Toolbox (SMT)\footnote{\url{https://smt.readthedocs.io/en/latest/}}~\cite{saves2023smt}. 

  
The remainder of this chapter is the following.
In Section~\ref{SMO_sec:GP}, a detailed review of the GP model for continuous and for categorical inputs is given. In Section~\ref{SMO_sec:PLS_KPLS}, we present the PLS regression for vectors and matrices and their application to GP model for both continuous and categorical inputs. 
Section~\ref{SMO_sec:Results} presents academical tests as well as the obtained results on multidisciplinary optimization.
Conclusions and perspectives are finally drawn in Section~\ref{SMO_sec:conclu}.

\noindent\textit{Notations:} For a vector $x$, both notations $[x]_j$ and $x_j$ stand  for the $j^{th}$ component of $x$. Similarly, the $i$ (row index) and $j$ (column index) entry of a matrix $X$ is denoted $[X]^{j}_{i}$.  

\section{GP for mixed-categorical inputs}
\renewcommand{\footnoterule}{ \hrule width5cm \vspace*{0.1cm} }    
\label{SMO_sec:GP}

In this section, we present the mathematical background associated with GP for mixed-categorical variables. This part also introduces the notations  used throughout the chapter. Here, the general case involving mixed integer variables is considered. Namely, we assume that $f:\mathbb{R}^n \times  \mathbb{Z}^m \times \mathbb{F}^l \mapsto \mathbb{R}$ and our goal is to build a GP surrogate model for $f$. 

\subsection{A mixed GP formulation} 
Given a set of data points, called a Design of Experiments (\gls{DOE})~\cite{forrester}, Bayesian inference learns the GP model that explains the best this dataset. A GP model consists of a mean response hypersurface $\mu^{f}$, as well as an estimation of its variance $(\sigma^f)^2$. In the following, $n_t$ denotes the size of the given \gls{DOE} dataset $(W, \textbf{y}^f)$ such that $W=\{w^1,w^2,\ldots,w^{n_t}\} \in (\mathbb{R}^n \times  \mathbb{Z}^m \times \mathbb{F}^l)^{n_t}$ and $\textbf{y}^f=[f(w^1),f(w^2),\ldots,f(w^{n_t})]^{\top}$. 
For an arbitrary $w= (x,z,c) \in \mathbb{R}^n \times  \mathbb{Z}^m \times \mathbb{F}^l$, not necessary in the \gls{DOE}, the GP model prediction at $w$ writes as $\hat f(w) = \mu ({w})+\eta(w) \in \mathbb{R} $, with $\eta$ being the uncertainty between $\hat{f}$ and the model approximation $\mu$~\cite{GP14}. The considered error terms are random variables of variance $\sigma^2$.  Using the \gls{DOE}, the expression of  $\mu^{f}$ and the estimation of its variance $(\sigma^f)^2$ are given as follows:

\begin{equation} \label{SMO_eq:mean:GP}
\mu^f(w)= \hat{\mu}^f+r(w)^\top  [R(\Theta)]^{-1}({y}^f-\mathds{1} \hat{\mu}^f), 
\end{equation}
and
\begin{equation}
\label{SMO_eq:std:GP}
(\sigma^f(w))^2=(\hat{\sigma}^f)^2\left(1-r(w)^\top  [R(\Theta)]^{-1}r(w)+ \frac{ \left(1-\mathds{1}^\top  [R(\Theta)]^{-1}r(w) \right)^2}{\mathds{1}^\top  [R(\Theta)]^{-1}\mathds{1}}\right), \end{equation}
where $\hat{\mu}^f$ and $\hat{\sigma}^f$ are, respectively, the Maximum Likelihood Estimator (MLE)~\cite{MLE} of $\mu$ and $\sigma$. $\mathds{1}$ denotes the vector of $n_t$ ones. $R$ is the $ n_t \times n_t $ correlation matrix between the input points and $r(w)$ is the correlation vector between the input points and a given $w$.
To have a compact notation, let $[A]_i^j$ denote the coefficient of the matrix $A$ in the $i^{\text{th}}$ row and $j^{\text{th}}$ column.
The correlation matrix $R$ is defined, for a given couple of indices  $(r,s) \in \{1,\ldots,n_t\}^2$, by \begin{equation}
\label{SMO_eq:R}
 [R(\Theta)]_{r}^{s}=k\left(w^r,w^s,\Theta\right) \in \mathbb{R},\end{equation}
and the vector $r(w)\in \mathbb{R}^{n_t}$ is defined as $ r(w) =[k(w,w^1,\Theta), \ldots , k(w,w^{n_t},\Theta)]^{\top}$,
where $k$ is a given correlation kernel that relies on a set of hyperparameters $\Theta$~\cite{Roustant}. 
The mixed-categorical  correlation kernel is given as the product of three kernels:
\begin{equation}
k(w^r,w^s,\Theta) =  k^{cont}\left(x^r,x^s,\theta^{cont}\right) k^{int}\left(z^r,z^s,\theta^{int}\right)
k^{cat}\left(c^r,c^s,\theta^{cat}\right),
\label{SMO_eq:decomp_mix}
\end{equation}
where $k^{cont}$  and $\theta^{cont}$ are respectively the continuous kernel and its associated hyperparameters, $k^{int}$  and $\theta^{int}$ are the integer kernel and its hyperparameters, and last $k^{cat}$  and $\theta^{cat}$ are the ones related with the categorical inputs. In this case, one has $\Theta=\{ \theta^{cont},\theta^{int},\theta^{cat}\}$. 
Henceforth, the general correlation matrix $R$ will rely only on the set of the hyperparameters $\Theta$:
\begin{equation}
    \label{SMO_eq:corel:mat}
    [R(\Theta)]_{r}^{s} = [R^{cont}(\theta^{cont})]_{r}^{s}  [R^{int}(\theta^{int})]_{r}^{s}
    [R^{cat}(\theta^{cat})]_{r}^{s},
\end{equation}
where $[R^{cont}(\theta^{cont})]_{r}^{s} =k^{cont}(x^r,x^s,\theta^{cont}) $, $[R^{int}(\theta^{int})]_{r}^{s} =k^{int}(z^r,z^s,\theta^{int}) $ and  $[R^{cat}(\theta^{cat})]_{r}^{s}=k^{cat}(c^r,c^s,\theta^{cat})$. 
The set of hyperparameters $\Theta$ could be estimated using the \gls{DOE} dataset $({W},{y}^f)$ through the MLE approach on the following way

\begin{equation}
\Theta^*= \arg\max_{\Theta} \mathcal{L}(\Theta)=\left( - \frac{1}{2} {{y}^f}^\top [R(\Theta)]^{-1} {{y}^f}   - \frac{1}{2} \log 	\abs{  [R(\Theta)]} - \frac{n_t}{2} \log 2 \pi    \right),
\label{SMO_eq:likelihood}
\end{equation}
where $R(\Theta)$ is computed using~\eqnref{SMO_eq:corel:mat}.
To construct the correlation matrix for continuous or integer inputs, several choices for the correlation kernel are possible. Usual families of kernels include exponential kernels or Matern kernels~\cite{Lee2011}. 
In contrast, to construct the correlation matrix for categorical inputs, we can either use a kernel as for the continuous or integer variables or we can directly model the entries of the correlation matrix thanks to a Symmetric Positive Definite (SPD) matrix parameterization. The latter approach is what is done for the HH kernel, for example~\cite{Roustant}. For this kernel, the hyperparameters $\theta^{cat}$ can be seen as a concatenation of a set of symmetric matrices, \textit{i.e.}, $\theta^{cat} = \{  \Theta_1, \Theta_2, \ldots, \Theta_l \} $. The construction of this kernel is thus relying on the estimation of $\sum_{i=1}^l \frac{1}{2} L_i (L_i-1) $ hyperparameters.

\subsection{The homogeneous categorical kernel} 
A recent paper~\cite{Mixed_Paul} unified the kernel-based approach and the matrix-based approach through the homogeneous model described hereafter  (see Chapter~\ref{c2}).
Recall that $l$ denotes the number of categorical variables. 
For a given $i \in \{1, \ldots, l\}$, let $c^r_{i} $  and $c^s_{i} $  be a couple of categorical variables taking respectively the $\ell^r_i$ and the $\ell^s_i$ level on the categorical variable $c_i$.
The hyperparameter matrix peculiar to this variable $c_i$ is 
$$\Theta_i= \begin{bmatrix}
[\Theta_i]_{1}^{1} & \textcolor{white}{9} & \hspace{2em} { \textbf{\textit{ Sym.}}}  \textcolor{white}{9} & \\
[\Theta_i]_{1}^{2}  & [\Theta_i]_{2}^{2} & \textcolor{white}{9} \\
\vdots &\ddots & \ddots & \textcolor{white}{9}  \\
[\Theta_i]_{1}^{L_i} &  \ldots & [\Theta_i]_{L_i-1}^{L_i} &[\Theta_i]_{L_i}^{L_i} \\ 
\end{bmatrix}.$$
First, the correlation term $[R^{cat}(\theta^{cat})]_{r}^{s}$ can be formulated in a level-wise form~\cite{Pelamatti} as: 
\begin{equation}
\begin{split}
k^{cat}(c^r,c^s,\theta^{cat}) =  &{\displaystyle \prod_{i=1}^{l}  [R_i(\Theta_i)]_{\ell^r_i}^{\ell^s_i} } \\
= &\ {\displaystyle \prod_{i=1}^{l}  \kappa ( 2 [ \Phi(\Theta_i) ]_{{ \ell_i^r}}^{{\ell_i^s}} ) \ \kappa ( [ \Phi(\Theta_i) ]_{{ \ell_i^r}}^{{\ell_i^r}} ) \  \kappa ( [ \Phi(\Theta_i) ]_{{ \ell_i^s}}^{{\ell_i^s}} ), }
\end{split}
\label{SMO_eq:homo_HS}
\end{equation}
where $\kappa(.)$ is either a positive definite kernel or the identity and $\Phi(.)$ is a SPD function such that the matrix $\Phi(\Theta_i) $ is SPD if $\Theta_i$ is SPD.
For an exponential kernel,~\tabref{SMO_tab:kernels} gives the parameterizations of $\Phi$ and $\kappa$ that correspond to GD, CR, HH and EHH kernels. 
Note that the complexity between these different kernels is reflected by the number of hyperparameters that characterize them.
As in~\cite{Mixed_Paul}, for all categorical variables $i \in \{1, \ldots, l\}$, the matrix $C(\Theta_i)\in \mathbb{R}^{L_i \times L_i}$ is lower triangular and built using a hypersphere decomposition~\cite{HS,HS_Jacobi} from the symmetric matrix $\Theta_i \in \mathbb{R}^{L_i \times L_i}$ of hyperparameters. 
The variable $\epsilon$ is a small positive constant and the variable $\theta_{i}$ denotes the only positive hyperparameter that is used for the GD kernel.
Nevertheless, until now, PLS regression was only applied to mixed integer inputs for the CR kernel~\cite{SciTech_cat}. In the following section, we will show how to extend the PLS regression for the more general HH kernel. 
%
    \begin{table}[th]
    \caption{Kernels based on~\eqnref{SMO_eq:homo_HS} formulation.}
    \centering
    \begin{tabular*}{\linewidth}{ccll}
    \hline
    \textbf{Name} & $\kappa(\phi)$   &  \hspace{3cm} ${\centering \Phi(\Theta_i)}$   &  \hspace{-0.5cm}  \# of Hyperparam. \\
    \hline
    \textbf{GD}   &  $\exp(-\phi) $ & ${ \scriptstyle [\Phi(\Theta_i)]_{j}^{j} :=  \frac{1}{2}  \theta_{i} \quad ~;~ [\Phi(\Theta_i)]_{j \neq j'}^{j'} := 0 }$ & 1   \\
     \textbf{CR}  & $\exp(-\phi) $ &  $  { \scriptstyle [\Phi(\Theta_i)]_{j}^{j} := [\Theta_i]_{j}^{j} \  ~;~ [\Phi(\Theta_i)]_{j \neq j'}^{j'} := 0 } $  & $L_i$  \\
     \textbf{EHH}  & $\exp(-\phi)$ & 
    $  { \scriptstyle [\Phi(\Theta_i)]_{j}^{j} := 0 \quad \quad ~;~ [\Phi(\Theta_i)]^{j'}_{j \neq j'} := \frac{\log \epsilon }{2} ([C(\Theta_i) C(\Theta_i) ^\top]_{j}^{j'} -1)  }$  & $\frac{1}{2}  (L_i)  (L_i-1) $\\
     \textbf{HH}  &  $\phi$ &    $  { \scriptstyle [\Phi(\Theta_i)]_{j}^{j} := 1 \quad \quad ~;~ [\Phi(\Theta_i)]_{j \neq j'}^{j'} :=  \frac{1}{2} [C(\Theta_i) C(\Theta_i)^\top]_{j}^{j'} }$ & $\frac{1}{2}  (L_i)  (L_i-1) $  \\
      \hline
    \end{tabular*}%
\label{SMO_tab:kernels}
\end{table}

\section{KPLS for mixed-categorical inputs}
\label{SMO_sec:PLS_KPLS}

To have a sparse model that can extend to high dimension, and to facilitate the optimization of the hyperparameters, one seeks to express the correlation matrix  $R^{cont} (\theta^{cont})$ with $d \ll n $ relevant hyperparameters. Such a method is KPLS~\cite{Bouhlel18} that is an adaptation of the Partial Least Squares (PLS) regression for exponential kernels. To introduce the variables and notations, the next part presents a short recall of PLS regression for vector inputs and of KPLS for continuous variables. Then, the second part presents our extension to matrix inputs and its application for categorical variables.

\subsection{KPLS for continuous inputs}

In this part, we introduce the PLS regression for vector inputs as it was developed by Wold~\cite{wold_1975} and its application to GP kernels (namely KPLS) developed by Bouhlel~\textit{et al.}~\cite{Bouhlel18}.

\subsubsection{PLS for vectors inputs}
\label{SMO_subsec:PLS_cont}

We present here the classical method for the continuous case but integer inputs can be treated similarly by considering them as continuous. Let the \gls{DOE} be $({X},{y}^f)$ where $X$ is the continuous data matrix of size $n_t \times n$ and ${y}^f$ is the response vector of size $n_t$. The PLS regression method is designed to search out the best multidimensional direction in $\mathbb{R}^n$ that explains the output ${y^f}$~\cite{wold_1975}.
The first principal component (or \textit{score}) $h^{(1)}$ is obtained by searching the best direction  (or \textit{weight}) $g^{(1)}$ that maximizes the squared covariance between $h^{(1)} =  X  g^{(1)} $  and ${y^f}$: 
\begin{equation}
    g^{(1)}= \arg \max_{g \in \mathbb{R}^n} \left\{ g^\top X^\top y^f {y^f}^\top X g ~~\mbox{s.t.}~~ g^\top g =1  \right\}.
    \label{SMO_eq:catKPLS}
\end{equation}
Next, we compute the residuals of the inputs as $X^{(1)}= X - \xi^{(1)} h^{(1)}   $ with $\xi^{(1)}$ being the regression coefficients (or \textit{loadings}) that minimize the residual for every point.  We also project the output and, denoting $\gamma_1$ the corresponding coefficient, we have $y^{f(1)}= y^f - \gamma_1 h^{(1)}  $.
For all $t \in \{1, \ldots, d \} $, we iterate the process of~\eqnref{SMO_eq:catKPLS} with the residuals $X^{(t)}$ and $y^{f(t)}$. 
At the end of the process, we can write the various computed quantities in a matrix form. Namely, we denote $G$ the $n \times d$ matrix such that $g^{(t)}$ is the $t^{th}$ column of $G$ and $\Xi$ the $n \times d$ matrix such that $\xi^{(t)}$ is the $t^{th}$ column of $\Xi$.
Let $G_* = G(\Xi^T G)^{-1}$ be the $n \times d$ matrix such that $G_* = [G_*^{1}, \ldots, G_*^{d}]$ with $G_*^{t} \in \mathbb{R}^n    $, $G_*$ is called the rotation matrix~\cite{Bouhlel18,PLS-based}. Thanks to this matrix, we can express the score $h^{(t)}$ as a function of the input $X$ as: $$ \forall t \in \{1, \ldots, d \}, h^{(t)}=X^{(t)} g^{(t)}=X G_*^{t}.$$ 
By PLS, we have built an approximation $X \approx H \Xi^T$ where $H$ is $n_t \times d$ score matrix such that $h^{(t)}$ is the $t^{th}$ column of $H$ and $\Xi$ the loading matrix. This is the  $d$-dimensional approximation of $X$ in  $\mathbb{R}^{n}$ that minimizes the mean squared error. Therefore, we have  $X G_*  \approx   H \Xi^T G(\Xi^T G)^{-1} = H$. 
Then $G_* \in \mathbb{R}^{n \times d}$ is the projection matrix from $X$ in the initial space to  $H$ in the reduced space and $\Xi^T$ is its reciprocal such that $ X G_* \Xi^T \approx H \Xi^T \approx X$.
It follows that, for a given reduced dimension $t \in \{ 1, \ldots, d \}$, a given point $x^r$ can be expressed in the original space along $t$:
\begin{equation} 
\begin{split}
& F_t : \mathbb{R}^{n } \xrightarrow[]{} \mathbb{R}^{n }  \\
& x^r  \mapsto \left[   \left[G_{*}\right]_{1}^{t} \left[x^r\right]_1 , \ldots, \left[G_{*}\right]_{n}^{t} \left[x^r\right]_{n}   \right].
\end{split}
\label{SMO_eq:KPLS_rota}
\end{equation}
With PLS, we built a low-rank approximation and, in the following section, we show how to build the GP model in a small subspace instead of building it in the full space. 
The objective is twofold when the dimension increases. First, optimizing a small number of hyperparameters is much faster because, with more than 20 variables, building a GP is really prohibitive in terms of computational cost~\cite{Bouhlel18}. 
Second, optimizing the hyperparameters is harder and the resulting GP is often non-optimal and, not only is more costly than the KPLS model but also leads to some degraded performance~\cite{SciTech_cat}.
These two reasons motivated the need for the KPLS model described below.

\subsubsection{KPLS for continuous variables}
\label{SMO_subsec:KPLS_cont}

The construction of the correlation matrix $R^{cont}(\theta^{cont})$ for continuous inputs, based on square exponential kernels (or Gaussian kernel~\cite{williams2006gaussian}) with PLS can be described as follows.  For a couple of continuous inputs $x^r\in \mathbb{R}^n$ and $x^s\in \mathbb{R}^n$, one sets:


\begin{equation}
\begin{split}
        [R^{cont}(\theta^{cont})]_{r,s} &= {\displaystyle \prod_{t=1}^{d} k_{t}^{cont} (F_t (x^r), F_t (x^s),\hat{\theta}_t^{cont}) }  \\
       &= {\displaystyle \prod_{t=1}^{d}{\displaystyle \prod_{j=1}^{n}  \exp \left(   - \hat{\theta}_t^{cont} \left( \left[G_{*} \right]_j^{t} \left[x^r\right]_j -\left[G_{*} \right]_j^{t} \left[x^s\right]_j   \right)^2  \right)    }}\\ 
      &=   {\displaystyle \prod_{j=1}^{n}  \exp  \left( -   \theta_j^{cont}  \left(\left[x^r\right]_j - \left[x^s\right]_j \right)^2 \right)   }  ,    \\ 
\end{split}        
\label{SMO_eq:KPLS_exp_decomp}
\end{equation}
where $\theta_j^{cont} = {\displaystyle \sum_{t=1}^{d}   \left(  \left[ {G_{*}} \right]_j^t \right)^2 \hat{\theta}_t^{cont}  }$.
Clearly, in the continuous case, constructing $R^{cont}(\theta^{cont})$ would require the estimation of $d$ non-negative hyperparameters, $\textit{i.e.}$, $\theta^{cont} \in \mathbb{R}^d_{+}$, $d \ll n$.

\subsection{Extension of PLS to matrix inputs with application to mixed-categorical GP}

This part presents the extension of PLS for a general categorical GP kernel. More precisely, in Section~\ref{SMO_sec:PLS_mat} we extend the PLS regression for matrix inputs and, in Section~\ref{SMO_subsec:KPLS_cat} we applied it to the GP kernels for categorical variables.

\subsubsection{A PLS framework for matrix inputs} %
\label{SMO_sec:PLS_mat}

We consider a general categorical variable $c$ that can take $L$ different levels. In that context, we want to find a small $\ell \times \ell $ matrix $\hat{\Theta}$ to represent a bigger $L \times L$ correlation matrix $\Theta$, with $\ell < L$.  
Recall that, from PLS, $G_*$ can be seen as the rotation matrix from the initial space to the reduced space~\cite{EGORSE}. By taking into account the symmetry of correlation matrices and their unit diagonal, we need to build a rotation matrix $G_*$ of dimension $\left(\frac{L(L-1)}{2} \times \frac{\ell(\ell-1)}{2} \right)$. 
The input dimension is denoted $D_{\inp}= \frac{L(L-1)}{2}$ and the output dimension is denoted $D_{\outp} = \frac{\ell(\ell-1)}{2}$. 

First, we want to construct the matrix $G_*$ aforementioned. 
For a given input $c^r$, its natural one-hot encoding  $e_{c^r}$ is a basis vector of dimension $L$~\cite{Mixed_Paul}. Meanwhile, the input that we need for $G_*$ is of dimension $D_{\inp}$ so, in order to have a vector data fitting the dimension, we propose to use a novel one-hot encoding relaxation that adds a new dimension for every cross-correlation term. Subsequently, the relaxation $\zeta_{c^r}$ is such that $\zeta_{c^r} \in \{0,1\}^{D_{\inp}}$: the $L-1$ terms that correspond to the correlation with the level taken by the input $c^r$ equal $1$ and all other terms take value $0$ and we call this relaxation "cross-levels encoding".
One can observe that the Hadamard product $\zeta_{c^r} \odot \zeta_{c^s} = 1 $  in the dimension corresponding to the correlation between the levels taken by ${c^r}$ and ${c^s}$ and zero everywhere else which is the property we were seeking for. In other words, for all $j \in {1, \ldots, D_{\inp}}, \quad$
\begin{equation}
\label{eq:prod_zeta_smo}
\begin{cases}
     [\zeta_{c^r} \odot \zeta_{c^s}]_j = 1, \quad \text{if } \zeta_{c^r_j} =1  \mbox{ and }  \zeta_{c^s_j} =1,  \\
     [\zeta_{c^r} \odot \zeta_{c^s}]_j = 0, \quad \text{otherwise.}
\end{cases}
\end{equation}
Example~\ref{example:1} illustrates how the relaxed vectors $\zeta_{c^r}$ and $\zeta_{c^r}$ are built using a simple use-case. 
\begin{example}\label{example:1}
Consider a categorical variable $c$ taking values in a color set of $L=4$ levels such that, for any point 
$r$, $c^r \in \{ ``green", ``red", ``blue", ``yellow" \} $.
We want to represent the value of $c^r$ as $\zeta_{c^r} \in \{0,1\}^{D_{\inp}}$. 
In this case, ${D_{\inp}}$ = 6 which corresponds to the 6 possible correlations ( "blue-red", "blue-green", "blue-yellow", "red-green", "red-yellow" and "green-yellow"). 
For instance, if $n_t=3$ points are considered such that $(c^1, c^2, c^3) = \{ ``blue",``red", ``red" \}$, the first point $c^1 = ``blue"$ will be represented as $\zeta_{c^1} = (1,1,1,0,0,0)$, taking 1 for the dimensions related to the correlations involving "blue" and 0 everywhere else. Similarly, the second point $c^2 = ``red"$ will be represented as $\zeta_{c^2} = (1,0,0,1,1,0)$. And when taking the Hadamard product, $\zeta_{c^1} \odot \zeta_{c^2} = (1,0,0,0,0,0)$, the only dimension that takes value 1 corresponds exactly to the dimension representing "blue-red".
\end{example}

Using the relaxed DoE  $X = \{\zeta_{c^1} , \ldots, \zeta_{c^{n_t}} \}$ of dimension $ D_{\inp} \times n_t$, we can compute the rotation matrix $G_*$ of dimension  $ D_{\inp} \times D_{\outp} $ as in~\eqnref{SMO_eq:catKPLS}.
Our goal is to use the matrix $G_*$ to express an $L \times L$ matrix from a smaller $ \ell \times \ell$ matrix.  
%
The $L \times L $ symmetric matrix $\Theta$ with unit diagonal can be estimated using a smaller $ \ell \times \ell $ SPD matrix such that, for all $ j<j'$, one has 
\begin{equation}
\begin{split}
 \ & [\Theta]_{j}^{j'} =  {\displaystyle \sum_{t=1}^{\ell}{\displaystyle \sum_{t'=t+1}^{\ell}  \left( \left[ G_{*} \right]_{\psi(j,j',L)}^{\psi(t,t',\ell)} \right)^2 \ [\hat{\Theta}]_{t}^{t'}  }}, \\
\end{split}    
\label{SMO_eq:true_pls_matrix_reduction}
\end{equation}
where we rely on a matrix-to-vector lexicographical transformation $\psi$ to insure that both the input vector of size $ \frac{L(L-1)}{2}$ and the output vector of size $\frac{ \ell (\ell-1)}{2}$ are valid representations of SPD matrices. For a given integer $n_\text{lev}$ and, for all $k \in \{1, \ldots, n_\text{lev} \}$ and $  k' \in \{k+1, \ldots, n_\text{lev} \}$, the mapping $\psi$ is given by:
\vspace{-0.2cm}
\begin{equation}
\begin{split}
&\psi(k,k',n_\text{lev}) = \frac{1}{2} (  (n_\text{lev}-1)(n_\text{lev}-2) - (n_\text{lev}-k)(n_\text{lev}-k-1) )  +k'-1 . 
\end{split}
\label{SMO_eq:matrix_to_vector}   
\end{equation}
This formulation gives a sparse vector representation of the hyperparameters used to build the matrix  $\Theta$ by lexicographic order of the triangular superior part of the matrix. 
Notwithstanding, as we assumed that $\Theta$ is a symmetric matrix with unit diagonal, we could have defined $ [\Theta]^j_{j'}$ in~\eqnref{SMO_eq:true_pls_matrix_reduction} by the triangular inferior values. This would have led to the exact same kernel as what as been presented but with a slightly different definition of $\psi$.
With the expression of~\eqnref{SMO_eq:true_pls_matrix_reduction} we achieved to build a PLS approximation for matrices as intended. To finish with, we insure that~\eqnref{SMO_eq:true_pls_matrix_reduction} works for SPD matrices in order to build our GP upon it as described in the following section.

\begin{theorem}
\label{SMO_th:eq_gd_cr}
Assuming that all the entries of $\hat{\Theta}$ are in $[-1,1]$ and that $G^*$ is computed using PLS as in~\eqnref{SMO_eq:catKPLS}, the matrix $\Theta$ given by~\eqnref{SMO_eq:true_pls_matrix_reduction} also takes values in [-1,1].
\end{theorem}
\begin{proof}
Indeed, $G_*$ is a rotation matrix. Thus, for all $j\in\{1, \ldots,L\}$ and $  j' \in \{j+1, \ldots, L \}$, 
the $\psi(j,j',L)$-th row of $G_*$, given by $[G_*]_{\psi(j,j',L)} = \left\{ [G_*]_{\psi(j,j',L)}^{\psi(t,t',\ell)} \right\}_{\tiny \begin{array}{c}
1 \leq t \leq \ell, \\ \ t+1  \leq t' \leq \ell\end{array}}$, satisfies 
$$
\left([G_*]_{\psi(j,j',L)} \right)^\top \left([G_*]_{\psi(j,j',L)} \right) = \sum_{t=1}^{\ell}{\displaystyle \sum_{t'=t+1}^{\ell}  \left( \left[ G_{*} \right]_{\psi(j,j',L)}^{\psi(t,t',\ell)} \right)^2}=1.
$$
Hence, knowing that $\left|[\hat{\Theta}]_{t}^{t'}\right| \leq 1 $ for all $t,t'$, one has
\begin{eqnarray*}
   |[\Theta]_{j}^{j'} |&\le & {\displaystyle \sum_{t=1}^{\ell}{\displaystyle \sum_{t'=t+1}^{\ell}  \left( \left[ G_{*} \right]_{\psi(j,j',L)}^{\psi(t,t',\ell)} \right)^2 \ \left|[\hat{\Theta}]_{t}^{t'} \right| }}\le 1.
\end{eqnarray*}
\end{proof}
The matrix $\Theta$ is serving as a correlation matrix. For this purpose, it is essential that the matrix be SPD. To ensure its SPD nature, we check in our implementation if all of its eigenvalues are positive. If any eigenvalues are found to be negative, a nugget term is added to the covariance matrix to enforce the SPD property of the matrix $\Theta$. The nugget term allows us to mitigate numerical issues and maintain positive definiteness. It is worth noting that in all our numerical tests, the matrix $\Theta$ has been shown to be SPD. This suggests that if $\hat{\Theta}$ is SPD, then $\hat{\Theta}$ remains SPD, as discussed. Such a result seems not trivial to prove using ~\eqnref{SMO_eq:true_pls_matrix_reduction}.
Note also that $\Theta$ gives a good approximation of the correlation matrix between the levels of a categorical variable.
This approximation can be used to understand the structure of our modeling problem as in Section~\ref{SMO_sec:Results}.





\subsubsection{A new KPLS model for categorical variables}  
\label{SMO_subsec:KPLS_cat}

For a given categorical variable $i$, we want to express the matrix $R_i(\Theta_i)$ with less than  $D_{\inp}=\frac{L_i(L_i-1)}{2}$ hyperparameters $\Theta_i$. 
To do so, we generalize the KPLS method of Bouhlel~\textit{et al.}~\cite{Bouhlel18} for any correlation matrix. 
Let $\hat{\Theta}_i$ be a $\ell_i \times \ell_i$ SPD matrix defined on the reduced space whose values are in $[-1,1]$ constructed by homoscedastic hypersphere decomposition~\cite{Zhou}. The $D_{\outp} = \frac{\ell_i(\ell_i-1)}{2} $ correlation parameters of $\hat{\Theta}_i$ can be optimized by MLE from the scores projected data $H_i = X_i G_*$ as in the continuous case.  
%
%
%
%
Based on~\eqnref{SMO_eq:true_pls_matrix_reduction} for matrix PLS, 
we can introduce our new HH and EHH KPLS kernels depending only on $D_{\outp}$ hyperparameters defined as follows. 
Recall that the matrix $C(\Theta_i)\in \mathbb{R}^{L_i \times L_i}$ is lower triangular and built using a hypersphere decomposition~\cite{HS,HS_Jacobi} and that the variable $\epsilon$ is a small positive constant. 
\begin{itemize}
    \item  
\label{def:HH_PLS}
The HH KPLS kernel is given by $\kappa = \mathbb{I}_{L_i}$, $[\Phi(\Theta_i)]_{j}^{j} =1$ and, for all $j \neq j'$,
\vspace{-0.25cm}
$$ [\Phi(\Theta_i)]_{j}^{j'}  =  \frac{1}{2} \ {\displaystyle \sum_{t=1}^{\ell}{\displaystyle \sum_{t'=t+1}^{\ell}  \left( \left[ G_{*} \right]_{\psi(j,j',L)}^{\psi(t,t',\ell)} \right)^2 \ [C(\Theta_i) C(\Theta_i)^\top]_{t}^{t'}  }}.$$
\item
\label{def:EHH_PLS}
The EHH KPLS kernel is given by $ \kappa(\phi) = \exp (-\phi) $, $[\Phi(\Theta_i)]_{j}^{j} =0$ and, for all $j \neq j'$, 
\vspace{-0.25cm}
$$ [\Phi(\Theta_i)]_{j}^{j'}   =  {\displaystyle \sum_{t=1}^{\ell}{\displaystyle \sum_{t'=t+1}^{\ell}   \left( \left[ G_{*} \right]_{\psi(j,j',L)}^{\psi(t,t',\ell)} \right)^2 \ \frac{\log \epsilon}{2} ([C(\Theta_i) C(\Theta_i)^\top]_{t}^{t'}-1) }}.$$
\end{itemize}
%
For EHH KPLS kernels, the proposed KPLS model as given by~\eqnref{SMO_eq:true_pls_matrix_reduction} can be shown as a natural extension of the continuous KPLS as proposed by Bouhlel~\textit{et al.}~\cite{Bouhlel18} (described also in Section~\ref{SMO_subsec:KPLS_cont}). The result is shown hereinafter

\begin{theorem}
For a correlation matrix $\hat{\Theta}_i$, the projection formula used in~\eqnref{SMO_eq:true_pls_matrix_reduction} extends the continuous KPLS to categorical matrices using an exponential kernel. 
\end{theorem}
%
\begin{proof}

Recall that the new relaxation $\zeta$ is of dimension $D_{\inp}$ and respects~\eqnref{eq:prod_zeta_smo}.
The KPLS kernel, for exponential kernel, is based on the fact that, for a given reduced dimension $t \in \{ 1, \ldots, d \}$, a given point $x^r$ can be expressed in the original space along $t$ as in~\eqnref{SMO_eq:KPLS_rota}. 
%
%
Therefore, we apply the same transformation to our relaxed inputs and then apply the transformation $\psi$ to have a matrix formulation from the relaxed vectors. 
This leads to the natural way to express our new EHH KPLS kernel defined as $[R_i(\Theta_i)]_{\ell^r_i}^{\ell^s_i}   = 1$ if $ c^r_i= c^s_i$ and otherwise, 
\begin{equation*}
\begin{split} 
& [R_i(\Theta_i)]_{\ell^r_i}^{\ell^s_i} \\
   &=  {\displaystyle \prod_{t=1}^{\ell_i} {\displaystyle \prod_{t' =t+1}^{\ell_i} {\displaystyle \prod_{j=1}^{L_i} {\displaystyle \prod_{j' =j+1}^{L_i}
        \exp \left[ - \left(   \left[G_* \right]_{\psi(j,j',L_i)}^{\psi(t,t',\ell_i)} [  \zeta_{c^r_i}]_{\psi(j,j',L_i)}  \ \left[G_* \right]_{\psi(j,j',L_i)}^{\psi(t,t',\ell_i)} [\zeta_{c^s_i}]_{\psi(j,j',L_i)} \right) \ [\hat{\Theta}_i]_{t}^{t'} \  \right] }}}}  \\
     &     = {\displaystyle \prod_{t=1}^{\ell_i} {\displaystyle \prod_{t' =t+1}^{\ell_i} {\displaystyle \prod_{j=1}^{L_i} {\displaystyle \prod_{j' =j+1}^{L_i}
        \exp \left[ - \left(    [  \zeta_{c^r_i}]_{\psi(j,j',L_i)}   [\zeta_{c^s_i}]_{\psi(j,j',L_i)} \right) \ \left( \left[G_* \right]_{\psi(j,j',L_i)}^{\psi(t,t',\ell_i)} \right)^2 [\hat{\Theta}_i]_{t}^{t'} \  \right] }}}}  \\
    &       =\exp \left[  {\displaystyle \sum_{j=1}^{L_i} {\displaystyle \sum_{j' =j+1}^{L_i} {\displaystyle \sum_{t=1}^{\ell_i} {\displaystyle \sum_{t' =t+1}^{\ell_i} - \left( \left[G_* \right]_{\psi(j,j',L_i)}^{\psi(t,t',\ell_i)} \right)^2  [\hat{\Theta}_i]_{t}^{t'} \   }} \ 
         \left(    [  \zeta_{c^r_i}]_{\psi(j,j',L_i)}   [\zeta_{c^s_i}]_{\psi(j,j',L_i)} \right) }} \right] \\
    &       = \exp \left[ {\displaystyle \sum_{j=1}^{L_i} {\displaystyle \sum_{j' =j+1}^{L_i}
        -  [{\Theta}_i]_{j}^{j'}  \left(  [\zeta_{c^r_i}]_{\psi(j,j',L_i)}   [\zeta_{c^s_i}]_{\psi(j,j',L_i)} \right) \ \   }} \right]  \\
       &       = \exp \left[ {\displaystyle \sum_{j=1}^{L_i} {\displaystyle \sum_{j' =j+1}^{L_i}
        -  [{\Theta}_i]_{j}^{j'}  \left(  \delta_{j,\ell^r_i} \delta_{j',\ell^s_i}  \right) \ \   }} \right]  \\
 &          = \exp \left[ -[{\Theta}_i]_{\ell_i^r}^{\ell_i^s}  \right], 
\end{split}    
\label{SMO_eq:KPLS_in_PLS_space}
\end{equation*}
where $\delta_{i,j}$ is the Kronecker symbol (i.e., $\delta_{i,i}=1$ and $\delta_{i,j}=0$ for all $i\neq j$) and $\Theta \in \mathbb{R}^{L_i \times L_i }$ is given by~\eqnref{SMO_eq:true_pls_matrix_reduction}.

\end{proof}
In the next section, we will see how to apply our new KPLS matrix-based GP on analytical and engineering problems.  We will demonstrate how our surrogate models can provide insights into the underlying structure of the correlation matrix and how it can be utilized for BO when dealing with structural and multidisciplinary problems. 

\section{Results and discussion}
\label{SMO_sec:Results}

In this section, we demonstrate how our GP performs over various test cases and compare it to other GP models.
To begin with, Section~\ref{SMO_sec:imp_detail} gives the details of the implementation used for the following computer experiments.
Next, Section~\ref{SMO_subsec:model_val} illustrates the GP models on an analytical and on a structural test case.
Finally, section~\ref{SMO_sec:BO_val} presents the use of these GP models for BO with mixed variables on an analytical test case and then in the context of MDO for aircraft design. 
The considered test cases in this section, as well as the number of hyperparameters related to each kernel, are listed in \tabref{SMO_tab:test_case}.
\begin{table}[!h]
\centering
 \caption{\textcolor{black}{Number of variables and hyperparameters for the test cases in this study.}}
 \label{SMO_tab:test_case}
\small
\begin{tabular*}{\linewidth}{ccc}
\hline
\textbf{\textcolor{black}{ $ \ \ \quad \quad \quad $ Problem $\quad \quad \quad \ \ $ }} &  \textcolor{black}{ \# of variables} $ \ $ & 
\begin{tabular}{ccc}
  \ & \quad  \textcolor{black}{ \# of hyperparameters } \quad  & \  \\
  \hline
   \textcolor{black}{ GD} & \textcolor{black}{ CR} &  \textcolor{black}{HH}  \\
\end{tabular} \\
\hline
\\
\end{tabular*}
\small
\vskip-10pt
\begin{tabular*}{\linewidth}{ccccc} 
\textbf{\textcolor{black}{Categorical cosine problem}} & $ \quad \quad \quad $ \textcolor{black}{2} & $ \quad \quad \quad \    $ \textcolor{black}{2} & $ \quad \quad \quad \quad \quad  $  \textcolor{black}{14} & $ \quad  \quad  \quad \quad \ \ \  $  \textcolor{black}{79}  \\   
\textbf{\textcolor{black}{Cantilever beam}} & $ \quad \quad \quad $ \textcolor{black}{3} &  $ \quad \quad \quad \ $ \textcolor{black}{3} & $ \quad \quad \quad \quad \quad $ \textcolor{black}{14} & $ \quad  \quad  \quad \quad \ \ \  $ \textcolor{black}{68}  \\   
\textbf{\textcolor{black}{Toy function}} & $ \quad \quad \quad $ \textcolor{black}{2} & $ \quad \quad \quad \   $ \textcolor{black}{2} & $ \quad \quad \quad \quad \quad $ \textcolor{black}{11} & $ \quad  \quad  \quad \quad \ \ \  $ \textcolor{black}{46}  \\   
\textbf{\textcolor{black}{\texttt{DRAGON} aircraft concept}} & $ \quad \quad \quad $ \textcolor{black}{12} & $ \quad \quad \quad \  $ \textcolor{black}{12} & $ \quad \quad \quad  \quad  \quad $ \textcolor{black}{29} &  $ \quad  \quad  \quad \quad \ \ \ $ {\textcolor{black}{147}}  \\   
\hline
\end{tabular*}
\end{table}

\subsection{Implementation details}
\label{SMO_sec:imp_detail}
Optimizing the likelihood with respect to the hyperparameters necessitates the use of an efficient gradient-free algorithm. In this study, we have employed COBYLA~\cite{COBYLA} from the Python library Scipy, which employs default termination criteria related to the trust region size. Since COBYLA is a local search algorithm, we have employed a multi-start technique for enhanced robustness. Our models and their implementation can be accessed in the SMT v2.0 toolbox~\cite{SMT2019, saves2023smt}. In SMT 2.0\footnote{\url{https://smt.readthedocs.io/en/latest/}}, the default number of starting points for COBYLA is set to 10, distributed evenly. We utilize a straightforward noiseless Kriging model with a constant prior. It is important to note that the absolute exponential kernel and the squared exponential kernel behave similarly for categorical variables. The correlation values range between $2.06e-9$ and $0.999999$ for both continuous and categorical hyperparameters. Consequently, we select the constant $\epsilon$ to correspond to a correlation value of $2.06e-9$. The DoE are generated through Latin Hypercube Sampling (LHS)~\cite{LHS}, and the validation sets consist of evenly spaced points.

For BO without constraint, we are using the EGO method of SMT 2.0 with the aforementioned GP models from the same toolbox. For BO under constraints, we are using the SEGOMOE method~\cite{bartoli:hal-02149236} with the mean criterion for the  metamodels of constraints~\cite{SEGO-UTB}. We are using the same GP models for both objective and constraints. The optimization of the WB2s infill criterion~\cite{bartoli:hal-02149236} is done using SNOPT~\cite{gill2005snopt}.
To compare BO with, we used a Multi-Objective Evolutionary Algorithm (MOEA)~\cite{Petrowski2017} named the Non-dominated Sorting Genetic Algorithm II (NSGA-II) \cite{nsga2} due to its low configuration effort and high performance. The \gls{NSGA2} algorithm used is the implementation from the toolbox pymoo~\cite{pymoo} with the default parameters (probability of crossover = 1, eta = 3). Pareto fronts are not relevant in our study as we are considering single-objective optimization. We note that, although NSGA-II is designed for multi-objective optimization, for the purpose of establishing a baseline reference (for comparison), we have used NSGA-II to solve our mono-objective optimization problem. In fact, to the best of our knowledge, NSGA-II is the only open-source optimization solver available for addressing mixed-variable constrained optimization.

In this chapter, all results are obtained using an Intel® Xeon® CPU E5-2650 v4 @ 2.20 GHz core and 128 GB of memory with a Broadwell-generation processor front-end and a compute node of a peak power of 844 GFlops. 

Note that when using KPLS, the GP models built based on the EHH and HH kernels demonstrate comparable performance, with a slight advantage for HH in our numerical tests. For this reason, and to enhance the readability of the numerical section, we have decided to report only results on the HH kernel. Nevertheless, the method related to the use of KPLS within the EHH kernel (i.e., EHH-PLS) is available in the SMT 2.0 toolbox.

\subsection{Surrogate modeling}
\label{SMO_subsec:model_val}

In this section, we validate our model on both a state-of-the-art modeling problem in Section~\ref{SMO_subsec:catcos} and on an structural cantilever beam problem in Section~\ref{SMO_subsec:beam}. More precisely, the matrix based PLS model is compared with literature model0s like  GD, CR or HH. This section shows that the PLS information can capture the shape of the correlation matrix between the various levels of a categorical variable.

\subsubsection{Analytic validation on a categorical cosine problem ($n= 1$, $m=0$, $l=1$ and $L_1=13$)}
\label{SMO_subsec:catcos}

In this section, we investigate the categorical cosine problem, as outlined in~\cite{Roustant}, to showcase the behavior of the proposed kernels. The black-box function, denoted as $f$, relies on both a continuous variable within the range of $[0,1]$ and on a categorical variable with 13 distinct levels. 
Consequently, the relaxed dimension (i.e., the number of hyperparameters) is 14 for the construction of a continuous GP with CR and the most general GP with our new relaxation is of dimension 79.
Appendix~\ref{subsec:cosine} provides a detailed description of this function.
A given point of the \gls{DOE} is denoted as $w= (x,c)$, where $x$ represents the continuous variable and $c$ represents the categorical variable. This work aims at modeling the interplay between the various variables (together with their respective levels) as well as their impact on the targeted function. Notably, the categorical variable exhibits two distinct groups of curves, each comprising a subset of the 13 levels. The first group encompasses levels 1 to 9, while the second group consists of levels 10 to 13. Within each group, we observe strong positive correlations, implying that variables within the same group exhibit a similar behavior. Conversely, strong negative correlations manifest between the two groups, indicating distinct behavior and characteristics between them.  
In this example, the number of relaxed dimensions for continuous relaxation is 14. A \gls{LHS} \gls{DOE} with 98 points ($14\times 7$, if 7 points per dimension are considered) is chosen to built the GP models. On this test case, the number of hyperparameters to optimize is therefore $2$ for GD and HH with 2D PLS, $14$ for CR and $79$ for HH as indicated in~\tabref{SMO_tab:resRoustant}. 

\begin{figure}[H]
\begin{center}
\vspace{-0.25cm}
        \subfloat[GD kernel (2 hyperparameters: 1 cat. and 1 cont.)]{
      \centering
		\includegraphics[   height=4.4cm, width=5cm]{images/corr_gower_R.jpg}
    }
\vspace{-0.25cm}
      \subfloat[HH with PLS kernel (2 hyperparameters: 1 cat. and 1 cont.)]{
      \centering
		\includegraphics[   height=4.4cm, width=5cm]{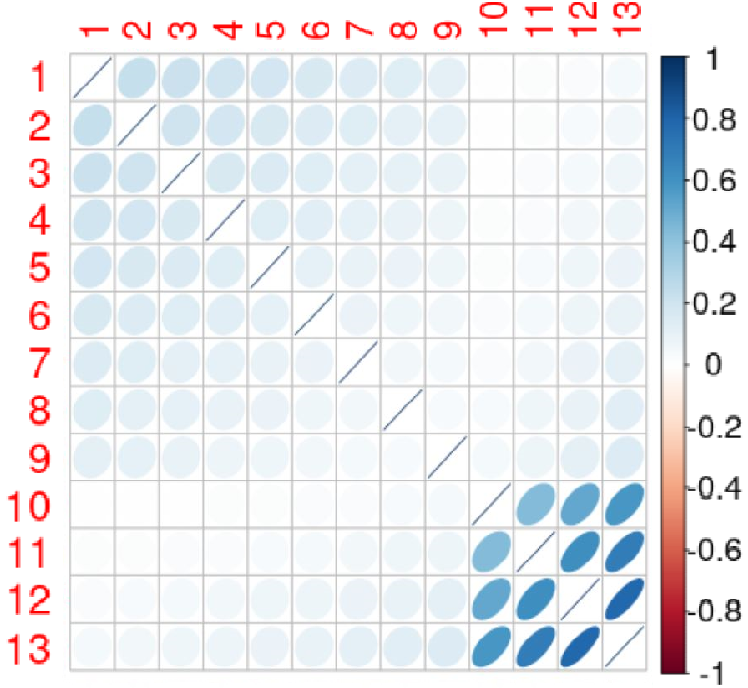}
    	\includegraphics[   height=4.4cm, width=8.6cm]{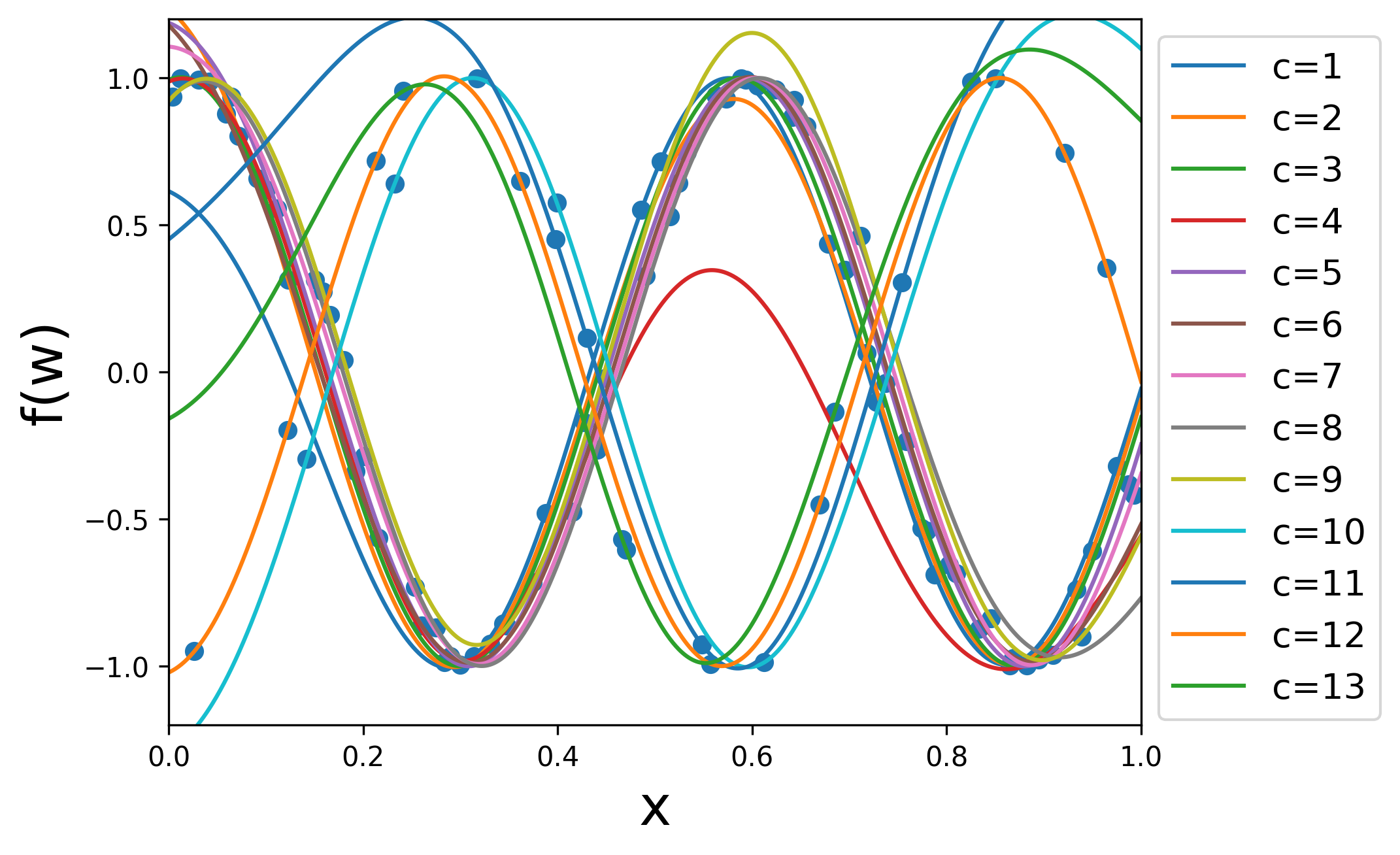}
    }
\vspace{-0.25cm}
      \subfloat[CR kernel (14 hyperparameters: 13 cat. and 1 cont.)]{
      \centering
		\includegraphics[   height=4.4cm, width=5cm]{images/corr_cont_R.jpg}
    	\includegraphics[   height=4.4cm, width=8.6cm]{images/CR_50_Roustant_curves.jpg}
    }

\vspace{-0.25cm}
    \subfloat[HH kernel (79 hyperparameters: 78 cat. and 1 cont.)]{
  \centering
	\includegraphics[   height=4.4cm, width=5cm]{images/corr_plot_HH_v2.JPG}
	\includegraphics[   height=4.4cm, width=8.6cm]{images/curves_hh_v2.JPG} 
     }
\caption{Correlation matrices and associated predictions on the cosine problem using a DoE of 98 points. }
\label{SMO_Roustant_comp}
\end{center} 
\end{figure}

For GD, CR, HH and HH with PLS, the associated mean posterior models are shown on~\figref{SMO_Roustant_comp} on the right and the estimated correlation matrices $R_i=R_1^{cat}$ are displayed on the left.
The latter matrices can be interpreted as such: 
for two given levels $\{\ell_r^1,\ell_s^1\}$, the correlation term $[R_1]_{\ell_r^1,\ell_s^1}$ is in blue for correlation value close to 1, in white for correlations close to 0 and in red for value close to -1; moreover the thinner the ellipse, the higher the correlation and we can see that the correlation between a level and itself is always 1. 
\vspace{-0.1cm}

At first glance, one can see, on~\figref{SMO_Roustant_comp}, that the predicted values remain properly within the interval $[-1,1]$ only with the HH kernel and looked badly estimated with GD. To quantify this assumption, we compute the Root Mean Square Error (RMSE) and Predictive Variance Adequacy (PVA)~\cite{PVA} for every model as 
\begin{equation*}
\label{SMO_eq:RMSE}
\mbox{RMSE}= \sqrt{\underset{i=1}{\overset{N}{\sum}} \frac{1}{N} \left( \hat{f}(w_i) - f(w_i) \right)^2} \quad \mbox{and} \quad 
\mbox{PVA}= \log \left( {\underset{i=1}{\overset{N}{\sum}} \frac{1}{N} \frac{\left( \hat{f}(w_i) - f(w_i) \right)^2}{ [\sigma^f(w_i)]^2 }  } \right),
\end{equation*}
where $N$ is the size of the validation set, $\hat{f}(w_i)$ and $\sigma^f(w_i) $ are the mean and standard deviation predictions of our GP model at a point $w_i$, $f(w_i)$ is the associated true value and the validation set consists of $N=13000$ evenly spaced points (see Appendix~\ref{subsec:cosine}). The values are reported in~\tabref{SMO_tab:resRoustant} and show that the PVA is constant indicating that the estimation of the variance is proportional to the RMSE. Nevertheless, the PVA is smaller for the EHH and HH kernels because the optimization has not converged yet after 887 seconds (79 parameters being hard to optimize in that case). However, a longer run gives a RMSE of 1.280 and a PVA of 21.95 for HH, which, once again is around the other PVA values. Note that the RMSE obtained with HH and EHH are significantly smaller than the errors obtained with all other methods, even with an incomplete optimization.

In~\tabref{SMO_tab:resRoustant}, we also included the results obtained for different sizes of the PLS reduced matrix, namely $2 \times 2,3 \times 3,4 \times 4$ and $5 \times 5$. As expected, these results show that the more complex the model, the smaller the RMSE. Moreover, with the same number of hyperparameters the HH 2D PLS outperforms GD in terms of accuracy and retrieves part of the correlation matrix structure. Also, the HH 2D PLS model can lead to almost similar performance as CR with significantly less hyperparameters. 
Indeed, PLS reduces slightly the accuracy of the model because KPLS is a simplified model optimized in a small subspace but it reduces the run time by a factor of around 10 on this categorical cosine test case. Both the CR and HH models have been built for illustration purposes but are intractable in time for real test cases and the latter is the reason for developing reduced order model in this chapter.
Notwithstanding, in terms of computational costs, computing PLS for matrices could be time-consuming (130 seconds for 2D PLS versus 2.5 seconds for GD). However, this cost is a fixed cost making our method of particular interest for large datasets and high number of dimensions.

We also include the metrics for CR+PLS but this method combines all hyperparameters and does a unique PLS in a continuous space, therefore, this method cannot be used to retrieve the correlation matrix and to study the categorical variable, it has no explainability. Nevertheless, in terms of predictive power, this method is associated with a computational time similar to GD for a better predictive power which makes it a good method for quick prediction on uncertain zones as in the context of BO. For the method developed in this chapter, we tested several PLS components but these components correspond to a small correlation matrix and therefore should be of the form $\frac{n(n+1)}{2}, n \in \mathbb{N}$. We computed the models with the first number of PLS components ($1,3,6,10$) as indicated in~\tabref{SMO_tab:resRoustant} but the predictive gain is not significant, especially given the increase in run time to build the model when the number of components increases. In particular, after 6 components (4x4 matrix), we have a cost comparable with the full model, and, even if these 6 hyperparameters are optimized completely, the prediction is still rough whereas the incomplete optimization of the HH model gives a significantly better prediction still. For that reason, in the following experiments and comparisons, we will stick to the 1 component (2x2 matrix) PLS to retrieve the matrix because this method is less expensive and almost as efficient as its variants with more parameters. As mentioned before, the HH run has not converged after 887 seconds but we used the same internal parameters in the SMT 2.0 software to compare the various models fairly. To be as exhaustive as possible, we also added the results for the EHH kernel, more details about the models without PLS are available in~\cite{Mixed_Paul}.

\begin{table}[!h]
\centering
 \caption{Kernel comparison for the cosine test case in terms of number of hyperparameters, time, RMSE and PVA metrics.}
 \label{SMO_tab:resRoustant}
\small
\begin{tabular}{cccccc}
\hline
  \textbf{Kernel} &  \# of Hyperparam. & run time (s) & RMSE & PVA   \\
  \hline 
  \textbf{GD} &  2 & 2.5 & 30.079&  21.99  \\   
 \hdashline 
  \textbf{CR} &  14 & 19& 22.347 & 23.04 \\
  \textbf{CR+PLS} &   2 & 2.8 & 26.376  & 21.87  \\

  \hdashline 
  \textbf{HH}  & 79 & 887 & 5.330   & 15.34  \\
 
 \textbf{HH+PLS(2x2)} & 2 & 130 & 26.087  & 21.86  \\
  \textbf{HH+PLS(3x3)} & 4 & 326 & 25.504  & 21.95  \\
  \textbf{HH+PLS(4x4)} & 7 & 819& 23.01 & 22.13  \\
  \textbf{HH+PLS(5x5)} & 11 & 1787& 23.04 & 22.13  \\
  \hdashline 
  \textbf{EHH}  & 79 & 959 & 6.858   & 15.46  \\

\hline
\end{tabular}
\end{table}

\subsubsection{Structural modelling: a cantilever beam bending problem ($n=2$, $ m=0$, $l=1$ and $L_1=12$) }
\label{SMO_subsec:beam}

\begin{figure}[H]
\begin{center}

	\subfloat[GD kernel (3 hyperparameters: 1 cat. and 2 cont.)]{
      \centering 
		\includegraphics[   height=5.5cm, width=5cm]{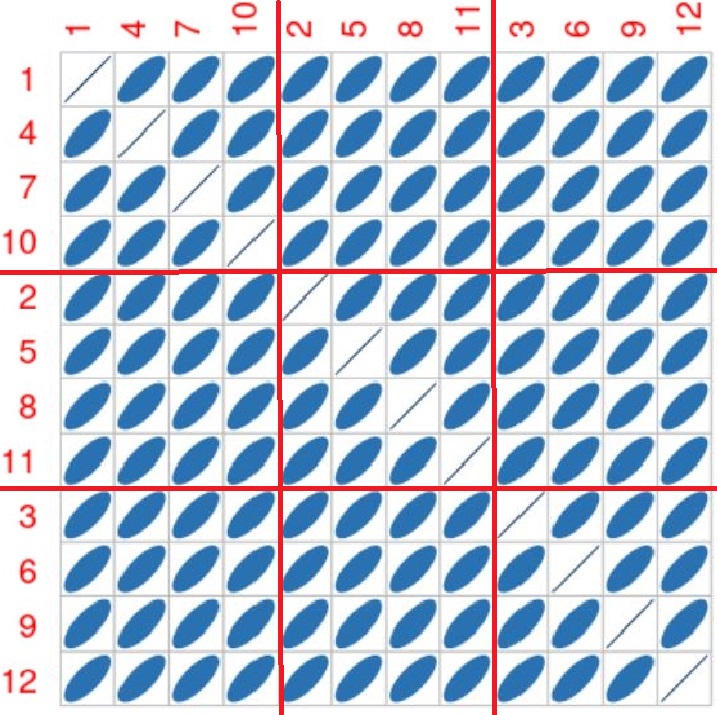} \label{SMO_corr_canti_gower}
     }  
              \hspace{.6 cm}
            \subfloat[HH with PLS kernel (3  hyperparameters: 1 cat. and 2 cont.)]{
      \centering
		\includegraphics[   height=5.5cm, width=5cm]{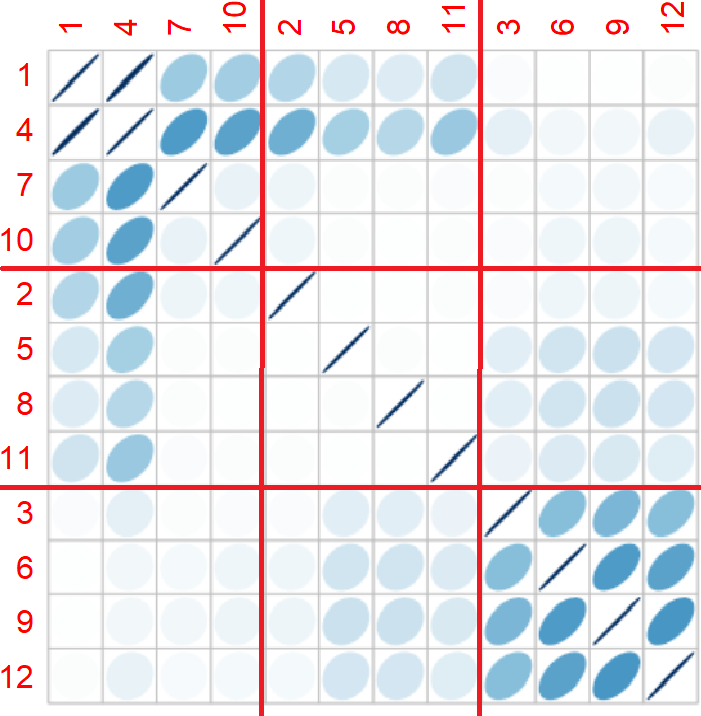}\label{SMO_corr_canti_hh_pls}
     }  
     \hspace{.1 cm}
      \centering
		\includegraphics[   height=5.5cm, width=1cm]{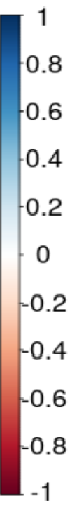}
     
        \subfloat[CR kernel (12 hyperparameters: 10 cat. and 2 cont.)]{
      \centering
		\includegraphics[   height=5.5cm, width=5cm]{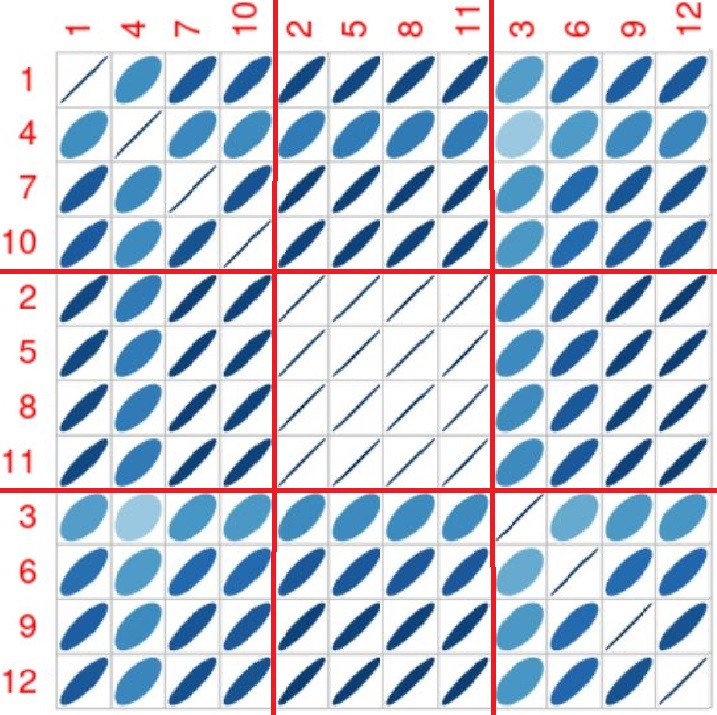} \label{SMO_corr_canti_cr}
     }  
        \hspace{.6 cm}
        \subfloat[HH kernel (66 hyperparameters: 64 cat. and 2 cont.)]{
      \centering
		\includegraphics[   height=5.5cm, width=5cm]{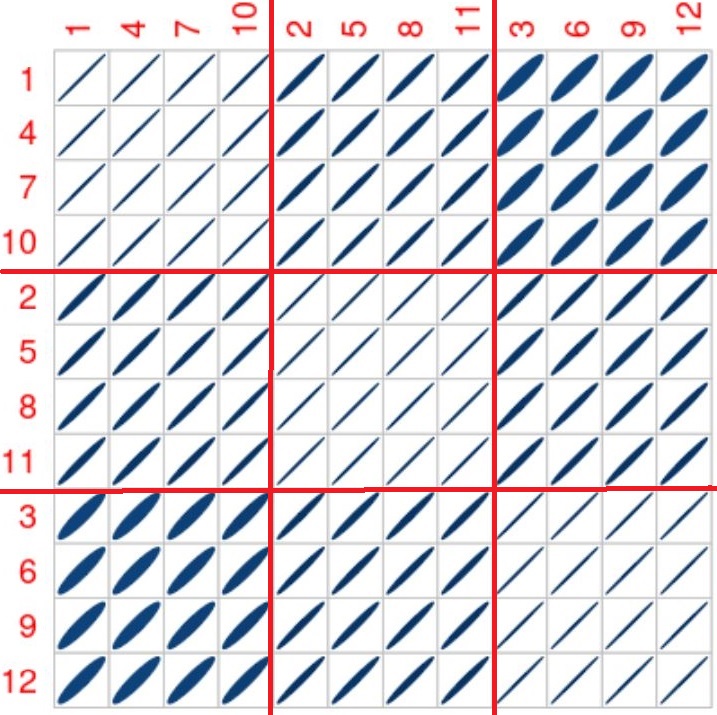}\label{SMO_corr_canti_ehh}
     }  
    \hspace{.1 cm}
      \centering
		\includegraphics[   height=5.5cm, width=1cm]{images/legend_pos_neg.jpg}
     
\caption{Correlation matrix $R_1^{cat}$  using different choices for $\Theta_1$ for the categorical variable $\tilde{I}$ from the cantilever beam problem.}
\label{SMO_corr_Cantilever}
\end{center} 
\end{figure}

A classic engineering problem frequently employed for model validation is the beam bending problem in its linear elasticity range~\cite{Roustant, Cheng2015TrustRB}. This problem serves as an illustrative example and involves a cantilever beam subjected to a load applied at its free end, denoted as $F$. The specific setup of the problem is depicted in~\figref{fig:beam_SMO}. In accordance with the findings presented in~\cite{ Cheng2015TrustRB}, the Young's modulus for the material is determined to be $E=200$ GPa, and a load of $F=50$ kN has been chosen for the analysis. Furthermore, following the methodology outlined in~\cite{Roustant}, a total of 12 potential cross-sections can be used for the beam. These cross-sections encompass four distinct shapes (square, circle, I and star), each with the possibility of being either full, thick, or hollow, as visually depicted in~\figref{fig:beam_SMO_shape}. For a given cross-section, which consists of a specific shape and thickness, its size is determined by the surface area denoted as $S$. Additionally, each cross-section is associated with a normalized moment of inertia $\tilde{I}$ around the neutral axis, representing a latent variable connected to the beam's shape~\cite{oune2021latent}. 
\begin{figure}[ht]
\centering

\vspace{-5pt}
\captionsetup{justification=raggedright,singlelinecheck=false}
\subfloat[Bending problem.]{
\begin{tikzpicture}
    \hspace{-3pt}
    \point{origin}{-0.75}{-0.25};
    \point{begin}{0}{0};
    \point{end}{5}{0};
    \point{end_bot}{4.99}{-0.9};
    \point{end_up}{5}{0.5};
    \beam{2}{begin}{end};
    \support{3}{begin}[-90];
    \load{1}{end}[90]   ;
    \notation{1}{end_up}{$F=50kN$};

     \draw[<->] (end) -- (end_bot) node[midway, right] {$\delta$} ;
     \draw[<->] (0,0.5) -- (5,0.5) node[midway, above] {L};
     
    \draw
      [-, ultra thick] (begin) .. controls (1.5, +.01) and (2.5, -.15) .. (4.93, -0.9)
      [-, ultra thick] (begin) .. controls (1.5, +.01) and (2.5, -.2) .. (4.85, -1.5)
      [-, ultra thick] (begin) .. controls (1.5, +.01) and (2.5, -.4)   .. (4.78, -1.9);
  \end{tikzpicture}
\label{fig:beam_SMO}   
}
\subfloat[Possible cross-section shapes.]{
\centering
\begin{tikzpicture}
\hspace{-50pt}


\tstar{0.25}{0.5}{6}{0}{thick,fill=yellow,xshift=+3.6cm,yshift= -1.2cm}
\tstar{0.14}{0.28}{6}{0}{thick,fill=white,xshift=+3.6cm,yshift= -1.2cm}

\tstar{0.25}{0.5}{6}{0}{thick,fill=yellow,xshift=+2.4cm,yshift= -1.2cm}
\tstar{0.08}{0.16}{6}{0}{thick,fill=white,xshift=+2.4cm,yshift= -1.2cm}

\tstar{0.25}{0.5}{6}{0}{thick,fill=yellow,xshift=+1.2cm,yshift= -1.2cm}

\fill[green,even odd rule] (3.6,0) circle (0.5) (3.6,0) circle (0.33);
\draw (3.6,0) circle (0.5) ;
\draw (3.6,0) circle (0.33) 
; 
\fill[green,even odd rule] (2.4,0) circle (0.5)(2.4,0) circle (0.17);
\draw (2.4,0) circle (0.5) ;
\draw (2.4,0) circle (0.17) 
; 
\fill[green,even odd rule] (1.2,0) circle (0.5) ;
\draw (1.2,0) circle (0.5) ;

\def\pos{-2.4}
\fill[blue,even odd rule]  (\pos-0.5,-1.7+1.2) -- (\pos-0.5,-0.7+1.2) -- (\pos+0.5,-0.7+1.2) -- (\pos+0.5,-1.7+1.2) -- cycle ;

\def\pos{-1.2}
\fill[blue,even odd rule]  (\pos-0.5,-1.7+1.2) -- (\pos-0.5,-0.7+1.2) -- (\pos+0.5,-0.7+1.2) -- (\pos+0.5,-1.7+1.2) -- cycle   (\pos-0.25,-1.45+1.2) -- (\pos-0.25,-0.95+1.2) -- (\pos+0.25,-0.95+1.2) -- (\pos+0.25,-1.45+1.2) -- cycle ;

\def\pos{0}
\fill[blue,even odd rule]  (\pos-0.5,-1.7+1.2) -- (\pos-0.5,-0.7+1.2) -- (\pos+0.5,-0.7+1.2) -- (\pos+0.5,-1.7+1.2) -- cycle   (\pos-0.35,-1.55+1.2) -- (\pos-0.35,-0.85+1.2) -- (\pos+0.35,-0.85+1.2) -- (\pos+0.35,-1.55+1.2) -- cycle ;

\def\pos{-2.4}
\fill[black] (\pos-0.24,-0.5-1.2) -- (\pos-0.24,0.5-1.2) -- (\pos+0.24,0.5-1.2)  -- (\pos+0.24,-0.5-1.2)   -- cycle ;
\fill[black] (\pos-0.5,-0.5-1.2) -- (\pos-0.5,-0.18-1.2) -- (\pos+0.5,-0.18-1.2)  -- (\pos+0.5,-0.5-1.2)   -- cycle ; 
\fill[black] (\pos-0.5,0.18-1.2) -- (\pos-0.5,0.5-1.2) -- (\pos+0.5,0.5-1.2)  -- (\pos+0.5,0.18-1.2)   -- cycle ;  

\def\pos{-1.2}
\fill[black] (\pos-0.19,-0.5-1.2) -- (\pos-0.19,0.5-1.2) -- (\pos+0.19,0.5-1.2)  -- (\pos+0.19,-0.5-1.2)   -- cycle ;
\fill[black] (\pos-0.5,-0.5-1.2) -- (\pos-0.5,-0.25-1.2) -- (\pos+0.5,-0.25-1.2)  -- (\pos+0.5,-0.5-1.2)   -- cycle ; 
\fill[black] (\pos-0.5,0.25-1.2) -- (\pos-0.5,0.5-1.2) -- (\pos+0.5,0.5-1.2)  -- (\pos+0.5,0.25-1.2)   -- cycle ;  

\def\pos{0}
\fill[black] (\pos-0.14,-0.5-1.2) -- (\pos-0.14,0.5-1.2) -- (\pos+0.14,0.5-1.2)  -- (\pos+0.14,-0.5-1.2)   -- cycle ;
\fill[black] (\pos-0.5,-0.5-1.2) -- (\pos-0.5,-0.32-1.2) -- (\pos+0.5,-0.32-1.2)  -- (\pos+0.5,-0.5-1.2)   -- cycle ; 
\fill[black] (\pos-0.5,0.32-1.2) -- (\pos-0.5,0.5-1.2) -- (\pos+0.5,0.5-1.2)  -- (\pos+0.5,0.32-1.2)   -- cycle ;  

\point{un}{-2.15}{-0.80};
\notation{1}{un}{\tiny 1};
\point{deux}{-2.15+1.2}{-0.80};
\notation{1}{deux}{\tiny 2};
\point{trois}{-2.15+2.4}{-0.80};
\notation{1}{trois}{\tiny 3};

\point{quatre}{-2.15+3.6}{-0.80};
\notation{1}{quatre}{\tiny 4};
\point{cinq}{-2.15+4.8}{-0.80};
\notation{1}{cinq}{\tiny 5};
\point{six}{-2.15+6}{-0.80};
\notation{1}{six}{\tiny 6};

\point{sept}{-2.15}{-0.80-1.2};
\notation{1}{sept}{\tiny 7};
\point{huit}{-2.15+1.2}{-0.80-1.2};
\notation{1}{huit}{\tiny 8};
\point{neuf}{-2.15+2.4}{-0.80-1.2};
\notation{1}{neuf}{\tiny 9};

\point{dix}{-2.15+3.6}{-0.80-1.2};
\notation{1}{dix}{\tiny 10};
\point{onze}{-2.15+4.8}{-0.80-1.2};
\notation{1}{onze}{\tiny 11};
\point{douze}{-2.15+6}{-0.80-1.2};
\notation{1}{douze}{\tiny 12};

  \end{tikzpicture}
  \label{fig:beam_SMO_shape}    
}

\captionsetup{justification=centering,singlelinecheck=false}
\caption{Cantilever beam problem~\cite[Figure 6]{Mixed_Paul}.}
\end{figure}
Hence, the problem at hand involves modeling with two continuous variables: the length $L$, ranging from 10 to 20 meters, and the surface area $S$, ranging from 1 to 2 square meters. Additionally, there is a categorical variable, $\tilde{I}$, with 12 levels representing the various cross-section options available. The tip deflection, at the free end, $\delta$ is given by $$ \delta = f( \tilde{I}, L,S) = \frac{F}{3E} \frac{L^3}{S^2\tilde{I}}. $$
As a result,  the relaxed dimension used to construct the GP model using the CR method is 14, while the relaxed dimension for the most general GP model employing the HH method is 68.
To compare our models, we draw a 98 point \gls{LHS} as training set and the validation set is a grid of $12\times30\times30=10800$ points. For the four models GD, CR, HH and HH with PLS, the correlation matrix associated to every model are drawn in~\figref{SMO_corr_Cantilever} showing the predicted correlations between the available cross-sections. We recall that these matrices can be interpreted as such:
for two given levels $\{\ell_r^1,\ell_s^1\}$, the correlation term $[R_1]_{\ell_r^1,\ell_s^1}$ is in blue for correlation value close to 1, in white for correlations close to 0 and the thinner the ellipse, the higher the correlation.

The models are summarized in~\tabref{SMO_tab:resCantilever} indicating the complexity of each model and the information that could be recovered from it. For the HH kernel, the indicated computational time corresponds to the duration required to fully converge all 68 hyperparameters. In fact, the computational cost and difficulty to optimize the likelihood in spaces of dimension superior to ten is the biggest limitations of HH and such exhaustive kernels. This increase in difficulty to converge and associated computational cost is one of the main motivations for our PLS method and for simpler models.
As expected, we have 3 groups of 4 shapes depending on their respective thickness (respectively, the full levels \{1,4,7,10\}, the medium levels \{2,5,8,11\}, and the hollow levels \{3,6,9,12\}). The more the thickness is similar, the higher the correlation: the thickness has more impact than the shape of the cross-section on the tip deflection. However, given the database, two points with similar $L$ and $S$ values will have similar output whatever the cross-section. The effect of the cross-section on the output is always the same (in the form of $\frac{1}{\tilde{I}}$) leading to an high correlation after maximizing the likelihood. 
%
%
\resizebox{1.0\textwidth}{!}{%
\begin{threeparttable}[!ht]
\centering
\caption{Results of the cantilever beam models.}
\begin{tabular}{ccccc}
\hline
\textbf{Categorical kernel} &
{Identified clusters} & 
\# hyps.  &$\ $ Time (s) & {RMSE (cm)} \\
\hline 
 \textbf{GD}   &
 -- &
 3 
 & 8
 & {1.3858}

\\  
\textbf{HH+PLS}   & 
\textbf{Hollow} cluster and \textbf{Full} cluster (partly) &
3 & 38
& {1.2989}

\\
\textbf{CR}  &
\textbf{Medium} cluster & 
14  
& 89
& {1.1604}

\\
\textbf{HH}   &
\textbf{Full}, \textbf{Medium} and \textbf{Hollow} & 
68  
& 2128 
& {0.1247}

\\
\hline
\end{tabular}
\label{SMO_tab:resCantilever}
\end{threeparttable}
}
For the GD model in~\figref{SMO_corr_canti_gower}, there is only one mean positive correlation value, therefore no structural information can be extracted from this unique value. On the contrary, the HH model is the most general one and can model every cross-correlation value independently from the others and, in~\figref{SMO_corr_canti_ehh}, we can distinguish the three groups of four shapes, as expected, because the shapes of the cross-sections are not significant in comparison with their thickness.  Concerning the CR model, its structure favors the medium group, that is well represented, but the thick and hollow shapes are not distinguished in~\figref{SMO_corr_canti_cr}. 
To finish with, the HH+PLS model in~\figref{SMO_corr_canti_hh_pls} captures well the hollow shapes correlations and managed to capture half the full shapes correlations. More precisely, HH+PLS captured, with only one categorical hyperparameter (3 in total), the correlation between the full square and the full circle cross-sections (indexed 1 and 4) but failed to capture the correlation with full I and full star shapes (indexed 7 and 10). Consequently, with only one categorical hyperparameter, our model performs really well and is able to reconstruct structure from the data, thus outperforming the GD model for the same computational cost.  \bigbreak

To conclude, this section showed the capability of our PLS model to capture structures in the data while using only a small number of hyperparameters. Our model eases the optimization of the likelihood function and reduces the computational cost associated with the GP surrogate model. Notably, our work was efficiently applied to a structural modeling problem. However, building the GP surrogate is only part of the total optimization cost, and, in the next section, we show how our GP can be used in the context of high-dimensional MDO for aircraft design.



 \subsection{Surrogate-based optimization: Bayesian optimization}
\label{SMO_sec:BO_val}

Efficient Global Optimization (EGO) is a well-known Bayesian optimizer that relies on GP to find out the optimum of an unconstrained black-box problem that can be evaluated a limited amount of times~\cite{Jones98}.
The workflow of EGO begins with building a first GP model based on an initial DoE, followed by employing an acquisition function to guide the selection of the next point that will be evaluated through the expensive black-box function. 
The most commonly used acquisition function is the expected improvement and, once a new point has been evaluated, the GP model is updated and the selection process repeats with the updated GP. At every step, a new model is built and a new point is evaluated until a maximal budget is reached. 
Hereinafter, we will use the GP models aforementioned to optimize expensive-to-evaluate black-box problems involving mixed integer variables using the EGO algorithm.
Moreover, EGO has been generalized to tackle constrained problems by Sasena \textit{et al.}~\cite{sasena2002exploration} with an algorithm called SEGO and used in the optimizer SEGOMOE~\cite{bartoli:hal-02149236}.

\subsubsection{Analytic validation on a mixed optimization problem ($n=1, \ m=0, \ l=1,$ and $L_1 = 10$)}
\label{SMO_sec:MI-BO}
\begin{figure}[H]
\begin{minipage}[b]{.6\linewidth}
\centering
\hspace{-1.25cm}
\subfloat[Convergence curves: medians of 20 runs.]{
\includegraphics[height=5cm,,width=7cm]
{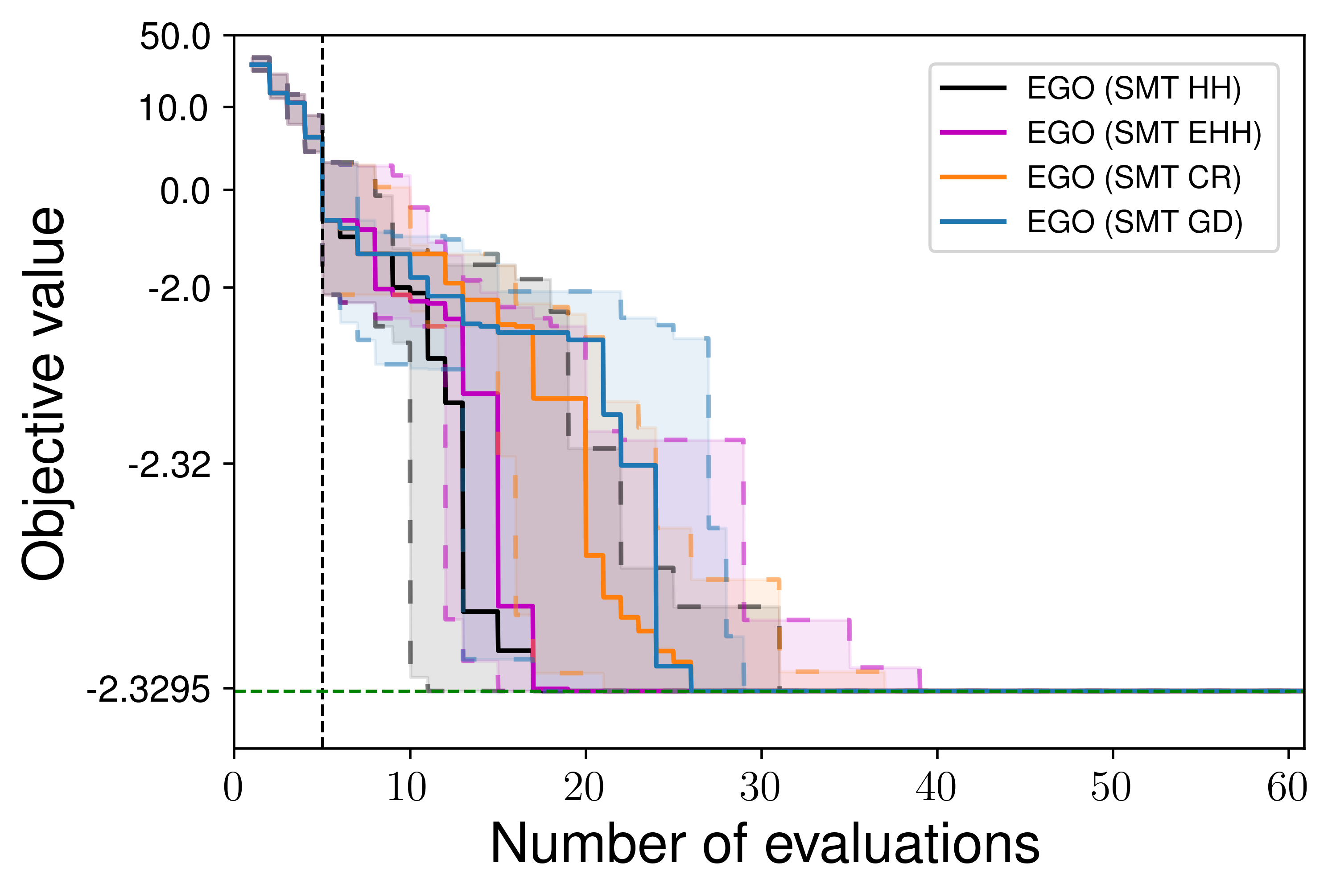}
\label{SMO_convmi}
}
\end{minipage}
\begin{minipage}[b]{.4\linewidth}
\centering 
\hspace{-1.5cm}
\subfloat[Boxplots after 25 evaluations.]{
\includegraphics[height=5cm,width=6.2cm]{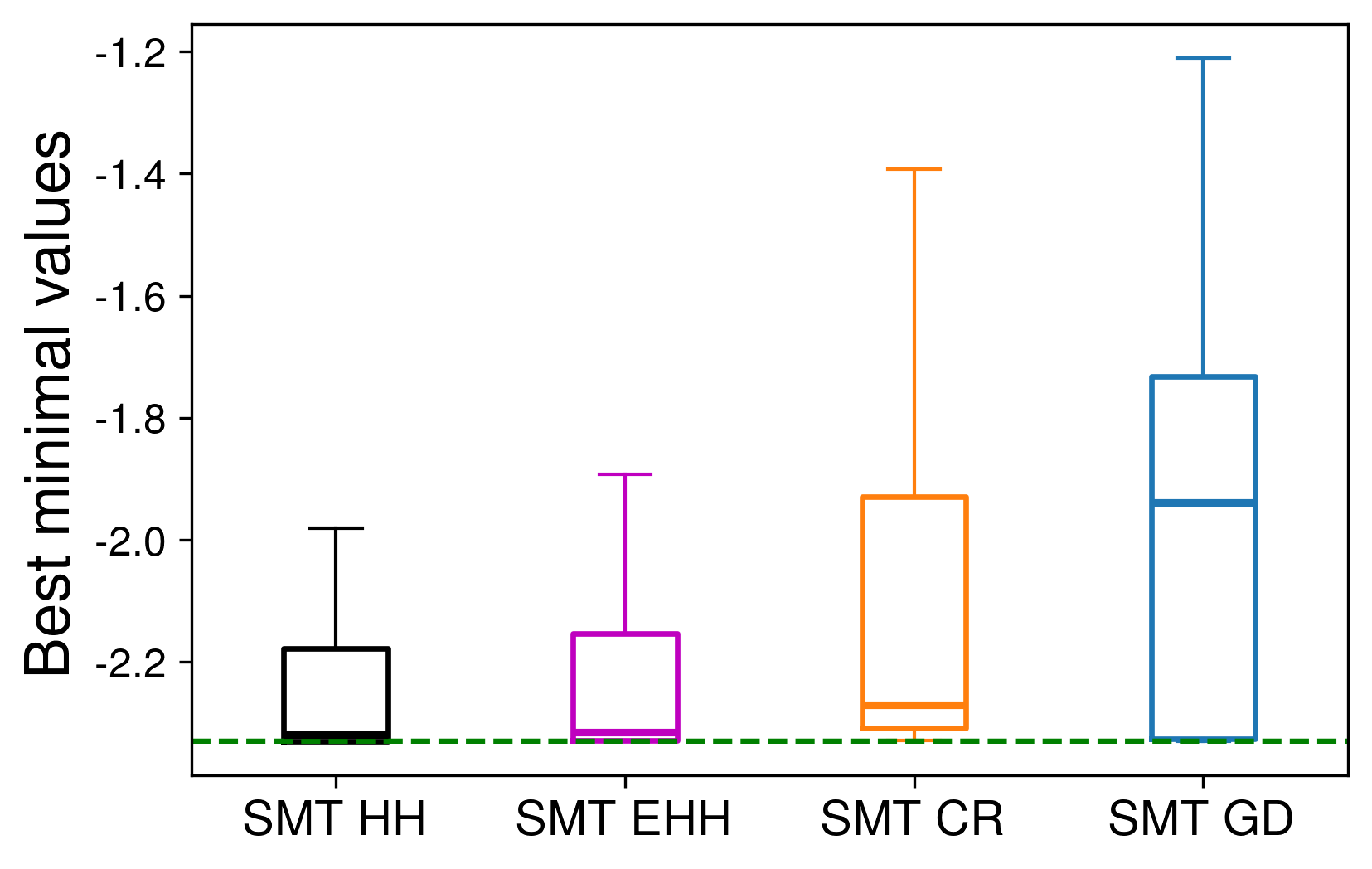}
\label{SMO_mini_mi}
}
\end{minipage}
\caption{Optimization results for the Toy function~\cite{CAT-EGO} for 20 DoE of 5 points.}
\label{SMO_res_optim_mi}
\end{figure}

\begin{figure}[H]
\begin{minipage}[b]{.6\linewidth}
\centering
\hspace{-1.25cm}
\subfloat[Convergence curves: medians of 20 runs.]{
\includegraphics[height=5cm,,width=7cm] {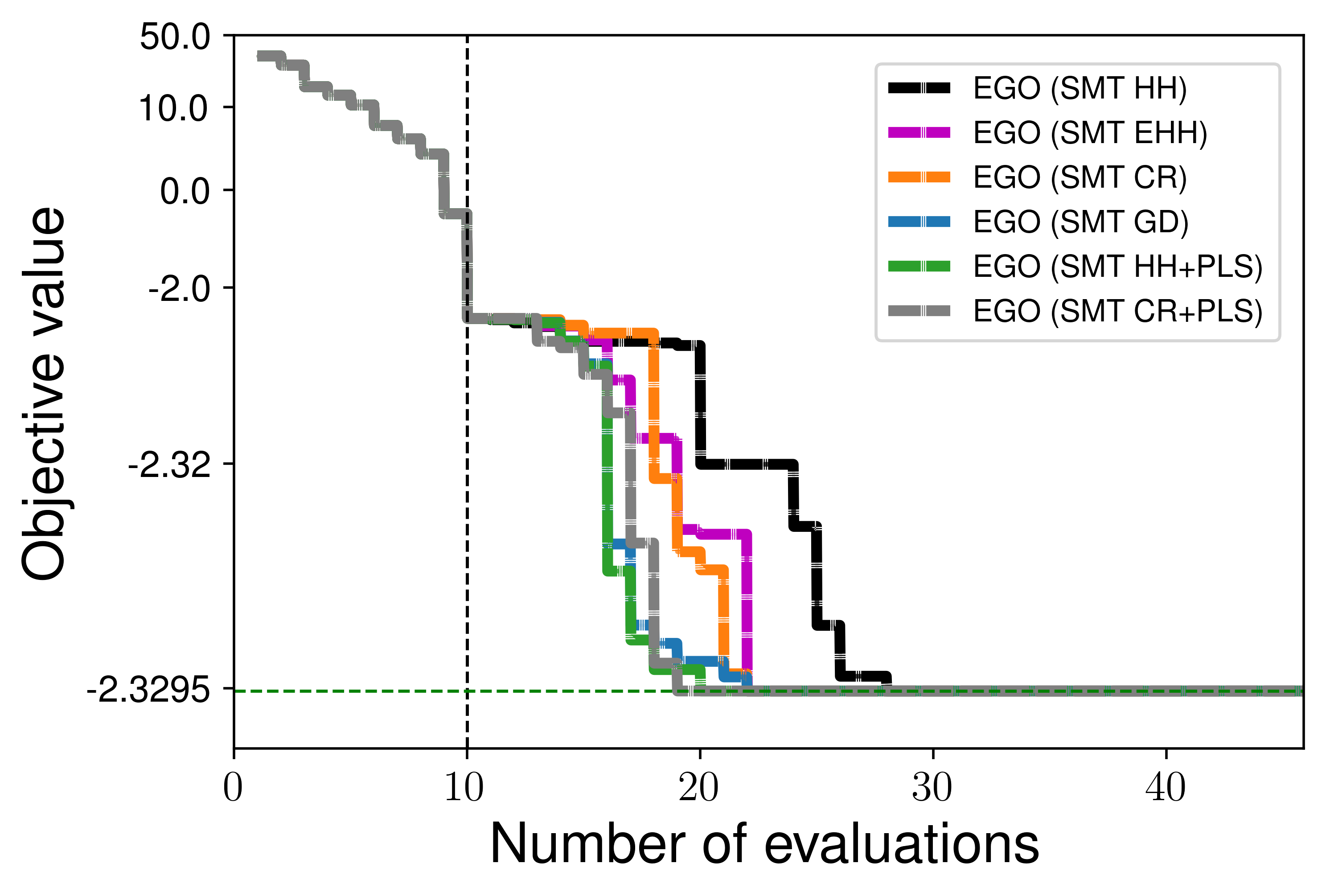}
\label{SMO_convmi2}
}
\end{minipage}
\begin{minipage}[b]{.4\linewidth}
\centering 
\hspace{-1.5cm}
\subfloat[Boxplots after 25 evaluations.]{
\includegraphics[height=5cm,width=6.2cm]{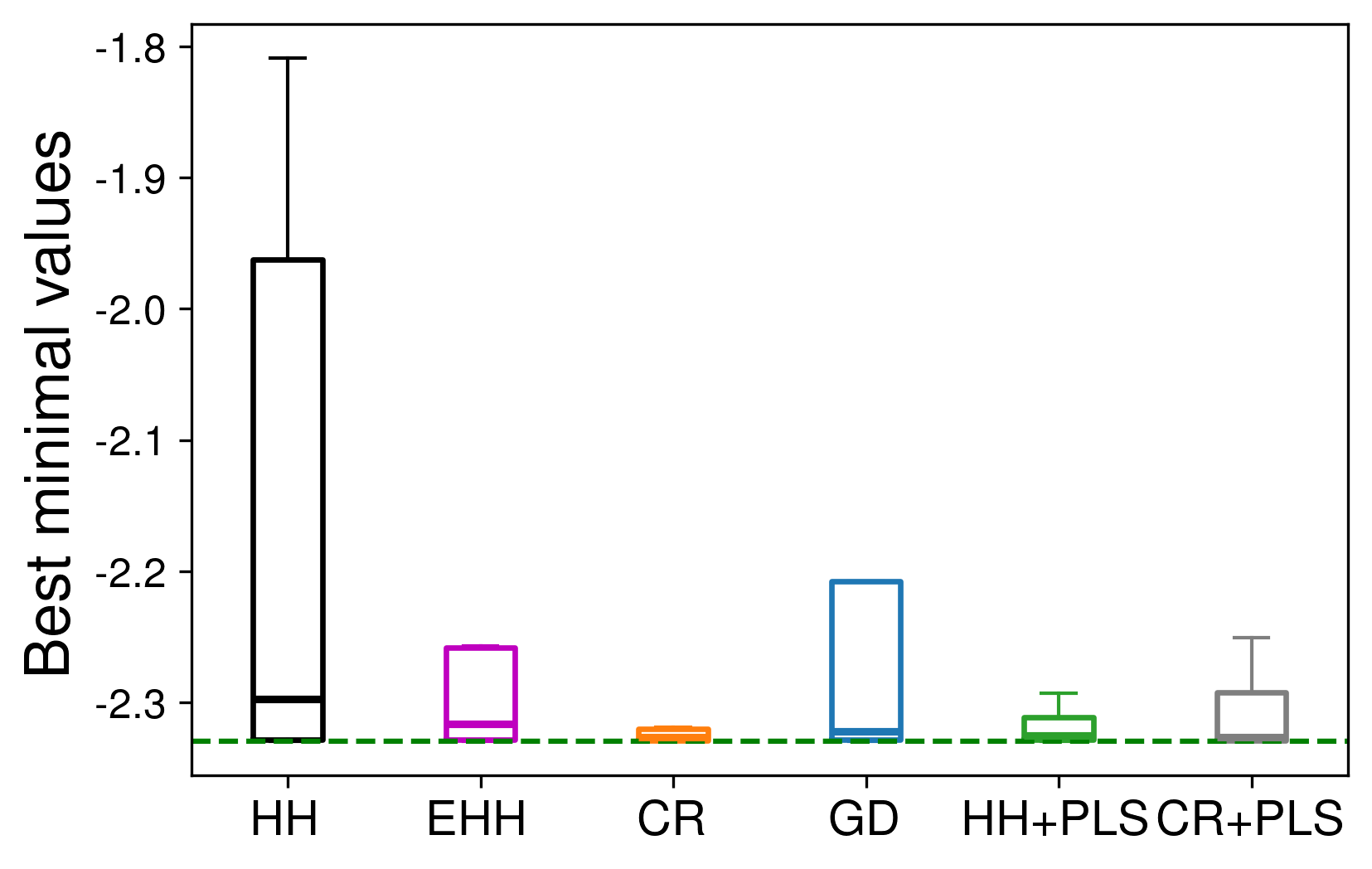}
\label{SMO_mini_mi2}
}
\end{minipage}
\caption{Optimization results for the Toy function~\cite{CAT-EGO} for 20 DoE of 10 points.}
\label{SMO_res_optim_mi2}
\end{figure}

The mixed test case that illustrates BO is a toy test case~\cite{CAT-EGO} detailed in Appendix~\ref{app:Toy}. 
This test case has two variables, one continuous and one categorical with 10 levels.
As a result, the relaxed dimension used to construct the GP model using the CR method is 11, while the relaxed dimension for the most general GP model employing the HH method is 46.
In~\figref{SMO_res_optim_mi} and~\figref{SMO_res_optim_mi2} six GP models are being compared. These six models are four classical GP models, namely GD, CR, EHH and HH and two PLS based models, namely our new model HH+PLS and the previously developed CR+PLS model~\cite{SciTech_cat}. 
In particular, in~\cite{SciTech_cat}, we used this CR+PLS method coupled with a criterion to choose automatically the number of PLS components that gives the best prediction for optimization. This adaptive PLS method for mixed integer has been applied to the MDO of \gls{DRAGON} as detailed in Section~\ref{SMO_subsec:res_AD}.
To assess the performance of our algorithm, we performed 20 runs with different initial DoE sampled by \gls{LHS}.
Every DoE consists of 5 points in~\figref{SMO_res_optim_mi} and of 10 points in~\figref{SMO_res_optim_mi2}. For both experiments, we chose a budget of 55 infill points.
Figure \ref{SMO_convmi} and \figref{SMO_convmi2} plot the convergence curves for the six methods. To visualize the data dispersion, the boxplots of the 20 best solutions after 25 evaluations are plotted in~\figref{SMO_mini_mi} and~\figref{SMO_mini_mi2}. The computational times for every method are indicated in~\tabref{SMO_tab:resToy} for a 5 point DoE and in~\tabref{SMO_tab:resToy2} for a 10 point DoE. 
We note that the overall computational cost is derived by the optimization cost related to the maximization of the infill criterion. In fact, such optimization is often related to the number of the design variables rather than the size of the DoE.
On this test case, our method gives the best results with the 5 point DoE in terms of median convergence speed and dispersion among the 20 DoE. For the 10 point DoE, our method is among the faster together with CR+PLS. However, even if the HH+PLS method has been shown to be efficient for solving this test case, it is still more costly than CR+PLS or GD because the computational cost associated to the reconstruction of the matrix of hyperparameters is significant. 
Nevertheless, it is a method based only on two hyperparameters (one categorical and one continuous) making it around 20 times easier to optimize than HH or EHH and 3 times faster for better performance. In any case, using a 5 point DoE is slightly more efficient than using a 10 point DoE because BO is known to perform better with a smaller DoE for a given budget of evaluations~\cite{le2021revisiting}. But this effect is not significant for methods that use PLS, as PLS  benefits greatly from the initial DoE information to find the most interesting search directions. This explains why PLS methods are performing better with a 10 point DoE than with a 5 point DoE.
\begin{table}[ht]    
\centering
 \caption{Results of the Toy problem optimization (5 point DoE and 55 infill points).}
\small
\small
\begin{tabular}{ccc}
\hline
\textbf{Kernel}& $ \ $ number of hyperparameters  & $\ $optimization duration (s) $\ $    \\
\hline 
\textbf{GD} & 2 &   315 \\
\hdashline 
\textbf{CR} & 11 &  503 \\
\textbf{CR+PLS} & 2 &  320 \\
\hdashline 
\textbf{HH}  & 46 &  1983 \\
\textbf{HH+PLS} & 2 &  646 \\
\hdashline 
\textbf{EHH} & 46 &  2086 \\
\hline
\end{tabular}
\label{SMO_tab:resToy}
\end{table}

\begin{table}[ht]    
\centering
 \caption{Results of the Toy problem optimization (10 point DoE and 55 infill points).}
\small
\small
\begin{tabular}{ccc}
\hline
\textbf{Kernel}& $ \ $ number of hyperparameters  & $\ $optimization duration (s) $\ $    \\
\hline 
\textbf{GD} & 2 &  314  \\
\hdashline 
\textbf{CR} & 11 & 479  \\
\textbf{CR+PLS} & 2 & 326  \\
\hdashline 
\textbf{HH}  & 46 &  2142 \\
\textbf{HH+PLS} & 2 & 662  \\
\hdashline 
\textbf{EHH} & 46 & 2079   \\
\hline
\end{tabular}
\label{SMO_tab:resToy2}
\end{table}

\subsubsection{Multidisciplinary design optimization for a green aircraft 
($n=10, \ m=0, \ l=2, L_1 = 17$ and $L_2 = 2$)}
\label{SMO_subsec:res_AD}
\begin{figure}[H]
\begin{centering}
\includegraphics[height=5cm]{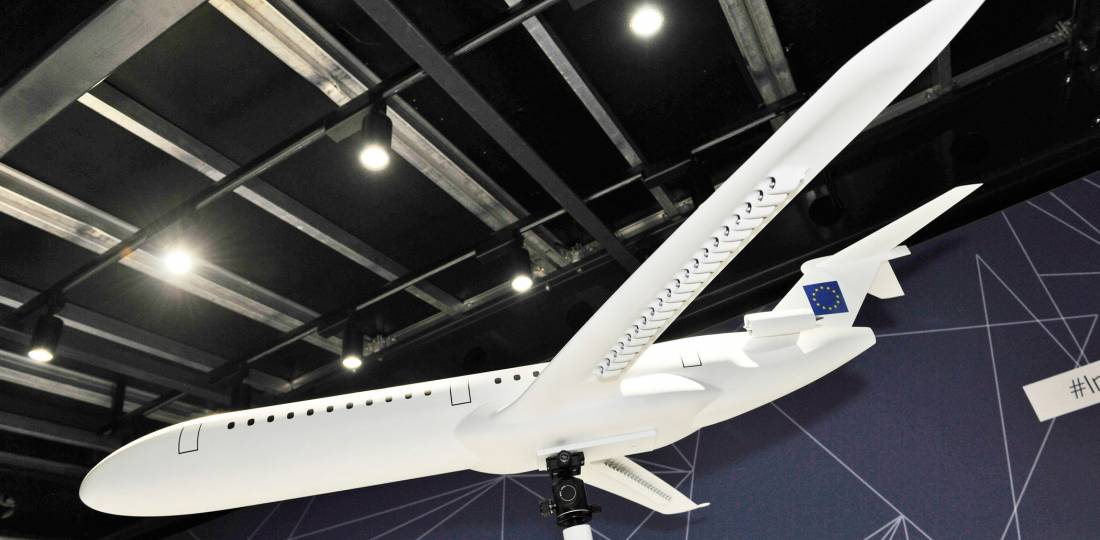}
   \caption{``\texttt{DRAGON}'' aircraft mock-up.}
        \label{SMO_Dragon2020}
\end{centering}
\end{figure}

For the core MDO application, we apply the Future Aircraft Sizing Tool with Overall Aircraft Design (FAST-OAD)~\cite{David_2021} on ``\texttt{DRAGON}'' 
(Distributed fans Research Aircraft with electric Generators by ONera), an innovative aircraft currently under development. The \gls{DRAGON} aircraft concept in~\figref{SMO_Dragon2020} has been introduced by ONERA in 2019~\cite{schmollgruber} within the scope of the European CleanSky 2 program\footnote{\href{SMO_https://www.cleansky.eu/technology-evaluator}{\color{blue}https://www.cleansky.eu/technology-evaluator}} which sets the objective of 30\% reduction of CO2 emissions by 2035 with respect to 2014 state-of-the-art. A first publication in SciTech 2019~\cite{schmollgruber} was followed by an up-to-date report in SciTech 2020~\cite{schmollgruber2}.
In response to this ambitious goal, ONERA introduced a concept for a distributed electric propulsion aircraft that makes significant strides in enhancing fuel efficiency by optimizing propulsive performance. This is realized by an increase in the bypass ratio through a strategic placement of numerous compact electric fans on the wing pressure side, as an alternative to the use of larger turbofans. This design decision effectively resolves the challenges associated with large under-wing turbofans and grants the aircraft the capability to operate at transonic speeds. Consequently, the primary design objective for the  \gls{DRAGON} revolves around accommodating a passenger capacity of 150 individuals and facilitating travel over a range of 2750 Nautical Miles at a speed of Mach 0.78.

\begin{figure}[H]

  \centering 
\includegraphics[  ,height=6cm]{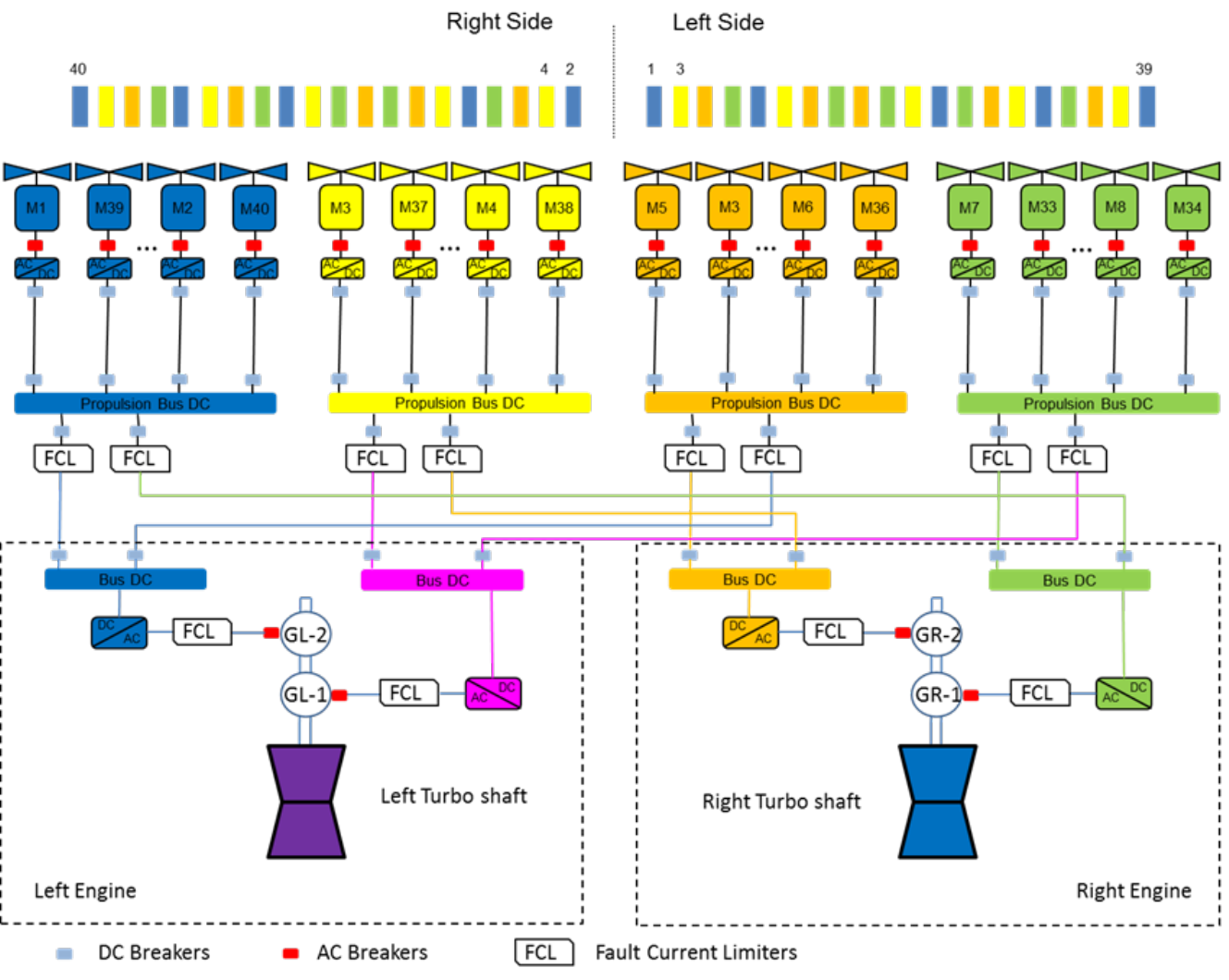}
 \caption{Turboelectric propulsive architecture.}
 \label{SMO_DragonArchitecture}
\label{SMO_Dragon}
\end{figure}

The integration of distributed propulsion in an aircraft introduces certain trade-offs. It necessitates the use of a turbo-electric propulsion system to provide the necessary power to drive the electric fans, which, in turn, contributes to increased intricacy and added weight. Typically, this power is generated onboard by coupling turboshafts to electric generators. The generated electrical power is subsequently transmitted to the electric fans through an electric architecture designed to ensure resilience in the face of potential single component failures. This safety feature is achieved through the deployment of redundant components, as illustrated  in~\figref{SMO_DragonArchitecture}.
The initial setup comprises two turboshafts, four generators, four propulsion buses with cross-feed capabilities, and 40 fans. This configuration was selected during the preliminary study phase due to its compliance with safety standards. Nevertheless, it was not specifically tailored for optimizing weight. Given that the turboelectric propulsion system significantly contributes to the overall weight of the aircraft, there is a specific interest in optimizing this system, especially concerning the number and type of individual components, each characterized by discrete or categorical values.
This time, as the evaluations are expensive, we are doing only 10 runs instead of 20.
Also, to have realistic results, the constraints violation will be forced to be less than $10^{-3}$. From now on, let MAC denote the Mean Average Chord, VT, the Vertical Tail, HT, the Horizontal Tail and TOFL, the Take-off Field Length. The optimizations are realized with SEGOMOE.

In~\cite{saves2021constrained,SciTech_cat}, the \gls{DRAGON} configuration has been already optimized to attain such goals. In particular, in~\cite{SciTech_cat}, the model has been updated based on the results obtained in~\cite{saves2021constrained} that display limitations based on the turboshaft layout with the turbogenerators at the rear of the fuselage so it would be advantageous to locate the turbogenerator below the wing. 
Nonetheless, adopting this approach would impose constraints on the available space allocated for the electric fans. Consequently, it would restrict the maximum achievable propulsive efficiency. To address the inherent trade-off between a lighter propulsion system and an enhanced propulsive efficiency, we integrated the layout as a categorical variable within the optimization problem of~\figref{SMO_tab:dragon}. 
Finally, five constraints are considered on this optimization problem. Among these constraints, the TOFL, climb duration, and the top of climb slope angle exert a significant influence on the design of the hybrid electric propulsion system. Additionally, a portion of the wing trailing edge near the wingtip must be kept unobstructed to accommodate ailerons, thus limiting the available space for the electric fans. Lastly, compliance with wingspan limits is mandated by airport regulations.

To know how optimizing the fuel mass will impact the aircraft design, we are considering the optimization problem described in~\tabref{SMO_tab:dragon}. We can now solve a constrained optimization problem with 10 continuous design variables and 2 categorical variables with 17 and 2 levels respectively, for a total of 12 design variables. For the optimization, this new problem is a hard test case involving 29 relaxed variables and 5 constraints. The definition of the architecture variable is given in~\tabref{SMO_tab:dragon_archi1} and the definition of the turboshaft layout is given in~\tabref{SMO_tab:dragon_archi2}.
For modeling electric architectures, it is more efficient to represent the architectural choices using two integer variables instead of one categorical variable. 
However, taking this approach expands the range of potential architectures beyond the initial 17 configurations. Yet, there are two important constraints to consider for these possible setups.
The first constraint relates to the electrical connections between components: ensuring a certified electric architecture is crucial, and figuring out how to connect, for example, 8 motors to 6 generators is not straightforward.
The second constraint is connected to the distributed propulsion system, especially the numerous propellers. Managing this system involves addressing a substantial number of potential failures in the electro-mechanical architecture as for both stability and redundancy, not all electric connections are allowed.
Consequently, to simplify the optimization problem and avoid introducing many constraints, the model uses a single categorical variable to represent the various feasible architectures.

Note that a simplified analysis has been done in a conference paper~\cite{SciTech_cat}, the latter was an optimization of the same aircraft configuration but with simpler methods, both HH and HH with PLS had never been tested before. 
The relaxed dimension used to construct the GP model using the CR method is 29 as indicated in~\tabref{SMO_tab:dragon}, while the relaxed dimension for the most general GP model employing the HH method is 137.

\resizebox{1.0\textwidth}{!}{%
    \hspace{-0.1cm}
\begin{threeparttable}[ht]
\centering
\vspace*{-0.3cm}
\caption{Definition of the ``\texttt{DRAGON}'' optimization problem.}
\begin{tabular}{lllrr}
& Function/variable & Nature & Quantity & Range\\
\hline
\hline
Minimize & Fuel mass & cont & 1 &\\
\hline
with respect to & \mbox{Fan operating pressure ratio} & cont & 1 & $\left[1.05, 1.3\right]$ \\  
     & \mbox{Wing aspect ratio} & cont & 1 &    $\left[8, 12\right]$ \\
    & \mbox{Angle for swept wing} & cont & 1 & $\left[15, 40\right]$  ($^\circ$) \\
     & \mbox{Wing taper ratio} & cont & 1 &    $\left[0.2, 0.5\right]$ \\
     & \mbox{HT aspect ratio} & cont & 1 &    $\left[3, 6\right]$ \\
    & \mbox{Angle for swept HT} & cont & 1 & $\left[20, 40\right]$  ($^\circ$) \\
     & \mbox{HT taper ratio} & cont & 1 &    $\left[0.3, 0.5\right]$ \\
 & \mbox{TOFL for sizing}  & cont &1 & $\left[1800, 2500\right]$ ($m$) \\
 & \mbox{Top of climb vertical speed for sizing} & cont & 1 & $\left[300, 800\right]$ ($ft/min$) \\
 & \mbox{Start of climb slope angle} & cont & 1 & $\left[0.075, 0.15\right]$ ($rad$) \\
 & \multicolumn{2}{l}{Total  continuous variables} & 10 & \\
 \hline
& \mbox{Architecture} & cat & 17 levels & \{1,2,3, \ldots,15,16,17\} \\
& \mbox{Turboshaft layout} & cat & 2 levels & \{1,2\} \\
 & \multicolumn{2}{l}{Total categorical variables} & 2 & \\
 \hline
  &   \multicolumn{2}{l}{\textbf{Total relaxed variables}} & {\textbf{29}} & \\
  \hline
subject to & Wing span \textless  \ 36   ($m$)  & cont & 1 \\
 & TOFL \textless  \ 2200 ($m$) & cont & 1 \\
 & Wing trailing edge occupied by fans  \textless  \ 14.4 ($m$) & cont & 1 \\
 & Climb duration \textless  \ 1740 ($s $) & cont & 1 \\
 & Top of climb slope \textgreater \ 0.0108 ($rad$) & cont & 1 \\
 & \multicolumn{2}{l}{\textbf{Total  constraints}} & {\textbf{5}} & \\
\hline
\end{tabular}

\label{SMO_tab:dragon}
\end{threeparttable}
}
\begin{table}[ht]
\centering
 \caption{Definition of the architecture variable and its 17 associated levels.}
\small
\begin{tabular}{ccc}
\hline
  \textbf{Architecture number} & number of motors & number of cores and generators\\
  \hline
  \textbf{1} & 8 &2 \\
  \textbf{2} & 12 &  2\\
  \textbf{3} & 16 &  2\\
  \textbf{4} &20 &2 \\
  \textbf{5} & 24 &  2\\
  \textbf{6} & 28 &  2\\
  \textbf{7} &32 & 2\\
  \textbf{8} & 36  & 2\\
  \textbf{9} & 40 &  2\\
  \textbf{10} & 8   & 4\\
  \textbf{11} & 16  & 4\\
  \textbf{12} & 24  & 4\\
  \textbf{13} & 32  & 4\\
  \textbf{14} & 40  & 4\\
  \textbf{15} & 12  & 6\\
  \textbf{16} & 24  & 6\\
  \textbf{17} & 36  & 6\\
\hline
\end{tabular}
\label{SMO_tab:dragon_archi1}
\end{table}
\begin{table}[H]
\centering
 \caption{Definition of the turboshaft layout variable and its 2 associated levels.}
\small
\begin{tabular}{cccccc}
\hline
  \textbf{Layout} & position & y ratio & tail & VT aspect ratio & VT taper ratio\\
  \hline 
  \textbf{1} & under wing &0.25 & without T-tail& 1.8 & 0.3 \\
  \textbf{2} & behind & 0.34 & with T-tail& 1.2 & 0.85\\
\hline
\end{tabular}
\label{SMO_tab:dragon_archi2}
\end{table}

To validate our method, we compare the 7 methods described in~\tabref{SMO_tab:dragon_meth} on the optimization of the \gls{DRAGON} aircraft concept with 5 and 10 points for the initial \gls{DOE} as before.

\begin{figure}[H]
\begin{minipage}[b]{.6\linewidth}
\centering
\hspace{-1.25cm}
\subfloat[Convergence curves: medians of 10 runs.]{
\includegraphics[height=5cm,,width=7cm]
{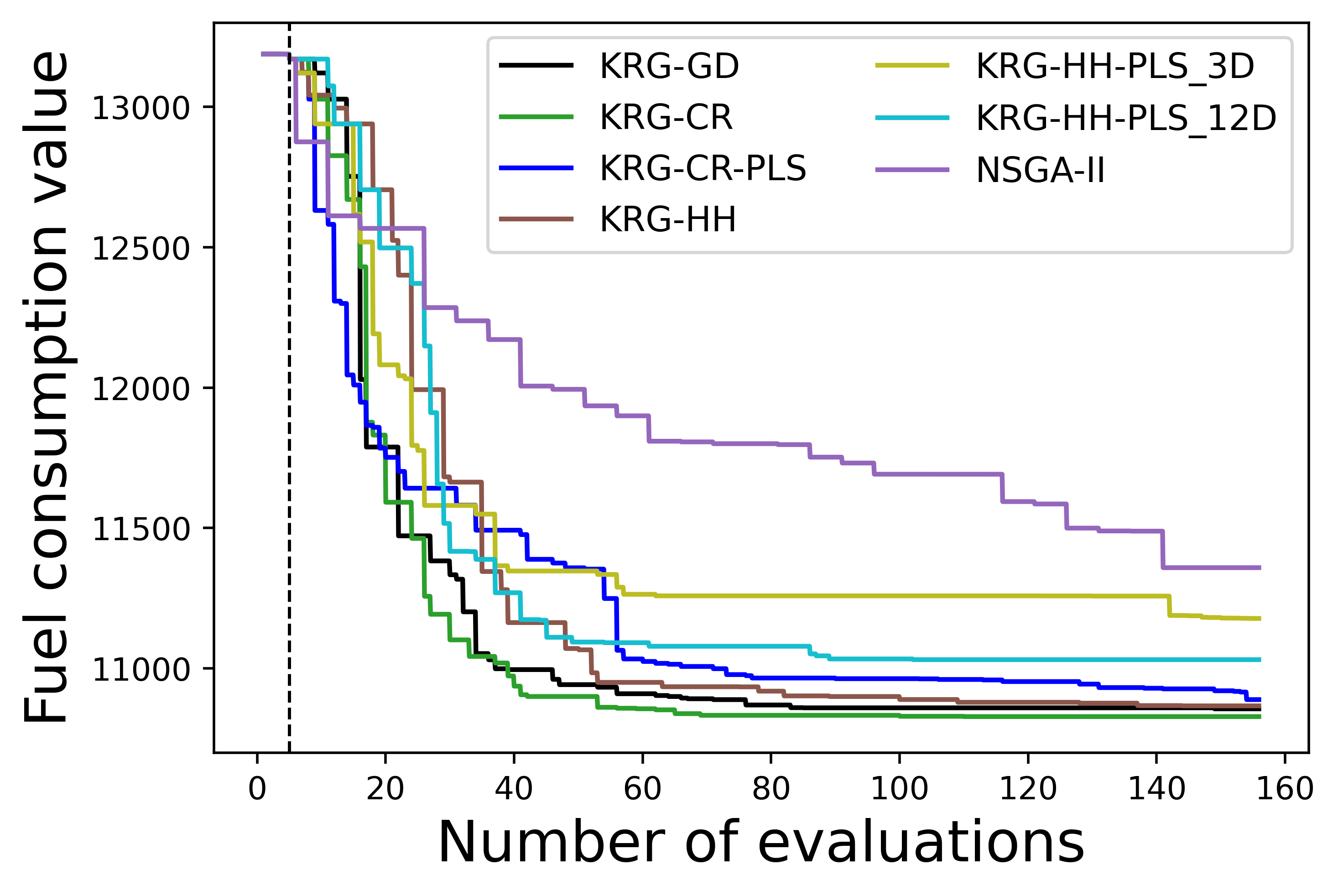}
\label{SMO_convmidragon5}
}
\end{minipage}
\begin{minipage}[b]{.4\linewidth}
\centering 
\hspace{-1.5cm}
\subfloat[Boxplots after 100 evaluations.]{
\includegraphics[height=5cm,width=6.2cm]{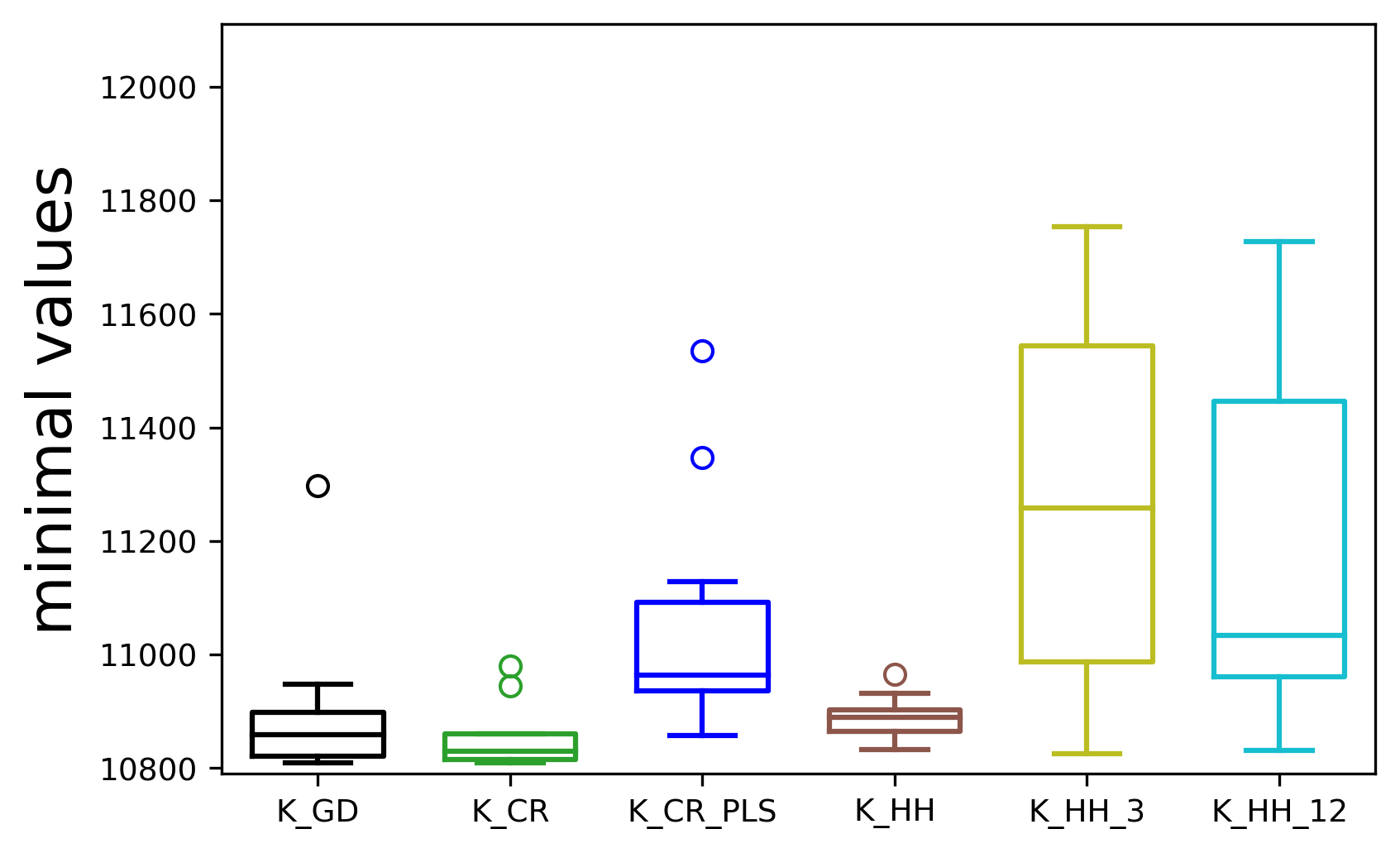}
\label{SMO_mini_midragon5}
}
\end{minipage}
\caption{Optimization results for the ``\texttt{DRAGON}'' aircraft~\cite{SciTech_cat}  for 10 DoE of 5 points.}
\label{SMO_res_optim_dragon5}
\end{figure}

\begin{figure}[H]
\begin{minipage}[b]{.6\linewidth}
\centering
\hspace{-1.25cm}
\subfloat[Convergence curves: medians of 10 runs.]{
\includegraphics[height=5cm,,width=7cm]
{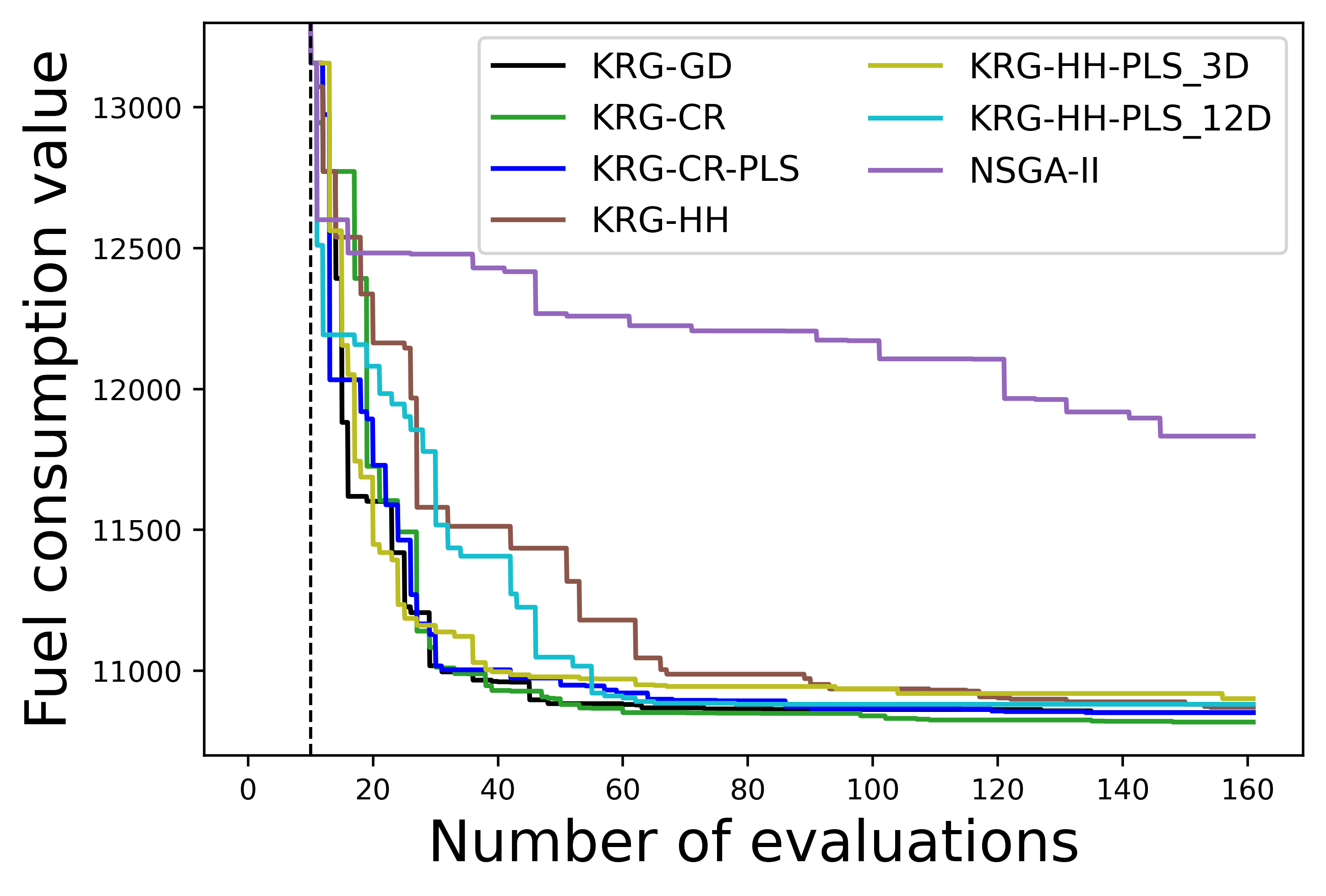}
\label{SMO_convmidragon10}
}
\end{minipage}
\begin{minipage}[b]{.4\linewidth}
\centering 
\hspace{-1.5cm}
\subfloat[Boxplots after 100 evaluations.]{
\includegraphics[height=5cm,width=6.2cm]{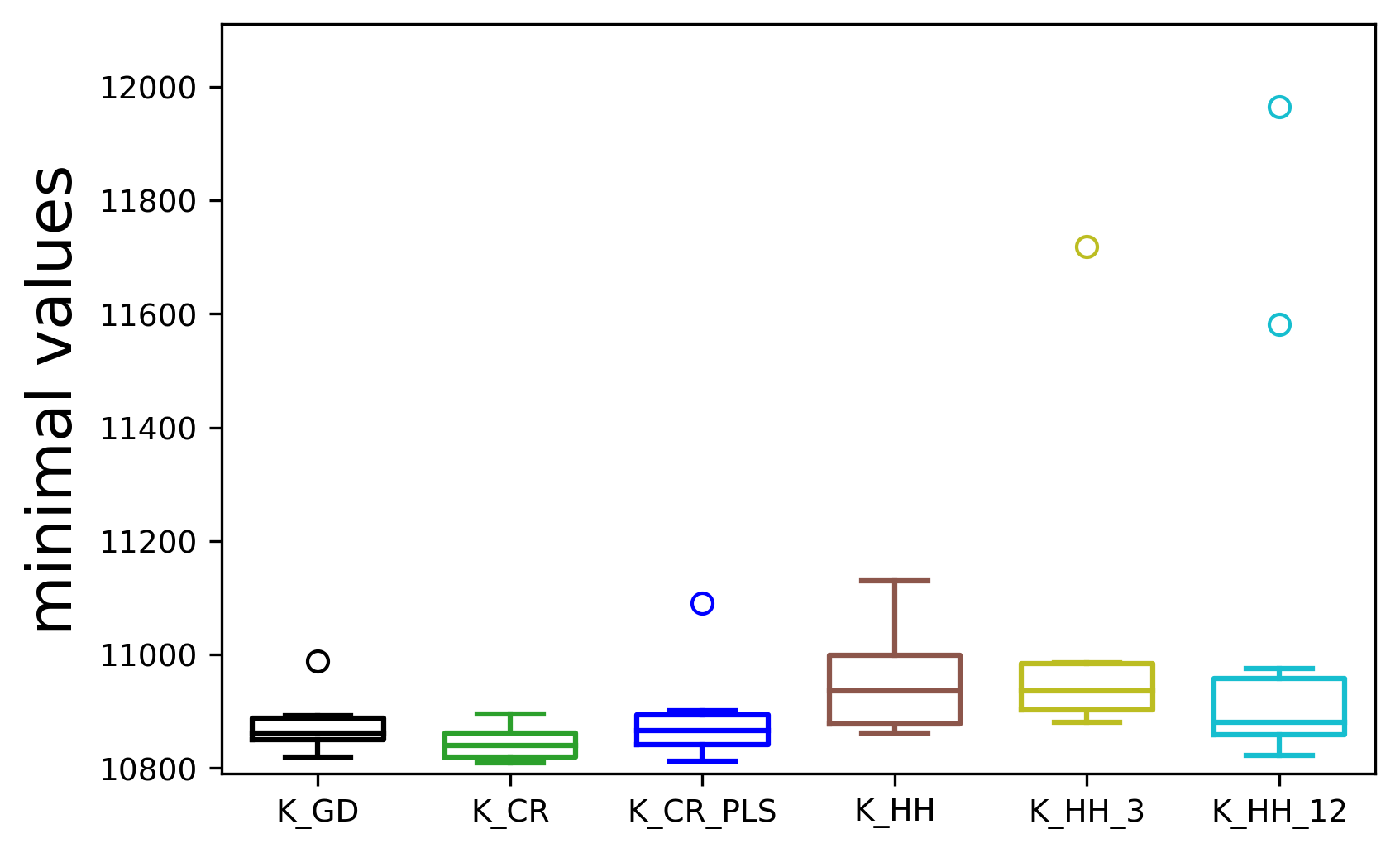}
\label{SMO_mini_midragon10}
}
\end{minipage}
\caption{Optimization results for the ``\texttt{DRAGON}'' aircraft~\cite{SciTech_cat}  for 10 DoE of 10 points.}
\label{SMO_res_optim_dragon10}
\end{figure}

\begin{table}[H]
\centering
 \caption{The various kernels compared on the MDO of ``\texttt{DRAGON}''. }
\small
\begin{tabular}{cccc}
\hline
    \textbf{Name} & \# of cat. params & \# of cont. params & Total \# of params\\
    \hline 
    \textbf{GD} & 2 & 10 & 12\\
     \hdashline
    \textbf{CR} & 19 & 10 & 29  \\
    \textbf{CR with PLS 3D} & Not applicable & Not applicable & 3 \\ 
    \hdashline
    \textbf{HH} & 137 & 10 & 147 \\
    \textbf{HH with PLS 3D} & 2 & 1 & 3 \\
    \textbf{HH with PLS 12D} & 2 & 10 & 12 \\
    \hdashline
    \textbf{NSGA-II} & Not applicable & Not applicable & Not applicable \\
    \hline
\end{tabular}
\label{SMO_tab:dragon_meth}
\end{table}   
As mentioned above, we are doing 10 runs for every method based on 10 starting DoE sampled by LHS to quantify the methods randomness. For every method and every starting DoE, we are running the method for a budget of 150 infill points, hence evaluating the black-box 155 times for the 5 point DoE and 160 times for the 10 point DoE. 
The results are given on~\figref{SMO_res_optim_dragon5} for the 5 point DoE and on~\figref{SMO_res_optim_dragon10} for the 10 point DoE.
More precisely,~\figref{SMO_convmidragon5} and~\figref{SMO_convmidragon10} display the convergence curves for the 7 methods and, to visualize the data dispersion, the boxplots of the 10 best solutions after 100 evaluations are shown in~\figref{SMO_mini_midragon5} and~\figref{SMO_mini_midragon10}. This computer experiments setup and figures are similar to what can be found in~\cite{SciTech_cat}.
We note that for this study, the computational cost of building the GP model is assumed to be negligible compared to the cost of evaluating the objective and the constraints at a given point. In fact, one simulation related to DRAGON is taking 2 to 5 minutes, which means half a day for a total of 150 simulations and around 60 days of computational time for running all the optimizations related to test case in this section.

These results confirm the previous analyses made on the analytic test cases. 
First, the methods without PLS (GD, CR and HH) converge slightly better with a 5 point DoE, although this effect is less significant in this case because 5 or 10 points are both small quantities compared to 150 iterations and also because the search space is larger than before.
Second, the methods with PLS (CR-PLS, HH-PLS\_3D and HH-PLS\_12D) greatly benefit from a bigger starting DoE because the more representative the initial data, the more relevant the computed principal components. 
We note that a small DoE could lead to have PLS methods stuck in irrelevant zones as it might not be well-posed. In fact, in the case where an important part of the design space is not featured in the DoE, the available data will not be able to learn on this zone and thus the approximate representative space (computed by the PLS) will be sub-optimal.
The results show that for an initial DoE of 5 points, the methods without PLS are both faster and more consistent than their equivalent with PLS whereas the opposite is observed with the  DoE of 10 points. The 10 point DoE results also display that the HH model is costly and too complicated to be used efficiently with small DoE and that simpler methods like CR-PLS or GD are most effective at the start of the optimization process.
Based on this observation, it seems that combining models is the overall best method to tackle whatever black-box optimization problem beforehand.
These results also display that BO is more fitted to tackle such problems than evolutionary algorithms and that a small DoE could lead to have PLS methods stuck in irrelevant zones which, once again favors the alternative method of combining models for future works.

In terms of aircraft design, the ideal configuration was determined with an estimated fuel consumption of 10,809 kilograms against 11,248 kilograms for the original reference configuration. This configuration corresponds to option 10, featuring fewer engines (8 in total), but incorporating 4 cores and electric generators.
The most advantageous layout positions the turbo-generators at the rear. This choice is influenced by the increased lever arm between the wing and the horizontal tail, which results from the maximum sweep angle applied to the horizontal tail. However, it is worth noting that the combination of high sweep and high aspect ratio is not adequately considered from a structural standpoint, leading to unrealistically heavy weights for the horizontal stabilizer. Despite this limitation, the optimization process yields a suitable trade-off based on the models used in FAST-OAD.  The optimum found in the previous study of~\cite{SciTech_cat} was 10,816 kilograms which show that our new algorithms are more efficient. Still the new aircraft configuration is really close to the previous one, the changes are on the wing taper ratio reduced from 0.235 to 0.22, the TOFL for sizing which is now at the lower bound of 1800 m and not at 1803 m and to finish with, the start of climb angle was slightly reduced from 0.104 rad to  0.1035 rad. 

\section{Conclusion}
\label{SMO_sec:conclu}



In this work, we proposed mixed-categorical metamodels based on GP for high-dimensional structural and multidisciplinary optimization. 
Our research was driven by the increasing complexity of engineering systems, which involve various disciplines and require optimization involving numerous design variables that could be either continuous, integer, and categorical. 
Our key findings center on the development of a more efficient approach for building surrogate models for large-scale mixed-categorical inputs. In~\cite{Mixed_Paul}, we introduced a mixed categorical kernel (EHH), a powerful tool for handling mixed-categorical variables which combine the matrix-based HH approach with the exponential kernel. However, we identified that the EHH and HH kernels effectiveness came at the cost of a significant increase in hyperparameters related to the GP surrogate model. To address this issue, we devised a novel approach by extending the partial least squares regression method as developed in~\cite{Bouhlel18} for continuous kernels, to reduce the number of hyperparameters while maintaining accuracy.

The significance of our research extends to both researchers and practitioners. For researchers, our work contributes to the evolving field of surrogate modeling for MDO. It offers a valuable solution to the challenge of high-dimensional mixed-categorical optimization, opening doors for further exploration in this domain.
Practitioners in engineering and optimization fields will find our findings beneficial as they provide a practical and efficient toolset for handling complex optimization problems. Our approach, implemented in the open-source software SMT~\cite{SMT2019,saves2023smt}, has been demonstrated effectively in structural and multidisciplinary applications, showcasing its real-world applicability.

Further works may include combining the several methods that now exist in the literature to have surrogate models that increase automatically in complexity when the size of the dataset increases along the optimization process. 
Also, the surrogate models can be coupled to any surrogate-based optimization algorithm. In particular, in~\cite{grapin_constrained_2022}, SEGOMOE has been extended to multi-objective optimization and we also consider extending high-dimensional GP models to both mixed and hierarchical variables to tackle technological choices and variable-size problems~\cite{saves2023smt,audet2022general}. 
Future performance benchmarks on industrial test cases should include comparisons with the Maximum Likelihood Estimation (MLE) approach for latent space identification, as demonstrated in latent map Gaussian process~\cite{oune2021latent}.
%
%
%

\recap{ 
\lettrine[lines=2, lhang=0.33, loversize=0.25, findent=1.5em]{I}{n} this chapter, we proposed a new mixed-categorical GP for high dimension based on PLS regression to reduce the number of hyperparameters the model relies on. In particular, the contributions of this chapter are listed below.
\begin{itemize}
    \item A new PLS regression for symmetric positive definite matrices has been developed. 
    \item The KPLS method has been extended to the general matrix based categorical EHH and HH kernels. 
    \item The HH with PLS and EHH with PLS kernels have been implemented into the open-source software SMT.
    \item The various categorical GP models have been compared in terms of approximation performance. They also showcase  their respective abilities to capture the structure of the correlation matrix and this work highlights the interest of PLS in that context as it can capture a lot of information while relying on few hyperparameters.
    \item The new KPLS models have been used to perform Bayesian optimization for multidisciplinary problems and, more precisely, to optimize the "\texttt{DRAGON}" aircraft concept. 
\end{itemize}

This chapter corresponds to the article under review:
{  \textit{Saves, P., Diouane, Y., Bartoli, N., Lefebvre, T., Morlier, J., “High-dimensional mixed-categorical Gaussian process with application to multidisciplinary design optimization for a green aircraft”, Structural and Multidisciplinary Optimization, 2024.}}
}

\chapter{\texttt{SMT 2.0}: A Surrogate Modeling Toolbox with a focus on Hierarchical and Mixed Variables Gaussian Processes} 


\chaptermark{\texttt{SMT 2.0} \MakeLowercase{with hierarchical and mixed variables} GP}
 \label{c4}

\setlength{\fboxrule}{0pt}
\hspace{4cm}  \noindent\fbox{%
     \parbox{0.75\textwidth}{%
        I am also a power, and my power is strong as long as I may set the strength of my words against that of the world. Such is my only consolation. I know that the relapses into despair will be many and deep, but the miracle of liberation leads me to a goal that makes me dizzy: a reason for living. $\textit{(translated)}$       
 \\
        \hrule \vspace{0.2cm}
     \hspace*{\fill} Vårt behov av tröst är omättligtn, Stig Dagerman}%
} 

\objectif{ 

\lettrine[lines=2, lhang=0.33, loversize=0.25, findent=1.5em]{T}{he} Surrogate Modeling Toolbox (SMT) is an open-source Python package that offers a collection of surrogate modeling methods, sampling techniques, and a set of sample problems.  
This chapter presents \texttt{SMT 2.0}, a major new release of SMT that introduces significant upgrades and new features to the toolbox whose objectives are described hereinafter.
\begin{itemize}
    \item To propose new Gaussian process models for mixed discrete and hierarchical variables handling. 
    \item To propose new derivatives capabilities for the models.
    \item To implement new surrogate models, new sampling methods and new applications.
    \item To document the open-source software and propose new tutorials.
\end{itemize}


%
}


\minitoc


\setcounter{section}{-1}
\section{Synthèse du chapitre en français}

Le chapitre "SMT 2.0: A Surrogate Modeling Toolbox with a focus on Hierarchical and Mixed Variables Gaussian Processes" présente les développements récents de la version 2.0 de SMT, une boîte à outils open-source de modèles de substitution en Python.

L'utilisation de modèles de substitution est devenue une technique courante pour réduire l'effort computationnel dans des tâches d'exploration de l'espace de conception, de quantification de l'incertitude ou d'optimisation de processus coûteux. Dans ce contexte, SMT 2.0 vise à regrouper les modèles de substitution de la littérature tout en étendant et élargissant leurs capacités de modélisation.
Pour ce faire, SMT 2.0 introduit notamment la prise en compte des variables hiérarchiques et mixtes dans les modèles car ces types de variable sont souvent présents dans les problèmes d'ingénierie. De plus, SMT 2.0 étend la prise en compte des dérivées des modèles, permet la construction de modèles en haute-dimension et fournit des estimations de dérivées pouvant être, par exemple, utilisées à des fins d'optimisation. 
Toutefois, la principale contribution qu'apporte SMT 2.0 est la construction de modèles de substitution impliquant des variables mixtes et pouvant prendre en compte des hiérarchies de variables, notamment dans les processus gaussiens (également appelés modèles de krigeage). Pour gérer les variables mixtes, SMT 2.0 implémente plusieurs noyaux de corrélation, tels que le noyau basée sur la distance de Gower (GD), le noyau basé sur l'hypersphère homoscédastique (HH) et sa version exponentielle (EHH). Ces noyaux permettent de construire des modèles de substitution plus ou moins précis et coûteux suivant le problème et l'usage recherché par l'utilisateur. 
Qui plus est, SMT 2.0 introduit d'autres améliorations telles que de nouvelles procédures d'échantillonnage, de nouveaux modèles de substitution, les dérivées des noyaux de corrélation, les dérivées de la variance du modèle de krigeage, un critère adaptatif pour les problèmes à haute dimension,... Ce chapitre présente aussi des applications de l'optimisation bayésienne avec des variables hiérarchiques et mixtes, ainsi que des applications à la conception d'aéronefs ou à la fusion de données.

Du point de vue logiciel, SMT 2.0 maintient et améliore la modularité et la généricité de la version précédente (version 1.3). Ce chapitre décrit l'organisation du code source de SMT 2.0, qui est divisé en sous-modules pour les méthodes d'échantillonnage, les problèmes de benchmarking et les modèles de substitution. Il souligne également les efforts fournis pour assurer une documentation complète et de haute qualité, ainsi que des tests automatiques pour garantir la stabilité du logiciel.

En résumé, SMT 2.0 est une boîte à outils pour les modèles de substitution prenant en charge les variables hiérarchiques et mixtes. Grâce à ses nouvelles fonctionnalités et à sa praticité d'utilisation, SMT 2.0 facilite la résolution de problèmes industriels et offre de nouvelles possibilités pour l'optimisation et la quantification d'incertitude des systèmes complexes.

\section{Motivation and significance}
\label{sec:introSMT}
With the increasing complexity and accuracy of numerical models, it has become more challenging to run complex simulations and computer codes~\cite{Mader_ad,Kennedy_Bayes}. 
As a consequence, surrogate models have been recognized as a key tool for engineering tasks such as design space exploration, uncertainty quantification, and optimization~\cite{Hwang2018b}. 
In practice, surrogate models are used to reduce the computational effort of these tasks by replacing expensive numerical simulations with closed-form approximations~\cite[Ch.~10]{Martins2021}.
To build such a model, we start by evaluating the original expensive simulation at a set of points through a Design of Experiments (\gls{DOE}).
Then, the corresponding evaluations are used to build the surrogate model according to the chosen approximation, such as Kriging, quadratic interpolation, or least squares regression.

The Surrogate Modeling Toolbox (SMT) is an open-source framework that provides functions to efficiently build surrogate models~\cite{SMT2019}.
Kriging models (also known as Gaussian processes) that take advantage of derivative information are one of SMT's key features~\cite{bouhlel2019gradient}.
Numerical experiments have shown that SMT achieved lower prediction error and computational cost than Scikit-learn~\cite{scikit-learn} and UQLab~\cite{UQLab} for a fixed number of points~\cite{ToolSMT}.
SMT has been applied to rocket engine coaxial-injector optimization~\cite{DL1}, aircraft engine consumption modeling~\cite{DGP1}, numerical integration~\cite{eliavs2020periodic}, multi-fidelity sensitivity analysis~\cite{drouet2023multi}, high-order robust finite elements methods~\cite{karban2021fem,kudela2022recent}, planning for photovoltaic solar energy~\cite{chen2020surrogate}, wind turbines design optimization~\cite{jasa2022effectively}, porous material optimization for a high pressure turbine vane~\cite{wang2023transpiration}, chemical process design~\cite{savage2020adaptive} and many other applications.
 
In systems engineering, architecture-level choices significantly influence the final system performance, and therefore, it is desirable to consider such choices in the early design phases~\cite{chan2022trying}. 
Architectural choices are parameterized with discrete design variables; examples include the selection of technologies, materials, component connections, and number of instantiated elements.
When design problems include both discrete variables and continuous variables, they are said to have \emph{mixed variables}.

When architectural choices lead to different sets of design variables, we have \emph{hierarchical} variables~\cite{Hutter,Architecture}.
For example, consider different aircraft propulsion architectures~\cite{fouda2022automated}.
A conventional gas turbine would not require a variable to represent a choice in the electrical power source, while hybrid of pure electric propulsion would require such a variable.
The relationship between the choices and the sets of variables can be represented by a hierarchy.

Handling hierarchical and mixed variables requires specialized  surrogate modeling techniques~\cite{Effectiveness}.
To address these needs, \texttt{SMT 2.0} is offering researchers and practitioners a collection of cutting-edge tools to build surrogate models with continuous, mixed and hierarchical variables.
The main objective of this chapter is to detail the new enhancements that have been added in this release compared to the original \texttt{SMT 0.2} release~\cite{SMT2019}.

There are two new major capabilities in \texttt{SMT 2.0}: the ability to build surrogate models involving mixed variables and the support for hierarchical variables within Kriging models.
To handle mixed variables in Kriging models, existing libraries such as BoTorch~\cite{balandat2020botorch}, Dakota~\cite{Dakota}, DiceKriging~\cite{DiceKriging}, LVGP~\cite{zhang2020latent}, Parmoo~\cite{parmoo}, and Spearmint~\cite{GMHL} implement simple mixed models by using either continuous relaxation (CR), also known as \emph{one-hot encoding}~\cite{GMHL}, or a Gower distance (GD) based correlation kernel~\cite{Gower}. 
KerGP~\cite{Roustant} (developed in R) implements more general kernels but there is no Python open-source toolbox that implements more general kernels to deal with mixed variables, such as the homoscedastic hypersphere (HH)~\textcolor{black}{\cite{Zhou}} and exponential homoscedastic hypersphere (EHH)~\textcolor{black}{\cite{Mixed_Paul}} kernels.
Such kernels require the tuning of a large number of hyperparameters but lead to more accurate Kriging surrogates than simpler mixed kernels~\cite{Pelamatti,Mixed_Paul}. 
\texttt{SMT 2.0} implements all these kernels (CR, GD, HH, and EHH) through a unified framework and implementation. 
%
%
\textcolor{black}{  To handle hierarchical variables, no library in the literature can  build peculiar surrogate models except \texttt{SMT 2.0}, which implements two Kriging methods for these variables. 
Notwithstanding, most softwares are compatible with a naïve strategy called the imputation method~\cite{Effectiveness}   
but this method lacks depth and depends on arbitrary choices.
This is why~\citet{Hutter} proposed a first kernel, called \texttt{Arc-Kernel} which in turn was generalized by~\citet{Horn_hier} with a new kernel called the \texttt{Wedge-Kernel}~\cite{DACE_hier}.
None of these kernels are available in any open-source modeling software. Furthermore, thanks to the framework introduced in~\citet{audet2022general}, our proposed kernels are sufficiently general so that all existing hierarchical kernels are included within it. Section 4 describes the two kernels implemented in \texttt{SMT 2.0} that are referred as \texttt{SMT Arc-Kernel} and \texttt{SMT Alg-Kernel}.
In particular, \texttt{Alg-Kernel} is a novel hierarchical kernel introduced in this chapter. } 
Table~\ref{tab:comparison} outlines the main features of the state-of-the-art modeling software that can handle hierarchical and mixed variables.
\begin{table}[H]
\caption{Comparison of software packages for hierarchical and mixed Kriging models. \checkmark = implemented. * = user-defined. 
}
\resizebox{\columnwidth}{!}{%
\begin{tabular}{l c c c c c c c c }
\hline
\textbf{Package} & \texttt{BOTorch} & \texttt{Dakota} & \texttt{DiceKriging} & \texttt{KerGP} & \texttt{LVGP} & \texttt{Parmoo}  & \texttt{Spearmint} & \texttt{SMT 2.0} \\
\hline
\textbf{Reference} & \cite{balandat2020botorch} & \cite{Dakota} & \cite{DiceKriging} & \cite{Roustant} &\cite{zhang2020latent}  &\cite{parmoo} &\cite{GMHL} & This chapter    \\
\textbf{License}   & MIT & EPL & GPL & GPL & GPL & BSD &  GNU & BSD  \\
\textbf{Language}  & Python & C & R & R & R & Python & Python &  Python \\

\textbf{Mixed var.} & \checkmark & \checkmark & \checkmark & \checkmark & \checkmark & \checkmark & \checkmark & \checkmark \\
\textit{GD kernel} & \checkmark & \checkmark & \checkmark & *  & &  &  & \checkmark \\
\textit{CR kernel} &  &  &  & & \checkmark  &  \checkmark  &  \checkmark & \checkmark\\
\textit{HH kernel} &  &  & & \checkmark  &  &  &  & \checkmark  \\
\textit{EHH kernel} &  &  &  & *  &  & &   & \checkmark \\
\textbf{Hierarchical var.} &  &  &  &  &  & &   & \checkmark \\
\hline
\end{tabular}
}
\label{tab:comparison}
\end{table}


\texttt{SMT 2.0} introduces other enhancements, such as additional sampling procedures, new surrogate models, new Kriging kernels (and their derivatives), Kriging variance derivatives, and an adaptive criterion for high-dimensional problems. 
\texttt{SMT 2.0} adds applications of Bayesian optimization (BO) with hierarchical and mixed variables or noisy co-Kriging that have been successfully applied to aircraft design~\cite{SciTech_cat}, data fusion~\cite{condearenzana}, and structural design~\cite{RaulAIAA}.
The \texttt{SMT 2.0} interface is more user-friendly and offers an improved and more detailed documentation for users and developers\footnote{\url{http://smt.readthedocs.io/en/latest}}. 
\texttt{SMT 2.0} is hosted publicly\footnote{\url{https://github.com/SMTorg/smt}} and can be directly imported within Python scripts.
It is released under the New BSD License and runs on Linux, MacOS, and Windows operating systems.
Regression tests are run automatically for each operating system whenever a change is committed to the repository.
In short, \texttt{SMT 2.0} builds on the strengths of the original SMT package while adding new features. On one hand, the emphasis on derivatives (including prediction, training and output derivatives) is maintained and improved in \texttt{SMT 2.0}. On the other hand, this new release includes support for hierarchical and mixed variables Kriging based models.
For the sake of reproducibility, an open-source notebook is available that gathers all the methods and results presented on this chapter\footnote{\url{https://github.com/SMTorg/smt/tree/master/tutorial/NotebookRunTestCases_Paper_SMT_v2.ipynb}}.

The remainder of the chapter is organized as follows.
First, we introduce the organization and the main implemented features of the release in Section~\ref{sec:description}. 
Then, we describe the mixed-variable Kriging model with an example in Section~\ref{sec:mixed}. 
Similarly, we describe and provide an example for a hierarchical-variable Kriging model in Section~\ref{sec:hier}.
The Bayesian optimization models and applications are described in  Section~\ref{sec:BO}.
Finally, we describe  the other relevant contributions in Section~\ref{sec:other} and conclude in Section~\ref{sec:concl}.

\section{\texttt{SMT 2.0}: an improved surrogate modeling toolbox}
\label{sec:description}

From a software point of view, \texttt{SMT 2.0} maintains and improves the modularity and generality of the original SMT version~\cite{SMT2019}. 
In this section, we describe the software as follows.
Section~\ref{sec:legacy} describes the legacy of \texttt{SMT 0.2}.
Then, Section~\ref{sec:organization} describes the organization of the repository.
Finally, Section~\ref{sec:newSMT2} shows the new capabilities implemented in the \texttt{SMT 2.0} update.

\subsection{Background on SMT former version: \texttt{SMT 0.2} }
\label{sec:legacy}  

SMT~\cite{SMT2019} is an open-source collaborative work originally developed by ONERA, NASA Glenn, ISAE-SUPAERO/ICA and the University of Michigan.
Now, both Polytechnique Montréal and the University of California San Diego are also contributors. 
\texttt{SMT 2.0} updates and extends the original SMT repository capabilities among which the original publication~\cite{SMT2019} focuses on different types of derivatives for surrogate models detailed hereafter.

\paragraph{A Python surrogate modeling framework with derivatives}
One of the original main motivations for SMT was derivative support. In fact, none of the existing packages for surrogate modeling such as Scikit-learn in Python~\cite{scikit-learn}, SUMO in Matlab~\cite{Gorissen2010} or GPML in Matlab and Octave~\cite{williams2006gaussian} focuses on derivatives. 
Three types of derivatives are distinguished: prediction derivatives, training derivatives, and output derivatives.
SMT also includes new models with derivatives such as Kriging with Partial Least Squares (KPLS)~\textcolor{black}{\cite{Bouhlel18}} and Regularized Minimal-energy Tensor-product Spline (RMTS)~\textcolor{black}{\cite{Hwang2018b}}. 
These developed derivatives were even used in a novel algorithm called Gradient-Enhanced Kriging with Partial Least Squares (GEKPLS)~\textcolor{black}{\cite{bouhlel2019gradient}} to use with adjoint methods, for example~\textcolor{black}{\cite{bouhlel2020}}. 

\paragraph{Software architecture, documentation, and automatic testing}
SMT is organized along three main sub-modules that implement a set of sampling techniques (\textbf{sampling\_methods}), benchmarking functions (\textbf{problems}), and surrogate modeling techniques (\textbf{surrogate\_models}). 
The toolbox documentation\footnote{\url{https://smt.readthedocs.org}} is created using reStructuredText and Sphinx, a documentation generation package for Python, with custom extensions.
Code snippets in the documentation pages are taken directly from actual tests in the source code and are automatically updated. 
The output from these code snippets and tables of options are generated dynamically by custom Sphinx extensions. 
This leads to high-quality documentation with minimal effort.
Along with user documentation, developer documentation is also provided to explain how to contribute to SMT.
This includes a list of API methods for the \textbf{SurrogateModel}, \textbf{SamplingMethod}, and \textbf{Problem} classes, that must be implemented to create a new surrogate modeling method, sampling technique, or benchmarking problem.
When a developer submits a pull request, it is merged only after passing the automated tests and receiving approval from at least one reviewer.
The repository on GitHub\footnote{\url{https://github.com/SMTorg/smt}} is linked to continuous integration tests (GitHub Actions) for Windows, Linux and MacOS, to a coverage test on coveralls.io and to a dependency version check for Python with DependaBot.
Various parts of the source code have been accelerated using Numba~\cite{numba}, an LLVM-based just-in-time (JIT) compiler for numpy-heavy Python code. Numba is applied to conventional Python code using function decorators, thereby minimizing its impact on the development process and not requiring an additional build step. 
For a mixed Kriging surrogate with 150 training points, a speedup of up to 80\% is observed, see Table~\ref{tab:numba_benchmark}.
The JIT compilation step only needs to be done once when installing or upgrading SMT and adds an overhead of approximately 24 seconds on a typical workstation.
\textcolor{black}{
In this chapter, all results are obtained using an Intel® Xeon® CPU E5-2650 v4 @ 2.20 GHz core and 128 GB of memory with a Broadwell-generation processor front-end and a compute node of a peak power of 844 GFlops.}

\begin{table}[H]
\caption{Impact of using Numba on training time of the hierarchical Goldstein problem. Speedup is calculated excluding the JIT compilation table, as this step is only needed once after SMT installation.}
\begin{center}
\resizebox{0.85\columnwidth}{!}{%
\begin{tabular}{lcccc}
\hline
\textbf{Training set} & Without Numba & Numba & Speedup & JIT overhead \\
\hline
\textbf{15 points} & 1.3 s & 1.1 s & 15\% & 24 s \\
\textbf{150 points} & 38 s & 7.4 s & 80\% & 23 s \\
\hline
\end{tabular}
}
\end{center}
\label{tab:numba_benchmark}
\end{table}

\subsection{Organization of \texttt{SMT 2.0}} 
\label{sec:organization}
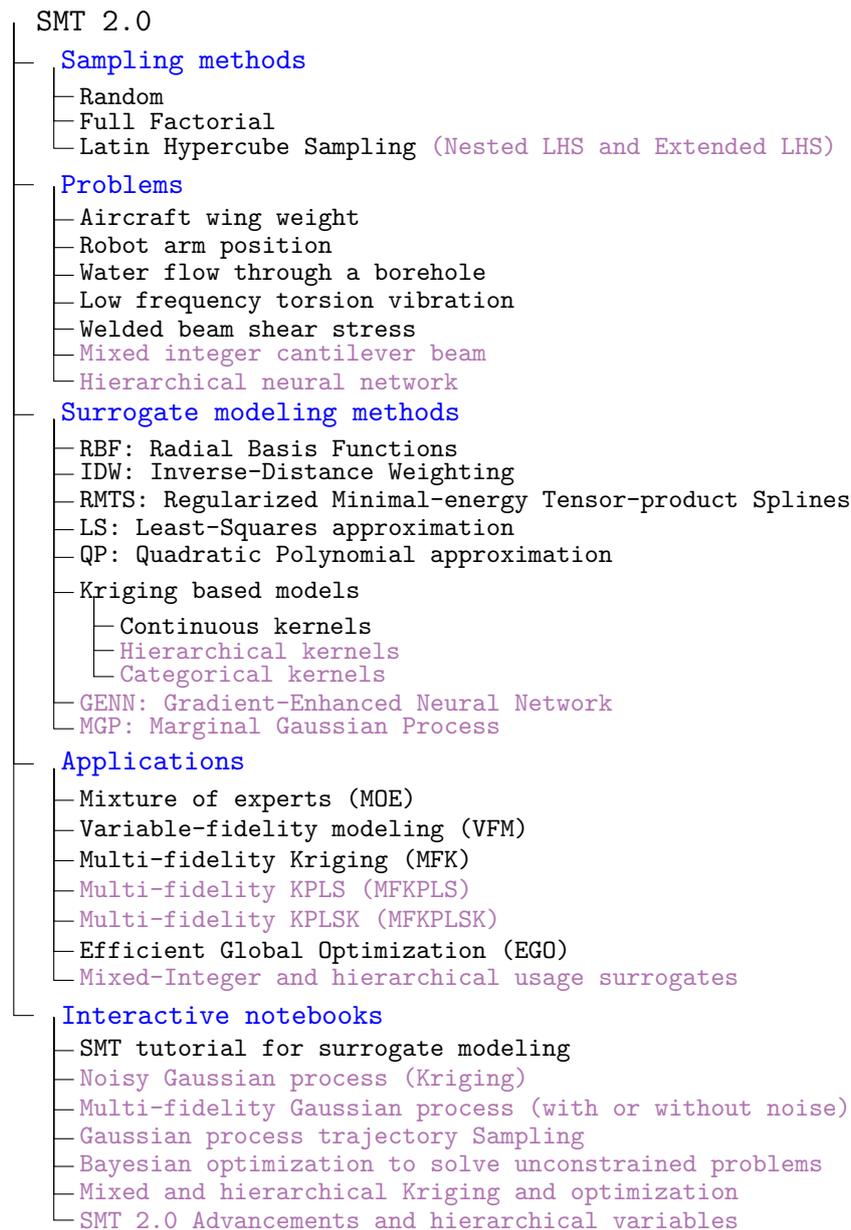
\begin{figure}[h!]
\centering
\vspace{3pt}
\begin{forest}
      for tree={
        font= \small \ttfamily,
        grow'=0,
        child anchor=west,
        parent anchor=south,
        anchor=west,
        calign=first,
        inner xsep=10pt,
        inner ysep = -3pt,
        edge path={
          \noexpand\path [draw, \forestoption{edge}]
          (!u.south west) +(7.5pt,0) |- (.child anchor)  \forestoption{edge label};
        },
        file/.style={edge path={\noexpand\path [draw, \forestoption{edge}]
          (!u.south west) +(7.5pt,0) |- (.child anchor) \forestoption{edge label};},
          inner xsep=2pt, inner ysep=-3pt ,font=\footnotesize \ttfamily
                     },
        before typesetting nodes={
          if n=1
            {insert before={[,phantom]}}
            {}
        },
        fit=band,
        before computing xy={l=15pt},
      }  
    [ \   \normalsize \texttt{SMT 2.0}
      [ \textcolor{blue}{Sampling methods}
        [Random, file
        ]      
        [Full Factorial, file
        ]
        [ Latin Hypercube Sampling \textcolor{Orchid}{(Nested LHS and Extended LHS)}, file
        ]
      ]
      [ \textcolor{blue}{Problems}
        [Aircraft wing weight, file
        ]      
        [Robot arm position, file
        ]  
        [Water flow through a borehole, file
        ]   
        [Low frequency torsion vibration, file
        ]    
        [Welded beam shear stress, file
        ]  
        [ \textcolor{Orchid}{Mixed integer cantilever beam}, file
        ]
        [ \textcolor{Orchid}{Hierarchical neural network}, file
        ]      
      ]
      [ \textcolor{blue}{Surrogate modeling methods}
        [ RBF: Radial Basis Functions, file
        ]
        [IDW: Inverse-Distance Weighting, file
        ]
        [RMTS: Regularized Minimal-energy Tensor-product Splines, file
        ]
         [LS: Least-Squares approximation, file
        ]      
         [QP: Quadratic Polynomial approximation,  file
        ]
        [Kriging based models,  file
            [Continuous kernels,  file]
            [\textcolor{Orchid}{Hierarchical kernels},  file]
            [\textcolor{Orchid}{Categorical kernels},  file]
        ]
        [\textcolor{Orchid}{GENN: Gradient-Enhanced Neural Network},  file ]
         [ \textcolor{Orchid}{MGP: Marginal Gaussian Process},  file ]
      ]
      [ \textcolor{blue}{Applications}
          [   Mixture of experts (MOE) , file]
          [  Variable-fidelity modeling (VFM) , file]
          [  Multi-fidelity Kriging (MFK) , file]
         [    \textcolor{Orchid}{ Multi-fidelity KPLS (MFKPLS)} , file]
         [    \textcolor{Orchid}{ Multi-fidelity KPLSK (MFKPLSK)} , file]
         [   Efficient Global Optimization (EGO) , file]
         [  \textcolor{Orchid}{ Mixed-Integer and hierarchical usage surrogates}, file]
      ]
        [ \textcolor{blue}{Interactive notebooks}
          [  SMT tutorial for surrogate modeling, file]
          [  \textcolor{Orchid}{  Noisy Gaussian process (Kriging) }, file]
          [ \textcolor{Orchid}{Multi-fidelity Gaussian process (with or without noise)} , file]
         [  \textcolor{Orchid}{ Gaussian process trajectory Sampling} , file]
         [ \textcolor{Orchid}{ Bayesian optimization to solve unconstrained problems} , file]
         [ \textcolor{Orchid}{ Mixed \& hierarchical Kriging and optimization}, file]
         [ \textcolor{Orchid}{SMT 2.0 Advancements and hierarchical variables}, file]
      ]
    ]
 \end{forest}
\caption{\label{fig:smt_codes} Functionalities of \texttt{SMT 2.0}. 
The new major features implemented in \texttt{SMT 2.0} compared to \texttt{SMT 0.2} are highlighted with the lavender color.}
\end{figure}

The main features of the open-source repository \texttt{SMT 2.0} are described in~\figref{fig:smt_codes}. 
More precisely, \texttt{Sampling Methods}, \texttt{Problems} and \texttt{Surrogate models} are kept from \texttt{SMT 0.2} and two new sections \texttt{Models applications} and \texttt{Interactive notebooks} have been added to the architecture of the code. 
These sections are highlighted in blue and detailed on~\figref{fig:smt_codes}. 
The new major features implemented in \texttt{SMT 2.0} are highlighted in lavender whereas the legacy features that were already in present in the original publication for \texttt{SMT 0.2}~\cite{SMT2019} are in black.

\subsection{New features within \texttt{SMT 2.0} } 
\label{sec:newSMT2}
The main objective of this new release is to enable Kriging surrogate models for use with both hierarchical and mixed variables.
Moreover, for each of these five sub-modules described in Section~\ref{sec:organization}, several improvements have been made between the original version and the \texttt{SMT 2.0} release.

\paragraph{Hierarchical and mixed design space}
A new design space definition class \texttt{DesignSpace} has been added that implements hierarchical and mixed functionalities. 
Design variables can either be continuous ({\footnotesize \texttt{FloatVariable}}), ordered ({\footnotesize \texttt{OrdinalVariable}}) or categorical ({\footnotesize \texttt{CategoricalVariable}}). 
The integer type ({\footnotesize \texttt{IntegerVariable}}) represents a special case of the ordered variable, specified by bounds (inclusive) rather than a list of possible values.
The hierarchical structure of the design space can be defined using {\footnotesize \texttt{declare\_decreed\_var}}: this function declares that a variable is a decreed variable that is activated when the associated meta variable takes one of a set of specified values, see Section~\ref{sec:hier} for background.
The \texttt{DesignSpace} class also implements mechanisms for sampling valid design vectors (i.e. design vectors that adhere to the hierarchical structure of the design space) using any of the below-mentioned samplers, for correcting and imputing design vectors, and for requesting which design variables are acting in a given design vector. 
Correction ensures that variables have valid values (\textit{e.g.} integers for discrete variables)~\cite{Effectiveness}, and imputation replaces non-acting variables by some default value ($0$ for discrete variables, mid-way between the bounds for continuous variables in \texttt{SMT 2.0})~\cite{Zaefferer}.

\paragraph{Sampling} 
SMT implements three methods for sampling.
The first one is a naïve approach, called \texttt{Random} that draws uniformly points along every dimension.
The second sampling method is called \texttt{Full Factorial} and draws a point for every cross combination of variables, to have an "exhaustive" design of experiments. 
The last one is the \texttt{Latin Hypercube Sampling} (LHS)~\cite{LHS} that draws a point in every Latin square parameterized by a certain criterion. 
For \gls{LHS}, a new criterion to manage the randomness has been implemented and the sampling method was adapted for multi-fidelity and mixed or hierarchical variables.
More details about the new sampling techniques are given in Section~\ref{sec:sampling}.

\paragraph{Problems} 
SMT implements two new engineering problems: a mixed variant of a cantilever beam described in Section~\ref{sec:mixed} and a hierarchical neural network described in Section~\ref{sec:hier}.

\paragraph{Surrogate models} 
In order to keep up with state-of-art, several releases done from the original version developed new options for the already existing surrogates.
In particular, compared to the original publication~\cite{SMT2019}, \texttt{SMT 2.0} adds gradient-enhanced neural networks~\cite{bouhlel2020} and marginal Gaussian process~\cite{MGP} models to the list of available surrogates. 
More details about the new models are given in  Section~\ref{sec:surrogates}.

\paragraph{Applications} 
Several applications have been added to the toolbox to demonstrate the surrogate models capabilities. 
The most relevant application is efficient global optimization (EGO), a Bayesian optimization algorithm~\cite{Jones2001JOGO,lafage2022egobox}. 
EGO optimizes expensive-to-evaluate black-box problems with a chosen surrogate model and a chosen optimization criterion~\cite{Jones98}.
The usage of EGO with hierarchical and mixed variables is described in Section~\ref{sec:BO}.

\paragraph{Interactive notebooks} 
These tutorials introduce and explain how to use the toolbox for different surrogate models and applications\footnote{\url{ https://github.com/SMTorg/smt/tree/master/tutorial}}.
Every tutorial is available both as a \texttt{.ipynb} file and directly on Google colab\footnote{\url{https://colab.research.google.com/github/SMTorg/smt/ }}. 
In particular, a hierarchical and mixed variables dedicated notebook is available to reproduce the results presented on this chapter\footnote{\url{https://github.com/SMTorg/smt/tree/master/tutorial/NotebookRunTestCases_Paper_SMT_v2.ipynb}}.

In the following, Section~\ref{sec:mixed} details the Kriging based surrogate models for mixed variables, and Section~\ref{sec:hier} presents our new Kriging surrogate for hierarchical variables.
Section~\ref{sec:BO} details the EGO application and the other new relevant features aforementioned are described succinctly in Section~\ref{sec:other}.

\section{Surrogate models with mixed variables in \texttt{SMT 2.0}}
\label{sec:mixed}

As mentioned in Section~\ref{sec:introSMT}, design variables can be either of continuous or discrete type, and a problem with both types is a mixed-variable problem.
Discrete variables can be ordinal or categorical.  
A discrete variable is \emph{ordinal} if there is an order relation within the set of possible values. 
An example of an ordinal design variable is the number of engines in an aircraft.
A possible set of values in this case could be ${2, 4, 8}$.
A discrete variable is \emph{categorical} if no order relation is known between the possible choices the variable can take. 
One example of a categorical variable is the color of a surface.
A possible example of a set of choices could be ${\text{blue}, \text{red}, \text{green}}$.
The possible choices are called the \emph{levels} of the variable.


Several methods have been proposed to address the recent increase interest in mixed Kriging based models~\cite{Pelamatti, Zhou, Deng, Roustant,GMHL,Gower,cuesta2021comparison,SciTech_cat}. 
The main difference from a continuous Kriging model is in the estimation of the categorical correlation matrix, which is critical to determine the mean and variance predictions.
As mentioned in Section~\ref{sec:introSMT}, approaches such as CR~\cite{GMHL,SciTech_cat}, continuous latent variables~\cite{cuesta2021comparison}, and GD~\cite{Gower} use a kernel-based method to estimate the correlation matrix.
Other methods estimate the correlation matrix by modeling the correlation entries directly~\cite{Pelamatti, Deng, Roustant}, such as HH~\cite{Zhou} and EHH~\cite{Mixed_Paul}. 
The HH correlation kernel is of particular interest because it generalizes simpler kernels such as CR and GD~\cite{Mixed_Paul}.
In \texttt{SMT 2.0}, the correlation kernel is an option that can be set to either CR (\texttt{CONT\_RELAX\_KERNEL}), GD ( {\texttt{GOWER\_KERNEL}), HH (\texttt{HOMO\_HSPHERE\_KERNEL}) or EHH (\texttt{EXP\_HOMO\_HSPHERE\_KERNEL}).

\subsection{Mixed Gaussian processes}

The continuous and ordinal variables are both treated similarly in \texttt{SMT 2.0} with a continuous kernel, where the ordinal values are converted to continuous through relaxation. 
For categorical variables, four models (GD, CR, EHH and HH) can be used in \texttt{SMT 2.0} if specified by the API. 
This is why we developed a unified mathematical formulation that allows a unique implementation for any model.   

Denote $l$ the number of categorical variables.
For a given  $i \in \{1, \ldots, l\}$, the $ i^{\text{th}}$ categorical variable is denoted $c_i$ and its number of levels is denoted $L_i$. 
The hyperparameter matrix peculiar to this variable $c_i$ is 
$$\Theta_i= \begin{bmatrix}
[\Theta_i]_{1,1} & \textcolor{white}{9} & \hspace{2em} { \textbf{\textit{ Sym.}}}  \textcolor{white}{9} & \\
[\Theta_i]_{1,2}  & [\Theta_i]_{2,2} & \textcolor{white}{9} \\
\vdots &\ddots & \ddots & \textcolor{white}{9}  \\
[\Theta_i]_{1,L_i} &  \ldots & [\Theta_i]_{L_i-1,L_i} &[\Theta_i]_{L_i,L_i} \\ 
\end{bmatrix}, $$
and the categorical parameters are defined as $\theta^{cat} = \{ \Theta_1 , \ldots, \Theta_l \}$.
For two given inputs in the DoE, for example, the $r^{\text{th}}$ and  $s^{\text{th}}$ points, let $c^r_{i} $ and $c^s_{i} $ be the associated categorical variables taking respectively the $\ell^i_r$ and the $\ell^i_s$ level on the categorical variable $c_i$. The categorical correlation kernel is defined by
\begin{equation}
\begin{split}
&k^{cat}(c^r,c^s,\theta^{cat}) = {\displaystyle \prod_{i=1}^{l}  \kappa (2 [ \Phi(\Theta_i) ]_{{ \ell_i^r},{\ell_i^s}} ) \  \kappa ( [ \Phi(\Theta_i) ]_{{ \ell_i^r},{\ell_i^r}} ) \  \kappa ( [ \Phi(\Theta_i) ]_{{ \ell_i^s},{\ell_i^s}} ) }
\end{split}
\end{equation}
where $\kappa$ is either a positive definite kernel or identity and $\Phi(.)$ is a symmetric positive definite (SPD) function such that the matrix $\Phi(\Theta_i) $ is SPD if $\Theta_i$ is SPD.
For an exponential kernel,~\tabref{tab:kernelsSMT} gives the parameterizations of $\Phi$ and $\kappa$ that correspond to GD, CR, HH, and EHH kernels. 
The complexity of these different kernels depends on the number of hyperparameters that characterizes them. 
As defined by \citet{Mixed_Paul}, for every categorical variable $i \in \{1, \ldots, l\}$, the matrix $C(\Theta_i)\in \mathbb{R}^{L_i \times L_i}$ is lower triangular and built using a hypersphere decomposition~\cite{HS,HS_Jacobi} from the symmetric matrix $\Theta_i \in \mathbb{R}^{L_i \times L_i}$ of hyperparameters. 
The variable $\epsilon$ is a small positive constant and the variable $\theta_{i}$ denotes the only positive hyperparameter that is used for the Gower distance kernel.

\begin{table}[H]
\caption{Categorical kernels implemented in \texttt{SMT 2.0}.}
\vspace{-0.5cm}
\begin{center}
\resizebox{\columnwidth}{!}{%
\begin{tabular}{cclc}
\hline
\textbf{Name} & $\kappa(\phi)$   &  \hspace{3cm} ${\centering \Phi(\Theta_i)}$   &  \textcolor{black}{\# of hyperparam.} \\
\hline
\texttt{SMT GD}   &  $\exp(\textcolor{black}{-}\phi) $ & ${ \displaystyle [\Phi(\Theta_i)]_{j,j} :=  \frac{1}{2}  \theta_{i} \quad  ~;~ [\Phi(\Theta_i)]_{j \neq j'} := 0 }$ & \textcolor{black}{1}   \\
\texttt{SMT CR}  & $\exp(\textcolor{black}{-}\phi) $ &  $  { \displaystyle [\Phi(\Theta_i)]_{j,j} := [\Theta_i]_{j,j} ~;~ [\Phi(\Theta_i)]_{j \neq j'} := 0 } $  & \textcolor{black}{$L_i$}  \\
\texttt{SMT EHH}  & $\exp(\textcolor{black}{-}\phi)$ & 
$  { \displaystyle [\Phi(\Theta_i)]_{j,j} := 0 \quad \quad ~;~ [\Phi(\Theta_i)]_{j \neq j'} := \frac{\log \epsilon }{2} ([C(\Theta_i) C(\Theta_i) ^\top]_{j,j'} -1)  }$  & \textcolor{black}{$\frac{1}{2}  (L_i)  (L_i-1) $}\\
\texttt{SMT HH}  &  $\phi$ &    $  { \displaystyle [\Phi(\Theta_i)]_{j,j} := 1 \quad \quad ~;~ [\Phi(\Theta_i)]_{j \neq j'} := \frac{1}{2} [C(\Theta_i) C(\Theta_i)^\top]_{j,j'} }$ & \textcolor{black}{$\frac{1}{2}  (L_i)  (L_i-1) $}  \\
\hline
\end{tabular}
}
\end{center}
\label{tab:kernelsSMT}
\end{table}

Another Kriging based model that can use mixed variables is Kriging with partial least squares (KPLS)~\cite{bouhlel_KPLSK}.
KPLS adapts Kriging to high dimensional problems by using a reduced number of hyperparameters thanks to a projection into a smaller space. 
Also, for a general surrogate, not necessarily Kriging, \texttt{SMT 2.0} uses continuous relaxation to allow whatever model to handle mixed variables. 
For example, we can use mixed variables with least squares (LS) or quadratic polynomial (QP) models.
We now illustrate the abilities of the toolbox in terms of mixed modeling over an engineering test case.

\subsection{An engineering design test-case}
\label{sec:beam}
A classic engineering problem commonly used for model validation is the beam bending problem~\cite{Roustant, Cheng2015TrustRB}.
This problem is illustrated on~\figref{fig:beam} and consists of a cantilever beam in its linear range loaded at its free end with a force $F$. 
As in~\citet{Cheng2015TrustRB}, the Young modulus is  $E=200$GPa and the chosen load is $F=50$kN.
Also, as in~\citet{Roustant}, 12 possible cross-sections can be used. These 12 sections consist of 4 possible shapes that can be either hollow, thick or full as illustrated in~\figref{fig:beam_shape}.

\begin{figure}[H]
\centering
\vspace{-5pt}
\captionsetup{justification=raggedright,singlelinecheck=false}
\subfloat[Bending problem.]{
\begin{tikzpicture}
    \hspace{-3pt}
    \point{origin}{-0.75}{-0.25};
    \point{begin}{0}{0};
    \point{end}{5}{0};
    \point{end_bot}{4.99}{-0.9};
    \point{end_up}{5}{0.5};
    \beam{2}{begin}{end};
    \support{3}{begin}[-90];
    \load{1}{end}[90]   ;
    \notation{1}{end_up}{$F=50kN$};

     \draw[<->] (end) -- (end_bot) node[midway, right] {$\delta$} ;
     \draw[<->] (0,0.5) -- (5,0.5) node[midway, above] {L};
     
    \draw
      [-, ultra thick] (begin) .. controls (1.5, +.01) and (2.5, -.15) .. (4.93, -0.9)
      [-, ultra thick] (begin) .. controls (1.5, +.01) and (2.5, -.2) .. (4.85, -1.5)
      [-, ultra thick] (begin) .. controls (1.5, +.01) and (2.5, -.4)   .. (4.78, -1.9);
  \end{tikzpicture}
\label{fig:beam}   
}
\subfloat[Possible cross-section shapes.]{
\centering
\begin{tikzpicture}
\hspace{-50pt}


\tstar{0.25}{0.5}{6}{0}{thick,fill=yellow,xshift=+3.6cm,yshift= -1.2cm}
\tstar{0.14}{0.28}{6}{0}{thick,fill=white,xshift=+3.6cm,yshift= -1.2cm}

\tstar{0.25}{0.5}{6}{0}{thick,fill=yellow,xshift=+2.4cm,yshift= -1.2cm}
\tstar{0.08}{0.16}{6}{0}{thick,fill=white,xshift=+2.4cm,yshift= -1.2cm}

\tstar{0.25}{0.5}{6}{0}{thick,fill=yellow,xshift=+1.2cm,yshift= -1.2cm}

\fill[green,even odd rule] (3.6,0) circle (0.5) (3.6,0) circle (0.33);
\draw (3.6,0) circle (0.5) ;
\draw (3.6,0) circle (0.33) 
; 
\fill[green,even odd rule] (2.4,0) circle (0.5)(2.4,0) circle (0.17);
\draw (2.4,0) circle (0.5) ;
\draw (2.4,0) circle (0.17) 
; 
\fill[green,even odd rule] (1.2,0) circle (0.5) ;
\draw (1.2,0) circle (0.5) ;

\def\pos{-2.4}
\fill[blue,even odd rule]  (\pos-0.5,-1.7+1.2) -- (\pos-0.5,-0.7+1.2) -- (\pos+0.5,-0.7+1.2) -- (\pos+0.5,-1.7+1.2) -- cycle ;

\def\pos{-1.2}
\fill[blue,even odd rule]  (\pos-0.5,-1.7+1.2) -- (\pos-0.5,-0.7+1.2) -- (\pos+0.5,-0.7+1.2) -- (\pos+0.5,-1.7+1.2) -- cycle   (\pos-0.25,-1.45+1.2) -- (\pos-0.25,-0.95+1.2) -- (\pos+0.25,-0.95+1.2) -- (\pos+0.25,-1.45+1.2) -- cycle ;

\def\pos{0}
\fill[blue,even odd rule]  (\pos-0.5,-1.7+1.2) -- (\pos-0.5,-0.7+1.2) -- (\pos+0.5,-0.7+1.2) -- (\pos+0.5,-1.7+1.2) -- cycle   (\pos-0.35,-1.55+1.2) -- (\pos-0.35,-0.85+1.2) -- (\pos+0.35,-0.85+1.2) -- (\pos+0.35,-1.55+1.2) -- cycle ;

\def\pos{-2.4}
\fill[black] (\pos-0.24,-0.5-1.2) -- (\pos-0.24,0.5-1.2) -- (\pos+0.24,0.5-1.2)  -- (\pos+0.24,-0.5-1.2)   -- cycle ;
\fill[black] (\pos-0.5,-0.5-1.2) -- (\pos-0.5,-0.18-1.2) -- (\pos+0.5,-0.18-1.2)  -- (\pos+0.5,-0.5-1.2)   -- cycle ; 
\fill[black] (\pos-0.5,0.18-1.2) -- (\pos-0.5,0.5-1.2) -- (\pos+0.5,0.5-1.2)  -- (\pos+0.5,0.18-1.2)   -- cycle ;  

\def\pos{-1.2}
\fill[black] (\pos-0.19,-0.5-1.2) -- (\pos-0.19,0.5-1.2) -- (\pos+0.19,0.5-1.2)  -- (\pos+0.19,-0.5-1.2)   -- cycle ;
\fill[black] (\pos-0.5,-0.5-1.2) -- (\pos-0.5,-0.25-1.2) -- (\pos+0.5,-0.25-1.2)  -- (\pos+0.5,-0.5-1.2)   -- cycle ; 
\fill[black] (\pos-0.5,0.25-1.2) -- (\pos-0.5,0.5-1.2) -- (\pos+0.5,0.5-1.2)  -- (\pos+0.5,0.25-1.2)   -- cycle ;  

\def\pos{0}
\fill[black] (\pos-0.14,-0.5-1.2) -- (\pos-0.14,0.5-1.2) -- (\pos+0.14,0.5-1.2)  -- (\pos+0.14,-0.5-1.2)   -- cycle ;
\fill[black] (\pos-0.5,-0.5-1.2) -- (\pos-0.5,-0.32-1.2) -- (\pos+0.5,-0.32-1.2)  -- (\pos+0.5,-0.5-1.2)   -- cycle ; 
\fill[black] (\pos-0.5,0.32-1.2) -- (\pos-0.5,0.5-1.2) -- (\pos+0.5,0.5-1.2)  -- (\pos+0.5,0.32-1.2)   -- cycle ;  

\point{un}{-2.15}{-0.8};
\notation{1}{un}{\tiny 1};
\point{deux}{-2.15+1.2}{-0.8};
\notation{1}{deux}{\tiny 2};
\point{trois}{-2.15+2.4}{-0.8};
\notation{1}{trois}{\tiny 3};

\point{quatre}{-2.15+3.6}{-0.8};
\notation{1}{quatre}{\tiny 4};
\point{cinq}{-2.15+4.8}{-0.8};
\notation{1}{cinq}{\tiny 5};
\point{six}{-2.15+6}{-0.8};
\notation{1}{six}{\tiny 6};

\point{sept}{-2.15}{-0.8-1.2};
\notation{1}{sept}{\tiny 7};
\point{huit}{-2.15+1.2}{-0.8-1.2};
\notation{1}{huit}{\tiny 8};
\point{neuf}{-2.15+2.4}{-0.8-1.2};
\notation{1}{neuf}{\tiny 9};

\point{dix}{-2.15+3.6}{-0.8-1.2};
\notation{1}{dix}{\tiny 10};
\point{onze}{-2.15+4.8}{-0.8-1.2};
\notation{1}{onze}{\tiny 11};
\point{douze}{-2.15+6}{-0.8-1.2};
\notation{1}{douze}{\tiny 12};

  \end{tikzpicture}
  \label{fig:beam_shape}    
}
\captionsetup{justification=centering,singlelinecheck=false}
\caption{Cantilever beam problem~\cite[Figure 6]{Mixed_Paul}.}
\end{figure}

To compare the mixed Kriging models of \texttt{SMT 2.0}, we draw a 98 point LHS as training set and the validation set is a grid of $12\times30\times30=10800$ points. 
For the four implemented methods, displacement error (computed with a root-mean-square error criterion), likelihood, number of hyperparameters and computational time for every model are shown in~\tabref{tab:resCantileverSMT}.
For the continuous variables, we use the square exponential kernel.
More details are found in~\cite{Mixed_Paul}.
As expected, the complex EHH and HH models lead to a lower displacement error and a higher likelihood value, but use more hyperparameters and increase the computational cost compared to GD and CR. 
On this test case, the kernel EHH is easier to optimize than HH but in general, they are similar in terms of performance. 
Also, by default \texttt{SMT 2.0} uses CR as it is known to be a good trade-off between complexity and performance~\cite{well-adapted_cont}.

\begin{table}[H]
\centering
\caption{Results of the cantilever beam models~\cite[Table 4]{Mixed_Paul}.}
\resizebox{\columnwidth}{!}{%
\begin{tabular}{ccccc}
\hline
\textbf{Categorical kernel} &  Displacement error (cm) & $ \ $ Likelihood  & \# of hyperparam.  
\\
\hline 
\texttt{SMT GD}   &1.3861 & 111.13&  3 
\\  
\texttt{SMT CR}   & 1.1671 & 155.32 & 14  
\\
\texttt{SMT EHH}   & {0.1613}   &
{236.25} & 68 
\\
\texttt{SMT HH}   & {0.2033} &
{235.66} & 68 
\\
\hline
\end{tabular}
}
\label{tab:resCantileverSMT}
\end{table}

\section{Surrogate models with hierarchical variables in \texttt{SMT 2.0}}
\label{sec:hier}

To introduce the newly developed Kriging model for hierarchical variables implemented in \texttt{SMT 2.0}, we present the general mathematical framework for hierarchical and mixed variables established by \citet{audet2022general}.
In \texttt{SMT 2.0}, two variants of our new method are implemented, namely \texttt{SMT Alg-Kernel} and \texttt{SMT Arc-Kernel}.
In particular, the \texttt{SMT Alg-Kernel} is a novel correlation kernel introduced in this chapter.

\subsection{The hierarchical variables framework}

A problem structure is classified as hierarchical when the sets of active variables depend on architectural choices.
This occurs frequently in industrial design problems.
In hierarchical problems, we can classify variables as neutral, meta (also known as dimensional) or decreed (also known as conditionally active) as detailed in~\citet{audet2022general}.
Neutral variables are the variables that are not affected by the hierarchy whereas the value assigned to meta variables determines which decreed variables are activated. For example, a meta variable could be the number of engines. If the number of engines changes, the number of decreed bypass ratios that every engine should specify also changes. However, the wing aspect ratio being neutral, it is not affected by this hierarchy. 

%

Problems involving hierarchical variables are generally dependant on discrete architectures and as such involve mixed variables. 
Hence, in addition to their role (neutral, meta or decreed), each variable also has a variable type amongst categorical, ordinal or continuous. 
For the sake of simplicity and because both continuous and ordinal variables are treated similarly~\cite{Mixed_Paul}, we chose to regroup them as \emph{quantitative variables}.
For instance, the neutral variables $x_{\neutral}$ may be partitioned into different variable types, such that $x_{\neutral} = (x_{\neutral}^{\cat},x_{\neutral}^{\quant})$ where $x_{\neutral}^{\cat}$ represents the categorical variables and $x_{\neutral}^{\quant}$ are the quantitative ones. 
The variable classification scheme in \texttt{SMT 2.0} is detailed in~\figref{fig:vars}.
%
\begin{figure}
\centering
\includegraphics[height=8.5cm,,width=14cm]{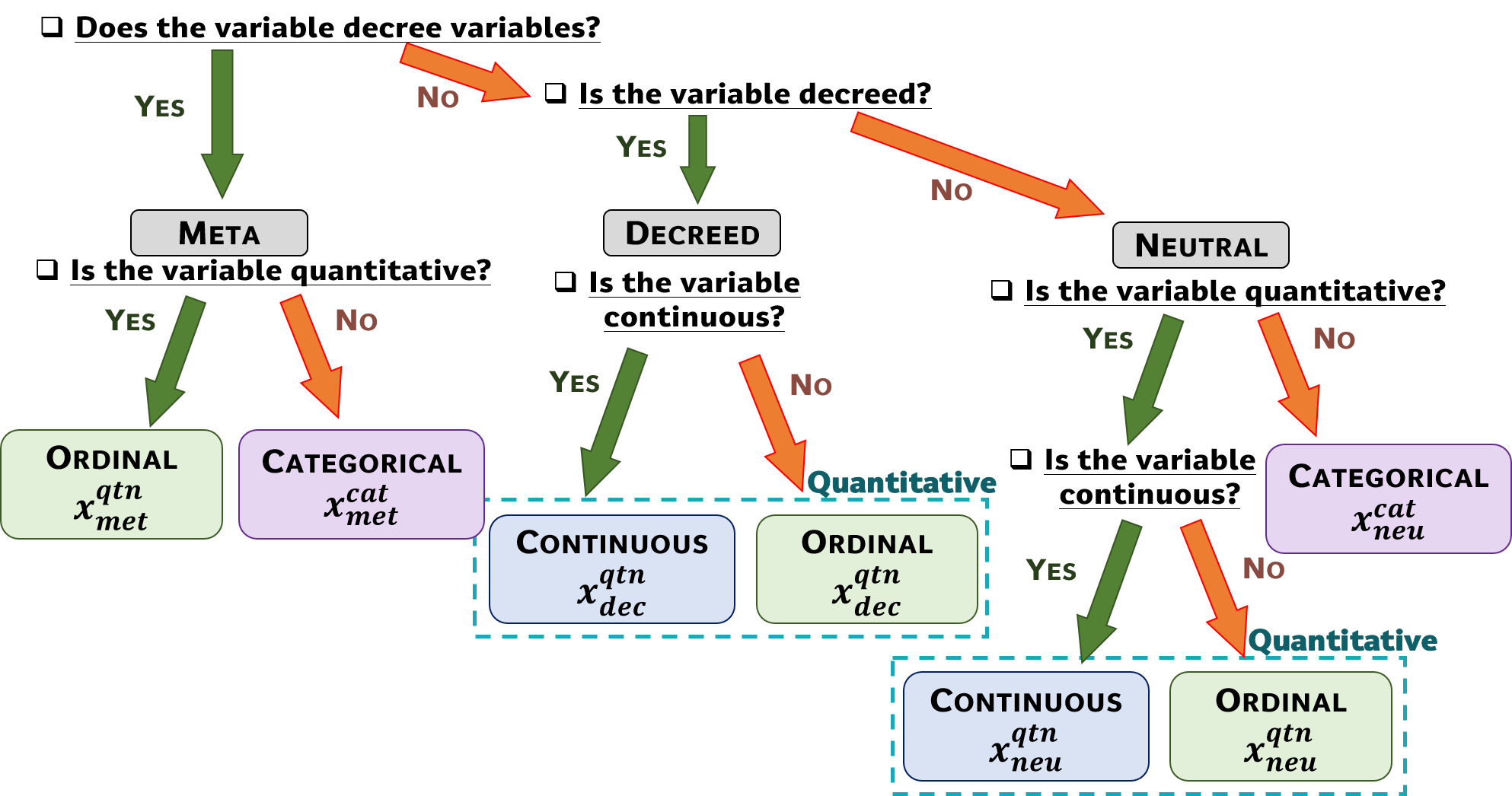}
\caption{Variables classification as used in \texttt{SMT 2.0}.
}
\label{fig:vars}
\end{figure}

To explain the framework and the new Kriging model, we illustrate the \textcolor{black}{inputs variables of the model using a classical machine learning problem related to the hyperparameters optimization of a fully-connected Multi-Layer Perceptron (MLP)~\cite{audet2022general} on~\figref{fig:MLP}}.
\textcolor{black}{In~\tabref{tab:hyp_NN}, we detail the input variables of the model related to the MLP problem (i.e., the hyperparameters of the neural network, together with their types and roles). To keep things clear and concise, the chosen problem is a simplification of the original problem developed by \citet{audet2022general}.}
\textcolor{black}{
Regarding the MLP problem of~\figref{fig:MLP} and following the classification scheme of~\figref{fig:vars}, we start by separating the input variables according to their role. In fact,
\begin{enumerate}
    \item changing the number of hidden layers modifies the number of inputs variables. Therefore, "\# of hidden layers" is a meta variable. 
    \item The number of neurons in the hidden layer number $k$ is either included or excluded. For example, the "\# of neurons in the $3^\text{rd}$ layer" would be excluded for an input that only has $2$ hidden layers. Therefore, "\# of neurons hidden layer $k$" are decreed variables.
    \item The "Learning rate", "Momentum", "Activation function" and "Batch size" are not affected  by the hierarchy choice. Therefore, they are neutral variables. 
\end{enumerate}
According to their types, the MLP input variables can be classified as follows:
\begin{enumerate}[start=4]
    \item The meta variable "\# of hidden layers" is an integer and, as such, is represented by the component $x^{\quant}_{\meta}$. 
    \item The decreed variables "\# of neurons hidden layer $k$" are integers and, as such, are represented by the component $x^{\quant}_{\decreed}$. 
      \item The "Learning rate", {"Momentum"},  "Activation function" and "Batch size" are, respectively, continuous, {for the first two} (every value between two bounds), categorical (qualitative between {three} choices) and integer (quantitative between 6 choices). Therefore, the "Activation function" {and the "Momentum" are} represented by the component $x^{\cat}_{\neutral}$. The "Learning rate" and the "Batch size" are represented by the component $x^{\quant}_{\neutral}$.
\end{enumerate}
}
\begin{table}[H]
\centering 
\caption{A detailed description of the variables in the MLP problem.}
\resizebox{\linewidth}{!}{%
\begin{tabular}{@{}llclcc@{}}
\toprule
\multicolumn{2}{l}{\textbf{MLP Hyperparameters}}    & Variable &  Domain     & Type   & Role  \\ \midrule
\multicolumn{2}{l}{\textbf{Learning rate}}                & $r$      & $[10^{-5}, 10^{-2}]$       & FLOAT &  NEUTRAL  \\
\multicolumn{2}{l}{\textbf{{Momentum}}}                & ${\alpha}$      & ${[0, 1]}$       & {FLOAT} &  {NEUTRAL}  \\
\multicolumn{2}{l}{\textbf{Activation function}}          & $a$      & $ \{$ReLU, Sigmoid, {Tanh}$\}$         &  ENUM & NEUTRAL \\
\multicolumn{2}{l}{\textbf{Batch size}}                    & $b$        & $ \{8, 16, \ldots, 128, 256 \} $     & ORD      &  NEUTRAL  \\  
\multicolumn{2}{l}{\textbf{\# of hidden layers}}          & $l$      & $\{1,2,3\} $    &  ORD       & META       \\
\multicolumn{2}{l}{\textbf{\# of neurons hidden layer $k$}} & $n_{k}$ & $ \{50, 51 , \ldots, 55 \} $      & ORD & DECREED \\
\bottomrule
\end{tabular}
}
\label{tab:hyp_NN}
\end{table}
\begin{figure}[H]
\resizebox{\columnwidth}{!}{
\centering
\begin{tikzpicture}


\node [draw, circle, minimum size=0.75cm] (N1c) at (0,0) {\large $\alpha_i^{(in)}$};
\node [draw, circle, minimum size=0.75cm, xshift=0cm, yshift=3cm, at=(N1c), label={\hspace{-0.5cm} \vspace{0.5cm} \Large Input layer}] (N1u) {\large $\alpha_1^{(in)}$};
\node [draw, circle, minimum size=0.75cm, xshift=0cm, yshift=-3cm, at=(N1c)] (N1b) {\large $\alpha_p^{(in)}$};

\node [at=($(N1c)!0.5!(N1b)$)] {\Huge \vdots} ;
\node [at=($(N1c)!0.5!(N1u)$)] {\Huge \vdots} ;


\node [draw, circle, minimum size=0.75cm, xshift=3.5cm, yshift=0cm, at=(N1c)] (N2c) {};
\node [draw, circle, minimum size=0.75cm, xshift=0cm, yshift=3cm, at=(N2c)] (N2u) {};
\node [draw, circle, minimum size=0.75cm, xshift=0cm, yshift=-3cm, at=(N2c)] (N2b) {};

\node [at=($(N2c)!0.5!(N2b)$)] {\Huge \vdots} ;
\node [at=($(N2c)!0.5!(N2u)$)] {\Huge \vdots} ;

\node [draw, circle, minimum size=0.75cm, xshift=2.5cm, yshift=0cm, at=(N2c)] (N3c) {};
\node [draw, circle, minimum size=0.75cm, xshift=0cm, yshift=3cm, at=(N3c)] (N3u) {};
\node [draw, circle, minimum size=0.75cm, xshift=0cm, yshift=-3cm, at=(N3c)] (N3b) {};

\node [at=($(N3c)!0.5!(N3b)$)] {\Huge \vdots} ;
\node [at=($(N3c)!0.5!(N3u)$)] {\Huge \vdots} ;

\node [draw, circle, minimum size=0.75cm, xshift=5cm, yshift=0cm, at=(N3c)] (N6c) {}; 

\node [draw, circle, minimum size=0.75cm, xshift=0cm, yshift=3cm, at=(N6c)] (N6u) {};

\node [draw, circle, minimum size=0.75cm, xshift=0cm, yshift=-3cm, at=(N6c)] (N6b) {};

\node [at=($(N6c)!0.5!(N6b)$)] {\Huge \vdots} ;
\node [at=($(N6c)!0.5!(N6u)$)] {\Huge \vdots} ;

\node [draw, circle, minimum size=0.75cm, xshift=2.5cm, yshift=0cm, at=(N6c)] (N7c) {}; 
\node [draw, fill=black!25, circle, minimum size=0.75cm, xshift=0cm, yshift=3cm, at=(N7c)] (N7u) {};
\node [draw, circle, minimum size=0.75cm, xshift=0cm, yshift=-3cm, at=(N7c)] (N7b) {};

\node [at=($(N7c)!0.43!(N7b)$)] {\Huge \vdots} ;
\node [at=($(N7c)!0.57!(N7u)$)] {\Huge \vdots} ;


\node [draw, circle, minimum size=1cm, xshift=3cm, yshift=0cm, at=(N7c), label={ \Large Output}] (N8c) {\large $\psi \left(\alpha^{(in)}\right )$}; 

\draw[ -{Stealth[length=2mm, width=1.5mm]} ] (N7c)--(N8c);
\draw[ -{Stealth[length=2mm, width=1.5mm]} ] (N7u)--(N8c);
\draw[ -{Stealth[length=2mm, width=1.5mm]}] (N7b)--(N8c);

\node [ circle, minimum size=0.75cm, at=($(N3c)!.33!(N6c)$)] (N4c) {};
\node [ circle, minimum size=0.75cm, at=($(N3u)!.33!(N6u)$)] (N4u) {};
\node [ circle, minimum size=0.75cm, at=($(N3b)!.33!(N6b)$)] (N4b) {};

\node [ circle, minimum size=1cm, at=($(N3c)!.66!(N6c)$)] (N5c) {}; 
\node [ circle, minimum size=1cm, at=($(N3u)!.66!(N6u)$)] (N5u) {};
\node [ circle, minimum size=1cm, at=($(N3b)!.66!(N6b)$)] (N5b) {};


\draw[ -{Stealth[length=2mm, width=1.5mm]} ] (N1c)--(N2c);
\draw[ -{Stealth[length=2mm, width=1.5mm]} ] (N1c)--(N2u);
\draw[ -{Stealth[length=2mm, width=1.5mm]} ] (N1c)--(N2b);

\draw[ -{Stealth[length=2mm, width=1.5mm]} ] (N1u)--(N2c);
\draw[ -{Stealth[length=2mm, width=1.5mm]} ] (N1u)--(N2u);
\draw[ -{Stealth[length=2mm, width=1.5mm]} ] (N1u)--(N2b);

\draw[ -{Stealth[length=2mm, width=1.5mm]} ] (N1b)--(N2c);
\draw[ -{Stealth[length=2mm, width=1.5mm]} ] (N1b)--(N2u);
\draw[ -{Stealth[length=2mm, width=1.5mm]} ] (N1b)--(N2b);

\draw[ -{Stealth[length=2mm, width=1.5mm]} ] (N2c)--(N3c);
\draw[ -{Stealth[length=2mm, width=1.5mm]} ] (N2c)--(N3u);
\draw[ -{Stealth[length=2mm, width=1.5mm]} ] (N2c)--(N3b);

\draw[ -{Stealth[length=2mm, width=1.5mm]} ] (N2u)--(N3c);
\draw[ -{Stealth[length=2mm, width=1.5mm]} ] (N2u)--(N3u);
\draw[ -{Stealth[length=2mm, width=1.5mm]} ] (N2u)--(N3b);

\draw[ -{Stealth[length=2mm, width=1.5mm]} ] (N2b)--(N3c);
\draw[ -{Stealth[length=2mm, width=1.5mm]} ] (N2b)--(N3u);
\draw[ -{Stealth[length=2mm, width=1.5mm]} ] (N2b)--(N3b);

\draw[ - ] (N3c) -- ($(N3c)!.66!(N4c)$); 
\draw[ - ] (N3c) -- ($(N3c)!.33!(N4u)$);
\draw[ - ] (N3c) -- ($(N3c)!.33!(N4b)$);

\draw[ - ] (N3u)-- ($(N3u)!.33!(N4c)$);
\draw[ - ] (N3u)-- ($(N3u)!.66!(N4u)$); 
\draw[ - ] (N3u)-- ($(N3u)!.25!(N4b)$); 

\draw[ - ] (N3b)-- ($(N3b)!.33!(N4c)$);
\draw[ - ] (N3b)-- ($(N3b)!.25!(N4u)$); 
\draw[ - ] (N3b)-- ($(N3b)!.66!(N4b)$); 

\node [at=($(N4c)!0.5!(N5c)$)] {\large \dots } ;
\node [at=($(N4u)!0.5!(N5u)$), rotate=160] {\large \dots } ;
\node [at=($(N4b)!0.5!(N5b)$), rotate=20] {\large \dots } ;

\draw[ -{Stealth[length=2mm, width=1.5mm]}  ] ($(N5c)!.5!(N6c)$)--(N6c); 
\draw[ -{Stealth[length=2mm, width=1.5mm]}  ] ($(N5c)!.66!(N6u)$)--(N6u);
\draw[ -{Stealth[length=2mm, width=1.5mm]}  ] ($(N5c)!.66!(N6b)$)--(N6b);

\draw[ -{Stealth[length=2mm, width=1.5mm]}  ] ($(N5u)!.75!(N6c)$)--(N6c);
\draw[ -{Stealth[length=2mm, width=1.5mm]}  ] ($(N5u)!.5!(N6u)$)--(N6u); 
\draw[ -{Stealth[length=2mm, width=1.5mm]}  ] ($(N5u)!.8!(N6b)$)--(N6b); 

\draw[ -{Stealth[length=2mm, width=1.5mm]}  ] ($(N5b)!.75!(N6c)$)--(N6c);
\draw[ -{Stealth[length=2mm, width=1.5mm]}  ] ($(N5b)!.8!(N6u)$)--(N6u); 
\draw[ -{Stealth[length=2mm, width=1.5mm]}  ] ($(N5b)!.5!(N6b)$)--(N6b); 

\draw[ -{Stealth[length=2mm, width=1.5mm]} ] (N6c)--(N7c);
\draw[ -{Stealth[length=2mm, width=1.5mm]} ] (N6c)--(N7u);
\draw[ -{Stealth[length=2mm, width=1.5mm]} ] (N6c)--(N7b);

\draw[ -{Stealth[length=2mm, width=1.5mm]} ] (N6u)--(N7c);
\draw[ -{Stealth[length=2mm, width=1.5mm]} ] (N6u)--(N7u);
\draw[ -{Stealth[length=2mm, width=1.5mm]} ] (N6u)--(N7b);

\draw[ -{Stealth[length=2mm, width=1.5mm]} ] (N6b)--(N7c);
\draw[ -{Stealth[length=2mm, width=1.5mm]} ] (N6b)--(N7u);
\draw[ -{Stealth[length=2mm, width=1.5mm]} ] (N6b)--(N7b);

\node [circle, minimum size=1cm] (T1u) at ($(N2u)-(0.8,-0.8)$)  {};
\node [circle, minimum size=1cm] (T1b) at ($(N2b)-(0.8,0.8)$)  {};
\node [circle, minimum size=1cm] (T2u) at ($(N7u)+(0.8,0.8)$)  {};
\node [circle, minimum size=1cm] (T2b) at ($(N7b)+(0.8,-0.8)$)  {};
\draw[dotted, line width= 0.4mm] ($(T1u)+(0,0)$)--($(T1b)+(0,0)$);
\draw[dotted, line width= 0.4mm] ($(T2u)+(0,0)$)--($(T2b)+(0,0)$);

\draw[dotted, line width= 0.4mm] ($(T1u)+(0,0)$)--($(T2u)+(0,0)$) node[midway, above, black, 
]{\Large Hidden layers $l$} ;

\draw[dotted, line width= 0.4mm] ($(T1b)+(0,0)$)--($(T2b)+(0,0)$);

\draw [decorate, decoration = {calligraphic brace,raise=15pt, mirror, amplitude=14pt}, line width = 0.4mm ] ($(N6u)-(0.5,-0.6)$) --  ($(N6b)-(0.5,0.6)$); 

\draw[->] ($(N5c)+(0,0)$) -- (7.5, -5) node[inner sep=8pt, at end, yshift=-0.5cm] {\Large $n_{l-1}$} ;

\node [draw, minimum size=0.85cm, line width=0.5mm, at=(N7u)] (N7U) {};

\node [draw, minimum height=4.5cm, minimum width=12.5cm, xshift=2.5cm, yshift=6cm, line width=0.3mm, at=(N7U)] (box) {}; 

\draw[->, inner sep=5pt] (N7U.north) .. controls +(down:0cm) and +(right:0cm) .. (box.south);

\node (boxC) [draw, circle, minimum size=0.85cm, at=(box.west), xshift=2.5cm] {$\alpha_{2}^{(l-1)}$}  ;

\node (boxCC) [ minimum size=0.85cm, at=(box.west), xshift=0.75cm] {\Large$n_{l-1}$}  ;

\node (boxT) [draw, circle, minimum size=0.85cm, at=(boxC), yshift=1.5cm] {$\alpha_{1}^{(l-1)}$}  ;
\node (boxB) [ circle, minimum size=0.85cm, at=(boxC), yshift=-1.5cm] { \vspace{-1cm} \Large \vdots }  ;
\draw [decorate, decoration = {calligraphic brace,raise=15pt, mirror, amplitude=14pt}, line width = 0.4mm ] ($(boxT)-(0.2,-0.6)$) --  ($(boxB)-(0.2,0.6)$) ;


\node (boxScal) [draw, minimum size=0.85cm, at=(boxC), xshift=3.25cm, label={Scalar product}] { 
$   \sum \limits_{i=1}^{n_{l-1}} w_i \alpha_i^{(l-1)} $
}  ;

\draw[->] (boxT) -- (boxScal) node[midway, above, black]{$w_{1}$};
\draw[->] (boxC) -- (boxScal) node[midway, above, black]{$w_{2}$};

\node (boxAct) [draw, minimum size=0.85cm, at=(boxScal), xshift=3.25cm, label={Activation function}] {
\begin{tikzpicture}
        \begin{axis}[axis lines=center, width=3.5cm, axis x line shift=1cm, height=3.5cm,
         yticklabels={,,},
         xticklabels={,,}
         ]
        \addplot[color=black]{1/(1+exp(-x))-2};
        \end{axis}
        \end{tikzpicture}
}  ;

\node (nodeRed) [draw, circle, minimum size=0.85cm, fill=black!25, at=(boxAct), xshift=2.5cm] {$\alpha_{1}^{(l)}$};

\draw [->] (boxAct) -- (nodeRed);

\end{tikzpicture}
}
\caption{The Multi-Layer Perceptron (MLP) problem (figure adapted from~\cite[Figure 1]{audet2022general}).}
\label{fig:MLP} 
\end{figure}
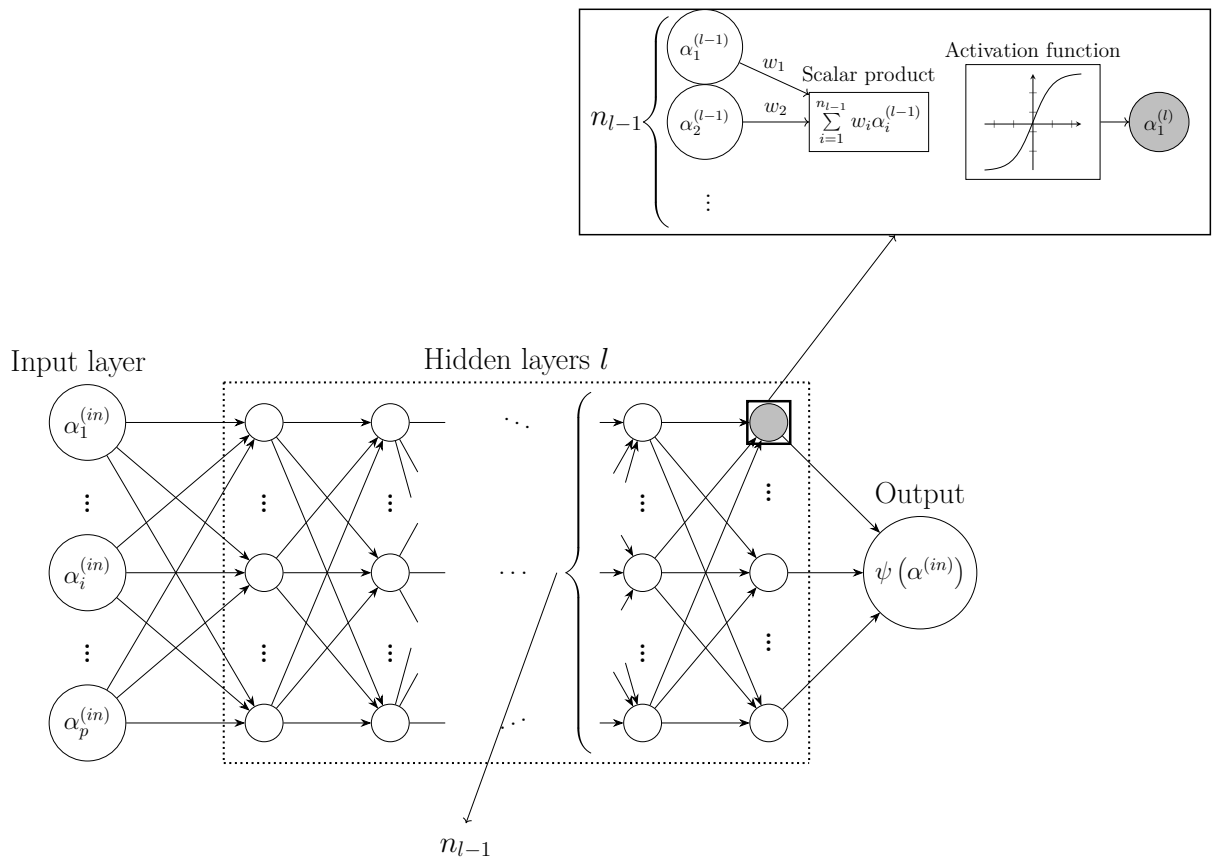

To model hierarchical variables, as proposed in~\cite{audet2022general}, we separate the input space $\mathcal{X}$ as $( \mathcal{X}_{\neutral},\mathcal{X}_{\meta},\mathcal{X}_{\decreed})$ where $\displaystyle \mathcal{X}_{\decreed}= \bigcup_{x_\meta \in \mathcal{X}_{\meta}} \mathcal{X}_{\acting}(x_\meta) $. 
Hence, for a given point $x \in \mathcal{X}$, one has $x=(x_{\neutral}, x_{\meta}, x_{\acting} (x_{\meta}) ),$ where $x_{\neutral} \in \mathcal{X}_{\neutral}$, $x_{\meta}\in \mathcal{X}_{\meta}$ and $x_{\acting} (u_{\meta})\in \mathcal{X}_{\acting} (u_{\meta})$ are defined as follows:
 \begin{itemize}
     \item The components $x_{\neutral} \in \mathcal{X}_{\neutral}$ gather all neutral variables that are not impacted by the meta variables but needed. 
     For example, in the MLP problem, $\mathcal{X}_{\neutral}$ gathers the possible learning rates, momentum, activation functions and batch sizes.
     Namely, from Table~\ref{tab:hyp_NN}, $ \mathcal{X}_{\neutral } = [10^{-5},10^{-2}] \times   [0,1] \times \{\mbox{ReLu}, \mbox{Sigmoid}, \mbox{Tanh} \} \times \{8,16,\ldots,256\} $.

     \item The components $x_{\meta}$ gather the meta (also known as dimensional) variables that determine the inclusion or exclusion of other variables. 
     For example, in the MLP problem, $\mathcal{X}_{\meta}$ represents the possible numbers of layers in the MLP. 
     Namely, from Table~\ref{tab:hyp_NN}, $ \mathcal{X}_{\meta} = \{1,2,3\}$.


     \item The components $x_{\acting} (x_{\meta})$, contain the decreed variables whose inclusion (decreed-included) or exclusion (decreed-excluded) is determined by the values of the meta components $x_{\meta}$. 
     For example, in the MLP problem, $\mathcal{X}_{\decreed}$ represents the number of neurons in the decreed layers.
     Namely, from Table~\ref{tab:hyp_NN}, $ \mathcal{X}_{\acting} (x_\meta=3) =  [50,55]^3$, $ \mathcal{X}_{\acting} (x_\meta=2) =  [50,55]^2$ and $ \mathcal{X}_{\acting} (x_\meta=1) =  [50,55]$. 
 \end{itemize}


\subsection{A Kriging model for hierarchical variables}

In this section, a new method to build a Kriging model with hierarchical variables is introduced based on the framework aforementioned. 
The proposed methods are included in \texttt{SMT 2.0}.
\subsubsection{Motivation and state-of-the-art}

Assuming that the decreed variables are quantitative,~\citet{Hutter} proposed several kernels for the hierarchical context. 
A classic approach, called the imputation method (\texttt{Imp-Kernel}) leads to an efficient paradigm in practice that can be easily integrated into a more general framework as proposed by~\citet{Effectiveness}. 
However, this kernel lacks depth and depends on arbitrary choices.
Therefore,~\citet{Hutter} also proposed a more general kernel, called \texttt{Arc-Kernel} and~\citet{Horn_hier} generalized this kernel even more and proposed a new formulation called the \texttt{Wedge-Kernel}~\cite{DACE_hier}.   
The drawbacks of these two methods are that they add some extra hyperparameters for every decreed dimension (respectively one extra hyperparameter for the \texttt{Arc-Kernel} and two hyperparameters for the \texttt{Wedge-Kernel}) and that they need a normalization according to the bounds. 
More recently,~\citet{pelamattihier} developed a hierarchical kernel for Bayesian optimization. 
However, our work is more general thanks to the framework introduced earlier~\cite{audet2022general} that considers variable-wise formulation and is more flexible as we are allowing sub-problems to be intersecting. 

In the following, we describe our new method to build a correlation kernel for hierarchical variables.
In particular, we introduce a new algebraic kernel called \texttt{Alg-Kernel} that behaves like the \texttt{Arc-Kernel} whilst correcting most of its drawbacks. 
In particular, our kernel does not add any hyperparameters, and the normalization is handled in a natural way. An extension for the decreed categorical variables is described in~\ref{app:catdec}. 

\subsubsection{A new hierarchical correlation kernel} 

For modeling purposes, we assume that the decreed space is quantitative, i.e., $ \mathcal{X}_{\decreed}= \mathcal{X}_{\decreed}^{\quant}$.
Let  $u \in \mathcal{X}$ be an input point partitioned as $u=(u_{\neutral}, u_\meta, u_{\acting}(u_{\meta}))$  and, similarly, $v\in \mathcal{X}$ is another input such that $v=(v_{\neutral}, v_\meta, v_{\acting}(v_{\meta}))$. 
The new kernel $k$ that we propose for hierarchical variables is given by
\begin{eqnarray} 
k(u, v) &= &k_{\neutral} (u_{\neutral},v_{\neutral}) \times  
k_\meta(u_\meta,v_\meta) \nonumber \\
& & \quad \quad \times \ 
k_{\meta,\decreed} ( [u_\meta, u_{\acting}(u_\meta)], [v_\meta, v_{\acting}(v_\meta)]), 
\label{eq:hier_ker}
\end{eqnarray}
where $k_{\neutral}$, $k_\meta$ and $k_{\meta,\decreed}$ are as follows:
 \begin{itemize}
\item  $k_{\neutral}$ represents the neutral kernel that encompasses both categorical and quantitative neutral variables, i.e., $k_{\neutral}$ can be decomposed into two parts $k_{\neutral} (u_{\neutral},v_{\neutral})= k^{\cat}(u_{\neutral}^{\cat},v_{\neutral}^{\cat})k^{\quant} (u_{\neutral}^{\quant},v_{\neutral}^{\quant}).$
The categorical kernel, denoted $k^{\cat}$, could be any Symmetric Positive Definite (SPD)~\cite{Mixed_Paul} mixed kernel (see Section~\ref{sec:mixed}). 
For the quantitative (integer or continuous) variables, a distance-based kernel is used. 
The  chosen quantitative kernel (Exponential, Matérn,...), always depends on a given distance $d$.
For example, the $n$-dimensional exponential kernel gives 

\begin{equation}  
 k^{\quant}(u^{\quant},v^{\quant}) = \displaystyle{ \prod^{n}_{i=1} \exp ( - d(u^{\quant}_i,v^{\quant}_i))}.
\end{equation}

    \item $k_{\meta}$ is the meta variables related kernel.
    It is also separated into two parts:
$k_{\meta} (u_{\meta},v_{\meta})= k^{\cat}(u_{\meta}^{\cat},v_{\meta}^{\cat})k^{\quant} (u_{\meta}^{\quant},v_{\meta}^{\quant})$ where the quantitative kernel is ordered and not continuous because meta variables take value in a finite set. 
\item  $k_{\meta,\decreed}$ is  an SPD kernel that models the correlations between the meta levels (all the possible subspaces) and the decreed variables. 
In what comes next, we detailed this kernel.
\end{itemize}

We can not separate the meta and decreed kernels as it would break the SPD property and the distance would be ill-defined as proven in~\ref{app:SPD}.

\subsubsection{Towards an algebraic meta-decreed kernel} 
Meta-decreed kernels like the imputation kernel or the \texttt{Arc-Kernel} were first proposed in~\cite{Zaefferer,Hutter} and the distance-based kernels such as \texttt{Arc-Kernel} or \texttt{Wedge-Kernel}~\cite{DACE_hier} were proven to be SPD. 
Nevertheless, to guarantee this SPD property, the same hyperparameters are used to model the correlations between the meta levels and between the decreed variables~\cite{Zaefferer}. 
For this reason, the \texttt{Arc-Kernel} includes additional continuous hyperparameters which makes the training of the GP models more expensive and introduces more numerical stability issues. 
In this context, we are proposing a new  algebraic meta-decreed kernel (denoted as \texttt{Alg-Kernel}) that enjoys similar properties as \texttt{Arc-Kernel} but without using additional continuous hyperparameters nor rescaling. 
In the \texttt{SMT 2.0} release, we included \texttt{Alg-Kernel} and a simpler version of \texttt{Arc-Kernel} that do not relies on additional hyperparameters.

Our proposed \texttt{Alg-Kernel} kernel is given by 
\begin{equation}
\begin{split}
 & k^{\text{alg}}_{\meta,\decreed} ( [u_\meta, u_{\acting}(u_\meta)], [v_\meta, v_{\acting}(v_\meta)]) \\&  \quad \quad = k^{\text{alg}}_\meta(u_\meta, v_\meta) \times   k^{\text{alg}}_{\decreed}(u_{\acting}(u_\meta),v_{\acting}(v_\meta)).
\end{split}
\end{equation}
Mathematically, we could consider that there is only one meta variable whose levels correspond to every possible included subspace. 
Let $I_{\text{sub}}$ denotes the components indices of possible subspaces, the subspaces parameterized by the meta component $u_\meta$ are defined as $\mathcal{X}_{\acting}(u_\meta=l), \ l \in I_{\text{sub}} $.
It follows that the fully extended continuous decreed space writes as $\mathcal{X}_{\decreed} = \bigcup_{l \in I_{\text{sub}}} \mathcal{X}_{\acting}(u_\meta=l)$ and $I_\decreed$ is the set of the associated indices. 
Let $I^{inter}_{u,v}$ denotes the set of components related to the space $ \mathcal{X}_{\acting} (u_\meta,v_\meta)$ containing the variables decreed-included in both  $ \mathcal{X}_{\acting} (u_\meta)$ and  $ \mathcal{X}_{\acting} (v_\meta) $.

Since the decreed variables are quantitative, one has 
\begin{equation}
\begin{split}
k^{\text{alg}}_{\decreed} (u_{\acting}(u_\meta),v_{\acting}(v_\meta)) 
&= k^{\quant} (u_{\acting}(u_\meta),v_{\acting}(v_\meta))\\
& = \prod_{i \in I_{u,v}^{inter}}   k^{\quant} (  [(u_{\acting}(u_\meta)]_i,[v_{\acting}(v_\meta)]_i )
\end{split}
\end{equation}
The construction of the quantitative kernel $k^{\quant}$ depends on a given distance denoted $d^{\text{alg}}$. 
The kernel $k^{\text{alg}}_\meta$ is an induced meta kernel that depends on the same distance $d^{\text{alg}}$ to preserve the SPD property of $k^{\text{alg}}_{\meta,\decreed}$. 
For every $i \in I_\decreed$, if $i \in I^{inter}_{u,v} $, the new algebraic distance is given by
\begin{equation}
d^{\text{alg}}( [u_{\acting} (u_\meta) ]_i , [v_{\acting} (v_\meta) ]_i )  = \left(\frac{2  | [u_{\acting} (u_\meta) ]_i  - [v_{\acting} (v_\meta) ]_i|}{ \sqrt{{ [u_{\acting} (u_\meta) ]_i }^2+1}\sqrt{{ [v_{\acting} (v_\meta) ]_i }^2+1}}\right)\theta_i,
\end{equation}
where $\theta_i \in \mathbb{R}^+$ is a continuous hyperparameter.
Otherwise, if $i \in I_\decreed$ but $i \notin I^{inter}_{u,v} $, there should be a non-zero residual distance between the two different subspaces $ \mathcal{X}_{\acting} (u_\meta)$ and $\mathcal{X}_{\acting} (v_\meta)$ to ensure the kernel SPD property. 
To have a residual not depending on the decreed values, our model considers that there is a unit distance
$$d^{\text{alg}}( [u_{\acting} (u_\meta) ]_i , [v_{\acting} (v_\meta) ]_i) = 1.0 \ \theta_i, \  \forall i \in I_\decreed \setminus I^{inter}_{u,v}. $$ 
%
The induced meta kernel $k^{\text{alg}}_{\meta}(u_{\meta},v_{\meta})$ to preserve the SPD property of $k^{\text{alg}}$ is defined as: 
\begin{equation}
k^{\text{alg}}_{\meta}(u_{\meta},v_{\meta}) 
=   \prod_{i \in I_\meta} \ k^{\quant}(1.0 \ \theta_i) .  
\label{eq:d_alg}
\end{equation}
%
The proof that our kernel is SPD is given in~\ref{app:SPD}. Not only our kernel of~\eqnref{eq:hier_ker} uses less hyperparameters than the \texttt{Arc-Kernel} (as we cut off its extra parameters) but it is also a more flexible kernel as it allows different kernels for meta and decreed variables.
Moreover, another advantage of our kernel is that it is numerically more stable thanks to the new non-stationary~\cite{hebbal2021bayesian} algebraic distance defined in~\eqnref{eq:d_alg}~\cite{wildberger2007rational}. 
Our proposed distance also does not need any rescaling by the bounds to have values between 0 and 1.
Moreover, this distance can be expressed for vectors inputs as given in~\ref{app:dalg}.

In what comes next, we will refer to the implementation of the kernels \texttt{Arc-Kernel} and \texttt{Alg-Kernel} by \texttt{SMT Arc-Kernel} and \texttt{SMT Alg-Kernel}. 
We note also that the implementation of \texttt{SMT Arc-Kernel} differs slightly from the original \texttt{Arc-Kernel} as we fixed some hyperparameters to 1 in order to avoid adding extra hyperparameters and use the formulation of~\eqnref{eq:hier_ker} and rescaling of the data. 

\subsubsection{Illustration on the MLP problem}

In this section, we illustrate the hierarchical \texttt{Arc-Kernel} on the MLP example. For that sake, we consider two design variables $u$ and $v$ such that  $u = (2.10^{-4}, \textcolor{black}{0.9}, \mbox{ReLU}, 16, 2, 55,51)$ and $v =  ( 5.10^{-3},\textcolor{black}{0.8}, \mbox{Sigmoid}, 64, 3, 50,54,53)$. 
Since the value of $u_\meta$ (i.e., the number of hidden layers) differs from one point to another (namely, $2$ for $u$  and $3$ for $v$), the associated variables $u_\acting(u_\meta)$ have  either 2 or 3 variables for the number of neurons in each layer (namely 55 and 51 for $u$, and 50, 54 and 53 for the point $v$).
In our case, \textcolor{black}{8} hyperparameters $([R_1]_{1,2}, \theta_1, \ldots, \textcolor{black}{\theta_7})$ will have to be optimized where $k$ is given by \eqnref{eq:hier_ker}.
These 7 hyperparameters can be described using our proposed framework as follows:
\begin{itemize}
    \item For the neutral components, we have $u_\neutral = (2.10^{-4}, \textcolor{black}{0.9} ,\mbox{ReLU}, 16)$ and $v_\neutral = (5.10^{-3},\textcolor{black}{0.8}, \mbox{Sigmoid}, 64)$. Therefore, for a categorical matrix kernel $R_1$  and a square exponential quantitative kernel,
\begin{equation*}
\centering
\begin{split}
k_{\neutral} (u_{\neutral},v_{\neutral}) &=  k^{\cat}(u_{\neutral}^{\cat},v_{\neutral}^{\cat})k^{\quant} (u_{\neutral}^{\quant},v_{\neutral}^{\quant}) \\ 
&= [R_1]_{1,2} \exp{[- \theta_1 (2.10^{-4} -5.10^{-3})^2]  }  \\
&\quad \ \ \textcolor{black}{ \exp{[- \theta_2 (0.9 - 0.8)^2]  } }  \exp{[ - \textcolor{black}{\theta_3} (16-64)^2 ]}.
\end{split}
\end{equation*}
The values $[R_1]_{1,2}$, $\theta_1$, \textcolor{black}{$\theta_2$} and $\theta_3$ need to be optimized. 
Here, $[R_1]_{1,2}$ is the correlation between "ReLU" and "Sigmoid".
\item 
For the meta components, we have $u_\meta = 2$ and $v_\meta = 3$.
Therefore, for a square exponential quantitative kernel,
\begin{equation*}
\centering
\begin{split}
k_{\meta} (u_{\meta},v_{\meta})&= k^{\cat}(u_{\meta}^{\cat},v_{\meta}^{\cat})k^{\quant} (u_{\meta}^{\quant},v_{\meta}^{\quant}) \\
&= \exp{[ - \textcolor{black}{\theta_4} (3-2)^2 ]}.
\end{split}
\end{equation*}
The value \textcolor{black}{$\theta_4$} needs to be optimized. 

\item  For the meta-decreed kernel, we have $[ u_\meta, u_\acting(u_\meta) ] = [2, (55,51)] $ and $[ v_\meta, v_\acting(v_\meta) ] = [3, (50,54,53)]$ which gives
\begin{eqnarray*}
 & & k^{\text{alg}}_{\meta,\decreed} ( [u_\meta, u_{\acting}(u_\meta)], [v_\meta, v_{\acting}(v_\meta)]) \\& &  \quad \quad = k^{\mathrm{alg}}_\meta(2,3) \ k^{\mathrm{alg}}_\decreed( (55,51), (50,54,53)).
\end{eqnarray*}
The distance $d^\text{alg} (51, 54) =  \left( \frac{  2 \times |51-54|}{ \sqrt{ 51^2+1 } \sqrt{ 54^2+1 }} \right) \textcolor{black}{\theta_6} = 2.178.10^{-3} \ \textcolor{black}{\theta_6}.$
In general, for surrogate models, and in particular in \texttt{SMT 2.0}, the input data are normalized.
With a unit normalization from $[50,55]$ to $[0,1]$, we would have  $d^\text{alg} (0.2, 0.8) =  \left( \frac{ 2 \times 0.6}{ \sqrt{ 0.2^2+1 } \sqrt{ 0.6^2+1 }} \right) \textcolor{black}{\theta_6} = 0.919 \ \textcolor{black}{\theta_6} .$ 
Similarly, we have, between 55 and 50, $d^\text{alg} (0, 1) =  1.414 \ \textcolor{black}{\theta_5}.$ 
Hence, for a square exponential quantitative kernel, one gets
\begin{eqnarray*}
 & & k^{\text{alg}}_{\meta,\decreed} ( [u_\meta, u_{\acting}(u_\meta)], [v_\meta, v_{\acting}(v_\meta)]) \\
& &  \quad \quad = \exp{[ - \textcolor{black}{\theta_7}]} \times \exp{[ - 1.414\ \textcolor{black}{\theta_5} ]} \times \exp{[ - 0.919\ \textcolor{black}{\theta_6} ]},
\end{eqnarray*} 
where the meta induced component is $k^{\text{alg}}_{\meta}(u_{\meta},v_{\meta}) =  \exp{[ - {\theta_7}] }$ because the decreed value $53$ in $v$ has nothing to be compared with in $u$ as in~\eqnref{eq:d_alg}.
The values  \textcolor{black}{$\theta_5$, $\theta_6$ and $\theta_7$} need to be optimized which complete the description of the hyperparameters.

We note that for the MLP problem, \texttt{Alg-Kernel} models use \textcolor{black}{10} hyperparameters whereas the \texttt{Arc-Kernel} would require \textcolor{black}{12} hyperparameters without the meta kernel ($\textcolor{black}{\theta_4}$) but with 3 extra decreed hyperparameters 
and the \texttt{Wedge-Kernel} would require \textcolor{black}{15} hyperparameters.
For deep learning applications, a more complex perceptron with up to $10$ hidden layers would require $\textcolor{black}{17}$ hyperparameters with \texttt{SMT 2.0} models against $\textcolor{black}{26}$ for \texttt{Arc-Kernel} and $\textcolor{black}{36}$ for \texttt{Wedge-Kernel}. 
The next section illustrates the interest of our method to build a surrogate model for this neural network engineering problem.

\end{itemize}

\subsection{A neural network test-case using \texttt{SMT 2.0}}
\label{sec:NN}
In this section, we apply our models to the hyperparameters optimization of a MLP problem aforementioned of~\figref{fig:MLP}. 
Within \texttt{SMT 2.0} an example illustrates this MLP problem. 
For the sake of showing the Kriging surrogate abilities, we implemented a dummy function with no significance to replace the real black-box that would require training a whole Neural Network (NN) with big data. 
This function requires a number of variables that depends on the value of the meta variable, i.e the number of hidden layers. 
To simplify, we have chosen only 1, 2 or 3 hidden layers and therefore, we have 3 decreed variables but deep neural networks could also be investigated as our model can tackle a few dozen variables. 
A test case (\emph{test\_hierarchical\_variables\_NN}) shows that our model \texttt{SMT Alg-Kernel} interpolates the data properly, checks that the data dimension is correct and also asserts that the inactive decreed variables have no influence over the prediction. 
In~\figref{fig:NN_hier} we illustrate the usage of Kriging surrogates with hierarchical and mixed variables based on the implementation of \texttt{SMT 2.0} for \emph{test\_hierarchical\_variables\_NN}. 

To compare the hierarchical models of \texttt{SMT 2.0} (\texttt{SMT Alg-Kernel} and \texttt{SMT Arc-Kernel}) with the state-of-the-art imputation method previously used on industrial application (\texttt{Imp-Kernel})~\cite{Effectiveness}, we draw a 99 point LHS (33 points by meta level) as a training set and the validation set is a LHS of $3\times 1000=3000$ points. 
For the \texttt{Imp-Kernel}, the decreed-excluded values are replaced by $52$ because the mean value $52.5$ is not an integer (by default, \texttt{SMT} rounds to the floor value). 
For the three methods, the precision (computed with a root-mean-square error RMSE criterion), the likelihood and the computational time are shown in~\tabref{tab:resNN}.
For this problem, we can see that \texttt{SMT Alg-kernel} gives better performance than the imputation method or \texttt{SMT Arc-kernel}. Also, as all methods use the same number of hyperparameters, they have similar time performances. 
A direct application of our modeling method is Bayesian optimization to perform quickly the hyperparameter optimization of a neural network~\cite{cho2020basic}. 

\begin{table}[H]
\centering
\caption{Results on the neural network model. }
\small
\resizebox{\columnwidth}{!}{
\begin{tabular}{cccc} 
\hline
\textbf{Hierarchical method}  &  Prediction error (RMSE) & $ \ $ Likelihood  & \# of hyperparam. 
\\
\hline 
\texttt{SMT Alg-kernel}  & \  {\textcolor{black}{ 3.7610}
} & {
\textcolor{black}{176.11}
} & \textcolor{black}{10}
\\   
\texttt{SMT Arc-kernel}  & \textcolor{black}{4.9208} 
& \textcolor{black}{162.01}
& \textcolor{black}{10} 
\\   
\texttt{Imp-kernel}  & \textcolor{black}{4.5455}
& \textcolor{black}{170.64}
& \textcolor{black}{10} 
\\
\hline
\end{tabular}
}
\label{tab:resNN}
\end{table}
\begin{figure}
\vspace{-42pt}
\begin{lstlisting}
from smt.sampling_methods import LHS
from smt.surrogate_models import KRG,MixIntKernelType,\
  MixHrcKernelType,DesignSpace,FloatVariable,\
  IntegerVariable,OrdinalVariable,CategoricalVariable
from smt.applications.mixed_integer import \
  MixedIntegerSamplingMethod,MixedIntegerKrigingModel
import f1_NN, f2_NN,f3_NN #dummy example
def test_hierarchical_variables_NN(self):
   def dummy_f(x):
       if x[0] == 1:
         y=(f1_NN(x[1],x[2],x[3],2**x[4],x[5]))
       elif x[0] == 2:
         y=(f2_NN(x[1],x[2],x[3],2**x[4],x[5],x[6]))
       elif x[0] == 3:
         y=(f3_NN(x[1],x[2],x[3],2**x[4],x[5],x[6],x[7]))
       return y
   # Define the mixed hierarchical design space
   ds = DesignSpace([
     IntegerVariable(1, 3), FloatVariable(1e-5, 1e-2),
     FloatVariable(0, 1),
     CategoricalVariable(["ReLU", "Sigmoid","Tanh"]),
     IntegerVariable(3, 8), IntegerVariable(50, 55),
     IntegerVariable(50, 55), IntegerVariable(50, 55),
   ])
   # activate x5 when x0 in [2, 3]; x6 when x0 == 3
   ds.declare_decreed_var(
     decreed_var=6, meta_var=0, meta_value=[2, 3])
   ds.declare_decreed_var(7, meta_var=0, meta_value=3)
   #Perform the mixed integer sampling 
   sampling = MixedIntegerSamplingMethod(
    LHS, ds, criterion="ese", random_state=42)
   Xt = sampling(100)
   Yt = dummy_f(Xt)
   #Build the surrogate
   sm = MixedIntegerKrigingModel(
     surrogate=KRG(design_space=ds, corr="abs_exp",
     categorical_kernel=MixIntKernelType.HOMO_HSPHERE,
     hierarchical_kernel=MixHrcKernelType.ALG_KERNEL)
   sm.set_training_values(Xt, Yt)
   sm.train()
   # Check prediction accuracy
   y_s = sm.predict_values(Xt)
   pred_RMSE = np.linalg.norm(y_s - Yt) / len(Yt)
   y_sv = sm.predict_variances(Xt)
   var_RMSE = np.linalg.norm(y_sv) / len(Yt)
   assert pred_RMSE < 1e-7
   assert var_RMSE < 1e-7
\end{lstlisting}
\caption{Example of usage of Hierarchical and Mixed Kriging surrogate.}
\label{fig:NN_hier}
\end{figure}

\section{Bayesian optimization within \texttt{SMT 2.0}}
\label{sec:BO}

Efficient global optimization (EGO) is a sequential Bayesian optimization algorithm designed to find the optimum of a black-box function that may be expensive to evaluate~\cite{Jones98}.
EGO starts by fitting a Kriging model to an initial DoE, and then uses an acquisition function to select the next point to evaluate.
The most used acquisition function is the expected improvement. 
Once a new point has been evaluated, the Kriging model is updated.
Successive updates increase the model accuracy over iterations.
This enrichment process repeats until a stopping criterion is met.

Because \texttt{SMT 2.0} implements Kriging models that handle mixed and hierarchical variables, we can use EGO to solve problems involving such design variables. 
Other Bayesian optimization algorithms often adopt approaches based on solving subproblems with continuous or non-hierarchical Kriging.
This subproblem approach is less efficient and scales poorly, but it can only solve simple problems. 
Several Bayesian optimization software packages can handle mixed or hierarchical variables with such a subproblem approach.
The packages include BoTorch~\cite{balandat2020botorch}, SMAC~\cite{SMAC3}, Trieste~\cite{picheny2023trieste}, HEBO~\cite{cowen-rivers_hebo_2020}, OpenBox~\cite{jiang2023openbox}, and Dragonfly~\cite{Dragonfly}. 

\subsection{A mixed optimization problem} 
\label{sec:MI-BO}

\figref{res_optim_mi_SMT} compares the four EGO methods implemented in \texttt{SMT 2.0}: \texttt{SMT GD}, \texttt{SMT CR}, \texttt{SMT EHH} and \texttt{SMT HH}. 
The mixed test case that illustrates Bayesian optimization is a toy test case~\cite{CAT-EGO} detailed in~\ref{app:Toy}. 
This test case has two variables, one continuous and one categorical with 10 levels.
To assess the performance of our algorithm, we performed 20 runs with different initial DoE sampled by LHS.
Every DoE consists of 5 points and we chose a budget of 55 infill points.
\figref{convmi} plots the convergence curves for the four methods. 
In particular, the median is the solid line, and the first and third quantiles are plotted in dotted lines. 
To visualize better the data dispersion, the boxplots of the 20 best solutions after 20 evaluations are plotted in~\figref{mini_mi}. 
As expected, the more a method is complex, the better the optimization. Both \texttt{SMT HH} and \texttt{SMT EHH} converged in around 18 evaluations whereas \texttt{SMT CR} and \texttt{SMT GD} take around 26 iterations as shown on~\figref{convmi}.
Also, the more complex the model, the closer the optimum is to the real value as shown on~\figref{mini_mi}. 

\begin{figure}[H]
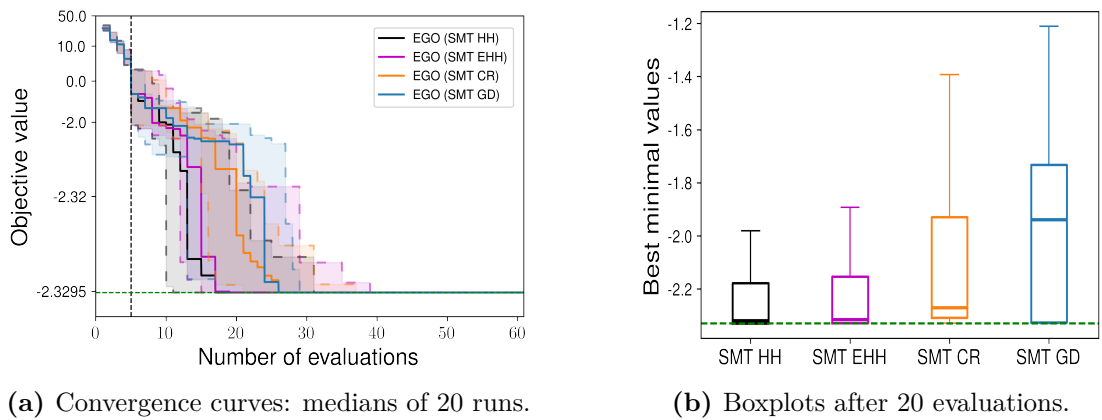

\begin{minipage}[b]{.6\linewidth}
\centering
\hspace{-1.25cm}
\subfloat[Convergence curves: medians of 20 runs.]{
\includegraphics[height=5cm,,width=7cm]
{images/miBO2.jpg}
\label{convmi}
}
\end{minipage}
\begin{minipage}[b]{.4\linewidth}
\centering 
\hspace{-1.5cm}
\subfloat[Boxplots after 20 evaluations.]{
\includegraphics[height=5cm,width=6.2cm]{images/minima_optim_SMT2_15it.jpg}
\label{mini_mi}
}
\end{minipage}
\caption{Optimization results for the Toy function~\cite{CAT-EGO}.}
\label{res_optim_mi_SMT}
\end{figure}

In~\figref{fig:Branin_mixed} we illustrate how to use EGO with mixed variables based on the implementation of \texttt{SMT 2.0}. 
The illustrated problem is a mixed variant of the Branin function~\cite{AMIEGO}.

\begin{figure}[t!]
\begin{lstlisting}
#Import the Mixed Integer API
from smt.surrogate_models import KRG,MixIntKernelType,\
  DesignSpace,FloatVariable,IntegerVariable
from smt.applications.mixed_integer import \
  MixedIntegerSamplingMethod as misamp
#Define the function
from smt.problems import Branin
fun = Branin(ndim=2)
#Define the mixed design space
design_space = DesignSpace([
  IntegerVariable(*fun.xlimits[0]),
  FloatVariable(*fun.xlimits[1]),
])
#Perform a mixed integer sampling with LHS
from smt.sampling_methods import LHS
smp = misamp(LHS,design_space,random_state=42)
xdoe = smp(10)
# Call the Bayesian optimizer
from smt.applications import EGO
criterion = "EI"  #'EI' or 'SBO' or 'LCB'
ego = EGO(xdoe=xdoe,
  n_iter=20,
  criterion="EI",
  random_state=42,
  surrogate=KRG(design_space=design_space,
   categorical_kernel=MixIntKernelType.GOWER))
x_opt, y_opt, _, _, _ = ego.optimize(fun=fun)
# Check if the result is correct
self.assertAlmostEqual(0.494, float(y_opt), delta=1)
\end{lstlisting}
\caption{Example of usage of mixed surrogates for Bayesian optimization.}
\label{fig:Branin_mixed}
\end{figure}

Note that a dedicated notebook is available to reproduce the results presented in this chapter and the mixed integer notebook also includes an extra mechanical application with composite materials~\cite{RaulAIAA}\footnote{\url{https://colab.research.google.com/github/SMTorg/smt/blob/master/tutorial/SMT_MixedInteger_application.ipynb} }. 

\subsection{A hierarchical optimization problem} 
\label{sec:HV-BO}
The hierarchical test case considered in this chapter to illustrate Bayesian optimization is a modified Goldstein function~\cite{pelamattihier} detailed in~\ref{app:Goldstein}.
The resulting optimization problem involves 11 variables: 5 are continuous, 4 are integer (ordinal) and 2 are categorical.
These variables consist in 6 neutral variables, 1 dimensional (or meta) variable and 4 decreed variables.
Depending on the meta variable values, 4 different sub-problems can be identified.
The optimization problem is given by:
\begin{equation}
\begin{split}
& \min  f( x^{\cat}_{\neutral}, x^{\quant}_{\neutral}, x^{\cat}_{m}, x^{\quant}_{\decreed}  ) \\
& \mbox{w.r.t.} \ \  
x^{\cat}_{\neutral} =w_2 \in \{ 0,1 \} \\
& \quad \quad \quad x^{\quant}_{\neutral} = (x_1,x_2,x_5,z_3,z_4) \in \{ 0,100 \}^3 \times \{ 0,1,2 \}^2   \\
& \quad \quad \quad  x^{\cat}_{m} = w_1 \in \{ 0,1,2,3 \} \\
& \quad \quad \quad x^{\quant}_{\decreed} = (x_3,x_4,z_1,z_2) \in \{ 0,100 \}^2 \times \{ 0,1,2 \}^2   
\end{split}
\end{equation}
Compared to the model choice of~\citet{pelamattihier}, we chose to model $x_5$ and $w_2$ as neutral variables even if $f$ does not depend on $x_5$ when $w_2=0$. 
Other modeling choices are kept; for example, $w_2$ is a so-called "binary variable" and not a categorical one~\cite{muller_so-mi_2013}. 
Similarly, we also keep the formulation of $w_1$ as a categorical variable but a better model would be to model it as a "cyclic variable"~\cite{tran:hal-03170761}.
The resulting problem is described in~\ref{app:Goldstein}.
To assess the performance of our algorithm, we performed 20 runs with different initial DoE sampled by LHS.
Every DoE consists of $n+1=12$ points and we chose a budget of $5n = 55$ infill points.
\textcolor{black}{To compare our method with a baseline}, we also tested the random search method thanks to the \texttt{expand\_lhs} new method~\cite{condearenzana} described in Section~\ref{sec:sampling} and we also optimized the Goldstein function using EGO with a classic Kriging model based on imputation method (\texttt{Imp-Kernel}). 
This method replaces the decreed-excluded variables by their mean values: $50$ or $1$ respectively for $(x_3,x_4)$ and $(z_1,z_2)$. 
\figref{convhier} plots the convergence curves for the four methods. 
In particular, the median is the solid line and the first and third quantiles are plotted in dotted lines. 
To visualize better the corresponding data dispersion, the boxplots of the 20 best solutions are plotted in~\figref{mini_hier}. 
The results in~\figref{res_optim} show that the hierarchical Kriging models of \texttt{SMT 2.0} lead to better results than the imputation method or the random search both in terms of final objective value and variance over the 20 runs and in term of convergence rate.
More precisely, \texttt{SMT Arc-Kernel} and \texttt{SMT Alg-Kernel} Kriging model gave the best EGO results and managed to converge correctly as shown in~\figref{mini_hier}. \textcolor{black}{More precisely, the 20 sampled DoEs led to similar performance and from one DoE, the method \texttt{SMT Alg-Kernel} managed to find the true minimum. However, this result has not been reproduced in other runs and is therefore not statistically significant. The variance between the runs is of similar magnitude regardless of the considered methods. }

\begin{figure}[ht]
\begin{minipage}[b]{.6\linewidth}
\centering
\hspace{-1.25cm}
\subfloat[Convergence curves: medians of 20 runs.]{
\includegraphics[height=5cm,,width=7cm]
{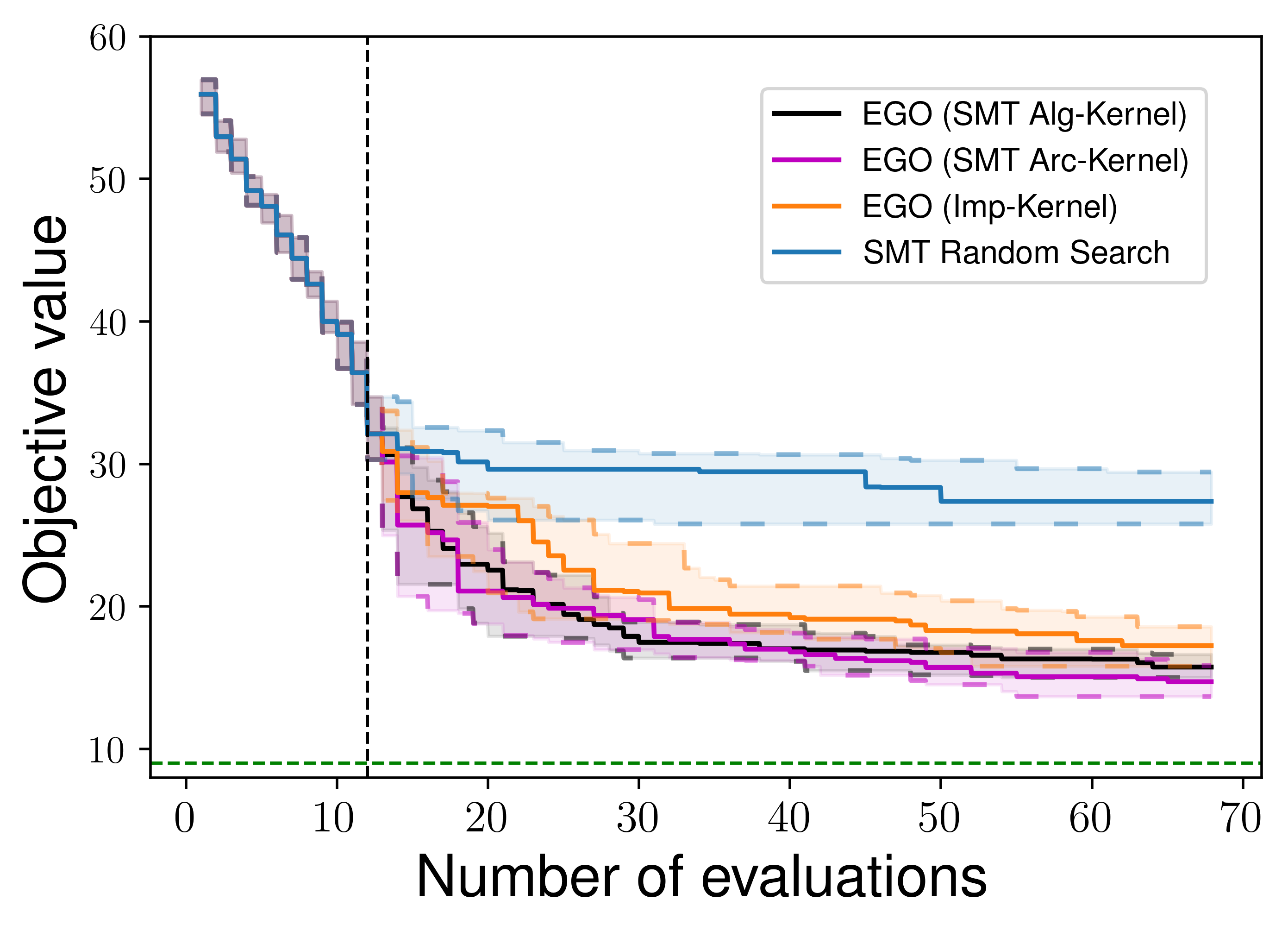}
\label{convhier}
}
\end{minipage}
\begin{minipage}[b]{.4\linewidth}
\centering 
\hspace{-1.5cm}
\subfloat[Boxplots after 67 iterations.]{
\includegraphics[height=5cm,width=6.2cm]{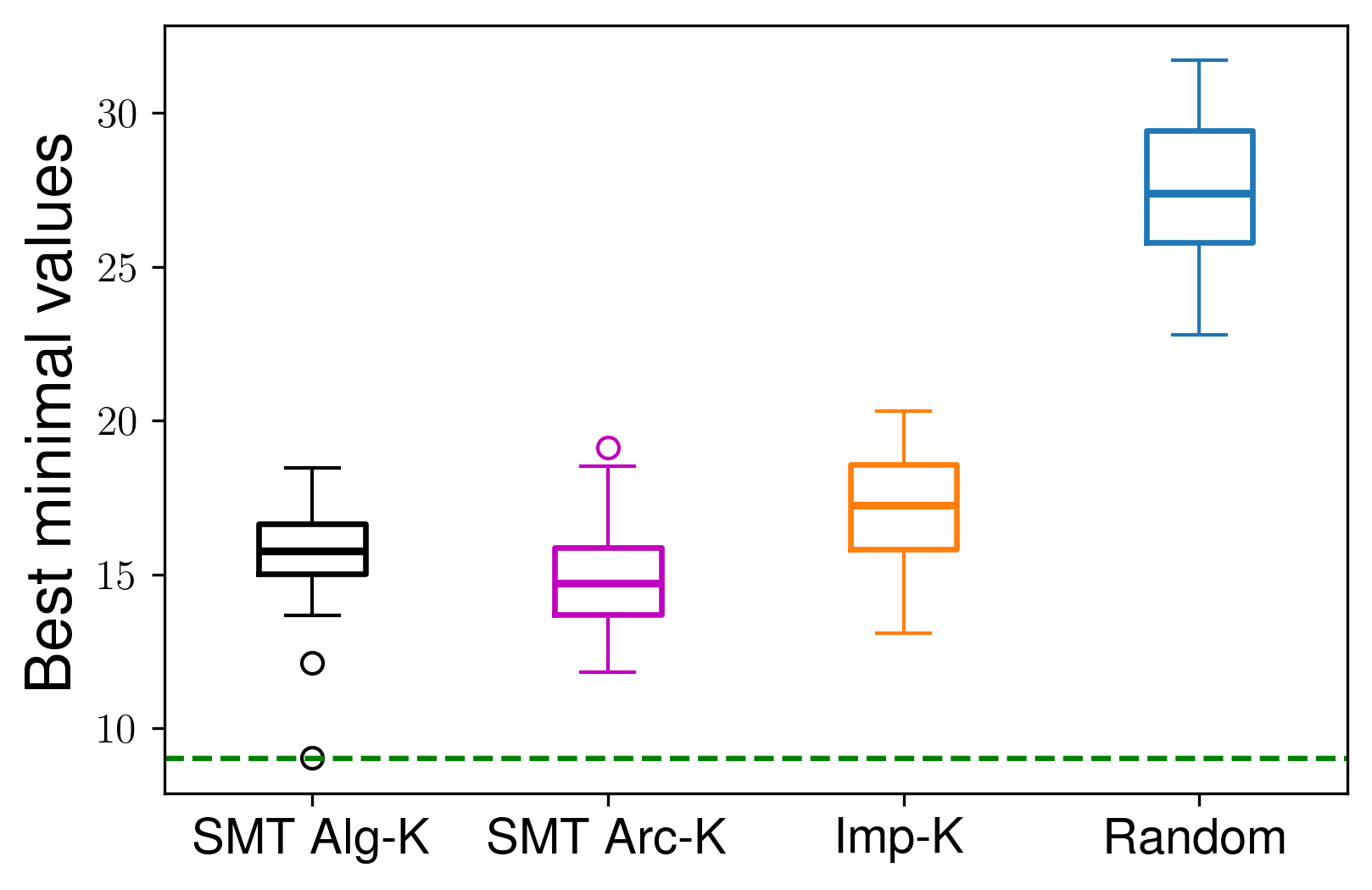}
\label{mini_hier}
}
\end{minipage}
\caption{Optimization results for the hierarchical Goldstein function.}
\label{res_optim}
\end{figure}

\section{Other relevant contributions in \texttt{SMT 2.0}}
\label{sec:other}
The new release \texttt{SMT 2.0} introduces several improvements besides Kriging for hierarchical and mixed variables. 
This section details the most important new contributions.  
Recall from Section~\ref{sec:organization} that five sub-modules are present in the code: \texttt{Sampling}, \texttt{Problems}, \texttt{Surrogate Models}, \texttt{Applications} and \texttt{Notebooks}.

\subsection{Contributions to \texttt{Sampling}}
\label{sec:sampling}
\paragraph{Pseudo-random Sampling}
The Latin Hypercube Sampling (LHS) is a stochastic sampling technique to generate quasi-random sampling distributions. 
It is among the most popular sampling method in computer experiments thanks to its simplicity and projection properties with high-dimensional problems. 
The LHS method uses the pyDOE package (Design Of Experiments for Python). 
Five criteria for the construction of LHS are implemented in SMT. 
The first four criteria (\texttt{center, maximin, centermaximin, correlation}) are the same as in pyDOE\footnote{\url{https://pythonhosted.org/pyDOE/index.html}}.
The last criterion \texttt{ese}, is implemented by the authors of SMT~\cite{LHS}. 
In \texttt{SMT 2.0} a new LHS method was developed for the Nested design of experiments (\texttt{NestedLHS})~\textcolor{black}{\cite{meliani2019multi}} to use with multi-fidelity surrogates.
A new mathematical method (\texttt{expand\_lhs})~\textcolor{black}{\cite{condearenzana}} was developed in \texttt{SMT 2.0} to increase the size of a design of experiments while maintaining the \texttt{ese} property.
Moreover, we proposed a sampling method for mixed variables, and the aforementioned LHS method was applied to hierarchical variables in~\figref{res_optim}.

\subsection{Contributions to \texttt{Surrogate models}}
\label{sec:surrogates}
\paragraph{New kernels and their derivatives for Kriging}
Kriging surrogates are based on hyperparameters and on a correlation kernel.
Four correlation kernels are now implemented in \texttt{SMT 2.0}~\cite{Lee2011}.
In SMT, these correlation functions are absolute exponential (\texttt{abs\_exp}), Gaussian (\texttt{squar\_exp}), Matern 5/2 (\texttt{matern52}) and Matern 3/2 (\texttt{matern32}).
In addition, the implementation of gradient and Hessian for each kernel makes it possible to calculate both the first and second derivatives of the GP likelihood with respect to the hyperparameters~\cite{SMT2019}. 

\paragraph{Variance derivatives for Kriging}
To perform uncertainty quantification for system analysis purposes, it could be interesting to know more about the variance derivatives of a model~\textcolor{black}{\cite{lopez, berthelin2022disciplinary, aerobest23cardoso}}. 
For that purpose and also to pursue the original publication about derivatives~\cite{SMT2019}, \texttt{SMT 2.0} extends the derivative support to Kriging variances and kernels.

\paragraph{Noisy Kriging}
In engineering and in big data contexts with real experiments, surrogate models for noisy data are of significant interest.
In particular, there is a growing need for techniques like noisy Kriging and noisy Multi-Fidelity Kriging (MFK) for data fusion~\cite{platt2022systematic}.
For that purpose, \texttt{SMT 2.0} has been designed to accommodate Kriging and MFK to noisy data including the option to incorporate heteroscedastic noise (using the \texttt{use\_het\_noise} option) and to account for different noise levels for each data source~\cite{condearenzana}.

\paragraph{Kriging with partial least squares}
Beside MGP, for high-dimensional problems, the toolbox implements Kriging with partial least squares (KPLS)~\cite{bouhlel_KPLSK} and its extension \gls{KPLSK}~\cite{Bouhlel18}. 
The PLS information is computed by projecting the data into a smaller space spanned by the principal components. 
By integrating this PLS information into the Kriging correlation matrix, the number of inputs can be scaled down, thereby reducing the number of hyperparameters required. 
The resulting number of hyperparameters $d_e$ is indeed much smaller than the original problem dimension $d$.
Recently, in \texttt{SMT 2.0}, we extended the KPLS method for multi-fidelity Kriging (MFKPLS and MFKPLSK)~\cite{meliani2019multi, MFKPLS,charayron2023towards}. 
We also proposed an automatic criterion to choose automatically the reduced dimension $d_e$ based on Wold's R criterion~\cite{wold_1975}.
This criterion has been applied to aircraft optimization with EGO where the number $d_e$ ($\texttt{n\_comp}$ in the code) for the model is automatically selected at every iteration~\cite{SciTech_cat}. 
\textcolor{black}{Special efforts have been made to accommodate KPLS for multi-fidelity and mixed integer data~\cite{MFKPLS,charayron2023towards}}.

\paragraph{Marginal Gaussian process} 
\texttt{SMT 2.0} implements Marginal Gaussian Process (MGP) surrogate models for high dimensional problems~\cite{EGORSE}. 
MGP are Gaussian processes taking into account hyperparameters uncertainty defined as a density probability function. 
Especially we suppose that the function to model $f : \Omega \mapsto \mathbb{R}$, where  $\Omega \subset \mathbb{R}^d$ and $d$ is the number of design variables, lies in a linear embedding $\mathcal{A}$ such as $\mathcal{A} = \{ u = Ax, x\in\Omega\},A \in \mathbb{R}^{d \times d_e}$ and $f(x)=f_{\mathcal{A}}(Ax)$ with $f(x)=f_{\mathcal{A}} : \mathcal{A} \mapsto \mathbb{R}$ and $d_e \ll d$.
Then, we must use a kernel $k(x,x')=k_{\mathcal{A}}(Ax,Ax')$ whose each component of the transfer matrix $A$ is an hyperparameter. 
Thus we have $d_e \times d$ hyperparameters to find. 
Note that $d_e$ is defined as $\texttt{n\_comp}$ in the code~\cite{MGP}. 

\paragraph{Gradient-enhanced neural network}
The new release \texttt{SMT 2.0} implements Gradient-Enhanced Neural Network (GENN) models~\cite{bouhlel2020}.
Gradient-Enhanced Neural Networks (GENN) are fully connected multi-layer perceptrons whose training process was modified to account for gradient information.
Specifically, the model is trained to minimize not only the prediction error of the response but also the prediction error of the partial derivatives: the chief benefit of gradient enhancement is better accuracy with fewer training points. 
Note that GENN applies to regression (single-output or multi-output), but not classification since there is no gradient in that case. 
The implementation is fully vectorized and uses ADAM optimization, mini-batch, and L2-norm regularization.
For example, GENN can be used to learn airfoil geometries from a database. 
This usage is documented in \texttt{SMT 2.0}\footnote{\url{https://smt.readthedocs.io/en/latest/_src_docs/examples/airfoil_parameters/learning_airfoil_parameters.html}}.

\subsection{Contributions to \texttt{Applications}}
\paragraph{Kriging trajectory and sampling}
Sampling a GP with high resolution is usually expensive due to the large dimension of the associated covariance matrix.
Several methods are proposed to draw samples of a GP on a given set of points.
To sample a conditioned GP, the classic method  consists in using a Cholesky decomposition (or eigendecomposition) of the conditioned covariance matrix of the process but some numerical computational errors can lead to non SPD matrix.
A more recent approach  based on Karhunen-Loève decomposition of the covariance kernel with the Nyström method has been proposed in~\cite{betz2014numerical} where the paths can be sampled by generating independent standard Normal distributed samples. 
The different methods are documented in the tutorial \emph{Gaussian Process Trajectory Sampling}~\cite{menz2021variance}.  

\paragraph{Parallel Bayesian optimization}
Due to the recent progress made in hardware configurations, it has been of high interest to perform parallel optimizations.
A parallel criterion called qEI~\cite{Ginsbourger2010} was developed to perform Efficient Global Optimization (EGO): the goal is to be able to run batch optimization. 
At each iteration of the algorithm, multiple new sampling points are extracted from the known ones.
These new sampling points are then evaluated using a parallel computing environment.
Five criteria are implemented in \texttt{SMT 2.0}: Kriging Believer ({\footnotesize  \texttt{KB}}), Kriging Believer Upper Bound ({\footnotesize \texttt{KBUB}}), Kriging Believer Lower Bound ({\footnotesize \texttt{KBLB}}), Kriging Believer Random Bound ({\footnotesize \texttt{KBRand}}) and Constant Liar ({\footnotesize \texttt{CLmin}})~\cite{roux2020efficient}.

\section{Conclusion}
\label{sec:concl}

\texttt{SMT 2.0} introduces significant upgrades to the Surrogate Modeling Toolbox.
This new release adds support for hierarchical and mixed variables and improves the surrogate models with a particular focus on Kriging (Gaussian process) models. 
\texttt{SMT 2.0} is distributed through an open-source license and is freely available online\footnote{ \url{https://github.com/SMTorg/SMT}}.
We provide documentation that caters to both users and potential developers\footnote{\url{https://smt.readthedocs.io/en/latest/}}. 
\texttt{SMT 2.0} enables users and developers collaborating on the same project to have a common surrogate modeling tool that facilitates the exchange of methods and reproducibility of results.

SMT has been widely used in aerospace and mechanical modeling applications. 
Moreover, the toolbox is general and can be useful for anyone who needs to use or develop surrogate modeling techniques, regardless of the targeted applications. 
SMT is currently the only open-source toolbox that can build hierarchical and mixed surrogate models.

\recap{ 
\lettrine[lines=2, lhang=0.33, loversize=0.25, findent=1.5em]{T}{his} chapter introduced SMT, an open-source Python package that offers a collection of surrogate modeling methods, sampling techniques, and a set of sample problems. In particular, this chapter focused on \texttt{SMT 2.0}, a major new release of SMT that introduced significant upgrades and new features to the toolbox. The main objectives fulfilled by this chapter are listed as follows. 
\begin{itemize}
    \item A new GP model has been developed for hierarchical variables handling based on a new distance and on a new kernel and has been validated on a neural network optimization problem.
    \item The use of derivatives has been strengthened because derivatives are of high interest for applications of surrogate modeling, such as uncertainty quantification or optimization.
    \item Several new GP models for mixed variables have been implemented. These models are the ones introduced and detailed in Chapter~\ref{c3} and Chapter~\ref{c4}.
    \item New surrogate models, new sampling methods and new applications have been developed and implemented in the \texttt{SMT 2.0} release. 
    \item A better documentation of the open-source software and new tutorials have been released to facilitate the usage for practitioners. 
    \item A special attention has been given to clean and speed up the code in order to reduce the computational times.
\end{itemize}


This chapter corresponds to the article
: {  \textit{Saves, P., Lafage, R., Bartoli, N., Diouane, Y., Bussemaker, J., Lefebvre, T., Morlier, J., Hwang, J. T., Martins, J. R. R. A., “SMT 2.0: A Surrogate Modeling Toolbox with Hierarchical and Mixed Variables Gaussian Processes”, Advances in Engineering Sofware, 2024.}}
}
}

\chapter{Applications to Bayesian optimization} \label{c5}

\vspace{-0.2cm}
\setlength{\fboxrule}{0pt}
\hspace{8cm} \noindent\fbox{%
    \parbox{0.5\textwidth}{%
  \hspace*{1.5cm}    Je suis allé au marché aux oiseaux  \\
 \hspace*{1.5cm}     Et j'ai acheté des oiseaux   \\
 \hspace*{1.5cm}     Pour toi \\
 \hspace*{1.5cm}     Mon amour...  \\
     \hrule \vspace{0.2cm}
     \hspace*{\fill} Paroles, Jacques Prévert}%
} 

\objectif{ 

\lettrine[lines=2, lhang=0.33, loversize=0.25, findent=1.5em]{T}{his} chapter presents practical applications of Bayesian optimization, that is an optimization technique relying on Gaussian process surrogate models. More precisely, the several  objectives of the chapter are listed below.
\begin{itemize}
    \item  To address the need of multi-objective optimization  for expensive-to-evaluate black-box problems featuring mixed hierarchical variables, many variables, high numbers of configurations and multimodal constraints through dedicated optimization algorithms.
    \item To investigate the application of Bayesian optimization to diverse analytical and engineering toy problems by demonstrating its adaptability in scenarios featuring high-dimensions, mixed variable types, constraints, and multi-objective considerations.
    \item To illustrate how Bayesian optimization could serve as an interesting tool for tackling complex optimization challenges across various domains and to emphasize the potential of Bayesian optimization for complex systems optimization, in particular for eco-design of aircraft configurations.  

\end{itemize}


}

\minitoc


\setcounter{section}{-1}
\section{Synthèse du chapitre en français}

Ce chapitre présente des applications concrètes de l'optimisation bayésienne basée sur des modèles de substitution. 
En particulier, les modèles utilisés dans ce travail sont construits à partir de processus gaussiens (\gls{GP} pour Gaussian Process).
Ce chapitre est principalement focalisé sur la conception optimale d'aéronefs plus respectueux de l'environnement ainsi que sur l'optimisation de systèmes complexes. 
Ces domaines d'application ont été motivés par les besoins substantiels d'algorithmes d'optimisation multi-objectif pour optimiser des fonctions boîtes noires coûteuses à évaluer. 
En particulier, les défis auxquels nous avons été confrontés comprenaient notamment, d'une part la gestion de variables hiérarchiques mixtes, et d'autre part la gestion d'un grand nombre de variables et de contraintes d'égalité ou d'inégalité multimodales dans le processus d'optimisation.

Pour commencer, nous avons validé et illustré notre approche algorithmique sur divers cas test analytiques. 
Pour ce faire, nous avons utilisé l'optimiseur \gls{SEGOMOE} et les processus gaussiens développés dans cette thèse et implémentés dans \gls{SMT} pour minimiser le nombre d'évaluations coûteuses de fonctions, même lorsque nous faisions face à un grand nombre de variables de conception. 
\gls{SEGOMOE} s'est avéré performant pour résoudre des problèmes pratiques d'optimisation avec 2 à 5 objectifs, tout en tenant compte de plusieurs contraintes. 
Après validation, nos algorithmes d'optimisation bayésienne ont pu être utilisés pour optimiser des systèmes concrets dans le cadre du projet européen AGILE 4.0 ce qui a permis de générer des résultats optimisés sur des problèmes multi-objectif avec des variables mixtes. 

\section{Introduction}

The previous chapters introduced \gls{GP} models in high-dimension with mixed integer and hierarchical variables. During the thesis, these \gls{GP} surrogate models have been extensively employed for Bayesian optimization and applied to several engineering problems, with a particular focus on applications within the field of aircraft design. 
Notwithstanding, we firstly validated these GP models optimization capabilities on various analytic problems featuring categorical, integer and continuous variables. These models were also validated and applied to multi-objective and constrained optimization problems. Leastwise, it appears essential to acknowledge the international and collaborative efforts that have played a major role, especially for the industrial applications and we would like to thanks all our coworkers for making this collaboration possible and fruitful. 
In summary, this chapter puts into practice the GP models that were developed in the previous chapters while showcasing several Bayesian optimization extensions to constrained, multi-objective and high-dimension optimization problems. 

Section~\ref{sec:bo} introduces the Bayesian optimization framework, while further details regarding the specific analytic and engineering test cases considered and optimized in our works can be found in Section~\ref{sec:anal_optim} and Section~\ref{sec:eng_optim}, respectively. To finish with, conclusions are drawn in Section~\ref{sec:conclusionc6}.

\section{Bayesian optimization}
\label{sec:bo}

The proposed optimization process is based on a sequential enrichment approach named Efficient Global Optimization (\gls{EGO})~\cite{Jones98} together with Super EGO (\gls{SEGO})~\cite{sasena2002exploration}, an evolution of~\gls{EGO} to handle constraints. 
This algorithm is a~\gls{BO} algorithm based on successive enrichment of a \gls{GP} model~\cite{williams2006gaussian} (also denoted by Kriging model~\cite{krige1951statistical}) as detailed in the following sections.

\subsection{Mono-objective efficient global optimization}
\label{sec:ego}

In this section, the objective is to solve a mono-objective unconstrained optimization problem of the form 
\begin{equation}
   \min_{x \in \Omega} \ f(x)    \label{eq:opt_prob_mono_uncons_cont}
\end{equation}
where $f:\mathbb{R}^n  \mapsto \mathbb{R}$ is the objective function and $\Omega \subset \mathbb{R}^n$ represents the bounded continuous design set for the $n$ continuous variables.  
The function $f$ is an expensive-to-evaluate simulation with no exploitable derivative information.
To solve Problem~\ref{eq:opt_prob_mono_uncons_cont}, the \gls{EGO} framework~\cite{Mockus,Jones98} builds a surrogate GP model of the objective function $f$ using a \gls{DOE} of $l$ sampled points in the design domain $\Omega$.
Namely, we assume that our unknown black-box function $f$ follows a GP with mean $\mu^{f}$ and standard deviation $\sigma^f$, $\textit{i.e.}$,
\begin{equation}
f \sim \hat{f}=
\mathcal{GP}
\left(\mu^{f}, [\sigma^f]^2\right). 
\label{eq:GP:f_concl}
\end{equation}
The optimal solution is estimated by iteratively enriching the GP via a search strategy that balances the exploration of the design space $\Omega$ and the minimization of the surrogate model of $f$.
The point that we will evaluate next is the one that gives the best improvement \textit{a priori} according to an acquisition function like the \gls{EI}~\cite{Jones98} defined over the model. The objective value at this new point will then be evaluated and incorporated into the next surrogate model. 
Namely, at each iteration of the EGO method, the search strategy requires solving a maximization sub-problem of the so-called acquisition function~\cite{Frazier, bartoli:hal-02149236,vazquez_convergence_2010,villemonteix2009informational} denoted $\alpha$.
The acquisition function being fully defined by the GP, the search strategy is computationally inexpensive and straightforward and does not require evaluating the computationally expensive function $f$.
The DoE is updated sequentially using the optimal solutions of the sub-problems and this process repeats until reaching a maximum number of evaluations (\textit{i.e.}, the maximal budget). 
The main steps of the EGO framework are summarized in Algorithm~\ref{algo:EGO}.

\smallskip
\begin{algorithm}[H]
\SetAlgoLined
{\textbf{Inputs:}}  Objective function $f$, initial sample $DoE^{(0)}$,  maximal budget \textit{max\_nb\_it}; 

\For{ $t = 0$ \textbf{to} \mbox{max\_nb\_it} - 1}{
\vspace{.2cm}
\begin{adjustwidth}{0pt}{40pt}
\begin{enumerate}
    \item   Build the surrogate model using the GP based on $DoE^{(t)}$\;
    \item Find $\bm{x}^{(t+1)}$, solution of the enrichment maximization sub-problem based on the acquisition function $\alpha^{(t)}$;
    \item Evaluate the objective function at $\bm{x}^{(t+1)}$\;
    \item Update the DoE\;
\end{enumerate}
\end{adjustwidth}
\smallskip
} 
{\textbf{Outputs:}}   The best point found in $DoE^{\textit{max\_nb\_it \ } }$; 
\caption{The Efficient Global Optimization framework.}
 \label{algo:EGO}
\end{algorithm}

\paragraph{The enrichment sub-problem}

The BO framework relies on the information provided by the GP (namely, $\hat\mu^{(t)}$ and $\hat\sigma^{(t)}$) to build the enrichment strategy. The latter is guided by the following maximization sub-problem:
\begin{equation}
    \max\limits_{\bm{x} \in \Omega} \alpha^{(t)}(\bm{x}),
    \label{eq:iner_opt}
\end{equation}
where $\alpha^{(t)}: \mathbb{R}^d \mapsto \mathbb{R}$ is the chosen acquisition function based on the GP at step $t$.
There are numerous acquisition functions in the literature~\cite{frazierTutorialBayesianOptimization2018, ShahriariTakingHumanOut2016, WangMaxvalueentropysearch2017, bartoli:hal-02149236} and the choice of $\alpha$ is essential for the enrichment process.
In particular, the \gls{EI}~\cite{Jones2001JOGO} acquisition function is the most commonly used in \gls{BO}.
Considering the $t^{\text{th}}$ iteration of the BO framework, the expression of $\alpha_{EI}^{(t)}$ depends on the predictions $\hat\mu^{f(t)}$ and $\hat\sigma^{f(t)}$.
For a given point $\bm{x}\in \Omega$, if  $\hat\sigma^{(t)}(\bm{x}) = 0$, then $\alpha_{EI}^{(t)}(\bm{x})=0$.
Otherwise,
\begin{equation}
    \small
    \alpha_{EI}^{(t)}(\bm{x})=\left(y_{min}^{(t)} - \hat{\mu}^{f(t)}(\bm{x})\right) \Phi \left( \frac{y_{min}^{(t)} - \hat{\mu}^{f(t)}(\bm{x})}{\hat\sigma^{f(t)}(\bm{x})} \right) + \hat\sigma^{(t)}(\bm{x}) \phi \left(  \frac{y_{min}^{(t)} - \hat{\mu}^{f(t)}(\bm{x})}{\hat\sigma^{f(t)}(\bm{x})} \right),
    \label{eq:EI}
\end{equation}
where the functions $\Phi$ and $\phi$ are, respectively, the cumulative distribution function and the probability density function of the standard normal distribution. $y_{min}^{(t)}$ is the current minimum given by $y_{min}^{(t)} = \min \bm{Y}^{(t)}$, with $ \bm{Y}^{(t)} = \{f \bm{(x}^i), \forall  \bm{x}^i \in DoE^{(t)} \} $.
Mono-objective EGO without constraint is implemented in \gls{SMT}~\cite{SMT2019,saves2023smt}.
However, within the EGO framework, it is possible to tackle problems with non linear constraints using different mechanisms for different computational costs~\cite{const,  frazierTutorialBayesianOptimization2018,Priem_thesis,ShahriariTakingHumanOut2016,sasena2002exploration}. Among these methods, \gls{SEGO} is of high interest for its ability to learn and take into account the uncertainties of the constraints as exposed in the next section.

\paragraph{Efficient Global Optimization with Random and Supervised Embeddings}

\begin{figure}[htb]
    \centering
    \includegraphics[width=0.9\textwidth,height=8cm]{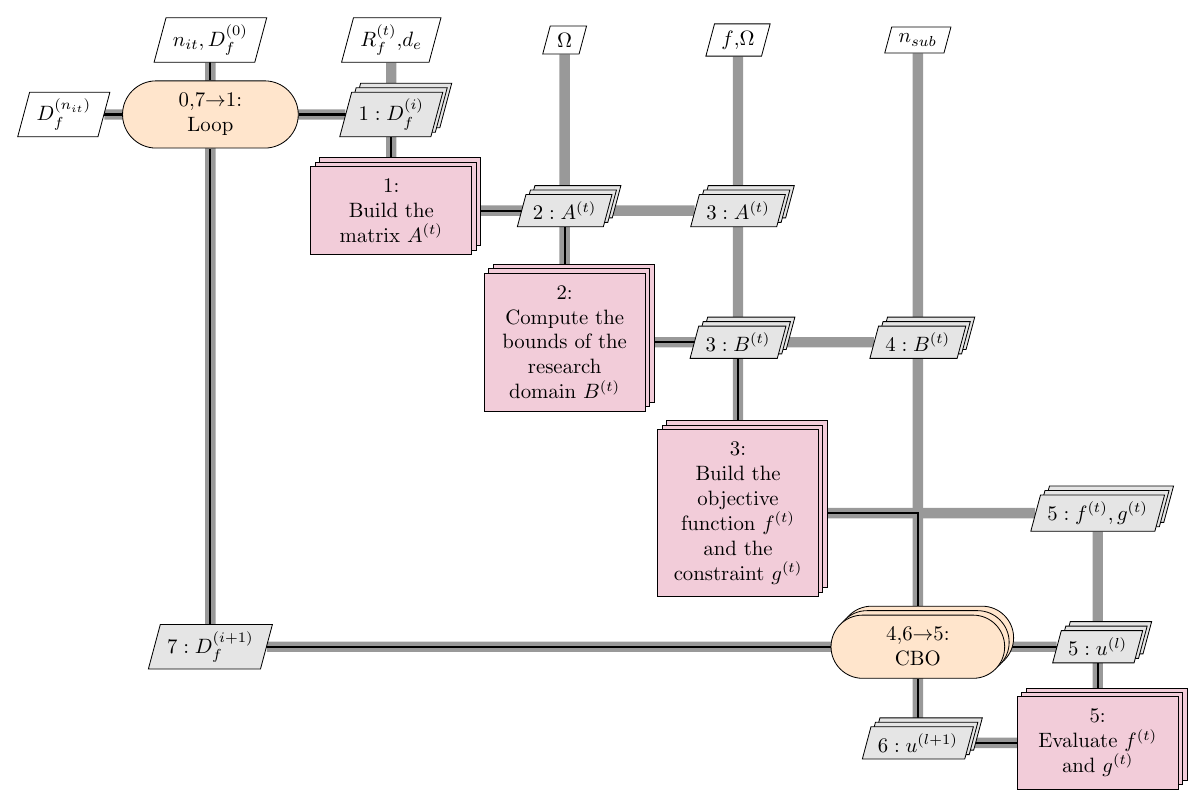}
    \caption{XDSM diagram of the EGORSE framework.}
    \label{fig:xdsm}
\end{figure}

\citet{EGORSE} introduces the algorithm \gls{EGORSE} described succinctly as follows. First, the objective function is supposed to depend only on the effective dimensions $d_e \ll n$ where it is assumed that it exists a function $f_{\As}: \mathbb{R}^{d_e} \mapsto \mathbb{R}$ such as $f_{\As}(\A\bm{x}) = f(\bm{x})$ with $\A \in \mathbb{R}^{d_e \times n}$, $\As = \left\{ \bm{u} = \A\bm{x} , \ \forall \bm{x} \in \Omega \right\}$ and\ $\Omega = [-1,1]^n$~\cite{wang2017batched,RREMBO}.
The idea is then to perform the optimization procedure in the reduced linear subspace $\As$ so that the number of hyper-parameters to estimate and dimension of the design space are reduced to $d_e$ instead of $n$.
This allows to build inexpensive GPs and will ease the acquisition function optimization.
Using a subspace (based on $\A$) for the optimization requires finding the effective dimension of the reduced design space $\mathcal{B} \subset \mathbb{R}^{d_e}$ as well as the backward application $\gamma: \mathcal{B} \mapsto \Omega$.
\gls{EGORSE} focuses on defining the optimization problem with a linear subspace and developing efficient procedures for constructing these embedding subspaces.
Most existing high dimension Bayesian optimization methods rely on random linear subspaces  meaning that no information is used to incorporate \textit{a priori} which may slow down the optimization process. 
In \gls{EGORSE}, a recursive search, with $T\in \mathbb{N}$ supervised reduction dimension methods, is performed to find supervised linear subspaces so that the most important search directions for space exploration are included. 
Using an initial DoE, one can use a \gls{PLS} regression \cite{hellandStructurePartialLeast1988} to build such linear embedding prior for the optimization process.
Furthermore, the new search design of the optimization problem (within the linear embedding subspace) is a necessary step.
Most methods rely on a classic optimization problem formulation that may limit the process performance due to very restricted new design space.
Here, once an appropriate linear subspace is found, the optimization problem is turned into a constrained optimization problem (CBO) to limit the computational cost of the algorithm~\cite{frazierTutorialBayesianOptimization2018,Priem_thesis,bartoli:hal-02149236,SEGO-UTB}.
Among the possibilities of coupling supervised and random embeddings, using both  Gaussian random and PLS transfer matrix $\A$ lead to the best results~\cite{Priem_thesis}. The  flow chart of the method is described by the \gls{XDSM} \cite{Lambe2012} diagram of Figure~\ref{fig:xdsm}.

\subsection{Efficient global optimization under constraints}
\label{sec:sego}

In this section, we consider that all the design variables are continuous. Namely, in this section, the design space will be restricted to $\Omega \subset \mathbb{R}^n$. 
The problem to tackle writes
\begin{equation}
   \min_{ x \in \Omega}  \left\{ f(x) ~~\mbox{s.t.}~~ g(w)\leq 0 ~\mbox{and}~ h(w)=0 \right\}
    \label{eq:opt_prob_cont_cons}
\end{equation}
where $f:\mathbb{R}^n \mapsto \mathbb{R}$ is the objective function, $g:\mathbb{R}^n  \mapsto \mathbb{R}^m$ gives the inequality constraints, and $h:\mathbb{R}^n  \mapsto \mathbb{R}^p$ returns the equality constraints.
$\Omega \subset \mathbb{R}^n$ represents the bounded continuous design set for the $n$ continuous variables.  The functions $f$, $g$, and $h$ are typically expensive-to-evaluate simulations with no exploitable derivative information.

To solve Problem~\ref{eq:opt_prob_cont_cons}, we are using \gls{SEGO}, an extension of the EGO framework that has been described in Section~\ref{sec:ego}. In fact, within EGO, at a given iteration $t$, a GP surrogate model is computed based on the current DoE to approximate the black-box $f$.  Henceforth, one wants to estimate the best new point to evaluate, as it is costly, by taking into account the model information to converge as fast as possible to the real optimum of the black-box. The next point that will be evaluated is the one that gives the best improvement \textit{a priori} according to the acquisition function such that the constraints are respected.
Once the point has been chosen inside the design space, the objective and constraints are evaluated at this point and incorporated into the next surrogate model. 
The hyperparameters that characterize and define the model are thus updated at each iteration until the stopping criterion is met. 
The Bayesian optimization process is made from these GPs in an iterative manner.  
To tackle constrained Bayesian optimization, EGO was extended to \gls{SEGO} by~\citet{sasena2002exploration}. \gls{SEGO} uses surrogate models of the constraints to give an estimation of the search space $\Omega_f$ in which the function $f$ is optimized through a given criterion. The latter was enhanced to tackle multi-modal and equality constraints with the Upper Trust Bound (UTB) criterion~\cite{SEGO-UTB,bartoli:hal-02149236,BRAC_AIAA}.
The acquisition function that we use is the \gls{WB2S} criterion (Watson and Barnes $2^{nd}$ criterion with scaling)~\cite{bartoli:hal-02149236} that is known to be more robust than the Expected improvement (\gls{EI}) criterion, especially in high dimension~\cite{jesus2021surrogate}. 
WB2s is smoother than the WB2 criterion~\cite{sasena2002exploration} and is less multimodal than EI. Algorithm~\ref{alg:optim} details the SEGO optimization procedure.

\smallskip
\begin{algorithm}[htb]
\SetAlgoLined
{\textbf{Inputs:}}  Initial DoE $\mathscr{D}_0$ and set $t=0$. Search space $ \Omega$. 
\\
\While{the stopping criterion is not satisfied}{
\vspace{.2cm}
\begin{adjustwidth}{0pt}{40pt}
\begin{enumerate}
    \item Build the surrogate model of the objective function to obtain the mean and standard deviation prediction at a given point: $(\hat{\mu}^f,\hat{\sigma}^f)$ from the DoE $\mathscr{D}_t$.
    \item Build the surrogate models for every constraints ($\hat{\mu}^{g_i},\hat{\sigma}^{g_i}$) and ($\hat{\mu}^{h_j},\hat{\sigma}^{h_j}$) from the DoE $\mathscr{D}_t$ to compute an estimation of the search space $\Omega_f$.
    \item Construct the acquisition function $\text{WB2s}(.) = \phi\left( \hat{\mu}^f(.),\hat{\sigma}^f(.) \right)$ from the objective model. 
    \item Maximize the acquisition function WB2s over $\Omega_f$: $x_{t}=  \underset{x \in \Omega_f}{\arg \max}\   \text{WB2s}(x)$. 
    \item Add $x_{t}, f(x_{t}), g(x_{t}), h(x_{t})$ to the DoE  $\mathscr{D}_{t+1}$. Increment $t$.
\end{enumerate}
\end{adjustwidth}
\smallskip
} 
{\textbf{Outputs:}}   The best point found in the final DoE; 
\caption{SEGO for continuous inputs.}
\label{alg:optim}
\end{algorithm}

\smallskip

To optimize expensive-to-evaluate black-box functions, we are using the \gls{SEGOMOE} algorithm~\cite{bartoli:hal-02149236} that combines SEGO with the Mixture Of Experts (MOE)~\cite{dimitri:SM02011, Liem2015}. 
The idea of MOE is to use an adaptive mixture of Kriging based models to tackle high dimension problems and heterogeneous functions. MOE approximates complex functions with heterogeneous behaviour by combining local surrogate models into a global one. In order to consider high dimension functions and to approximate objective functions and constraints, SEGOMOE uses adapted local Kriging-based models~\cite{Bouhlel18, bouhlel_KPLSK}. Moreover, some recent developments have been made in SEGOMOE to take into account non linear constraints~\cite{BRAC_AIAA}, mixed integer variables~\cite{SciTech_cat} and multi-objective applications~\cite{grapin_constrained_2022}. 
The general \gls{SEGOMOE} algorithm is developed in Python by ONERA and ISAE-SUPAERO and its performance has been validated and proven on different analytical and industrial test cases~\cite{bartoli:hal-02149236, AGILE_PIAS_MDO, Effectiveness}.

\subsection{Multi-objective efficient global optimization under constraints}
\label{sec:segomoo}

In this section, we consider multiple objectives and constraints with mixed integer variables. Namely, the problem to solve is given by
\begin{equation}
   \min_{w=(x,z,c) \in \Omega \times S \times \mathbb{F}^l}  \left\{ f(w) ~~\mbox{s.t.}~~ g(w)\leq 0 ~\mbox{and}~ h(w)=0 \right\}
    \label{eq:opt_prob_multi_cont}
\end{equation}

$f:\mathbb{R}^n \times  \mathbb{Z}^m \times \mathbb{F}^l \mapsto \mathbb{R}^q$ are the objective functions, $g:\mathbb{R}^n \times  \mathbb{Z}^m \times \mathbb{F}^l \mapsto \mathbb{R}^m$ gives the inequality constraints, and $h:\mathbb{R}^n \times  \mathbb{Z}^m \times \mathbb{F}^l \mapsto \mathbb{R}^p$ returns the equality constraints.
 $\Omega \subset \mathbb{R}^n$ represents the bounded continuous design set for the $n$ continuous variables.  $S \subset \mathbb{Z}^m$ represents the bounded integer set where $L_1, . . . , L_m$ are the numbers of levels of the $m$ quantitative integer variables on which we can define an order relation and $ \mathbb{F}^l = \{1, \ldots, L_1\} \times \{1, \ldots, L_2\} \times  \ldots \times \{1, \ldots, L_l\}$ is the design space for $l$ categorical qualitative variables with their respective  $L_1, . . . , L_l$ levels.
The functions $f$, $g$, and $h$ are typically expensive-to-evaluate simulations with no exploitable derivative information.

Some continuous, integer and categorical variables are involved within the different application cases that correspond to the problem~\ref{eq:opt_prob_multi_cont}, so the GP surrogate models have to be adapted to deal with~\cite{rufatoacreating,CompBO,deshwal2021bayesian,daxberger_mixed-variable_2020,oh2018bock,sheikh2022bayesian}.
For optimization, we use the continuous relaxation introduced by~\cite{GMHL} that relies on a one-hot encoding strategy~\cite{one-hot} to transform integer and categorical inputs into continuous ones.  In fact, the design  space  $ \Omega \times S \times \mathbb{F}^l $ is relaxed to a continuous space $ \Omega'$  constructed on the following way~\cite{SciTech_cat}:
\begin{itemize}
\item $\forall i \in \{1, \ldots,\ell\}$, the integer variable $z_i$ is relaxed within its bounds and treated as continuous. 
\item   $\forall j \in \{1, \ldots,l\}$, we use a relaxed one-hot encoding~\cite{one-hot} for the categorical variable $c_j$ (and its $L_j$ 
 associated levels) and add $L_j$ new continuous dimensions into $ \Omega'$. 
\end{itemize}
Therefore, we get, after relaxation, a new design space $\Omega'\subseteq \mathbb{R}^{d'}$ where $d'= d+\ell+ \sum_{j=1}^l L_j >d+\ell+l$. 
The initial Design of Experiments can be relaxed and expressed in the relaxed space $\Omega'$. Each costly function $f_i(\bm{x})$ is approximated by a GP characterized by its mean $\mu^{f_i} : $
$ \mathbb{R}^{d'} \to \mathbb{R}$ and its standard deviation $\sigma^{f_i} : \mathbb{R}^{d'} \to \mathbb{R}$ as
    $$\hat{f}_i(\bm{x})\sim \mathcal{GP}(\mu^{f_i}(\bm{x}),[\sigma^{f_i}(\bm{x})]^2) \quad i=1,\ldots, q.$$
For multi-objective, we assume that the $q$ components $f_i$ of $\bm{f}$ are independent and we define $\hat{\bm{f}}$ as the surrogate model of the objectives vector as 
$$\hat{\bm{f}}(x)\sim \mathcal{GP}(\hat{\y}(\bm{x}),\Sigma(\bm{x})) $$
where $\hat{\y}(\bm{x}): \mathbb{R}^{d'} \to \mathbb{R}^q$ is the GP prediction vector given by $[\mu^{f_1}(\bm{x}),\ldots \mu^{f_q}(\bm{x})]$ and $\Sigma(\bm{x})$ is a diagonal matrix whose diagonal vector is given by $ [\sigma^{f_i}(\bm{x})]^2, \forall i=1,\ldots, q$.
As before, the initial multi-objective optimization problem under constraints is replaced by an infill problem based on an acquisition function defined as 
\begin{equation}\label{eq:opt_acquisition_problem}
\left\lbrace
\begin{array}{l}
 \displaystyle \max_{{x}\in \mathbb{R}^{d'}} \quad \alpha^{\mbox{reg}}_{\bm{f}}(\bm{x}) \\
\mbox{s.t.}  \quad \hat{g}_{j}(\bm{x}) \leq 0\\
\quad  \quad \   \hat{h}_{k}(\bm{x}) = 0,\\
\end{array}\right.
\end{equation}
where $\hat{g_j}(\bm{x})$ and $\hat{h}_k(\bm{x})$ corresponds to the mean prediction of the GP constraint models and the regularized acquisition function~\cite{grapin_constrained_2022} is defined by
\begin{eqnarray} \label{eq:alpha:reg}
\alpha^{\mbox{reg}}_{\bm{f}}(\bm{x}) &:=& \gamma\;\alpha_{\bm{f}}(\bm{x}) -\psi(\mu_{\bm{f}}(\bm{x}))
\end{eqnarray}
where $\alpha_{\bm{f}}(\bm{x})$ is a standard scalar multi-objective acquisition function (EHVI, PI, MPI, \dots) depending on $\hat{\y}(\bm{x})$ and $\Sigma(\bm{x})$,
and $\gamma$ is a constant parameter. The function $\psi: \mathbb{R}^{d'} \to \mathbb{R}$ is a scalarization operator. Different choices exist for the function $\psi$ and two options are investigated in~\citet{grapin_constrained_2022}. Namely, for a given $\hat{\y}(\bm{x})\in \mathbb{R}^{q}$, we consider
\begin{eqnarray*}
\begin{cases}
& $$(\mbox{reg}=\max):\psi(\hat{\y}(\bm{x})) :=  \max_{i \leq q}\hyi$$,  \\
& $$(\mbox{reg}=\mbox{sum}) : \psi(\hat{\y}(\bm{x})):= \sum_{i=1}^{q} \hyi$$.  
\end{cases}
\end{eqnarray*}

\begin{figure}[!htb]
    \centering
    \includegraphics[scale=0.8]{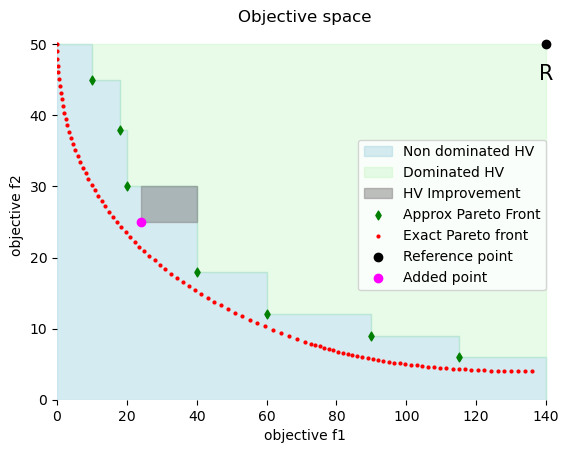}
    \caption{Hypervolume improvement: the hypervolume indicator of the non-dominated set (green points)
corresponds to the area dominated by it, up to $R$ (reference point in black).
The gray rectangle is the hypervolume improvement brought by the new added
point in magenta.}
    \label{fig:HVI}
\end{figure}

Due to potentially conflicting objectives, the solution of the optimization is not unique but a range of solutions is proposed.
The trade-off between these solutions is characterized by the notion of dominance:
a solution $\bm{x}$ is said to dominate another solution $\bm{x}'$ and denoted by $\bm{x} \preceq \bm{x}'$ if
$$f_i(\bm{x}) \le f_i(\bm{x}') \quad \forall i=1,\ldots, q.$$
The set of solutions representing optimal trade-offs is referred as Pareto set (PS) and the corresponding image of the PS in the objective space is known as the Pareto front~(PF), \textit{i.e}, 
$$\mbox{PF}:=\{\bm{f}(\bm{x})| \bm{x} \in \text{PS} \}.$$
Figure~\ref{fig:HVI} illustrates a Pareto front (red points) relative to two objectives. An approximation to this PF is given by the scattered green dots and the associated dominated hypervolume is given by the green area. 
Concerning the different BO criteria relative to the hypervolume improvement, an illustration is proposed in~\figref{fig:HVI}. The Expected Hyper-Volume Improvement (EHVI~\cite{zitzler2003performance,emmerich2006single}), the Probability of Improvement (PI~\cite{Jones2001JOGO}), and the Minimum of Probability of Improvement (MPI~\cite{rahat2017alternative}) are some multi-objective extensions of the well known Expected Improvement (\gls{EI}~\cite{Jones98}).
As seen in~\figref{fig:HVI}, the idea is to measure how much hypervolume improvement (grey area) could be achieved by evaluating a new point (magenta point) while considering the prediction and uncertainty of the GPs. Also, these criteria differ by how much they favor well-spread solutions. For instance, with EHVI, the hypervolume increase is small when adding a new value close to an already observed datum in the objective space.
To solve~\eqnref{eq:opt_acquisition_problem}, any optimization algorithm capable of considering non linear constraints can be used. This optimizer could either be derivative-free such as COBYLA (Constrained Optimization BY Linear Approximation~\cite{COBYLA}) or gradient-based such as SLSQP (Sequential Least Squares Programming~\cite{kraft1988software}) or SNOPT (Sparse Nonlinear
OPTimizer~\cite{gill2005snopt}) together with a multistart strategy.
This adaptive process is repeated until the total budget is reached. The feasible points of the final database represent the known Pareto optimal points. Nevertheless, as the set of known points has been enriched sequentially to increase the hypervolume, the final database can be used to build GP models for objectives and constraints as a post-processing step. These final GP models can be coupled to an evolutionary algorithm to deal with the multi-objective constrained problem and retrieve the approximated PF. The well known \gls{NSGA2} algorithm~\cite{nsga2} of the pymoo~\cite{pymoo} toolbox\footnote{\url{https://pymoo.org/}} is used for that purpose  and ultimately provides the predicted PF based on GPs with almost no additional computational cost.

In the end, the proposed strategy provides two outputs: the PF database and the predicted PF together with their associated respective PS. Comparing the proximity between these two PF is a good criterion to know if some additional enrichment points are needed or if the accuracy is sufficient.  
Algorithm~\ref{algo:SEGOMOE_Multiobj} details the \gls{SEGOMOE} optimization procedure. 

\smallskip
\begin{algorithm}[htb]
\SetAlgoLined
{\textbf{Inputs:}}  Initial DoE $\mathscr{D}_0$ and set $t=0$;   

\While{the stopping criterion is not satisfied}{
\vspace{.2cm}
\begin{adjustwidth}{0pt}{40pt}
\begin{enumerate}
    \item  Relax continuously integer and categorical input variables to a real bounded space $\Omega'$ of dimension $d'= d+\ell+ \sum_{j=1}^l   L_j$. Namely, we continuously relax the mixed categorical DoE to a continuous DoE $\mathscr{D}_t$ using the relaxation procedure;
    \item Build the GP model for each objective function $f_i(\bm{x})$ and each constraint $g_j(\bm{x})$, $h_k(\bm{x})$ related to the continuous DoE with PLS to reduce the number of the hyperparameters and compute an estimation of the search space $\Omega_f$;
    \item Build the acquisition function $\alpha^{\mbox{reg}}_{\bm{f}}(\bm{x})$;
    \item Maximize the acquisition function within the feasible domain  $\Omega_f$: $$\bm{x}_{t}:=  \underset{\bm{x} \in \Omega_f}{\arg \max}  \  \alpha^{\mbox{reg}}_{\bm{f}}(\bm{x})$$
    \item Add $\bm{x}_{t}, \bm{f}(\bm{x}_{t}) $ and $ g_1(\bm{x}_{t}), \ldots g_m(\bm{x}_{t}) , h_1(\bm{x}_{t}) \ldots h_p(\bm{x}_{t}) $ to the DoE  $\mathscr{D}_{t+1}$. Increment $t$;
\end{enumerate}
\end{adjustwidth}
\smallskip
} 
{\textbf{Post-process}}: Use the final database to build GP models for $f_i(\bm{x})$, $g_j(\bm{x})$ and $h_k(\bm{x})$. Then, apply NSGA-II algorithm to construct the predicted PF;\\
{\textbf{Outputs:}}  The PF database and the predicted PF; 
\caption{SEGOMOE for constrained multi-objective and mixed-integer problems.}
 \label{algo:SEGOMOE_Multiobj}
\end{algorithm}

\section{Analytic optimization problems}
\label{sec:anal_optim}

The analytic problems over which we have validated the Bayesian optimization and our proposed extensions are listed in~\tabref{tab:recap_anal} and come from several state-of-the-art papers~\cite{Pelamatti, Pelamatti2020, Roustant, CAT-EGO, AMIEGO, Gower, Vanaret, vanaret2024interval, pelamattihier, ZDT, binh1997mobes,zapotecas2018review,osyczka1995new}. These test cases have been selected to validate at least one of these aspects: high-dimension, mixed integer variables, constraints and multi-objective. The considered validation test cases are detailed  in~\tabref{tab:recap_anal} indicating, for every considered test case, its name, its number of objectives, constraints and variables, the nature of these variables and the reference of the paper in which we used this test case. 
Some of these test cases have already been presented in the previous chapters. For instance, the "Cosine curves" problem have been introduced in Chapter~\ref{c2}, the "Toy set 1" problems have been introduced in Chapter~\ref{c3} and the "Hierarchical Goldstein" have been introduced in Chapter~\ref{c4}. 

The "Modified Branin 10" and "Modified Branin 100" have been used in~\citet{EGORSE} to illustrate the efficiency of the \gls{EGORSE} algorithm to tackle really high dimension. 
The "Branin 1", "Branin 2", "Goldstein 1", "Goldstein 2" and "Toy set 2" have been used in~\cite{SciTech_cat} to validate the mixed integer GPs that use Gower distance, continuous relaxation and continuous relaxation coupled to partial least squares in the context of Bayesian optimization with mixed variables and constraints. This also applies to "Branin 3", "Branin 4" and "Branin 5" used for the same purpose in~\cite{saves2021constrained}.
To finish with, in~\citet{grapin_constrained_2022}, we relied on several test cases to validate multi-objective optimization. The ZDT test cases (ZDT 1, ZDT 2 and ZDT 3) have been declined into 3 test cases each with respectively 2, 5 and 10 variables for a total of 9 test cases. In addition the test cases "BNH","TNK" and "OSY" have been exploited to validate multi-objective Bayesian optimization under constraints.

\begin{table}[htb]
\centering
\caption{Definition of the analytical optimization problems.}
\resizebox{\linewidth}{!}{%
\small

\begin{tabular}{ccccc}
\hline
\textbf{Name} & \# of objs &  \# of cons  & \# of vars [\textbf{cont}, \textbf{int}, \textbf{cat}] & Reference \\
\hline
\textbf{Modified Branin 10} & 1 & --& [10,0,0]  & ~\cite{EGORSE} \\
\textbf{Modified Branin 100} & 1 & --& [100,0,0]  & ~\cite{EGORSE} \\
\textbf{Hierarchical Goldstein} & 1& -- & [5,4,1(4 levels)] & ~\cite{saves2023smt} \\
\textbf{Toy set 1} & 1& -- & [1,0,1(10 levels)] & ~\cite{saves2023smt} \\
\textbf{Toy set 2} & 1& -- & [1,0,1(10 levels)] &  ~\cite{SciTech_cat} \\
\textbf{Branin 1} & 1& -- & [2,0,2(2,2 levels)] &  ~\cite{SciTech_cat} \\
\textbf{Branin 2} & 1& -- & [10,0,2(2,2 levels)] &  ~\cite{SciTech_cat} \\
\textbf{Branin 3} & 1& 1 & [2,0,2(2,2 levels)] & ~\cite{saves2021constrained}\\
\textbf{Branin 4} & 1& 1 & [10,0,2(2,2 levels)] & ~\cite{saves2021constrained}\\
\textbf{Branin 5} & 1& -- & [1,1,0] & ~\cite{saves2021constrained}\\
\textbf{Cosine curves} & 1& -- & [1,0,1(13 levels)] &  ~\cite{SciTech_cat} \\  
\textbf{Goldstein 1} & 1& -- & [2,0,2(3,3 levels)] &  ~\cite{SciTech_cat} \\
\textbf{Goldstein 2} & 1& 1 & [2,0,2(3,3 levels)] &  ~\cite{SciTech_cat} \\
\textbf{ZDT 1-2d}  & 2 & -- & [2,0,0] & ~\cite{grapin_constrained_2022} \\
\textbf{ZDT 1-5d}  & 2 & -- & [5,0,0] & ~\cite{grapin_constrained_2022} \\
\textbf{ZDT 1-10d}  & 2 & -- & [10,0,0] & ~\cite{grapin_constrained_2022} \\
\textbf{ZDT 2-2d}  & 2 & -- & [2,0,0] & ~\cite{grapin_constrained_2022} \\
\textbf{ZDT 2-5d}  & 2 & -- & [5,0,0] & ~\cite{grapin_constrained_2022} \\
\textbf{ZDT 2-10d}  & 2 & -- & [10,0,0] & ~\cite{grapin_constrained_2022} \\
\textbf{ZDT 3-2d}  & 2 & -- & [2,0,0] & ~\cite{grapin_constrained_2022} \\
\textbf{ZDT 3-5d}  & 2 & -- & [5,0,0] & ~\cite{grapin_constrained_2022} \\
\textbf{ZDT 3-10d}  & 2 & -- & [10,0,0] & ~\cite{grapin_constrained_2022} \\
\textbf{BNH}  & 2 & 2 & [2,0,0] & ~\cite{grapin_constrained_2022} \\
\textbf{TNK}  & 2 & 2 & [2,0,0] & ~\cite{grapin_constrained_2022} \\
\textbf{OSY}  & 2 & 6 & [6,0,0] & ~\cite{grapin_constrained_2022} \\
\hline
\end{tabular}
}
\label{tab:recap_anal}
\end{table}

\subsection{High-dimension validation results}

\begin{figure}[H]
    \centering
    \subfloat[MB\_10 iteration convergence plot. \label{fig:ALL:mb10_nb}]{\includegraphics[width=0.50\textwidth]{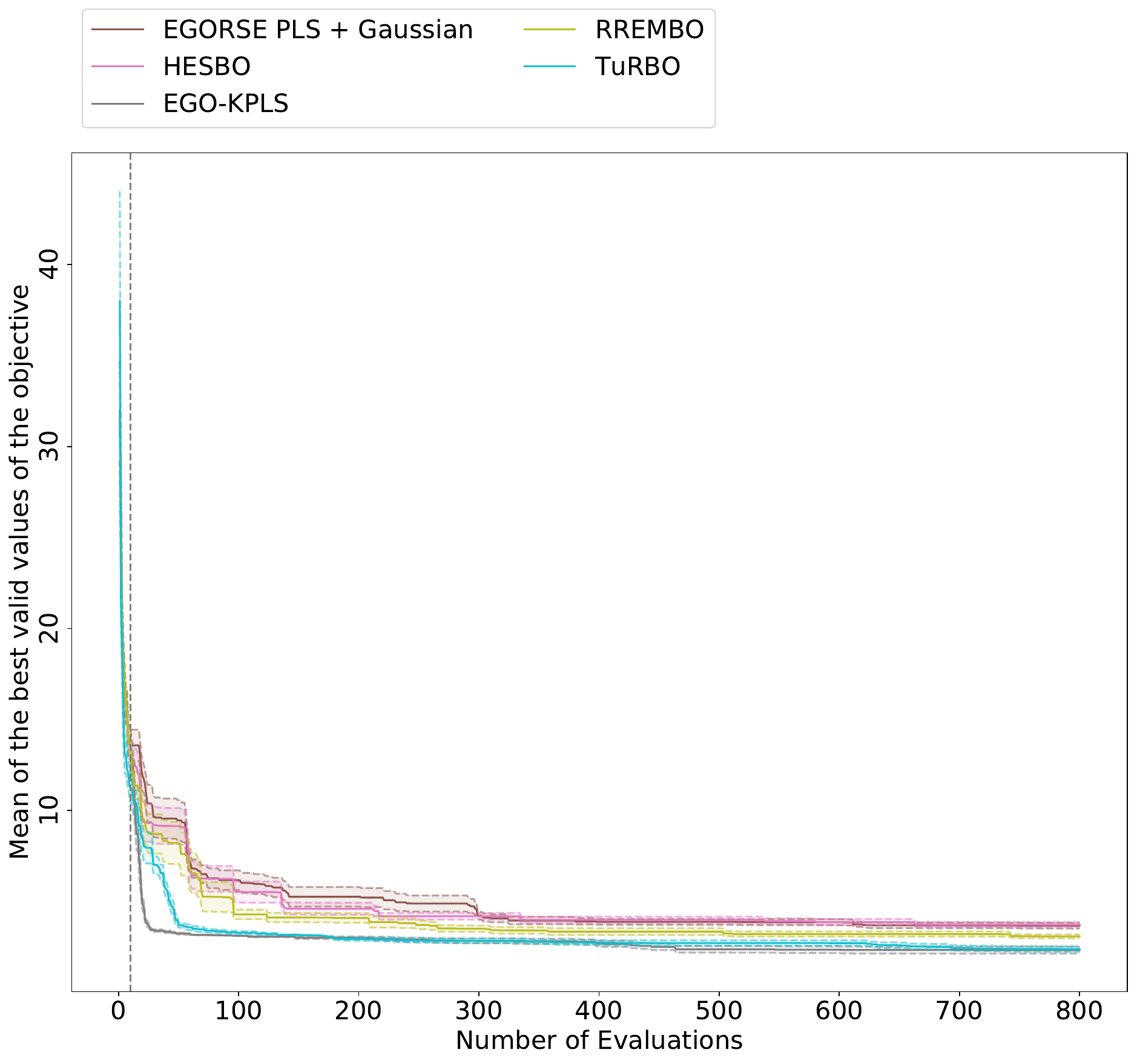}}
    \subfloat[MB\_100 iteration convergence plot. \label{fig:ALL:mb100_nb}]{\includegraphics[width=0.50\textwidth]{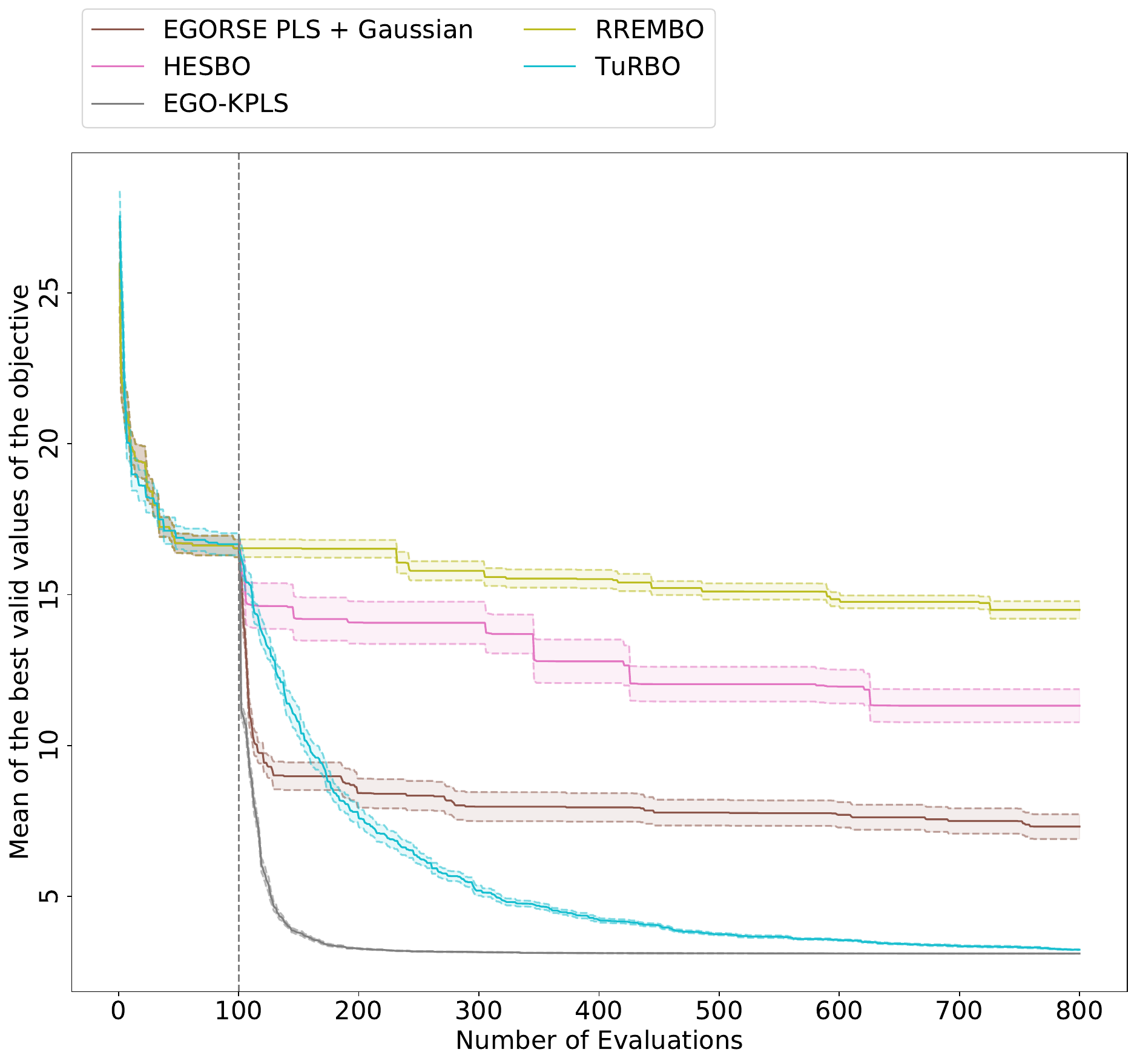}} \\
    \subfloat[MB\_10 time convergence plot. \label{fig:ALL:mb10_time}]{\includegraphics[width=0.50\textwidth]{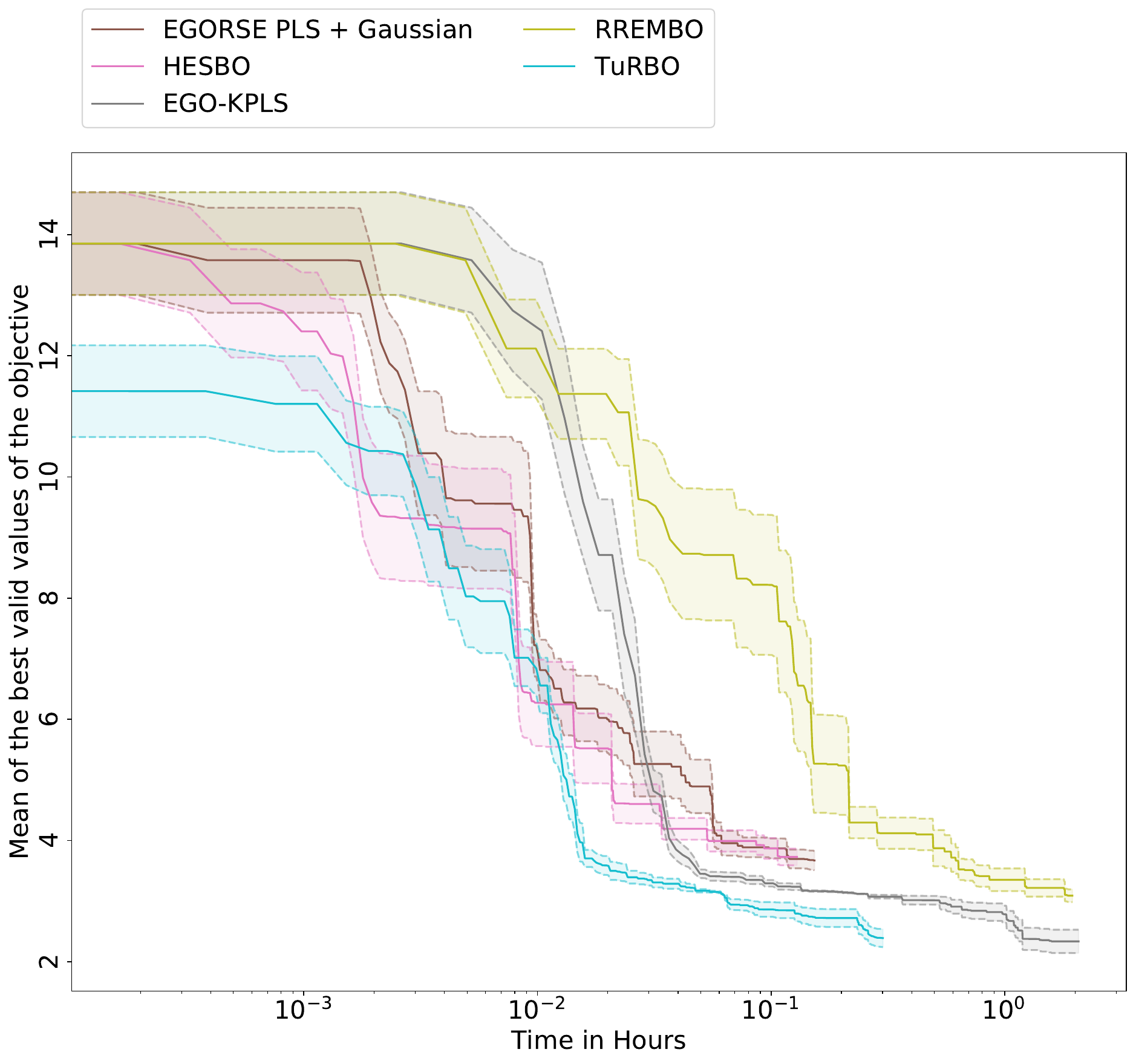}}
    \subfloat[MB\_100 time convergence plot. \label{fig:ALL:mb100_time}]{\includegraphics[width=0.50\textwidth]{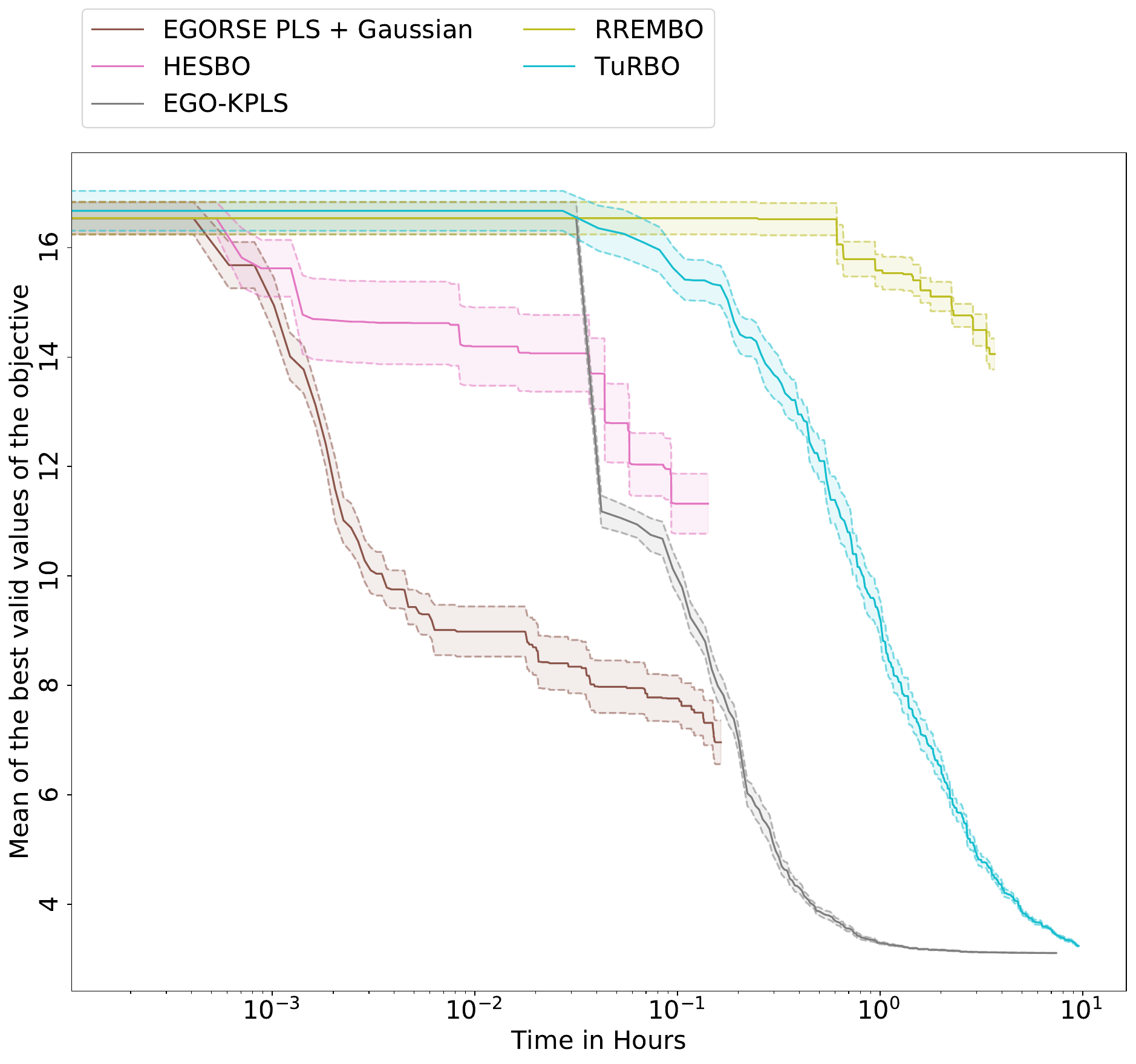}}
    \caption{Iteration and time convergence plots for 5 algorithms on the MB\_10 and MB\_100 problems. The grey vertical line shows the size of the initial DoE.}
    \label{fig:ALL}
\end{figure}

The considered class of problems is an adjustment of the Modified Branin (MB) problem~\cite{ParrInfillsamplingcriteria2012} whose number of design variables is artificially increased.
This problem is commonly used in the literature~\cite{RREMBO,wang2017batched,Hesbo} and it is  defined as follows:
\begin{equation}
    \min_{\bm{u} \in \Omega_1}{f_1(\bm{u})},
\end{equation}
where $\Omega_1 = [-5,10] \times [0,15]$ and
\begin{equation}
    f_1(\bm{u}) = \left[ \left( u_2 - \frac{5.1u_1^2}{4\pi^2} + \frac{5u_1}{\pi} - 6 \right)^2 +\left( 10 - \frac{10}{8\pi} \right) \cos{(u_1) + 1} \right] + \frac{5u_1 + 25}{15}.
\end{equation}
The modified version of the Branin problem is selected because it count three local minima including a global one.
The value of the global optimum is about ${\text{MB\_$d$}}_{min} = 1.1$.
Furthermore, the problem is normalized to have $\bm{u} \in [-1,1]^2$.
To artificially increase the number of design variables, a random matrix $\A_d \in \mathbb{R}^{2 \times d}$ is generated such that for all $\bm{x} \in [-1,1]^d$, $\A_d \bm{x} = \bm{u}$ belongs to $[-1,1]^2$.
An objective function MB\_$d$, where $d$ is the number of design variables, is defined such that $\text{MB\_$d$}(\bm{x}) = f_1(\A_d\bm{x})$.
Eventually, we solve the following optimization problem:
\begin{equation}
    \min\limits_{x \in [-1,1]^d} \text{MB\_$d$}(\bm{x}) = \min\limits_{x \in [-1,1]^d} f_1(\A_d\bm{x}).
\end{equation}
In the following, numerical experiments are conducted on two functions of respective dimension 10 and 100. These two test functions are denoted $\text{MB\_$10$}$ and $\text{MB\_$100$}$.
 {EGORSE} is compared to the following state-of-the-art algorithms:
\begin{itemize}
    \item  {TuRBO}~\cite{Turbo}: a Bayesian algorithm using confidence regions to favor the exploitation of the DoE data.
    Tests are performed with the TuRBO\footnote{https://github.com/uber-research/TuRBO} Python toolbox.
    \item    {EGO-KPLS}~\cite{SMT2019}: a Bayesian optimization method relying on the reduction of the number of GP hyper-parameters.
    This allows to speed up the GP building. 
    The SEGOMOE~\cite{bartoli:hal-02149236} Python toolbox is used.
    All the hyper-parameters of this algorithm are the default ones.
    The number of principal components for the KPLS model is set to two.
    \item   {RREMBO}~\cite{RREMBO}: a Bayesian optimization method using the random Gaussian transfer matrix to reduce the number of dimensions of the optimization problem.
      {RREMBO}\footnote{https://github.com/mbinois/RRembo} implementation of this algorithm is used.
    The parameters are also set by default.
    \item {HESBO}~\cite{Hesbo}: a Bayesian optimization algorithm using hash tables to generate the transfer matrix. 
    We use the  {HESBO}\footnote{https://github.com/aminnayebi/HesBO} Python toolbox with the default parameters.
\end{itemize}
For  {\gls{EGORSE}}, the version showing the best performance in term of convergence speed and robustness is selected, \textit{i.e.} EGORSE PLS + Gaussian with an initial DoE of $d$ points~\cite{EGORSE}.
To achieve this comparison, $10$ optimizations for each problem and  for each studied method are completed to analyze the statistical behavior of these BO algorithms.  A comparison between  {\gls{EGORSE}} and the four studied algorithms is performed in term of robustness, convergence speed both in CPU time and in number of iterations.
Figure~\ref{fig:ALL} provides the iteration convergence plots and the time convergence plots. 
Figure~\ref{fig:ALL:mb10_nb} shows that   {TuRBO} and   {EGO-KPLS} are converging the fastest and with a low standard deviation.
Moreover, the convergence plots of   {EGORSE},   {RREMBO} and   {HESBO} are hardly distinguishable. 
Figure~\ref{fig:ALL:mb100_nb} displays that   {EGO-KPLS} converges the fastest to the lowest values with a low standard deviation.
  {TuRBO} is also providing good performance even if it converges slower than   {EGO-KPLS}. 
\noindent Regarding the three methods using dimension reduction procedure,   {EGORSE} is converging to the lowest value with a relatively low standard deviation.
The good performance of   {EGO-KPLS} and   {TuRBO} is certainly due to the ability of these algorithms to search all over $\Omega$, which is not the case for other methods.
However, when the dimension of $\Omega$ increases, the ability to search all over $\Omega$ becomes a drawback. In fact, a complete search in $\Omega$ is intractable in time in this case.
Figures~\ref{fig:ALL:mb10_time} and~\ref{fig:ALL:mb100_time} depict the convergence CPU time necessary to obtain the regarded value.
First, the RREMBO, TuRBO and EGO-KPLS complete the optimization procedure in more than 8 hours on the MB\_100 problem against an hour on the MB\_10 problem.
This suggests that   {RREMBO},   {TuRBO} and   {EGO-KPLS} are intractable in time for larger problems.
Then, one can easily see that  {EGORSE} is converging the fastest in CPU time compared to the other algorithms on the MB\_100 problem.
In fact, the computation time needed to find the enrichment point is much lower than the one for   {TuRBO},   {EGO-KPLS} and   {RREMBO} which was sought in the definition of the enrichment sub-problem.
Finally,   {EGORSE} is converging to a lower value than   {HESBO} in a similar amount of time.
Thus,   {EGORSE} seems more interesting to solve high dimension problems than the studied algorithms.
Note that only   {HESBO} and   {EGORSE} are able to perform an optimization procedure on high dimension problems.
    
To sum up, a comparison of  {EGORSE} with other algorithms has been carried out and  it pointed out that   {TuRBO},   {EGO-KPLS} and   {RREMBO} are intractable in time for high dimension problems.  Furthermore,  {EGORSE} has appeared to be the most suitable to solve these high dimension problems efficiently.

\subsection{Mixed integer validation results}

For 10 different test cases (3 constrained and 7 unconstrained) we are considering 20 runs with different DoE sampled by \gls{LHS} for a total of 200 instances. These 10 test cases are  “Branin 1”, “Branin 2”,“Branin 3”, “Branin 4”, "Branin 5”, “Goldstein 1”, “Goldstein 2”, “Toy set 1", “Toy set 2” and "Cosine curves"~\cite{Pelamatti,Pelamatti2020,Roustant,CAT-EGO,AMIEGO,Gower,Vanaret,vanaret2024interval}. In order to compare our methods with the state-of-the-art, the following optimization algorithms are used: Bandit-BO, \gls{NSGA2}, SEGO with Kriging, SEGO with Gower distance and SEGO coupling Kriging and PLS.

The Bandit-BO implementation used is the one by Nguyen \textit{et al.}~\cite{Bandit-BO}, with neither parallelization nor batch evaluations.  
The \gls{NSGA2}~\cite{nsga2} algorithm used is the implementation from the toolbox pymoo~\cite{pymoo} with the default parameters (probability of crossover = 1, eta = 3). Fronts are not relevant in our study as we are considering single-objective optimization in this section.
The optimization with SEGO is made from \gls{SEGOMOE}~\cite{bartoli:hal-02149236} for both constrained and unconstrained  cases.
For SEGO using Gower distance~\cite{Gower} (denoted by SEGO+GD), we are considering the implementation of the Surrogate Modeling Toolbox (SMT)~\cite{SMT2019}, the open-source Python toolbox in which the computations associated to the present work have been done.
For SEGO using Kriging (denoted by SEGO+KRG), we also use the implementation from the toolbox SMT.  For the constrained analytical test cases, we are using the $UTB$ criterion~\cite{SEGO-UTB} in \gls{SEGOMOE}.
Some adaptions have been done to Bandit-BO and \gls{NSGA2} to consider both integer and categorical variables. As \gls{NSGA2} can only consider integer variables, categorical variables are treated as integer ones. 
Contrarily, Bandit-BO can treat only categorical variables so integer variables are treated as categorical ones.
Bandit-BO creates a GP model for each arm, so it requires at least $2\times N_c$ initial points, $N_c$ being the number of categorical possibilities for the problem inputs. For unconstrained test cases, these $2$ points by categorical possibility are sampled randomly. 
For constrained optimization, as we cannot use Bandit-BO, we use a continuous relaxed \gls{LHS}, then project the sampled points to obtain the mixed integer \gls{DOE}.

\begin{figure}[H]
 \subfloat[Data profiles for a tolerance of 2\%. \label{2DP}]{\includegraphics[width=0.50\textwidth]{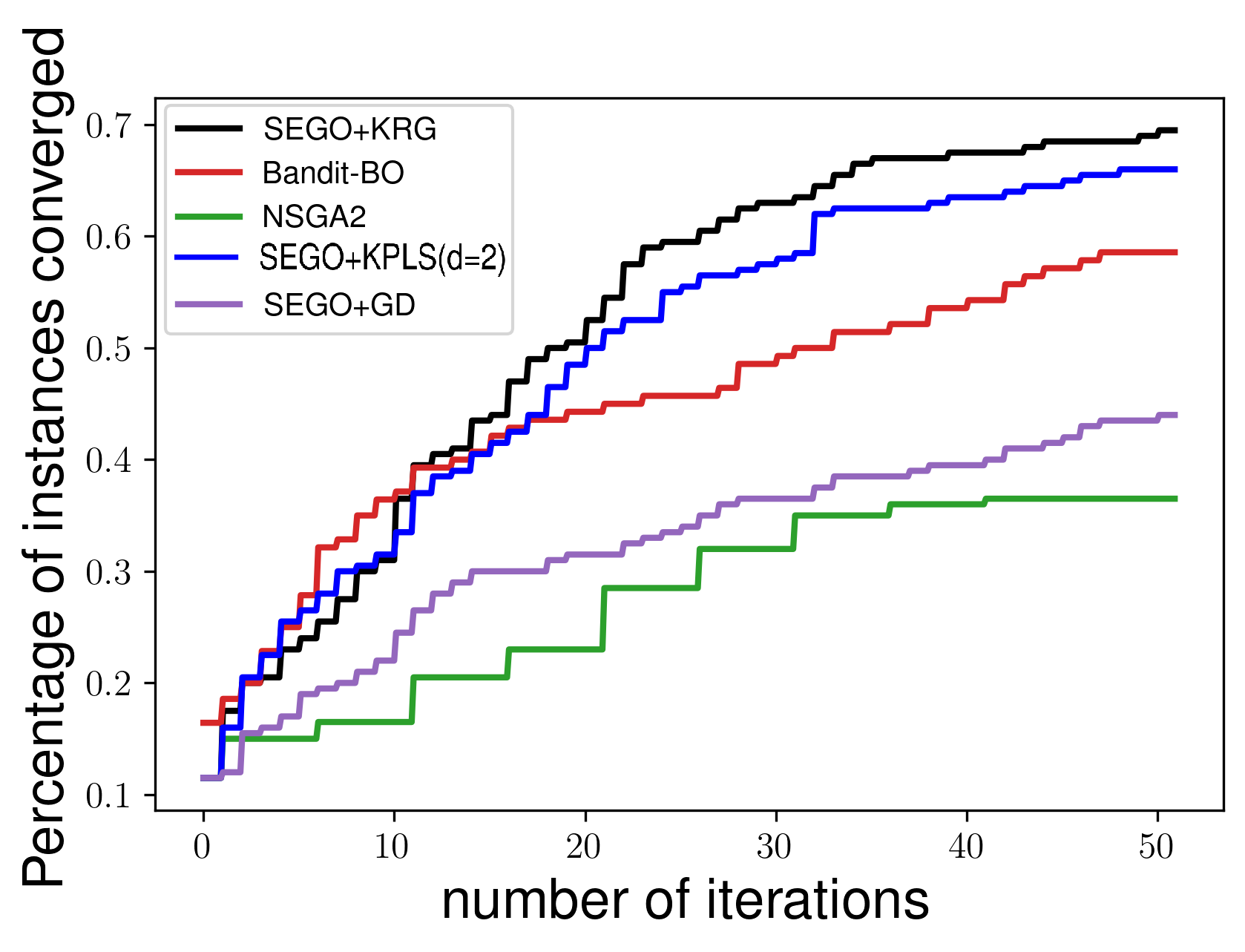}}
 \subfloat[Data profiles for a tolerance of 0.5\%. \label{0.5DP}]{\includegraphics[width=0.50\textwidth]{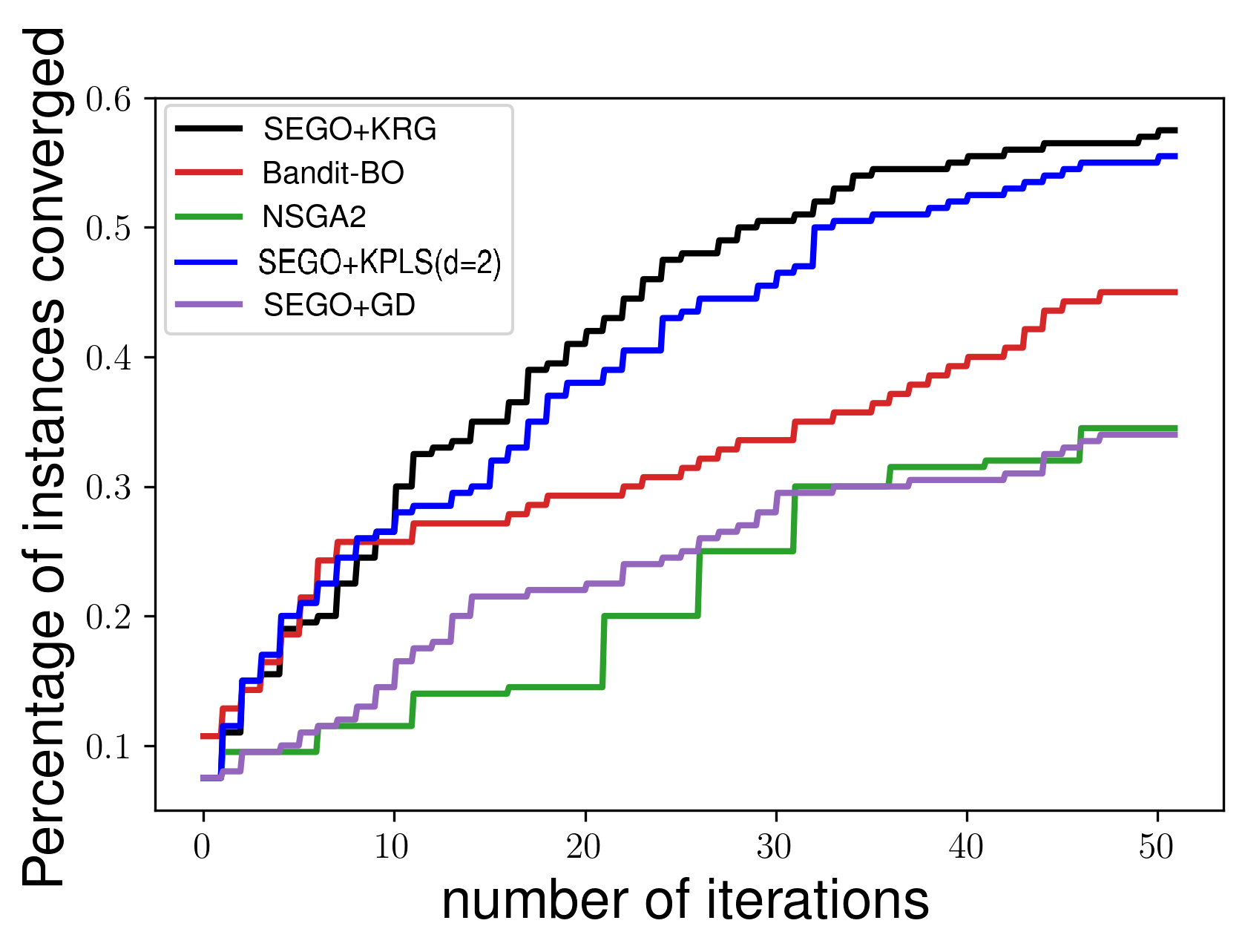}}
   \caption{Data profiles generated using 10 analytical test cases.}
     \label{DP_benchmark}
\end{figure}

For Kriging, Kriging with PLS or Gower distance, the hyperparameters are optimized with COBYLA~\cite{COBYLA} and the chosen model regression is constant. When optimizing with SEGO, the acquisition function is maximized using ISRES (Improved Stochastic Ranking Evolutionary Strategy~\cite{ISRESimp}) to find some interesting starting points and SNOPT (Sparse Nonlinear OPTimizer~\cite{gill2005snopt}) to finalize the process based on these starting points.
The squared exponential kernel is the only kernel considered for the \gls{BO} experiments.
In order to compute some statistical data (median and variance), we are doing 20 repetitions of the optimization process for a given method and an initial DoE size. We consider that a constraint is respected if the constraint violation is smaller than the threshold value $10^{-4}$.   The resulting percentages of instances that have converged after a given budget for every method are plotted on~\figref{DP_benchmark}. For unconstrained test cases, in order to compare with Bandit-BO, the size of the initial DoE is given by $\min \{5, 2\times N_c \}$ where $N_c$ is the number of categorical possibilities and for constrained test cases, we took 5 points for the initial DoE size. This allow us to compare Bandit-BO, SEGO, SEGO with Gower Distance, SEGO with KPLS and \gls{NSGA2}. Some tests are in dimension 2, so in order to compare with the same number of hyperparameters for all tests, we had to choose between 1 or 2 principal components for KPLS. As the number of points increases, the projected points are really closed to each other, so, to insure a certain stability, we choose KPLS with 2 principal components and noise evaluation, denoted by KPLS(d=2) in the following.
For constrained case, we keep only the inputs that give a constraint violation smaller than  $10^{-4}$.
A problem is considered solved if the error to the known solution is smaller than 2\% on~\figref{2DP} and smaller than 0.5\% on~\figref{0.5DP}. The mean error after 50 iterations can be found in~\tabref{tab:meanerror}. SEGO with PLS gives the smaller errors on constrained test cases   but, for unconstrained ones, SEGO-KRG performs the best. As we can see on the data profiles, SEGO-KRG and SEGO+KPLS(d=2) are similar and outperform the three other methods. However, SEGO with PLS is less efficient on unconstrained test cases and so SEGO-KRG gives slightly better results overall.

\begin{table}[H]
   \caption{Mean errors of each method after 50 iterations.}
   \begin{center}
   \resizebox{\columnwidth}{!}{%
      \begin{tabular}{*{6}{c}}
       \hline
       \textbf{ERRORS} & NSGA-II & Bandit-BO & SEGO+KRG & SEGO+KPLS(d=2) &SEGO+GD  \\
      \hline
     \textbf{Unconstrained test cases} & 14.29\% & 5.84\% &  2.27\% & 5.74\% &  8.41\%  \\
     \textbf{Constrained test cases} &   58.27\% & - &27.61\% & 25.41\% &  47.83\% \\
        \hline
      \end{tabular}
      }
   \end{center}
   \label{tab:meanerror}
\end{table}

\subsection{Multi-objective validation results}

The chosen test suite is constituted of three unconstrained and three constrained problems.  
First, the group of unconstrained problems is composed  of the three problems known as ZDT problems~\cite{ZDT}. Those problems are bi-objective with $d$ design variables. 
In this section, we are testing $d=2$, $d=5$ and $d=10$, leading to a total of $9$ unconstrained problems. The second group of test cases is composed of three constrained multi-objective problems. The first one is the Binh and Korn  problem (BNH)~\cite{binh1997mobes} with 2 objectives,  2 variables and 2 inequality constraints. The second problem is the Tanaka one (TNK)~\cite{zapotecas2018review} with 2 objectives,  2 design variables and 2 inequality  constraints. The third constrained problem is the Osyczka and Kundu problem (OSY)~\cite{osyczka1995new}  with 2 objectives,  6 design variables and 6 inequality constraints. 

In the context of multi-objective optimization, the efficiency of the tested methods will be evaluated based on the following two criteria:
\begin{itemize}
    \item the proximity between the obtained not dominated points and the associated explicit Pareto front;
    \item the distribution of the obtained points in the objective space to cover as largely as possible the Pareto-optimal possibilities.
\end{itemize}
Different performance indicators exist in the literature. In this context, the property of Pareto compliance is of high interest. Namely, a weakly Pareto compliant indicator $I$ is defined such that for two sets of points $A$ and $B$, if the elements of $A$ dominate those of $B$ then $I(A) \leq I(B)$. One efficient way to estimate the compliance indicator is given by inverted generational distance plus (${IGD}^{+}$)~\cite{ishibuchi2015modified}. The ${IGD}^{+}$ indicator is defined as follows, for a given set of points $A=\left \{a_1,a_2,\ldots,a_{|A|}\right \}$ and a reference set $Z=\left \{z_1,z_2,\ldots,z_{|Z|}\right \}$ with true values from the optimal Pareto front, the ${IGD}^{+}$ value of the set $A$ is
\begin{equation}
{IGD}^{+}(A) =  \; \frac{1}{|Z|} \; \bigg( \sum_{i=1}^{|Z|} {d_i^{+}}(z_i)^2 \bigg)^{1/2} \;\mbox{ with } \; d_i^{+}(z_i) = \max \{ \min_{a_i \in A}a_i - z_i, 0\} 
\end{equation}
where the notation $|Z|$ defines the cardinal of the set of points $Z$. In our comparison tests, the smaller is the value of ${IGD}^{+}$ the better is the tested method.
We note that evaluating the ${IGD}^{+}$ indicator requires the knowledge of the true optimal Pareto front to build the reference set $Z$. For that reason, we tested analytical multi-objective problems, for which we know explicitly the optimal Pareto fronts.
Due to stochastic effects caused by the starting points at each enrichment step and by the genetic algorithm, the ${IGD}^{+}$ convergence plots are averaged over 10 runs.

\begin{figure}[H]
\subfloat[]{
\includegraphics[width=0.33\textwidth]{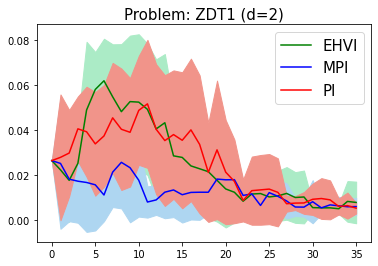}
\includegraphics[width=0.33\textwidth]{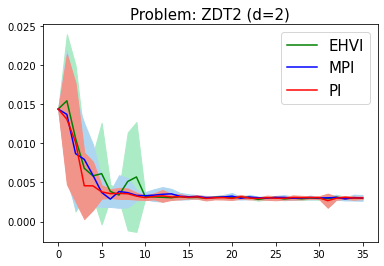}
\includegraphics[width=0.33\textwidth]{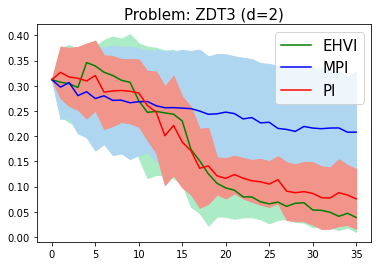}
}

\vspace{-1cm}
\subfloat[]{
\includegraphics[width=0.33\textwidth]{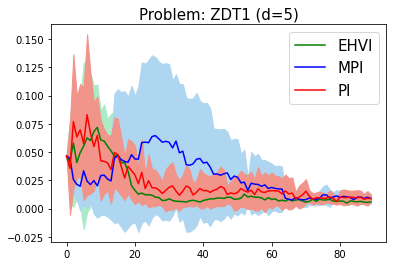}
\includegraphics[width=0.33\textwidth]{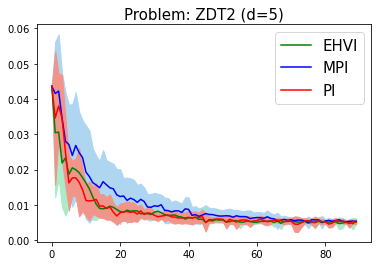}
\includegraphics[width=0.33\textwidth]{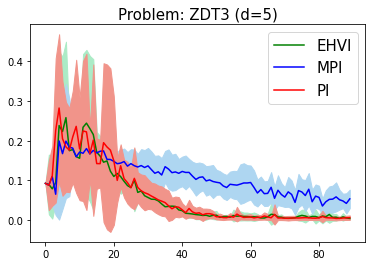}
}

\captionsetup[subfigure]{labelformat=customa}
\vspace{-1cm}
\subfloat[Unconstrained problems (from left to right: ZDT1, ZDT2 and ZDT3).\label{fig:uncons-IGD}]{
\includegraphics[width=0.33\textwidth]{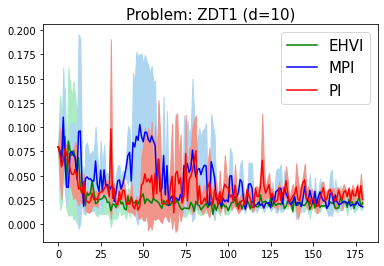}
\includegraphics[width=0.33\textwidth]{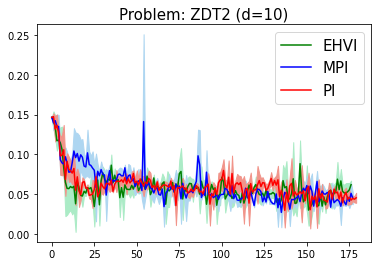}
\includegraphics[width=0.33\textwidth]{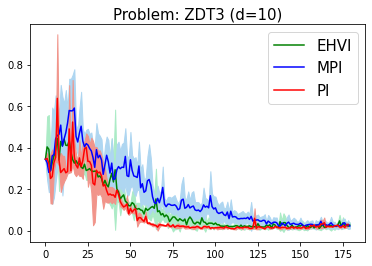}
}

\captionsetup[subfigure]{labelformat=customb}
\subfloat[Constrained problems (from left to right: BNH, TNK and OSY).\label{fig:cons-IGD}]{
\includegraphics[width=0.33\textwidth]{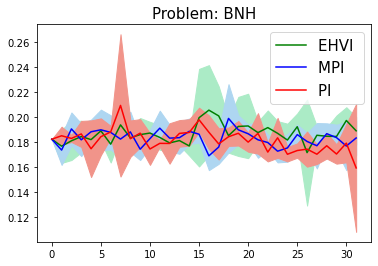}
\includegraphics[width=0.33\textwidth]{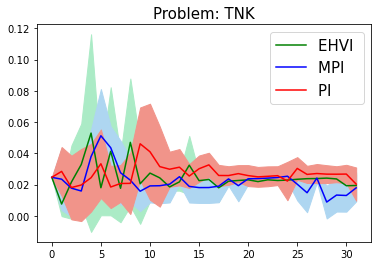}
\includegraphics[width=0.33\textwidth]{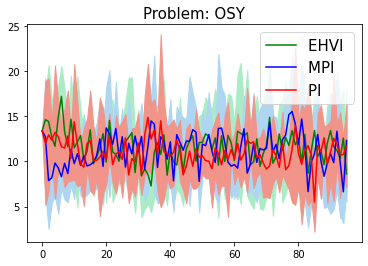}
}

\caption{Obtained convergence plots (\textit{i.e.}, ${IGD}^{+}$ values across iterations): a comparison of the acquisition functions $EHVI$, $PI$ and $MPI$ within SEGOMOE using ${IGD}^{+}$.}
\label{fig:IGD}
\end{figure}

\begin{figure}[H]

\subfloat[]{
\includegraphics[width=0.33\textwidth,clip=True]{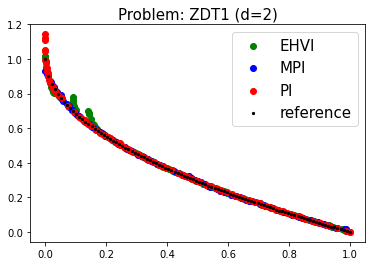}
\includegraphics[width=0.33\textwidth,clip=True]{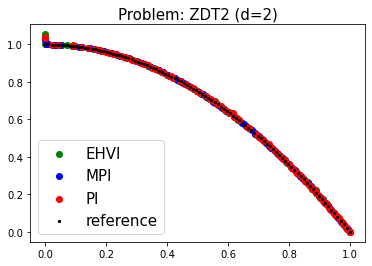}
\includegraphics[width=0.33\textwidth,clip=True]{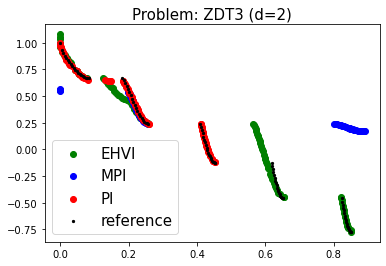}
}

\vspace{-1cm}
\subfloat[]{
\includegraphics[width=0.33\textwidth,clip=True]{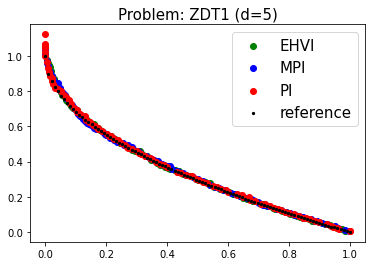}
\includegraphics[width=0.33\textwidth]{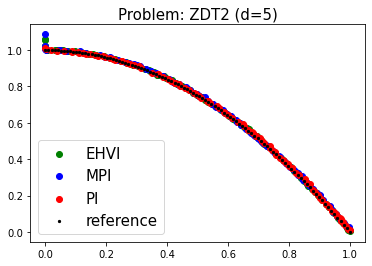}
\includegraphics[width=0.33\textwidth]{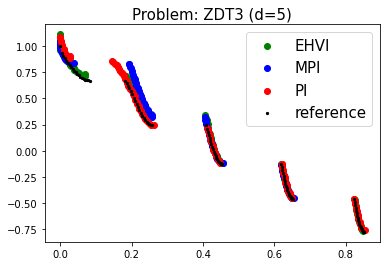}
}

\captionsetup[subfigure]{labelformat=customa}
\vspace{-1cm}
\subfloat[Unconstrained problems (from left to right: ZDT1, ZDT2 and ZDT3).\label{fig:uncons-pf}]{
\includegraphics[width=0.33\textwidth]{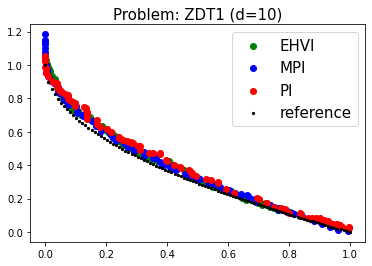}
\includegraphics[width=0.33\textwidth]{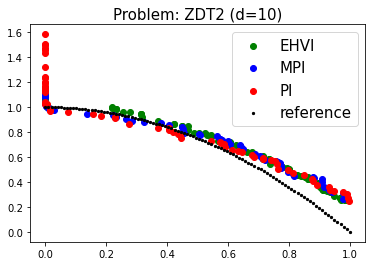}
\includegraphics[width=0.33\textwidth]{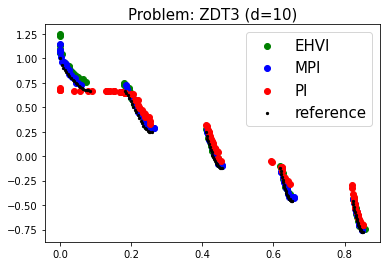}
}
\captionsetup[subfigure]{labelformat=customb}
\subfloat[Constrained problems (from left to right: BNH, TNK and OSY).\label{fig:cons-pf}]{
\includegraphics[width=0.33\textwidth]{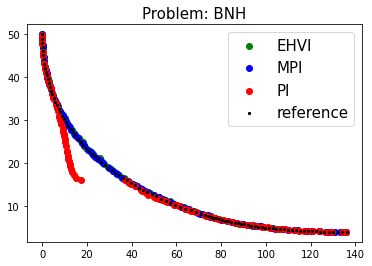} 
\includegraphics[width=0.33\textwidth]{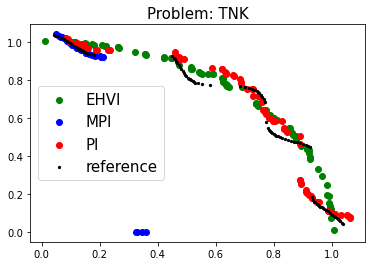}
\includegraphics[width=0.33\textwidth]{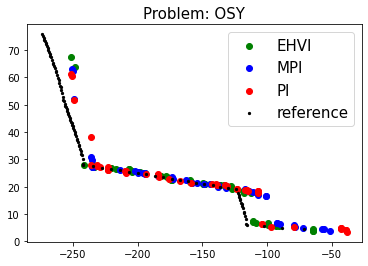}
}
  
\caption{A comparison of the obtained Pareto fronts obtained using $20d$ points.}  \label{fig:pf}
\end{figure}

The GP surrogate models of the objectives and constraints are computed using the SMT toolbox~\cite{SMT2019}. If no initial known points of the problem are provided, a \gls{LHS}~\cite{LHS} sampling method is used to draw a space-filling sample. 
Let $d$ be the dimension of the design space and $c$ the number of constraints, the sampling will be made of $2d+ 2c+1$ points and $20d$ evaluations of the problem will be computed in total, i.e $20d - (2d+2c+1)$ iterations lead to the presented results. 
Figure~\ref{fig:IGD} illustrates the convergence of the ${IGD}^{+}$ indicator through the iterations for the tested problems. For each problem, different acquisition criteria are compared: $EHVI$, $PI$ and $MPI$ whose mean value and associated dispersion are displayed. The associated Pareto fronts obtained at the final iteration are given in Fig.~\ref{fig:pf} together with the true Pareto front.

\captionsetup[subfigure]{labelformat=empty}
\captionsetup[subfigure]{labelformat=parens}

\captionsetup[subfigure]{labelformat=empty}
\captionsetup[subfigure]{labelformat=parens}
For unconstrained problems, see Fig.~\ref{fig:uncons-IGD}, one can see that \gls{SEGOMOE} performs well for the three tested acquisition functions, with comparable convergence plots for ${IGD}^{+}$. $EHVI$ is slightly better compared to $PI$ and $MPI$ for the unconstrained test cases. For the three constrained problems (BNH, TNK and OSY), convergence plots are given in Fig.~\ref{fig:cons-IGD}. 
One can see that the problems TNK ans OSY are very challenging for \gls{SEGOMOE} compare to the BNH problem. Similarly to the unconstrained test cases, we obtain comparable performance for the acquisition functions $EHVI$, $PI$ and $MPI$. 
Last, we note that for all the tested problems and for every infill criterion, \gls{SEGOMOE} converged to a good approximation of the explicit Pareto front. The dispersion over the ${IGD}^{+}$ scores becomes lower over the iterations, showing the consistency of each method. 
Even if the curse of dimensionality is a phenomenon leading to an exponential increase of the research space, our Bayesian algortihm do not lose much in efficiency but the chosen budget of iterations is proportional to the dimension of the problems. We can mention than $PI$ and $EHVI$ infill criteria are most of the time the ones giving the best results. The estimation of the $EHVI$ acquisition function is computationally expensive  compared to $PI$ which is the easiest to compute among the tested infill criteria.

These preliminary results are promising and some more realistic applications are considered in the next section.

\section{Engineering optimization problems}
\label{sec:eng_optim}
 
We applied our Bayesian optimization algorithm to varied engineering problems of high interest as decribed in~\tabref{tab:recap_eng}. In particular, the "Rover trajectory", "Cantilever beam", "Neural network" and "CEntral Reference Aircraft System (CeRAS)" optimization problems are known and tackled in various works related to this thesis for illustration purposes. 
More precisely, the "Mixed cantilever beam" and "Hierarchical neural network" problems were not optimized and have been used to illustrate the functioning of different GP models in Chapter~\ref{c2}, Chapter~\ref{c3} and Chapter~\ref{c4}. The rover trajectory problem was optimized to display the ability of \gls{EGORSE}~\cite{EGORSE} to solve really high dimension problems and the CeRAS test case was meant to illustrate a novel optimization method that relies on an adaptive strategy for choosing automatically the number of reduced dimension when using GP with continuous relaxation and PLS for mixed integer problems~\cite{SciTech_cat,PLS_Li_2002}. Also, another variant of the CeRAS test case was used to validate a new infill criterion for multi-objective Bayesian optimization~\cite{grapin_constrained_2022}.
%
Notwithstanding, the \gls{DRAGON} aircraft concept, the "aircraft upgrade" retrofitting design problem, the optimization of a "family of business jets" and the "supply chain and manufacturing optimization for producing aircraft" are meant for engineering purposes and direct applications in the context of the \gls{EU} project AGILE 4.0.  
Being real-life problems, these problems have numerous mixed integer variables together with several objectives and constraints.
\begin{table}[H]
\centering
\caption{Definition of the engineering optimization problems.}
\resizebox{\linewidth}{!}{%
\small

\begin{tabular}{ccccc}
\hline
\textbf{Name} & \# of objs &  \# of cons  & \# of vars [\textbf{cont}, \textbf{int}, \textbf{cat}] & Reference \\
\hline
\textbf{Rover\_600} & 1 & --& [600,0,0]  & ~\cite{EGORSE} \\
\textbf{Mixed cantilever beam} & 1 & --& [2,0,1(12 levels)]  & ~\cite{Mixed_Paul} \\
\textbf{Neural network} & 1& -- & [1,3,1(2 levels)] & ~\cite{saves2023smt} \\

\texttt{CeRAS} & 1& 2 & [6,2,2(2,2 levels)] & ~\cite{saves2021constrained}\\
\textbf{Bi-objective} \texttt{CeRAS} & 2 & 3 & [5,0,0] & ~\cite{grapin_constrained_2022} \\
\texttt{DRAGON} & 1& 5 & [10,0,2(17,2 levels) & ~\cite{SciTech_cat} \\
\textbf{Airframe upgrade design }& 4 &4  & [3,0,1(4 levels)] & ~\cite{bartoli2023Agile} \\
\textbf{Family of aircraft} & 2 &2  & [9,0,10(2 levels each)] & ~\cite{bartoli2023Agile} \\
\textbf{Production of aircraft} & 5 &2  & [0,0,8(21,21,21,21,6,5,4,5 levels)] & ~\cite{bartoli2023Agile} \\
\hline
\end{tabular}
}

\label{tab:recap_eng}
\end{table}
%
%
%
%
%
Based on~\tabref{tab:recap_eng}, the \gls{DRAGON} test case and the "Family of aircraft" problems have been chosen to illustrate our Bayesian optimization algorithm in what follows. 
In fact, the \gls{DRAGON} test case is a mono-objective test case but features five constraints and both continuous and categorical variables. We detail the optimization results for this particular aircraft in Section~\ref{sec:dragonc5}.
The "Family of aircraft" optimization problem is also representative of our optimization abilities because it requires the optimization of five objectives with over 100 millions possibilities
and under two constraints. The optimization results for this problem are detailed in Section~\ref{sec:familyc5}.
Concerning the other test cases presented in~\tabref{tab:recap_eng}, more optimization results are presented in the corresponding reference papers and detailed descriptions of these problems are given in Appendix~\ref{app:eng_cases}.

\subsection{"\texttt{DRAGON}" configuration optimization problem}
\label{sec:dragonc5}

This section details the analysis and optimization related to the \gls{DRAGON} aircraft concept. 
In particular, the optimization has been done using a novel algorithm based on \gls{KPLS} with a varying number of components automatically chosen along the optimization. Moreover, the analysis has been done using a software named \gls{FAST-OAD}.

\paragraph{The automatic PLS procedure}
\gls{KPLS} is an efficient method that can tackle high-dimensional problems by reducing the number of effective dimensions to a small number. This method is often used in industry and research~\cite{bhosekar2018advances,zuhal2021dimensionality,li2021improved}. The number of principal components being a key point for PLS, we propose a strategy to choose this number in an adaptive way during the iterations of the optimization process. As a rule of the thumb, taking 2 to 5 active components is, in general, efficient for most problems. However, when we do not have any $\textit{a-priori}$ information, an intuitive idea is to learn the number of active dimensions directly from the GP surrogates.
Following this idea, \citet{Bouhlel18} proposed to use the leave-one-out strategy to find the best number of hyperparameters for KPLS.
However, this strategy is not efficient because, at every iteration, it implies to  compute a large number of reduced size surrogate models to select the best one. As the size of the DoE increases, both the model computation and the leave-one-out criterion become more and more costly. 
In~\cite{PLS_Li_2002}, it has been proposed to use the adjusted Wold’s R criterion~\cite{WoldPLS} for dimensionality reduction purposes. Wold’s R criterion is based on cross-validation over $k$-folds \cite{PLS_Li_2002} and consists in measuring  the ratio of the PRedicted Error Sum of Squares (PRESS)~\cite{PRESS}. For a given number of components $d$, the $R$ criterion is defined as 
\begin{equation}\label{eqPRESS}
R(d)=\frac{\text{PRESS}(d+1)}{\text{PRESS}(d)}
\end{equation}
and should be smaller than a given threshold. This treshold is typically 1, which means we add components until the approximate error stops decreasing. However, a threshold of 0.9 or 0.95 could be used in order to be more selective and add fewer components. The construction of the KPLS models with an automatically chosen number of components is given in~\algref{alg:woldR}.

\smallskip
\begin{algorithm}[H]
\SetAlgoLined
{\textbf{Inputs:}} A DoE of size $n_t$, i.e, $\mathscr{D}= \{(x_1,y_1), \ldots, (x_{n_t},y_{n_t}) \} $ associated with a given function $\Psi$ ($\Psi(x_i)=y_i, \forall i=1,\ldots, n_t$). The maximal and minimal number of components $d_{\max}$ and $d_{\min}$, respectively. A threshold $\sigma  \in [0,1]$. A chosen number of folds $K$. $\forall k \in \{1, \ldots, K \}$ denote the corresponding known $n_k$ points: $\{(x_1,y_1), \ldots, (x_{n_k},y_{n_k}) \}$.

\textbf{Initialization:} 
\begin{itemize}
 \item For all subsamples $k \in \{1, \ldots, K \}$, build a KPLS model with $d_{\min} $ components over all data except the inputs in the subsample $k$. 
    \item For all $ i \in 1, \ldots, n_k$, let $\hat{y}_{i,-k}$ be the prediction of the KPLS model build without the points in the $k^{th}$ fold at the corresponding $x_i$. This fold is therefore used as a validation set to compute the $K$-fold cross validation PRESS as $ \text{PRESS}(d_{\min}) = {\displaystyle \sum_{k=1}^{K}{\displaystyle \sum_{i=1}^{n_k} \left(  y_i - \hat{y}_{i,-k}     \right)^2  }} $.
    \item Set $d_{\Psi} \xleftarrow{} d_{\min} $ .
\end{itemize}

\While{$d_{\Psi} \ < \  d_{\max}$}{
\vspace{.2cm}
\begin{adjustwidth}{0pt}{40pt}
\begin{enumerate}
    \item Divide  the DoE $\mathscr{D}$ into $K$ subsamples. 
  
    \item For all subsamples $k \in \{1, \ldots, K \}$, build a KPLS model with $d_{\Psi}$ components over all data except the inputs in the subsample $k$. 
    \item For all $ i \in 1, \ldots, n_k$, let $\hat{y}_{i,-k}$ be the prediction of the KPLS model build without the points in the $k^{th}$ fold at the corresponding $x_i$. This fold is therefore used as a validation set to compute the $K$-fold cross validation PRESS as $ \text{PRESS}(d_{\Psi}+1) = {\displaystyle \sum_{k=1}^{K}{\displaystyle \sum_{i=1}^{n_k} \left(  y_i - \hat{y}_{i,-k}     \right)^2  }} $. Set $R(d_\Psi)=\frac{\text{PRESS}(d_{\Psi}+1)}{\text{PRESS}(d_{\Psi})}$.
    \item 
     \uIf{$R(d_\Psi) \geq \sigma $}{
   \quad  \textbf{STOP};
   }
   \item Increment $d_{\Psi}$, \textit{i.e.}, $d_{\Psi} \xleftarrow{} d_{\Psi}+1$  .
  
 \end{enumerate}

\end{adjustwidth}
\smallskip
} 
\textbf{Outputs:} KPLS model with $ d_{\Psi}$ components;
   
\caption{
Adaptive dimension reduction for KPLS models.}
\label{alg:woldR}
\end{algorithm}

\smallskip
Note that when $d_{\min} =d_{\max}$,~\algref{alg:woldR} is equivalent to KPLS with a constant number of dimensions. The code implementation of our proposed method has been released in the toolbox SMT\footnote{\url{https://smt.readthedocs.io/en/latest/}}~\cite{saves2023smt}.

\paragraph{FAST-OAD}
\gls{FAST-OAD} represents a significant advance in the field of aircraft design and optimization, with particular relevance to the design and preliminary design phases.
Developed in collaboration between ONERA and ISAE-SUPAERO, \gls{FAST-OAD} is a sophisticated software package designed to meet the multiple challenges relative to aircraft design. Its main objective is to facilitate the creation of optimized aircraft configurations that align with the \gls{TLAR} within the sizing possibilities. One of the attributes of \gls{FAST-OAD} is its inherent multidisciplinary approach, enabling it to accommodate a wide range of aircraft configurations and architectures. This approach reflects the intrinsic multidisciplinary of modern aircraft design, covering multiple aspects such as aerodynamics, propulsion, structures and so on. By doing so, \gls{FAST-OAD} overcomes the limitations of traditional, discipline-centric design tools by providing an integrated platform for comprehensive design exploration. \gls{FAST-OAD} core functionality is underpinned by a point mass approach, which plays a central role to estimate essential aircraft design parameters such as fuel consumption or energy requirements. This approach enables rapid sizing and optimization, with execution times measured in seconds and minutes respectively, making it highly effective for aeronautical engineering applications~\cite{David_2021}. 
The latest open-source version of the \gls{FAST-OAD} software is available on Github\footnote{\href{https://github.com/fast-aircraft-design/FAST-OAD}{https://github.com/fast-aircraft-design/FAST-OAD}}.

\medskip

\begin{figure}[H]
 \subfloat[]{
   \centering
\centerline{
   \includegraphics[   height=4cm]{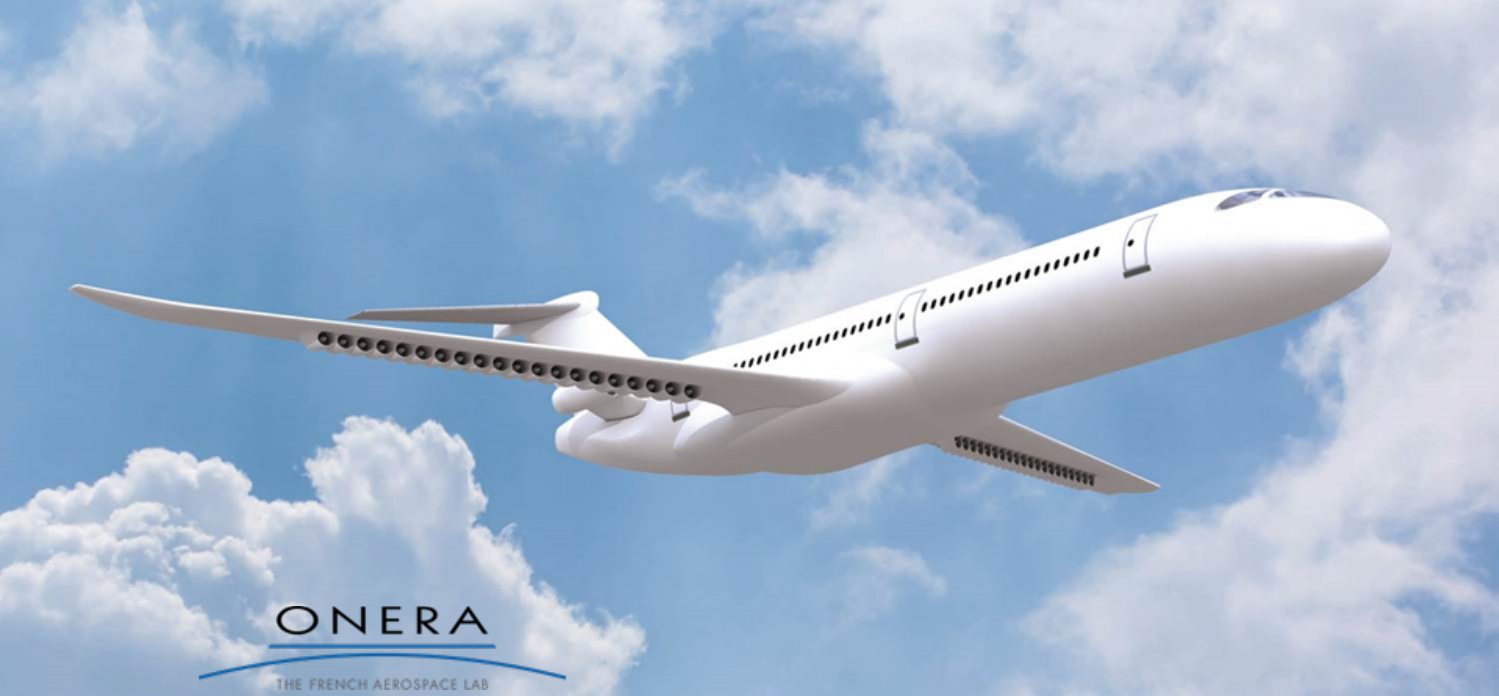}
}
}

\vspace{-0.8cm}
 \subfloat[]{
 \includegraphics[   height=5cm, width = 0.5\textwidth]{images/onera-dragon-21321165_wagner.jpg}
 \includegraphics[   height=5cm, width = 0.5\textwidth]{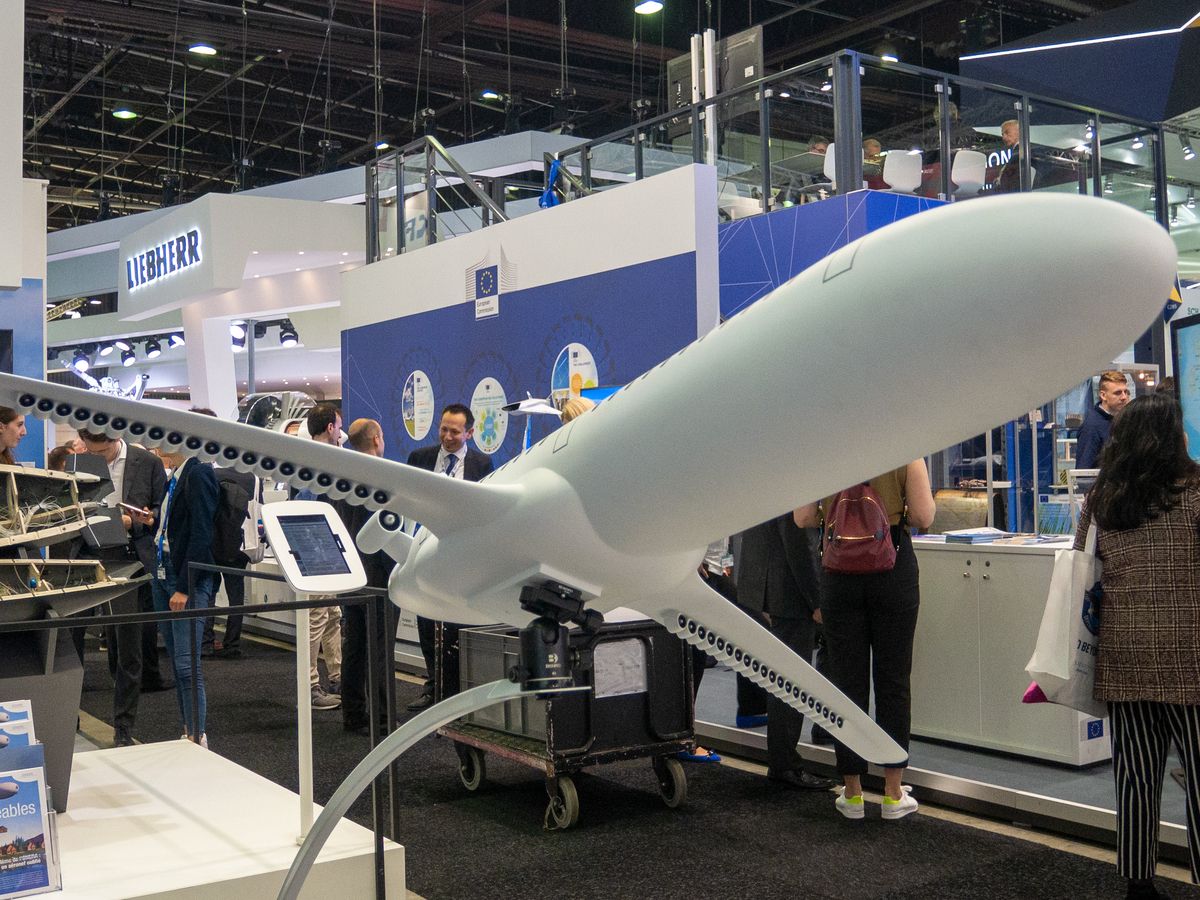}
 }
   \caption{``\texttt{DRAGON}'' aircraft mock-up.}
   \label{Dragon2020c6}
\end{figure}

The \gls{DRAGON} aircraft concept in~\figref{Dragon2020c6} has been introduced by ONERA in 2019~\cite{schmollgruber} within the scope of the European CleanSky 2 program\footnote{\href{https://www.cleansky.eu/technology-evaluator}{https://www.cleansky.eu/technology-evaluator}} 
which sets the objective of 30\% reduction of CO2 emissions by 2035 with respect to 2014 state-of-the-art. A first publication in SciTech 2019~\cite{schmollgruber} was followed by an up-to-date report in SciTech 2020~\cite{schmollgruber2}.
ONERA responded to this objective by proposing a distributed electric propulsion aircraft that improves the aircraft fuel consumption essentially by increasing the propulsive efficiency. Such efficiency increase is obtained through improvement of the bypass ratio by distributing a large number of small electric fans on the pressure side on the wing rather than having large diameter turbofans.
This design choice avoids the problems associated with large under-wing turbofans and at the same time allows the aircraft to travel at transonic speed. Thus the design mission set for \gls{DRAGON} is 150 passengers over 2750 Nautical Miles at Mach 0.78.

\captionsetup[subfigure]{labelformat=empty}

\captionsetup[subfigure]{labelformat=parens}

\begin{figure}[!htb]
\vspace*{-0.6cm}
  \centering 
\includegraphics[   ,height=6cm]{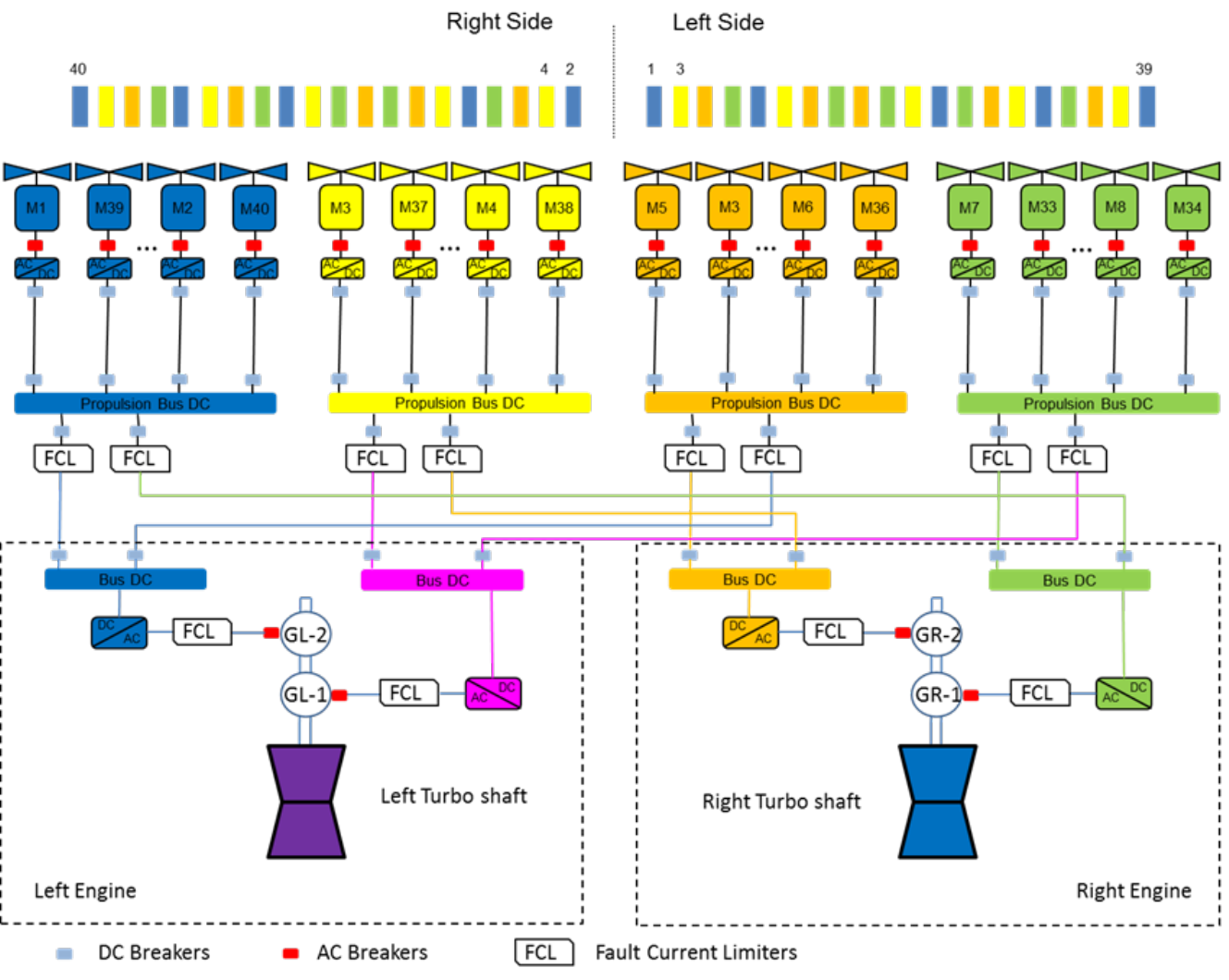}
 \caption{Turboelectric propulsive architecture.}
 \label{DragonArchitecturec6}
\end{figure}

The employment of a distributed propulsion comes at a certain cost; a turbo-electric propulsive chain is necessary to power the electric fans which brings additional complexity and weight. Typically, turboshafts coupled to electric generators are generating the electrical power on board of the aircraft. The power is then carried to the electric fans through an electric architecture sized to ensure robustness to single component failure. This safety feature is obtained with redundant components as depicted in~\figref{DragonArchitecturec6}.
The baseline configuration is two turboshafts, four generators, four propulsion buses with cross-feed and forty fans. This configuration was selected for the initial study as it satisfies the safety criterion. However it was not designed to optimize aircraft weight. The turboelectric propulsive chain being an important weight penalty, it is of particular interest to optimize the chain and particularly the number and type of each component, characterized by  some  discrete or particular values.

In Aerobest~\cite{saves2021constrained}, we solved a constrained optimization problem with 8 continuous design variables and 4 integer ones, for a total of 12 design variables. For both \gls{DRAGON} and “CeRAS”, we found that the best configuration with SEGO and Kriging with and without PLS regression were almost the same while PLS reduced the number of variables from 12 to 4. These first results were promising but no categorical variable was considered.
Here, the motivation is to add some of these variables in order to increase the number of considered possibilities.
Moreover, the \gls{MDA} fixed-point algorithm is now replaced to be more flexible and general as the \gls{DRAGON} test case is treated with Overall Aircraft Design method in \gls{FAST-OAD}~\cite{David_2021}.

In~\cite{saves2021constrained}, the observation of the optimization results from an aircraft design point of view led to surprising conclusions, ultimately resulting in questioning the validity of the models utilized for the design. Among the problems raised in this previous work, the sizing rules of the propulsive chain, based on the default architecture, proved to be insufficiently general to explore more diversified and valid architectures. The results indicated that an architecture comprising 40 motors connected to 6 cores of equal power would be optimal, although it is violating the design rule stating that each motor should be connected to a single core.

In consequence, we propose an up-to-date optimization problem for \gls{DRAGON} with the most recent enhancements.
The consideration of the variables related to architecture was revised to take full advantage of mixed variables optimization. Three configurations with different number of electric components were considered, each with their own sizing rules. The default configuration is conserved for the analysis and two new configurations, one with low distribution, the other with high distribution, were created and analysed to establish the sizing rules. The choice of architecture is hence represented by a categorical variable.
The number of motors was to remain an optimization variable as it is an important driver of the propulsive and aerodynamic efficiency of the aircraft but it is constrained by the type of architecture. In this analysis, the number of motors could only be multiple of 4, 8 or 12 depending of the selected architecture. The solution consisted in expending the levels of the categorical variable representing the architecture and assign a specific and valid number of motors to each level.
Additionally, it was noticed that the layout of \gls{DRAGON}, with the turbogenerators at the rear of the fuselage was disadvantageous for two reasons:
\begin{itemize}
    \item The added weight at the rear forces the wing to be more aft, reducing the level arm with the tail surfaces. To maintain the static stability, the tail surfaces have to be increased, adding weight and friction drag penalty.
    \item The electric cables running from the rear generators to the center part of the wing are heavy, accounting for $\sim$10\% of the total propulsion system weight.
\end{itemize}
For these two reasons, it would be advantageous to locate the turbogenerator below the wing. However, this would restrict the space available for the electric fans, which would in turn limit the maximum propulsive efficiency. To make the trade-off between lighter propulsion and better propulsive efficiency the layout was included into the optimization problem as a categorical variable.
Five constraints need to be accounted for when sizing the aircraft.
The firsts three constraints consist of maximum limits for the take off field length, climb duration and top of climb slope angle that are strong drivers for the sizing of the hybrid electric propulsion. The two last constraints consist of a portion of the wing trailing edge at the wingtip that should be left free for the ailerons (hence limiting the space available for the electric fans) and, finally, a wingspan limit is imposed by airport regulation.

To know how optimizing the fuel mass will impact the aircraft design, we are considering the new optimization problem described in~\tabref{tab:dragonc6}. We can now solve a constrained optimization problem with 10 continuous design variables and 2 categorical variables with 17 and 2 levels respectively, for a total of 12 design variables. For the optimization, this new problem is a hard test case involving 29 relaxed variables and 5 constraints. The definition of the architecture variable is given in~\tabref{tab:dragon_archi1c6} and the definition of the turboshaft layout is given in~\tabref{tab:dragon_archi2c6}. 

\begin{table}[H]
\centering
 \caption{Definition of the ``\texttt{DRAGON}'' optimization problem.}
\small

\resizebox{\columnwidth}{!}{%
\small

\begin{tabular}{lllrr}
 & Function/variable & Nature & Quantity & Range\\
\hline
\hline
Minimize & Fuel mass & cont & 1 &\\
 & \multicolumn{2}{l}{\bf Total objectives} & {\bf 1} & \\
\hline
with respect to & \mbox{Fan operating pressure ratio} & cont & 1 & $\left[1.05, 1.3\right]$ \\  
     & \mbox{Wing aspect ratio} & cont & 1 &    $\left[8, 12\right]$ \\
    & \mbox{Angle for swept wing} & cont & 1 & $\left[15, 40\right]$  ($^\circ$) \\
     & \mbox{Wing taper ratio} & cont & 1 &    $\left[0.2, 0.5\right]$ \\
     & \mbox{HT aspect ratio} & cont & 1 &    $\left[3, 6\right]$ \\
    & \mbox{Angle for swept HT} & cont & 1 & $\left[20, 40\right]$  ($^\circ$) \\
     & \mbox{HT taper ratio} & cont & 1 &    $\left[0.3, 0.5\right]$ \\
 & \mbox{TOFL for sizing}  & cont &1 & $\left[1800., 2500.\right]$ ($m$) \\
 & \mbox{Top of climb vertical speed for sizing} & cont & 1 & $\left[300., 800.\right]$($ft/min$) \\
 & \mbox{Start of climb slope angle} & cont & 1 & $\left[0.075., 0.15.\right]$($rad$) \\

 & \multicolumn{2}{l}{Total  continuous variables} & 10 & \\
 \cline{2-5}
& \mbox{Architecture} & cat & 17 levels & \{1,2,3, \ldots,15,16,17\} \\
& \mbox{Turboshaft layout} & cat & 2 levels & \{1,2\} \\

 & \multicolumn{2}{l}{Total categorical variables} & 2 & \\
 \cline{2-5}

  &   \multicolumn{2}{l}{\textbf{Total relaxed variables}} & {\textbf{29}} & \\
  \hline
  
subject to & Wing span \textless  \ 36   ($m$)  & cont & 1 \\
 & TOFL \textless  \ 2200 ($m$) & cont & 1 \\
 & Wing trailing edge occupied by fans  \textless  \ 14.4 ($m$) & cont & 1 \\
 & Climb duration \textless  \ 1740 ($s $) & cont & 1 \\
 & Top of climb slope \textgreater \ 0.0108 ($rad$) & cont & 1 \\

 & \multicolumn{2}{l}{\textbf{Total  constraints}} & {\textbf{5}} & \\
\hline
\end{tabular}
}
\label{tab:dragonc6}
\end{table}

\begin{table}[H]
\centering
\vspace*{-0.3cm}

 \caption{Definition of the architecture variable and its 17 associated levels.}
\small

\resizebox{0.85\columnwidth}{!}{%
\small

\begin{tabular}{cccc}
\hline
  \textbf{Architecture number} & number of motors & number of cores & number of generators\\
  \hline
  \textbf{1} & 8 &2 & 2\\
  \textbf{2} & 12 & 2 & 2\\
  \textbf{3} & 16 & 2 & 2\\
  \textbf{4} &20 &2 & 2\\
  \textbf{5} & 24 & 2 & 2\\
  \textbf{6} & 28 & 2 & 2\\
  \textbf{7} &32 & 2 & 2\\
  \textbf{8} & 36  & 2 & 2\\
  \textbf{9} & 40 & 2 & 2\\
 
  \textbf{10} & 8  &4 & 4\\
  \textbf{11} & 16 &4 & 4\\
  \textbf{12} & 24 &4 & 4\\
  \textbf{13} & 32 &4 & 4\\
  \textbf{14} & 40 &4 & 4\\

  \textbf{15} & 12 &6 & 6\\
  \textbf{16} & 24 &6 & 6\\
  \textbf{17} & 36 &6 & 6\\

\hline
\end{tabular}
}
\label{tab:dragon_archi1c6}
\end{table}

\begin{figure}[H]
\subfloat[Convergence curves.\label{DRAGONcurves}]{
\includegraphics[   height=5cm, width=0.52\textwidth]{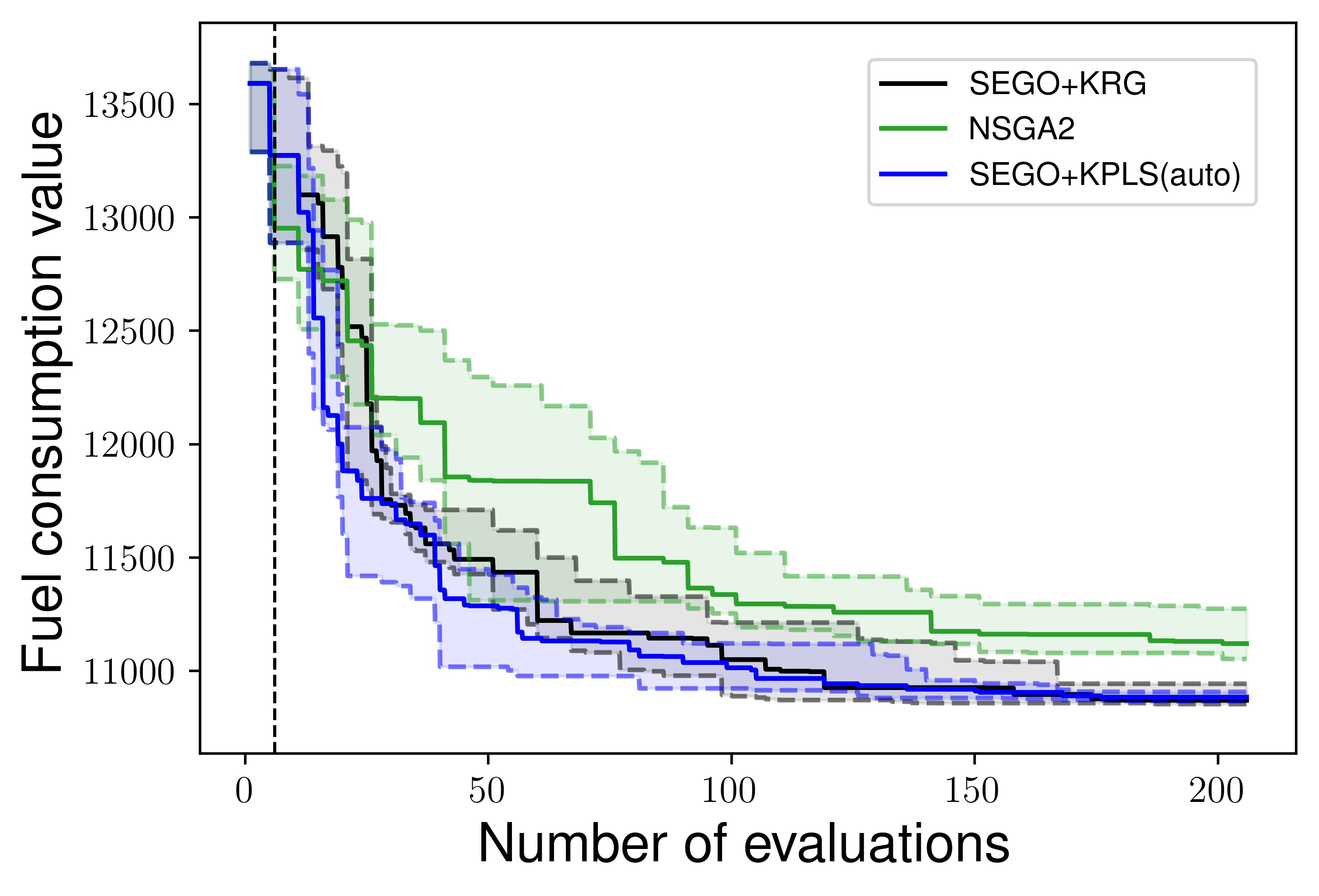}
}
\subfloat[Boxplots. \label{minima_DRAGON}]{
\includegraphics[   height=5cm, width=0.48\textwidth]{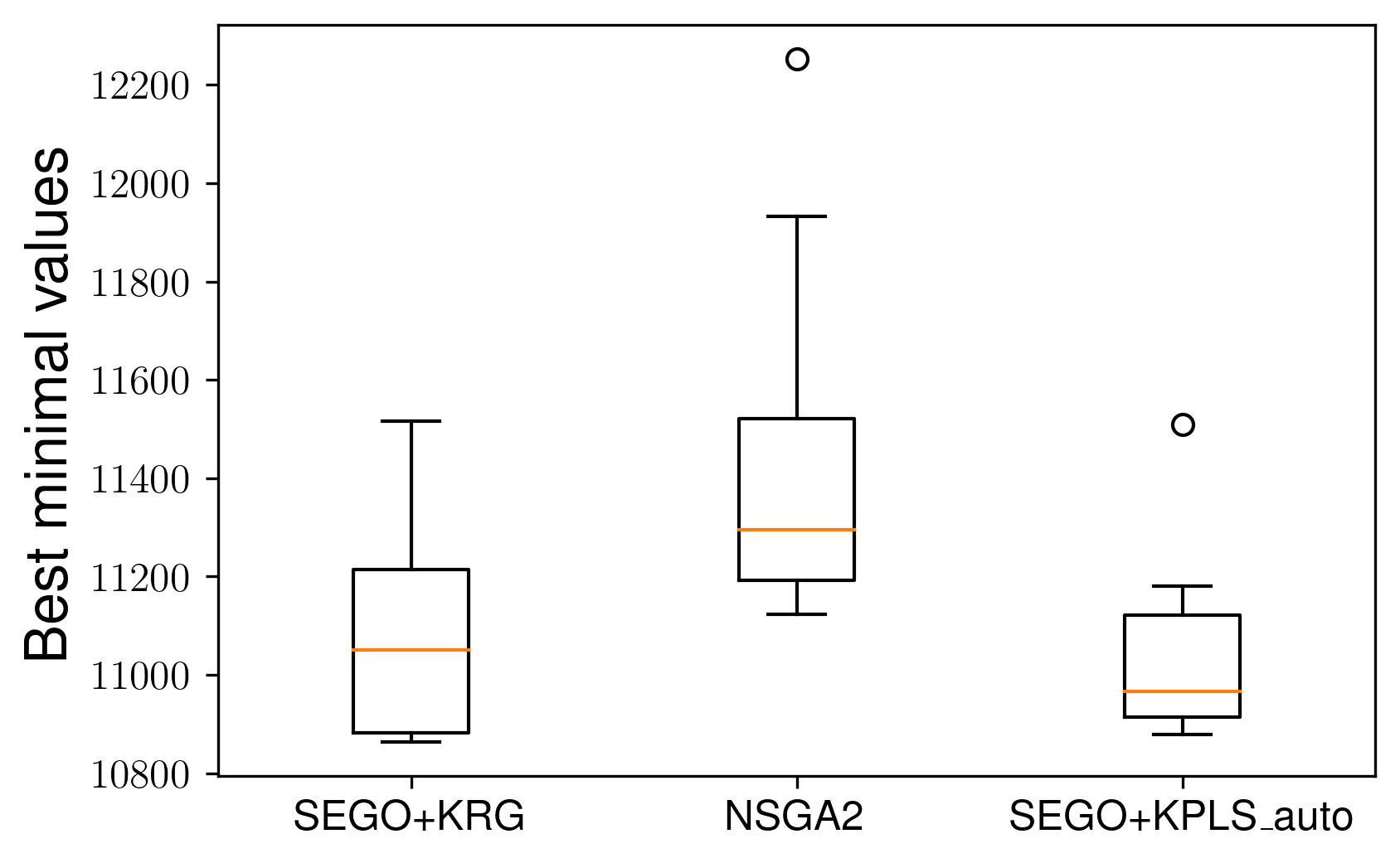}
}
\caption{``\texttt{DRAGON}'' optimization results using a DoE of 5 points over 10 runs. The Boxplots are generated, after 100 iterations, using the 10 best points.}
\label{DRAGON}
\end{figure}

\begin{figure}[H]
   
      \centering
	\includegraphics[height=5.5cm,width=12cm]{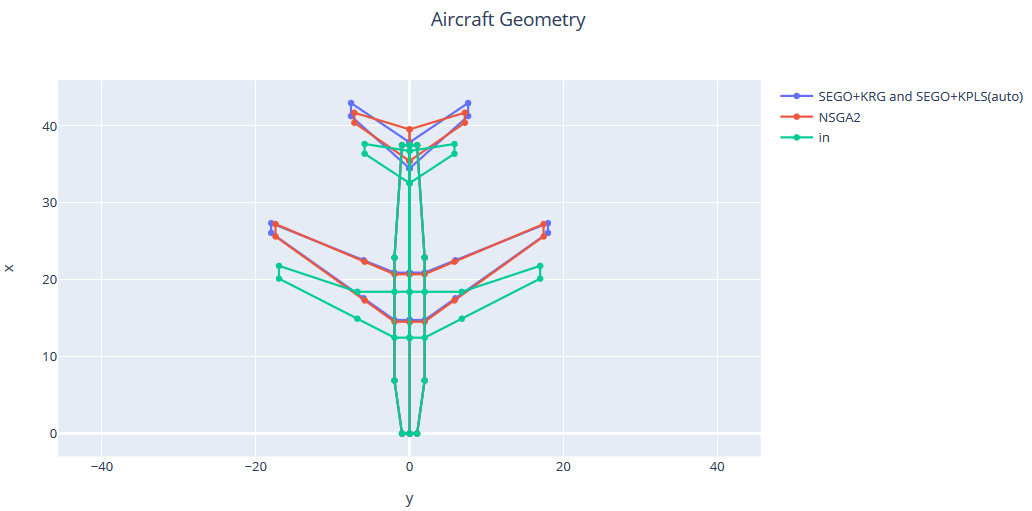}
	\caption{ “\texttt{DRAGON}'' best configuration geometry. Comparisons between the initial configuration, the NSGA2, the SEGO+KRG and the SEGO+KPLS(auto) results.}
     \label{DRAGON_geo}

\end{figure}

\begin{table}[H]
\centering
\vspace*{-0.3cm}

 \caption{Definition of the turboshaft layout variable and its 2 associated levels.}
\small

\resizebox{0.85\columnwidth}{!}{%
\small

\begin{tabular}{cccccc}
\hline
  \textbf{Layout} & position & y ratio & tail & VT aspect ratio & VT taper ratio\\
  \hline 
  \textbf{1} & under wing &0.25 & without T-tail& 1.8 & 0.3 \\
  \textbf{2} & behind & 0.34 & with T-tail& 1.2 & 0.85\\
 
\hline
\end{tabular}
}
\label{tab:dragon_archi2c6}
\end{table}
For the core \gls{MDO} application, we apply \gls{FAST-OAD} on \gls{DRAGON}, an innovative aircraft configuration under development.  
As the evaluations are expensive, we are doing only 10 runs instead of 20. For each run, we draw a random starting DoE of 5 points. 
Since we have no equality constraint and a small budget, we chose to tighten the treatment of constraints by considering only the average prediction of the models, instead of using the $UTB$ criterion which would broaden the treatment of constraints~\cite{SEGO-UTB}.
Also, to have realistic results, the constraints violation will be forced to be less than $10^{-3}$. On the following, let MAC denote the Mean Average Chord, VT be the Vertical Tail, HT be the Horizontal Tail and TOFL be the Take-off Field Length. 
To validate our method, we are comparing two variants of SEGO with automatic KPLS or with its more expensive version using Kriging, and NSGA2 in~\figref{DRAGON} where 10 runs are performed in order to plot the mean and the associated quartiles. The best method after 100 iterations is the proposed one involving automatic \gls{PLS} regression as shown in the boxplots in~\figref{minima_DRAGON}. In~\figref{DRAGONcurves}, after 200 iterations, we still find that SEGO is better than NSGA2 and that the PLS technique helps for the convergence.
We find that the best configurations found with SEGO and Kriging with and without PLS regression are almost the same. The description of this best configuration is given in~\tabref{tab:dragon_best} and the best geometries for NSGA2, SEGO+KRG and SEGO+KPLS(auto) are given in~\figref{DRAGON_geo}. SEGO+KRG and SEGO+KPLS(auto) are grouped as they give similar results.

\begin{table}[H]
   \caption{``\texttt{DRAGON}'' Optimal aircraft configuration.}
   \begin{center}
   \resizebox{0.7\columnwidth}{!}{%
    \small
      \begin{tabular}{*{3}{c}}
       \hline
        \textbf{Name} & Nature & Value \\
      \hline
      \textbf{Fuel mass} & cont &  10816 kg  \\
      \textbf{Wing span}  &  cont & 36 m \\
      \textbf{TOFL}  &  cont & 1722.7 m           \\
      \textbf{Wing trailing edge occupied by fan} &  cont & 10.65 m \\
      \textbf{Climb duration} & cont &  1735.3 s \\ 
       \textbf{Top of climb slope}  &  cont & 0.0108 rad \\
        \hdashline
     \textbf{\mbox{Fan operating pressure ratio}} & cont & 1.09\\
      \textbf{\mbox{Wing aspect ratio}} & cont & 10.9\\
      \textbf{\mbox{Angle for swept wing}} & cont & 32.2$^\circ$ \\
      \textbf{\mbox{Wing taper ratio}} & cont & 0.235\\
      \textbf{\mbox{HT aspect ratio}} & cont & 6\\
      \textbf{\mbox{Angle for swept HT}} & cont & 40$^\circ$ \\
      \textbf{\mbox{HT taper ratio}} & cont & 0.3 \\
      \textbf{\mbox{TOFL for sizing}} & cont & 1803 m \\
      \textbf{\mbox{Top of climb vertical speed for sizing}} & cont & 494 ft/min \\ 
      \textbf{\mbox{Start of climb slope angle}} & cont & 0.104 rad \\ 
      \textbf{\mbox{Architecture}} & cat & 10 \\
      \textbf{\mbox{Turboshaft layout}} & cat & 2 \\

      \hline
      \end{tabular}
    }
   \end{center}
   \label{tab:dragon_best}
\end{table}

The optimal configuration was found with Kriging for 10816 kg of fuel and KPLS found a similar point at 10819 kg. The best configuration found is the configuration 10, with the smaller number of 8 motors but with 4 cores and electric generators. 
On the point of view of the aircraft design, the number of motors was restricted to the minimum available and the two architectures allowing 8 motors (1 and 10) are performing almost identically, showing a low influence of this variable over the number of motors. In this case, it is the climb constraints that are sizing for the whole propulsion chain. Therefore the architecture sizing laws are not at play and no difference can be seen.

Despite an important space still available at the wing trailing edge, a small number of motors is looked for by the optimizer to lower the wetted surface area of the fans. This is probably oversimplified as the integration of the fans at the trailing edge is coarsely modeled. The fact that this direction is selected during the optimization tells us that this model should be updated to further improve the design.
The most favorable layout is found to be with the turbo-generators at the rear. The level arm between the wing and the horizontal tail being actually larger due to the maximum sweep angle employed for the horizontal tail. 
In particular, from a structural point of view, the combination of high sweep and high aspect ratio is too badly taken into account leading to unrealistic weight for the horizontal stabiliser, which probably favours too much this layout. Nevertheless, the trade-off found by optimization is correct given the models used in \gls{FAST-OAD}.

\subsection{Family of aircraft optimization problem}
\label{sec:familyc5}

This application case is part of the upgrade-driven stream and aims at designing a family of three business jet aircraft considering commonality options in the MDO workflow~\cite{bussemaker2022ac7}. The commonality options enable sharing one or more aircraft components between the aircraft in the family: by sharing components, design, certification, production, and maintenance costs are reduced. However, operating costs might be increased due to the use of components not designed for the typical flight conditions.

\begin{figure}[!htb]
      \centering
      \includegraphics[width=\linewidth]{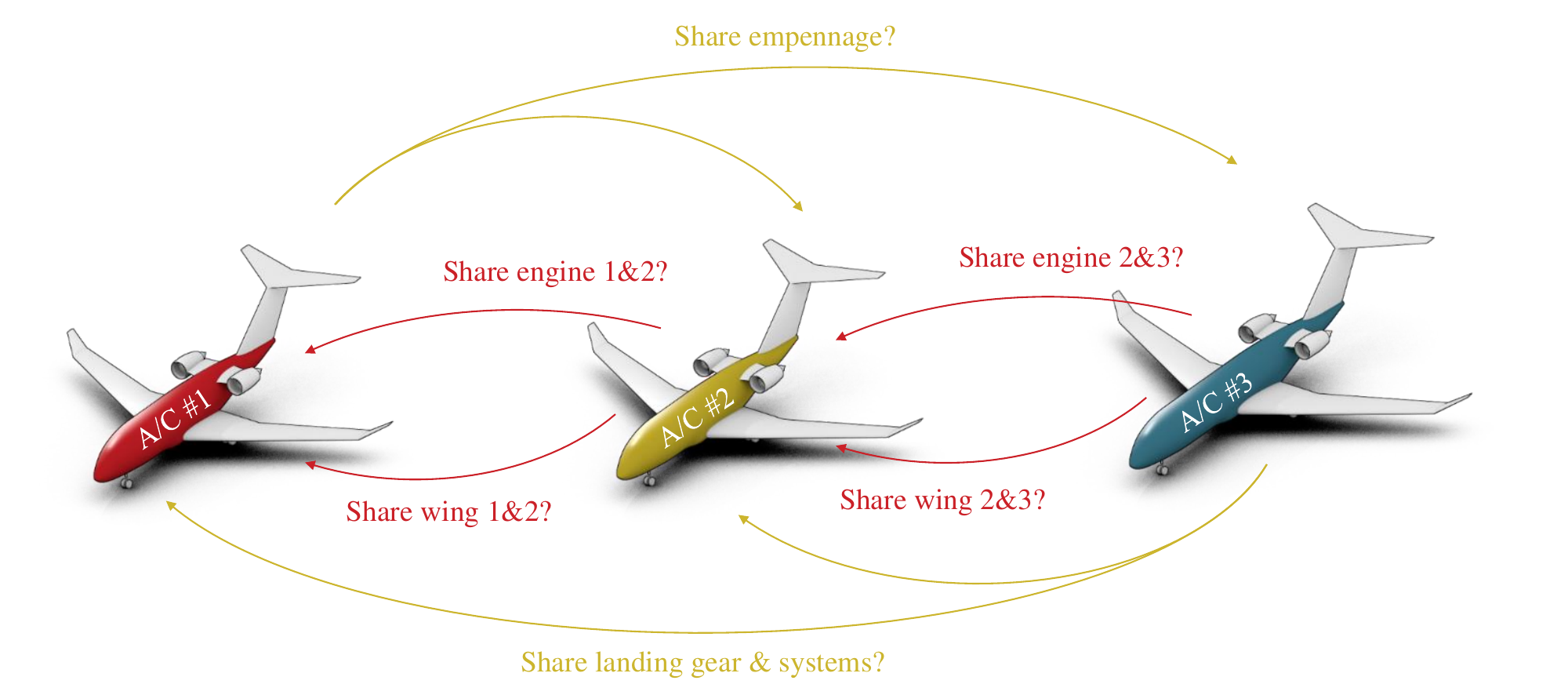}
   \caption{Visualization of the aircraft family including commonality sharing decisions.}
   \label{fig:ac7_sharing}
\end{figure}

Shareable components include wings, empennage (horizontal and vertical tail), engines, on-board systems (OBS), and landing gears, shown in Figure~\ref{fig:ac7_sharing}. Next to the commonality choices there are three design variables per wing: leading edge sweep, rear span location (determining flap size), and thickness-to-chord ratio. These three design variables per wing are inactive if the associated wing is a shared wing. Each aircraft family is optimized for two objectives: direct operating costs (DOC), representing the impact on performance, and Original Equipment Manufacturer (OEM) non-recurring costs, representing the potential benefit of commonality. Table~\ref{tab:pb_ac7} presents the optimization problem in more details.

\begin{table}[H]
    \centering
   \caption{Definition of the aircraft family design problem.}
   \small
      \begin{tabular}{lllrr}
 & Function/variable & Nature & Quantity & Range\\
\hline
\hline
Minimize & Direct Operating Costs & cont & 1 &\\
Minimize &  OEM Non-Recurring Costs& cont & 1 &\\
 & \multicolumn{2}{l}{\bf Total objectives} & {\bf 2} & \\
\hline
with respect to   &  Leading edge sweep* & cont & 3 & $\left[ 30.0, 42.0 \right]$ ($^\circ$)\\
 & Rear span location* & cont & 3 & $\left[ 0.72, 0.82 \right]$ ($\% chord$)\\
 & Wing thickness/cord ratio*  & cont & 3 & $\left[0.06, 0.11\right]$ \\
 & \multicolumn{3}{l}{* only active if the corresponding wing is \textbf{not} shared} \\
 & \multicolumn{2}{l}{Total continuous variables} & { 9} & \\ \cline{2-5}
& Engine commonality  1$\And$2 & cat  &  2 levels & sharing yes/no\\
& Engine commonality 2$\And$3 & cat  &  2 levels & sharing yes/no \\
& Wing commonality 1$\And$2 & cat  &  2 levels & sharing yes/no\\
& Wing commonality 2$\And$3 & cat  &  2 levels & sharing yes/no\\
& Landing gear commonality 1$\And$2 & cat  &  2 levels & sharing yes/no\\
& Landing gear commonality 2$\And$3 & cat  &  2 levels & sharing yes/no\\
& OBS commonality 1$\And$2 & cat  &  2 levels & sharing yes/no\\
& OBS commonality 2$\And$3 & cat  &  2 levels & sharing yes/no\\
& Empennage commonality 1$\And$3 & cat  &  2 levels & sharing yes/no\\
& Empennage commonality 3$\And$2 & cat  &  2 levels & sharing yes/no\\
& \multicolumn{2}{l}{ Total categorical variables} & { 10} & \\ \cline{2-5}
 & \multicolumn{2}{l}{\bf Total relaxed variables} & {\bf 29} & \\
\hline
subject to & 
\multicolumn{2}{l}{ Balanced  Field Length $\le 1524$ ($m$) } & 1 & \\
&\multicolumn{2}{l}{ Landing Field Length $ \le 762$ ($m$) }  & 1 & \\
 & \multicolumn{2}{l}{\bf Total constraints} & {\bf 2} & \\
\hline
\end{tabular}
   \label{tab:pb_ac7}
\end{table}

\begin{figure}[!htb]
  \centering
  \includegraphics[width=\textwidth]{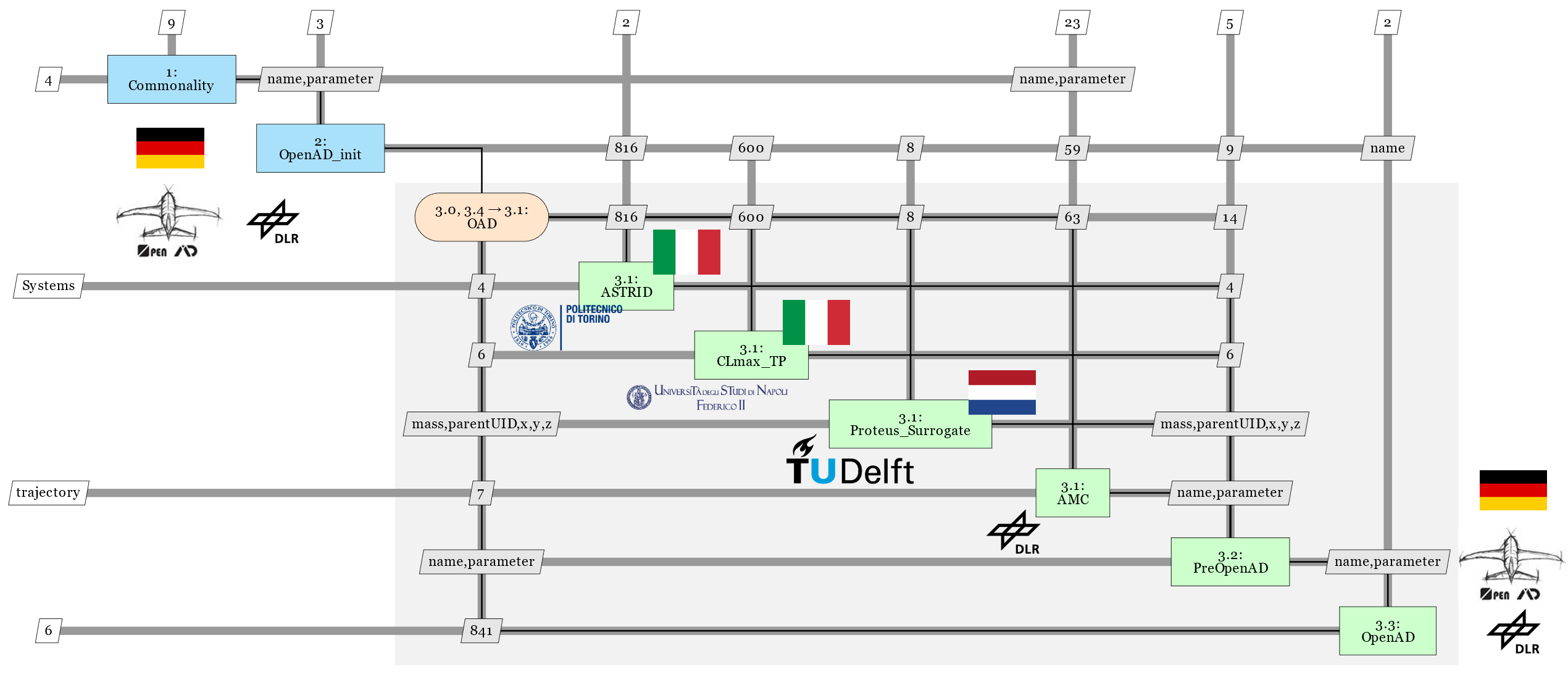}
  \caption{MDA for family aircraft design.}
  \label{fig:MDA_AC7}
\end{figure}

The design space is modeled using the Architecture Design Space Graph (ADSG) method~\cite{Architecture} implemented in the ADORE tool~\cite{bussemaker2022adore} developed during the AGILE 4.0 project\footnote{\url{https://www.agile4.eu/}}. MultiLinQ~\cite{agile4is2023}, also developed for AGILE 4.0, is used to couple the generate architecture instances to the central data schema used in the MDO workflow~\cite{bussemaker2022ac7}. The MDO workflow consists of two levels: an aircraft-level workflow (shown in Figure~\ref{fig:MDA_AC7}) applying commonality decisions and sizing one aircraft at a time, and a family-level workflow converging the three aircraft-level workflows and performing cost calculation on the family-level. For a more detailed overview of the design space model, coupling to the MDO workflow and implementation of the MDO workflow for the family design application case, the reader is referred to~\cite{bussemaker2022ac7,agile4is2023}.

For this application case, \gls{SEGOMOE} is used as optimization algorithm due to the need to handle hierarchical, mixed-discrete design variables, and the need to minimize the number of function evaluations as one family evaluation can take up to 2 hours. Hierarchical design variables were handled using the imputation method, where inactive variables are replaced by the mean value to prevent duplicate design vectors~\cite{Effectiveness}. \gls{SEGOMOE} was accessed through an ask-tell API implemented in WhatsOpt\footnote{\url{https://github.com/whatsopt/WhatsOpt-Doc}} running on a server at ONERA's premises as described in~\citet{lafage2019whatsopt}. For a more detailed description of how \gls{SEGOMOE} was coupled to ADORE and the MDO workflow, the reader is referred to~\cite{agile4is2023}.

\begin{figure}[!htb]
      \centering
      \includegraphics[width=\linewidth]{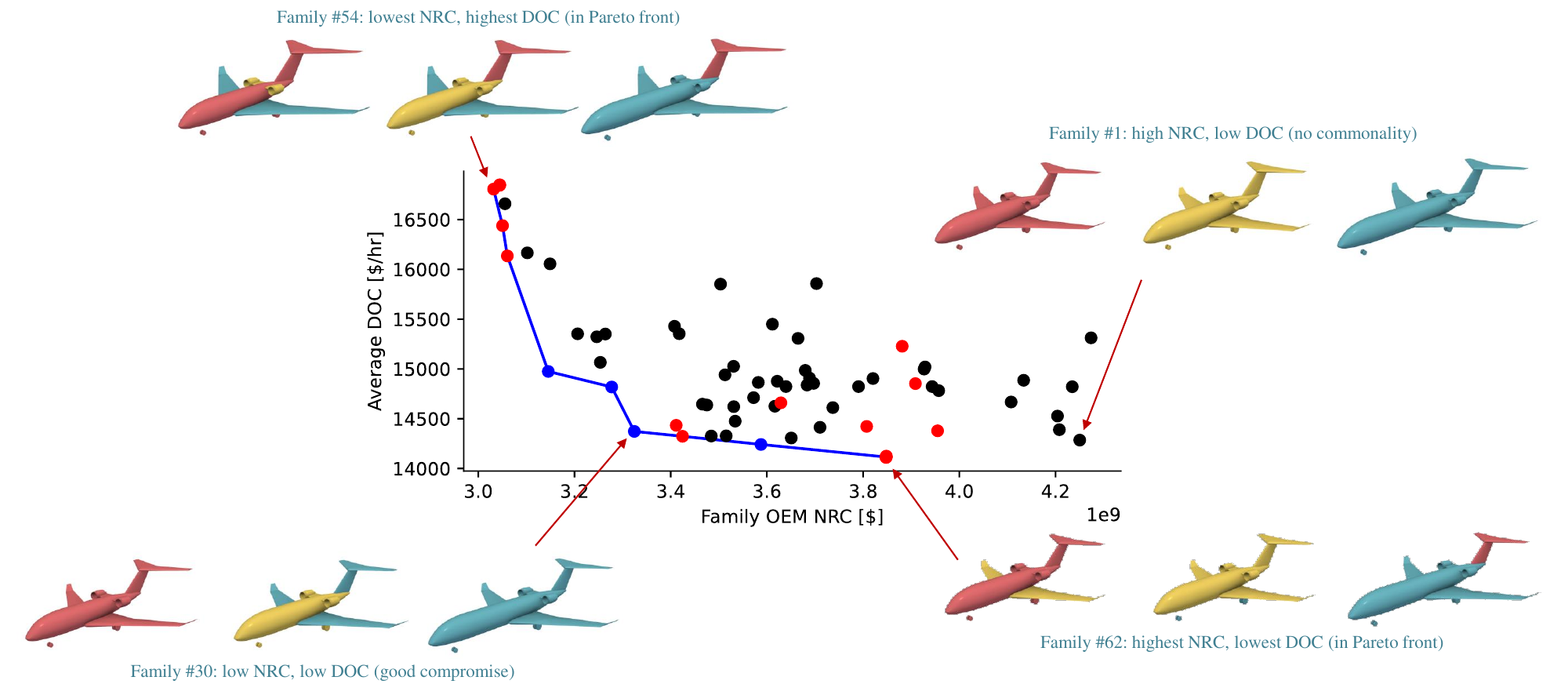}
   \caption{Results of the aircraft family design, showing the Pareto front for minimization of OEM non-recurring cost and DOC and several families. The initial DOE is given by the 50 black dots. Infill points generated by SEGOMOE are shown in red. Colors correspond to the originating aircraft as defined in Figure~\ref{fig:ac7_sharing}.}
   \label{fig:ac7_results}
\end{figure}

First, a 50 point DOE was executed to create the initial database of design points for \gls{SEGOMOE} and to verify the correct behavior of the MDO workflow. Then, \gls{SEGOMOE} was used to generate an addition of 18 infill points to explore the design space and extend the Pareto front. Figure~\ref{fig:ac7_results} shows the main Pareto front, with infill points shown in red. Several example families are shown with colors indicating the originating family member (colors are defined in Figure~\ref{fig:ac7_sharing}). 

As can be seen, when no components are shared (family  $\#1$) the OEM non-recurring costs are high whereas the operating costs are low, because all aircraft components are used at the operating points they are designed for. Introducing more component sharing reduces OEM non-recurring cost: family $\#62$ achieves the lowest operating costs at a reduced non-recurring cost. Family $\#54$ represents the opposite extreme: the lowest non-recurring cost coupled with the highest operating costs, as achieved by a high number of shared components. Family $\#30$ represents a good compromise, at a reduced non-recurring cost compared to family $\#62$, with only a moderate increase in DOC. 

\section{Conclusion}
\label{sec:conclusionc6}

To conclude, this chapter displayed practical applications of Bayesian optimization based on Gaussian process surrogate models.
In particular for optimization of eco-friendly aircraft designs and complex systems because, in that context, there is a pressing need for multi-objective optimization algorithms to deal with expensive-to-evaluate black-box functions. 
In fact, these problems often feature mixed hierarchical variables, numerous variables and multimodal equality or inequality constraints.
To begin with, we validated and illustrated our algorithm over various analytical test cases. 
Then, building upon this foundation, our work within the AGILE 4.0 project has yielded successful outcomes in addressing multi-objective problems with mixed integer variables.
Using the Efficient Global Optimization (\gls{EGO}) optimizer \gls{SEGOMOE}, we managed to keep the number of costly function evaluations low, even when dealing with a high number of design variables. 
\gls{SEGOMOE} has demonstrated its ability to solve practical problems involving 2 to 5 objective functions with several constraints and taking into account mixed integer variables. 
To apply \gls{EGO} to even more complex multidisciplinary processes, several developments need to be addressed. 
Future work will focus on extending the handling of hierarchical and mixed discrete variables while improving computational efficiency for larger databases.
Specifically, with regards to mixed integer variables, the intention is to integrate recent advancements done within the \gls{SMT} toolbox concerning mixed correlation kernels~\cite{Mixed_Paul} or hierarchical kernels~\cite{saves2023smt} for Gaussian processes into \gls{BO} and validate their effectiveness. 
These studies will be conducted within the context of the Horizon Europe COLOSSUS project\footnote{\url{https://colossus-sos-project.eu/}}. 
The COLOSSUS \gls{EU} project aims to develop a system-of-systems design methodology that facilitates the integrated optimization of aircraft, operations, and business models. 
This methodology will be applied to intermodal transport and wildfire-fighting scenarios.

\recap{ 

\lettrine[lines=2, lhang=0.33, loversize=0.25, findent=1.5em]{B}{ayesian} optimization is a powerful tool to optimize expensive-to-evaluate black-box functions and, in this chapter 
we have presented practical applications of Bayesian optimization based on Gaussian process (or EGO). In particular, we have extended and applied EGO to several configurations of aircraft design, notably for ecological purposes. This chapter responded to many problems described as follows. 
\begin{itemize}
    \item  The needs of multi-objective EGO for expensive-to-evaluate black-box problems featuring mixed hierarchical variables, many variables, high numbers of configurations and multimodal constraints have been fulfilled in the SEGOMOE software.    %
    \item The interest of SEGOMOE has been proven and illustrated on diverse analytical and engineering toy problems demonstrating its adaptability in scenarios featuring high-dimensions, mixed variable types, constraints, and multi-objective considerations.
    \item  It has been shown how SEGOMOE could serve as an interesting tool for tackling complex optimization challenges across various domains. In particular, SEGOMOE has been applied to many complex systems optimization, in particular for eco-design of aircraft configurations in the context of the AGILE 4.0 project. 
\end{itemize}


}

\chapter{Conclusions and perspectives} \label{c7}

\vspace{-0.2cm}
\setlength{\fboxrule}{0pt}
\hspace{6cm} \noindent\fbox{%
    \parbox{0.6\textwidth}{%
\hspace*{1.25cm}     Ich lebe grad, da das Jahrhundert geht. \\
\hspace*{1.25cm}     Man fühlt den Wind von einem großen Blatt,  \\
\hspace*{1.25cm}     das Gott und du und ich beschrieben hat \\
\hspace*{1.25cm}    und das sich hoch in fremden Händen dreht. \\
 \vspace{-0.3cm} \\
\hspace*{1.25cm}    Man fühlt den Glanz von einer neuen Seite, \\
\hspace*{1.25cm}    auf der noch Alles werden kann. \\
 \vspace{-0.3cm} \\
\hspace*{1.25cm}    Die stillen Kräfte prüfen ihre Breite \\
\hspace*{1.25cm}    und sehn einander dunkel an. 
     \\
     \hrule \vspace{0.2cm}
     \hspace*{\fill} Das Stunden-Buch, Rainer Maria Rilke}%
} 

\objectif{ 
\lettrine[lines=2, lhang=0.33, loversize=0.25, findent=1.5em]{T}{his} chapter concludes with the results and outcomes of this Ph.D. thesis and summarizes the main contributions and findings.
It underscores the significance of this research for both researchers and engineers, highlighting its implications for future studies. Additionally, recommendations for further researches are provided to underline the limitations of this study and put into perspective what remains to be done to continue expanding upon this work.
To finish with, practical recommendations are offered for industrial applications, with a particular emphasis on aircraft design and complex systems optimization.}










\setcounter{section}{-1}
\section{Synthèse du chapitre en français}

Ce chapitre conclut cette thèse de doctorat, dans laquelle une méthode d'optimisation a été développée pour traiter des simulations coûteuses à évaluer en grande dimension impliquant des variables mixtes. 
Cette méthode s'appuie sur un algorithme d'optimisation bayésienne capable de prendre en compte un grand nombre de variables mixtes avec plusieurs objectifs et contraintes. 
L'objectif principal de ce document était de répondre au besoin d'algorithme d'optimisation efficace pour optimiser des simulations en éco-conception d'avions malgré leurs coûts de calcul.  
La méthode centrale utilisée dans cette thèse est l'optimisation globale efficace (\gls{EGO} pour Efficient Global Optimization), un algorithme d'optimisation bayésienne basé sur des modèles de substitution. 
\gls{EGO} est particulièrement adapté aux problèmes sur lesquels peu de données et d'informations sont disponibles, ce qui est souvent le cas pour la conception d'aéronefs en raison du coût élevé des simulations. 
Néanmoins, l'adaptation d'EGO aux variables mixtes constituait un défi important, abordé dans cette thèse.

Le Chapitre~\ref{c2} a posé les bases de ce travail en permettant aux modèles de processus gaussiens (\gls{GP} pour Gaussian Process) de gérer les variables mixtes, combinant variables continues, entières et catégorielles. 
De plus, pour attaquer le problème de la haute dimension de la conception aéronautique, le Chapitre~\ref{c3} a étendu ces \gls{GP} mixtes à des espaces de grande dimension en se basant sur la régression par moindres carrés partiels.
De plus, la conception d'aéronefs implique souvent des variables hiérarchiques, où certaines variables influencent d'autres variables. Cette hiérarchie de variables a été prise en compte dans le Chapitre~\ref{c4}, dans lequel le modèle \gls{GP} a été étendu pour traiter les variables hiérarchiques. 
Dans le Chapitre~\ref{c4}, tous les modèles \gls{GP} ont été implémentés dans \gls{SMT}, une toolbox Python open-source pour permettre aux chercheurs et industriels de se servir facilement de ces travaux.
Pour finir, dans le Chapitre~\ref{c5}, toutes ces avancées méthodologiques ont été mises en œuvre pour optimiser des systèmes complexes dans le contexte de l'éco-conception d'aéronefs. De nouvelles méthodes et algorithmes ont été développés pour adapter le cadre EGO à des problèmes de contraintes multimodales, multi-objectif et de haute dimension, tout en prenant en compte des variables hiérarchiques mixtes.

Les résultats et les conclusions de cette thèse sont ici présentées. Ils démontrent que l'utilisation de modèles de substitution, en particulier les \gls{GP}, peut considérablement accélérer le processus d'optimisation pour la conception d'aéronefs respectueux de l'environnement, tout en réduisant le temps de calcul nécessaire pour les simulations coûteuses.
De plus, cette recherche met en évidence l'adaptabilité des techniques d'optimisation bayésienne pour résoudre des problèmes complexes, comme par exemples des problèmes multidisciplinaires  de conception, y compris ceux impliquant un grand nombre de variables mixtes.
Ces conclusions soulignent le potentiel de cette recherche pour l'industrie aérospatiale, où elle peut faciliter le développement de concepts innovants.
Enfin, cette thèse peut s'avérer utile pour des ingénieurs ou des chercheurs. Pour les chercheurs, elle offre des contributions significatives à la communauté des processus gaussiens et de l'optimisation bayésienne, ouvrant de nouvelles perspectives et fournissant des preuves empiriques de leur utilité. Pour les utilisateurs, en particulier ceux travaillant dans la conception aéronautique ou sur des systèmes complexes, cette recherche offre des outils pratiques et applicables qui peuvent informer la prise de décision et proposer des solutions pour améliorer des concepts d'ingénierie.

En conclusion, cette thèse marque une étape importante dans le domaine de l'optimisation en ingénierie, en particulier pour la conception d'aéronefs plus respectueux de l'environnement. Elle montre comment les modèles de substitution et l'optimisation bayésienne peuvent être adaptés avec succès pour relever les défis complexes de ce domaine, ouvrant ainsi la voie à des avancées futures dans la conception d'aéronefs plus efficaces et plus durables.

Pour étendre l'application de nos algorithmes à des processus multidisciplinaires encore plus complexes, nous avons identifié plusieurs axes de développement futurs. Ainsi, une prochaine étape portera sur la prise en charge des variables hiérarchiques et des variables discrètes mixtes grâce aux récents développements de \gls{GP} se basant sur ces variables.  
Nous pourrions également améliorer l'efficacité en terme de temps de calcul pour des bases de données plus volumineuses comme celles issues de données expérimentales, par exemple en soufflerie. 
Ainsi, des travaux de recherche seront menés dans le cadre du projet européen COLOSSUS, qui vise à développer une méthodologie de conception de système de systèmes pour l'optimisation intégrée des aéronefs, des opérations et des modèles économiques. Cette méthodologie sera appliquée à des scénarios de transport intermodaux et de lutte contre les incendies de forêt, appliquant encore une fois l'optimisation bayésienne dans des contextes de grande envergure et d'applications réelles.

\section{Thesis summary}


This thesis introduced a method for optimizing expensive-to-evaluate black-box simulations as quickly as possible. This method is based on a Bayesian optimization algorithm that can handle mixed categorical variables for high-dimensional problems with multiple objectives and constraints. 
This work has been motivated by the pressing need for optimization in eco-design of aircraft. In particular such systems are highly multidisciplinary and their Multidisciplinary Design Optimization (\gls{MDO}) need to be extended to account for the most recent industrial advances towards greener configurations.

To attain such goals, we used a surrogate-based algorithm called Efficient Global Optimization (\gls{EGO}) that is a Bayesian derivative-free optimization method tailored to account for sparse data context. 
In particular, this method relies on a Gaussian Process (\gls{GP}) metamodel so we needed to extend \gls{GP} to mixed variables as done in Chapter~\ref{c2}. 
Moreover, we needed to take into account a high number of design variables in the aircraft optimization problems, therefore, we extended to     high dimension the mixed \gls{GP} in Chapter~\ref{c3}.
Moreover, variables hierarchy often occurs in variable-size problems or within problems featuring technological choices. For example, if a variable corresponds to the number of engines, the number of motors will vary accordingly. Similarly, if the chosen motor is electric, a battery needs to be selected whereas no battery is needed for a thermal motor. In consequence, some variables influence some others, leading to a hierarchy of variables.
To take into account such hierarchy in the optimization problems, we extended the \gls{GP} model to account for hierarchical variables in Chapter~\ref{c4} and we implemented all these models in the \gls{SMT} open-source software. 
To finish with, in Chapter~\ref{c5}, we addressed our main goal, that is to optimize complex systems for eco-design of aircraft. 
To do so, we proposed new algorithms adapting the aforementioned \gls{EGO} framework to multimodal constraints, multi-objective and high dimension optimization while considering mixed hierarchical variables thanks to the previously developed \gls{GP}. All these algorithms have been implemented in the \gls{SEGOMOE} software.

\section{Results and findings}

To begin with, in Chapter~\ref{c2}, we compared several existing methods for \gls{GP} with mixed variables and we unified the distance-based approaches and the matrix-based approaches thanks to a novel unified framework. Moreover, we developed a new kernel, named EHH, that combines the well known exponential continuous kernel and the homoscedastic hypersphere kernel.
To do so, we constructed an homogeneous model that both unified and generalized the GD, CR, EHH and HH kernels. Moreover, this work gave new insights on these kernels based on a new expression in terms of hyperparameters instead of expressing the distance-based kernels in terms of input space dimensions. 
In the end, this framework led to a classification of kernels directly related to their respective number of hyperparameters. 
As Chapter~\ref{c2} connected matrix-based kernels and exponential kernels, in Chapter~\ref{c3}, we extended the high dimension \gls{KPLS} kernel developed for continuous exponential kernels to mixed variables. To do so, we proposed an extension of the \gls{PLS} from vectors to matrices and we showed the capability of our model to capture the structure in the data and to predict with very few hyperparameters the correlations between the various levels of a given categorical variable.
The implementation of these mixed variables \gls{GP} in the open-source software \gls{SMT} is detailed in Chapter~\ref{c4}. Moreover, Chapter~\ref{c4} detailed the extension of the mixed \gls{GP} to tackle hierarchical variables.
To finish with, in Chapter~\ref{c5}, we generalized Bayesian optimization with constraints with as much as 104 mixed variables and up to 5 objectives and we demonstrated how to optimize mixed hierarchical problems with these algorithms. This work led to optimized aircraft production and greener configurations. 

Our study has demonstrated that using surrogate models, especially \gls{GP}, can significantly accelerate the optimization process for environmentally friendly aircraft design while simultaneously reducing the computational time associated with costly simulations. Through extensive empirical testing and practical applications, we have shown the efficiency of our approaches to obtain aircraft configurations that consume less fuel and emit less carbon dioxide. 

Our research has also highlighted the adaptability of Bayesian optimization techniques to effectively address complex challenges in aircraft design, encompassing mixed-variable, high-dimensional, and multidisciplinary optimization problems. These findings underscore the potential of our research within the aerospace industry, as it facilitates the development of novel and innovate concepts over which the knowledge is limited, especially knowing that new aircraft computational simulations are costly and rough.

Furthermore, the significance of our work is relevant for both researchers and practitioners alike as we have shown both theoretical and practical aspects in this work.  
For researchers, our study offers substantial contributions to the Gaussian process and Bayesian optimization community, providing novel perspectives, open-source contributions and empirical evidence that our works can be used and extended for many applications and incorporated in many frameworks.
For practitioners, particularly those interested in aircraft design or complex systems, our research can be tested freely and applied straightforwardly. This work provides practical insights that can inform decision-making processes and drive enhancements in engineering domains.

\section{Practical recommendations for Bayesian optimization}
\label{sec:recommandationsc6}

For any practitioner that seeks to optimize a general, expensive-to-evaluate, unknown black-box problem, the following is recommended.

Concerning the optimization algorithms, we have seen throughout this work that Bayesian optimization is the most fitted method to respond to such problems even if several other methods exists. In particular, combining different approaches seems a promising path for future developments. For example trust-region efficient global optimization~\cite{TREGO}, surrogate-based evolutionary algorithm~\cite{sbea1}, or deep Gaussian process Bayesian optimization~\cite{hebbal2023deep} are of high interest. 

For the Bayesian optimization infill criterion, I suggest using WB2s~\cite{bartoli:hal-02149236} for mono-objective optimization and WB2s with probability of improvement and max regularization for multi-objective optimization~\cite{grapin_constrained_2022}. 
Note that, for the multi-objective multi-fidelity case, a criterion has been proposed in~\cite{MFKPLS}. 
To account for the constraints, the best method, in a general case, is to directly use the mean predictions of the constraints GP as constraints for the infill criterion optimization~\cite{BRAC_AIAA} but if the problem features equality or multimodals constraints, it would be better to use the UTB criterion~\cite{SEGO-UTB}. However, if the problem features more than 10 constraints, I would recommend using expected feasible improvement~\cite{Jones98} but more complex and up-to-date methods are available in the literature~\cite{hern2016general,picheny2014stepwise,picheny2016bayesian,pelamatti2023coupling}.

Concerning the GP model, this thesis highlighted the good performances of the KPLS model for both continuous and mixed variables. 
As such, it is recommended to use CR with PLS for mixed variables and classical GP with PLS for continuous variables. 
Regarding the choice of the reduced number of dimensions, this thesis proposes an adaptive criterion~\cite{SciTech_cat} that led to better performance than a fixed number of components strategy, even with a smaller average number of dimensions.
However, if the one-hot encoded dimension is smaller than 10, a classical GP without PLS is recommended. 
For the hierarchical variables handling, I would recommend using the Alg-Kernel developed in~\cite{saves2023smt} if the one-hot encoded dimension is smaller than 150. 
If the one-hot encoded dimension is higher than 150 and if the ratio between the size of the declared full space and the smaller valid space 
(known as the imputation ratio)~\cite{crawley2015system} is smaller than 5, I recommend using a mixed-integer model with imputation or correcting operator~\cite{agile4is2023}. 
Note that future works are required to investigate GP models for both high dimension ($\geq$150) and high imputation ratio ($\geq$5) as described in the following section.

\section{Limitations and perspectives}




This work covered many aspects and, as such, was limited in many respects, most of which have been highlighted in previous chapters.
In consequence, and despite the significant progress made in this research, we acknowledge its limitations and discern opportunities for future exploration. This section outlines these limitations and provides valuable recommendations for the direction of future research. 

First of all, the effectiveness of surrogate models, though substantial, remains subject to the quality of the available data and to the complexity of the design space~\cite{viana2021surrogate}. When the physics, the derivatives, the behavior, or any information is available, it is greatly recommended to take it into account. In particular, in aircraft design, we often receive data from different sources like wind tunnel, rough simulations, precise simulations with less data, various simulations of single disciplines, lower fidelity models,... 
To tackle this problem, R. Charayron~\cite{MFKPLS,charayron2023towards} proposed a multi-fidelity \gls{GP} surrogate model that integrates the mixed integer variables handling introduced in Chapter~\ref{c2}.
Furthermore, experimental data could be numerous, and in such context, the \gls{GP} model has to account for big data: in \gls{SMT}, a sparse Gaussian process~\cite{NIPS2016_7250eb93,pmlr-v206-moss23a} has been implemented for such purpose but this type of model has to be extended to handle mixed variables.
Another aspect worth mentioning eluded in this work is that we based our work on the assumption that we seek to optimize a unique expensive-to-evaluate black-box. However, in \gls{MDO}, we can optimize every discipline in parallel and not sequentially. This is what has been done in Efficient Global Multidisciplinary Design Optimization (EGMDO)~\cite{aerobest23cardoso,berthelin2022disciplinary} and, more generally, we restrict this work to the MultiDisciplinary Feasible (MDF) approach~\cite{delbecq2020benchmarking} but other monolitic architectures such as Independent Disciplinary Feasible (IDF) need to be investigate~\cite{Lucas_b}.
In consequence, our work is limited to one fidelity and everything as been treated as if there was only one discipline to optimize. Future works may include mixed hierarchical multi-fidelity surrogate models and mixed hierarchical EGMDO. 
Another limitation of our work is the restriction of surrogate modeling to \gls{GP}. This choice has been motivated by many reasons: the ability of \gls{GP} to provides good prediction even with very few data, its interpolation property that predicts correctly known configurations, and mostly because \gls{GP} provides an uncertainty quantification, required to perform \gls{BO}~\cite{frazierTutorialBayesianOptimization2018}. Notwithstanding other surrogate models are of high interest and need to be investigated like random forest~\cite{SMAC,SMAC3}, tree-structured Parzen estimator~\cite{TPE,Ozaki2022} or deep Gaussian process~\cite{hebbal2021bayesian}. The latter has shown good performance for multi-fidelity data and multi-objective optimization and could easily extend to mixed variables~\cite{hebbal2023deep}.

Similarly, we restrict ourselves to kernel engineering in a spatial continuous space but other approaches are also relevant for high-dimensional \gls{GP} with mixed variables such as latent variables~\cite{zhang2020latent,LVGP2}, manifold learning~\cite{yousefpour2023unsupervised,gautier2022fully}, trust regions~\cite{maus2022local,zhang2023kriging}, graph kernels~\cite{ok2020graph}, additives GP~\cite{Deng,pmlrv119} and many others. 
Concerning the latent variables approach~\cite{zhang2020latent}, it has not been investigated because of three major limitations. First, this method assumes a continuous order relation between the levels of the categorical variables, which is something we seek to be rid of to be as general as possible. Second, such an order relation and use of continuous kernel imply positive correlation, which is not always the case for categorical variables. Third, such a method embeds the categorical correlations in smaller subspaces but keeps the number of hyperparameters associated with the GP model high ( $l$ different 2D continuous latent subspaces but $ \Sigma_{i=1}^l \left( 2 L_i -3 \right)$ parameters to be tuned). More recently, the latent variable approach have been generalized and the first limitation has been removed thanks to a one-hot encoding~\cite{oune2021latent}. Another improvement brought by~\cite{oune2021latent} is that it requires a unique 2D latent space but this model still requires $ \Sigma_{i=1}^l \left(  2 L_i \right)$ parameters to be tuned. We also restrict ourselves to homoscedastic models but Compound Symmetry (CS)~\cite{Roustant} could be included as an intermediate model between one-hot encoding and Gower distance as CS has been proven to generalize GD. 

For high dimensional GP with mixed variables, future performance benchmarks on industrial test cases should include comparisons with the Maximum Likelihood Estimation (MLE) approach for latent space identification, as demonstrated in latent map Gaussian process~\cite{oune2021latent}, and, to be completely exhaustive, the EHH+PLS method, implemented in SMT could also be included in future tests. For hierarchical variables, it is important to acknowledge other approaches than the level-wise kernel method that we have followed like the subspaces approaches~\cite{DACE_hier, 10.3150/18-BEJ1049,pelamattihier,Talbi} or an extension of the latent map approach based on either classical imputation~\cite{Horn_hier} or random imputation~\cite{oune2021latent}.
However, our work is more general thanks to the framework introduced in Chapter~\ref{c4} that considers variable-wise formulation and is more flexible as we are allowing sub-problems to be intersecting. Moreover, a subspace approach can not be used when the number of possible subspaces is above 1,000 because it requires data point in every subspace and shared decreed components have to be defined multiple times which is indeed something to be avoided.

%
%

Another limitation is linked to the unavailability of derivatives. This limitation is twofold: first, the gradient-enhanced Kriging has not been extended to mixed variables, and second, the derivative of the likelihood with respect to the hyperparameters could not be computed. 
More precisely, the derivative of the likelihood is given by
$$ \frac{\partial \mathcal{L}(\Theta)}{\partial \theta_i } := \frac{1}{2} \text{Tr} \left( \left( [R(\Theta)]^{-1} \textbf{y}^f (\textbf{y}^f)^\top ([R(\Theta)]^{-1})^\top - [R(\Theta)]^{-1} \right) \frac{\partial [R(\Theta)]}{\partial \theta_i }   \right).$$ As such, to be able to compute this derivative, it is mandatory to be able to compute the derivative of the matrix kernel with respect to the hyperparameters $\frac{\partial [R(\Theta)]}{\partial \theta_i }$. For a continuous exponential kernel and for an hyperparameter $\theta_i$ that corresponds to a continuous or ordinal variable, $ \frac{\partial [R(\Theta)]}{\partial \theta_i } = -\textbf{d}_\textbf{i} [R(\Theta)] $ where $ \textbf{d}_\textbf{i} $ correspond to the distance between the $i^{\text{th}}$ components of the $n_t$ input points in the \gls{DOE}. For a matrix kernel like HH or EHH, and for an hyperparameter $\theta_i$  that corresponds to a categorical variable, $ \frac{\partial [R(\Theta)]}{\partial \theta_i } =  \frac{1}{\theta_i} [A] [R(\Theta)] $ where $[A]$ is a $n_t \times n_t$ matrix such that $ [A]_{r,s} = 1 $ if $\theta_i$ corresponds to the correlation between the level $\ell^r$ taken by $x^r$ and $\ell^s$ taken by $x^s$ and  $ [A]_{r,s} = 0 $ otherwise. Note that, for categorical hyperparameters $\theta_i$, $\frac{\partial^2 [R(\Theta)]}{\partial \theta_i^2 } =0$.
Consequently, gradient-based methods and Hessian-based methods may not be well-adapted for categorical variables handling. 
These limitations that have not been addressed in this thesis, therefore, the likelihood optimization relied on gradient-free algorithms like COBYLA~\cite{COBYLA}. Future works could consider enabling derivatives with mixed variables or replacing COBYLA with other approaches such as decomposition strategies or mesh adaptive direct search~\cite{Nomad18, barjhoux2018mixed,gamot2023hidden,Talbi,CAT-EGO}.

We could also consider extending gradient-enhanced Kriging with partial least squares to handle mixed variables thanks to the developments of Chapter~\ref{c3} in which we extended Kriging with partial least squares to handle mixed variables. Still, \gls{KPLS} suffers from many drawbacks and is not the only possible method to reduce the number of hyperparameters. For example, methods based on variables selection through metrics~\cite{spagnol_global_2019, cocchi2018chemometric,zhang2023indicator} or features extraction based on topology~\cite{oh2018bock, moriconi2020high} are seemingly less opaque. Future works could consider coupling different approaches, similarly at what has been done in EGORSE~\cite{EGORSE} in which random~\cite{RREMBO, kirschner2019adaptive} and supervised embeddings~\cite{MGP,Bouhlel18} have been successfully combined. However EGORSE combined the embeddings for the infill search during Bayesian optimization but no for building the surrogate model and future works may include the generalization of this method to handle really high dimension even with mixed variables. Forthcoming research will also consider extending high dimension surrogate models to tackle problems with strong hierarchy and hundreds of millions of configurations as in the "guidance, navigation, and control" problem.  
%
%
Concerning the \gls{SEGOMOE} optimization software for constrained problems, we assumed that every constraint is known~\cite{le2023taxonomy}. However, computational codes like multidisciplinary sizing codes can crash or cause unexpected errors~\cite{menz2023learning} and therefore, in that case, there is a need to account for unpredictable input points through "hidden constraints". In that context, \citet{tfaily2023efficient} developed a method to account for such constraints that need to be integrated in our algorithms still.  
Moreover, our works have been realized in Python, in both the SMT toolbox~\cite{SMT2019,saves2023smt} for the \gls{GP} surrogates and in SEGOMOE toolbox~\cite{bartoli:hal-02149236} for the multi-objective \gls{BO} under constraints. These softwares have also been interfaced within SBArchOpt by J.~\citet{bussemaker2023sbarchopt} to compare our algorithms with state-of-the-art on several architecture optimization problems including hidden constraints, mixed variables, millions of hierarchical architectures and multi-objective problems~\cite{agile4is2023}. Still Python provides an user-friendly environment for developers, researchers and practitioners but may be limited in terms of performance. For that reason R.~\citet{lafage2022egobox} proposed egobox, a Rust toolbox for efficient global optimization and future works may include enabling more functionalities in such efficient toolboxes. 

To finish with, the specific challenges posed by certain types of mixed-variable problems may require deeper investigation and tailored solutions. We plan to further investigate the integration of hierarchical kernels and mixed correlation kernels for Gaussian processes and its application  for Bayesian optimization. These advancements, as seen in the SMT toolbox, have the potential to enhance the effectiveness of our optimization approaches for high-dimensional mixed hierarchical problems. Moreover, the generalization of our findings to other complex engineering systems beyond aircraft design will be considered, with a particular focus on hyperparameters optimization for neural networks with E. Hallé-Hannan~\cite{audet2022general,schede2022survey}. As such, the principles and methods developed in this thesis are not limited to aircraft design. They can be adapted and extended to address sustainability challenges in other multidisciplinary engineering fields, such as automotive, renewable energy, and many more~\cite{Martins2021,xu2023using,yamawaki2021decomposition,huang2023mixed,gossard2022bayesian}.

This work will continue within the framework of the Horizon Europe COLOSSUS project, where we aim to apply our methodologies to the integrated optimization of aircraft, operations, and system-of-systems problems, with application to firefighting drones. To conclude this Ph.D. thesis, we highlight that this work has been an opportunity to strengthen numerous collaborations between ONERA, ISAE-SUPAERO, Polytechnique Montréal, the European Union and particularly DLR. It also displays the need to bridge the gap between academic research and industrial application and was an opportunity to assess the growing interest towards machine learning methods in engineering, particularly via the application of artificial intelligence.

\newpage
\let\oldcm=\chaptermark
\renewcommand{\chaptermark}[1]{\markboth{}{#1}}
\chaptermark{ LIST OF PUBLICATIONS }
 \phantomsection
\begin{spacing}{1.5}
\chapter*{\centering \large\textbf{LIST OF PUBLICATIONS}}
\addcontentsline{toc}{section}{\bf LIST OF PUBLICATIONS}
\noindent 

Parts of the work presented in this documents have been published in either journals or conferences. They are described hereinafter.

\textbf{List of awards:}

\begin{enumerate}

\item Published: 2023 – Saves, P., Bartoli, N., Diouane, Y., Lefebvre, T., Morlier, J., David, C., Nguyen Van, E., and Defoort, S., “AIAA Best MDO Paper Award 2022. Works presentation and award reception”. in: AIAA SciTech 2023, 2023, National Harbor (USA).
\end{enumerate}

\textbf{List of published research in Journals:} 

\begin{enumerate}

\item Published: 2023 – Saves, P., Diouane, Y., Bartoli, N., Lefebvre, T., Morlier, J., “A mixed-categorical correlation kernel for Gaussian process”, Neurocomputing, 2023, https://doi.org/10.1016/j.neucom.2023.126472

\item Published: 2024 – Saves, P., Lafage, R., Bartoli, N., Diouane, Y., Bussemaker, J., Lefebvre, T., Morlier, J., Hwang, J. T., Martins, J. R. R. A., “SMT 2.0: A Surrogate Modeling Toolbox with Hierarchical and Mixed Variables Gaussian Processes”, Advances in Engineering Software, 2024, https://doi.org/10.1016/j.advengsoft.2023.103571

\item Under Review: 2024 – Saves, P., Diouane, Y., Bartoli, N., Lefebvre, T., Morlier, J., “High-dimensional mixed-categorical Gaussian processes with application to multidisciplinary design optimization for a green aircraft”, Structural and Multidisciplinary Optimization, 2024

\item Under Review: 2024 – Bussemaker, J.H., Saves, P., Bartoli, N., Lefebvre, T., Lafage, R., Nagel, B., “System Architecture Optimization Strategies: Dealing with Expensive Hierarchical Multi-Objective Problems”, Structural and Multidisciplinary Optimization, 2024


\end{enumerate}

\textbf{List of conferences publications with proceedings:}

\begin{enumerate}

\item Published:  2021 –  Saves, P., Bartoli, N., Diouane, Y., Lefebvre, T., Morlier, J., David, C., Nguyen Van, E., and Defoort, S., “Constrained Bayesian optimization over mixed categorical variables, with application to aircraft design”. in: ECCOMAS Aerobest, 2021, Lisboa (Portugal), https://hal.archives-ouvertes.fr/hal-03346341/ 

\item Published: 2022 – Saves, P., Bartoli, N., Diouane, Y., Lefebvre, T., Morlier, J., David, C., Nguyen Van, E., and Defoort, S., “Bayesian optimization for mixed variables using an adaptive dimension reduction process: applications to aircraft design”. in: AIAA SciTech 2022, 2022, San Diego (USA). https://hal.archives-ouvertes.fr/hal-03514915/  

\item Published:  2022 –  Saves, P., Diouane, Y., Bartoli, N., Lefebvre, T., Morlier, J., “A general square exponential kernel to handle mixed-categorical variables for Gaussian process”. in: AIAA Aviation 2022, 2022, Chicago (USA), https://hal.archives-ouvertes.fr/hal-03700850/ 

\item Published:  2022 – Grapin, R., Diouane, Y., Morlier, J., Bartoli, N., Lefebvre T., Saves, P., Bussemaker, J.H., “Regularized Infill Criteria for Multi-objective Bayesian Optimization with Application to Aircraft Design”. in: AIAA Aviation 2022, 2022, Chicago (USA),  https://hal-mines-albi.archives-ouvertes.fr/hal-03753674/ 

\item Published: 2023 – Priem, R., Diouane, Y., Bartoli, N., Dubreuil, S., Saves, P., “High-dimensional efficient global opti-mization using both random and supervised embeddings”. In: Aviation 2023, 2023, San Diego (USA), https://hal.science/hal-04123951/ 

\item Published: 2023 –  Bartoli, N., Lefebvre, T., Lafage, R., Saves, P., Diouane Y., Morlier  J., Bussemaker, J.H., Donelli G., Gomes de Mello, J.M., Mandorino, M., Della Vecchia, P., “Multi Objective Bayesian Optimization with mixed integer variables for Aeronautical applications”. In: Aerobest 2023, 2023, Lisboa  (Portugal), https://hal.science/hal-04170287 

\item Under Review: 2024 – Bussemaker, J.H., Saves, P., Bartoli, N., Lefebvre T., Nagel B., “Surrogate-Based Optimization of System Architectures Subject to Hidden Constraints”, In: Aviation 2024, 2024, Las Vegas (USA)

\end{enumerate}

\textbf{List of communications and other publications:}

\begin{enumerate}
\item Published: 2020 – Saves, P., “On adapting the Super-Efficient Global Optimization solver to handle mixed-variables, with applications in aircraft design”. Master’s thesis. ISAE–SUPAERO - ONERA, 2020. 

\item Published: 2021 – Saves, P., “High dimensional multidisciplinary design optimization with Gaussian process for eco-design aircraft”.  MASCOT-NUM21 PdD day

\item Published: 2021 – Saves, P.,  “Présentation de la Toolbox SMT - L'I.A. pour les modèles de substitution”,   In: Sur-Day ISAE-SUPAERO  

\item Published: 2021 –  Saves, P., Bartoli, N., Diouane, Y., Lefebvre, T., and Morlier, J., “Enhanced kriging models within a bayesian optimization framework, to handle both continuous and categorical inputs”. In: SIAM CSE21, 2021, Fort Worth (USA)

\item Published: 2022 – Saves, P., Diouane, Y., Bartoli, N., Lefebvre, T., Morlier, J., “Towards a unified Gaussian process kernel-based correlation matrix representation for mixed-categorical variables”. in: JOPT 2022, 2022, Montréal (Canada) 

\item Published: 2022 – Saves, P., Diouane, Y., Bartoli, N., Lefebvre, T., Morlier, J., “New advances in Multidisciplinary Design Optimization with Gaussian Process for Eco-design Aircraft”. in: 3rd European MDO Workshop – Toward green aviation, 2022, Paris 

\item Published: 2023 – Saves, P., “High dimensional multidisciplinary design optimization with Gaussian process for eco-design aircraft”.  MASCOT-NUM 2023

\end{enumerate}

\end{spacing}
\let\chaptermark=\oldcm

\thispagestyle{plain}
\begin{appendices}
    \appendixpage
    \noappendicestocpagenum
    \addappheadtotoc
    \chapter{Test cases appendix}

\section{Categorical cosine case}
\label{subsec:cosine}
\begin{figure}[H]
\centering
\includegraphics[scale=.12]{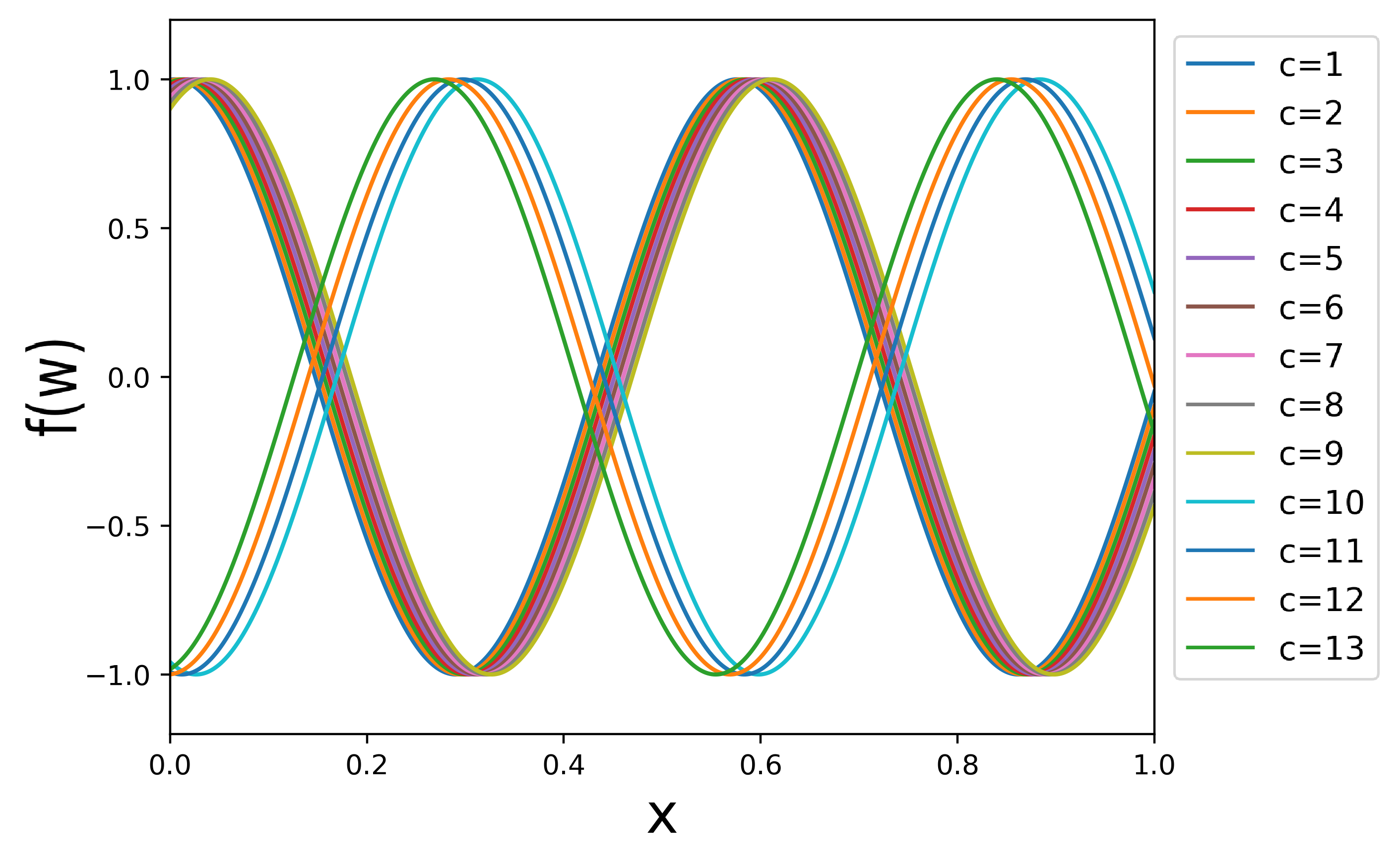}
\caption{Landscape of the cosine test case from~\cite{Roustant}.}
\label{fig:Roustant_ref}    
\end{figure}

This test case has one categorical variable with 13 levels and one continuous variable in $[0,1]$~\cite{Roustant}.
Let $w= (x,c )$ be a given point with  $x$ being the continuous variable and $c$ being the categorical variable, $c \in \{1, \ldots, 13\}$.

\begin{equation*}
\begin{split}
f(w) &= \cos \left( \frac{7 \pi}{2} x + \left( 0.4 \pi  + \frac{\pi }{15} c  \right) - \frac{c}{20} \right) , ~~~\mbox{if c $\in\{1,\ldots,9\}$ }  \\
f(w) &= \cos \left( \frac{7 \pi}{2} x  - \frac{c}{20} \right) , ~~~\mbox{if c $\in\{10,\ldots,13\}$ }  \\
\end{split}
\end{equation*}
The reference landscapes of the objective function (with respect to the categorical choices) are drawn on~\figref{fig:Roustant_ref}.

The DoE is given by a LHS of 98 points.
Our validation set is a evenly spaced grid of 1000 points in $x$ ranging  for every of the 13 categorical levels  for a total of 13000 points.


\section{Toy test function}
\label{app:Toy}
This Appendix gives the detail of the toy function of Section~\ref{sec:MI-BO}~\footnote{\url{https://github.com/jbussemaker/SBArchOpt}}.
First, we recall the optimization problem:
\begin{equation}
\begin{split}
& \min  f( x^{\cat}, x^{\quant}) \\
& \mbox{w.r.t.} \ \  x^{\cat} = c_1 \in \{ 0,1,2,3,4,5,6,7,8,9 \} \\
& \quad \quad \quad x^{\quant} = x_1 \in [ 0,1 ] \\
\end{split}
\end{equation}
The toy function $f$ is defined as
\begin{equation}
\begin{split}
   f({x_1}, {c_1})  =&   \mathds{1}_{c_1=0}  \ \cos(3.6 \pi(x-2)) +x -1 \\
     + &\mathds{1}_{c_1=1} \ 2 \cos(1.1 \pi \exp(x)) - \frac{x}{2} +2  \\
     + &\mathds{1}_{c_1=2}  \  \cos ( 2 \pi x)  + \frac{1}{2}x \\
  +&\mathds{1}_{c_1=3}  \  x ( \cos(3.4 \pi (x-1)) - \frac{x-1}{2})\\
  +& \mathds{1}_{c_1=4}  \   - \frac{x^2}{2} \\
 + &\mathds{1}_{c_1=5}  \  2 \cos(0.25 \pi \exp( -x^4))^2 - \frac{x}{2} +1 \\ 
 +& \mathds{1}_{c_1=6}   \ x \cos(3.4 \pi x ) - \frac{x}{2} +1 \\ 
 +&\mathds{1}_{c_1=7}   \  - x  (\cos(3.5 \pi x ) + \frac{x}{2}) +2 \\ 
 + &\mathds{1}_{c_1=8}   \ - \frac{x^5}{2} +1 \\ 
 +& \mathds{1}_{c_1=9}  \  - \cos (2.5 \pi x)^2 \sqrt{x} - 0.5 \ln (x+0.5)  - 1.3\\ 
\end{split}
\end{equation}

\section{Hierarchical Goldstein test function}
\label{app:Goldstein}
This Appendix gives the detail of the hierarchical Goldstein problem of Section~\ref{sec:HV-BO}~\footnote{\url{https://github.com/jbussemaker/SBArchOpt}}.
First, we recall the optimization problem:
\begin{equation}
\begin{split}
& \min  f( x^{\cat}_{\neutral}, x^{\quant}_{\neutral}, x^{\cat}_{m}, x^{\quant}_{\decreed}  ) \\
& \mbox{w.r.t.} \ \
x^{\cat}_{\neutral} =w_2 \in \{ 0,1 \} \\
& \quad \quad \quad x^{\quant}_{\neutral} = (x_1,x_2,x_5,z_3,z_4) \in [ 0,100 ]^3 \times \{ 0,1,2 \}^2  \\
& \quad \quad \quad   x^{\cat}_{m} = w_1 \in \{ 0,1,2,3 \} \\
& \quad \quad \quad 
 x^{\quant}_{\decreed} = (x_3,x_4,z_1,z_2) \in [ 0,100 ]^2 \times \{ 0,1,2 \}^2   
\end{split}
\end{equation}
The hierarchical and mixed function $f$ is defined as a hierarchical function that depends on $f_0$,  $f_1$, $f_2$ and $Gold_\text{cont}$ as describes in the following.
\begin{equation}
\begin{split}
   &f({x_1}, {x_2}, {x_3}, {x_4}, {z_1}, {z_2}, {z_3}, {z_4}, {x_5}, {w_1}, {w_2})  = \\
    &  \quad \  \mathds{1}_{w_1=0} f_0( {x_1}, {x_2},  {z_1}, {z_2}, {z_3}, {z_4}, {x_5}, {w_2}) \\
    & + \mathds{1}_{w_1=1} f_1 ({x_1}, {x_2}, {x_3}, {z_2}, {z_3}, {z_4}, {x_5}, {w_2})  \\
    & + \mathds{1}_{w_1=2} f_2 ({x_1}, {x_2}, {x_4}, {z_1}, {z_3}, {z_4}, {x_5}, {w_2})  \\
   & + \mathds{1}_{w_1=3}  Gold_{\text{cont}} ({x_1}, {x_2}, {x_3}, {x_4}, {z_3}, {z_4}, {x_5},  {w_2}).
\end{split}
\end{equation}
Then, the functions $f_0$, $f_1$ and $f_2$ are defined as mixed variants of $Gold_\text{cont}$ as such
\begin{equation}
\begin{split}
   &f_0( {x_1}, {x_2},  {z_1}, {z_2}, {z_3}, {z_4}, {x_5}, {w_2})  = \\
       & \mathds{1}_{z_2=0} \big(  \mathds{1}_{z_1=0}  Gold_{\text{cont}} ({x_1}, {x_2}, 20,20, {z_3}, {z_4}, {x_5},  {w_2})  \\
      & \quad   + \mathds{1}_{z_1=1} Gold_{\text{cont}} ({x_1}, {x_2}, 50,20, {z_3}, {z_4}, {x_5},  {w_2})   \\
      & \quad  + \mathds{1}_{z_1=2}Gold_{\text{cont}} ({x_1}, {x_2}, 80,20, {z_3}, {z_4}, {x_5},  {w_2})   \big) \\
       & \mathds{1}_{z_2=1} \big(  \mathds{1}_{z_1=0}  Gold_{\text{cont}} ({x_1}, {x_2}, 20,50, {z_3}, {z_4}, {x_5},  {w_2})  \\
      & \quad   + \mathds{1}_{z_1=1} Gold_{\text{cont}} ({x_1}, {x_2}, 50,50, {z_3}, {z_4}, {x_5},  {w_2})   \\
      & \quad  + \mathds{1}_{z_1=2}Gold_{\text{cont}} ({x_1}, {x_2}, 80,50, {z_3}, {z_4}, {x_5},  {w_2})   \big) \\
      & \mathds{1}_{z_2=2} \big(  \mathds{1}_{z_1=0}  Gold_{\text{cont}} ({x_1}, {x_2}, 20,80, {z_3}, {z_4}, {x_5},  {w_2})  \\
      & \quad   + \mathds{1}_{z_1=1} Gold_{\text{cont}} ({x_1}, {x_2}, 50,80, {z_3}, {z_4}, {x_5},  {w_2})   \\
      & \quad  + \mathds{1}_{z_1=2}Gold_{\text{cont}} ({x_1}, {x_2}, 80,80, {z_3}, {z_4}, {x_5},  {w_2})   \big) \\
  \end{split}
  \end{equation}
  \begin{equation*}
  \begin{split}
   &f_1 ({x_1}, {x_2}, {x_3}, {z_2}, {z_3}, {z_4}, {x_5}, {w_2})  = \\
       &  \mathds{1}_{z_2=0}  Gold_{\text{cont}} ({x_1}, {x_2}, {x_3},20, {z_3}, {z_4}, {x_5},  {w_2})  \\
      & \quad   + \mathds{1}_{z_2=1} Gold_{\text{cont}} ({x_1}, {x_2}, {x_3},50, {z_3}, {z_4}, {x_5},  {w_2})   \\
      & \quad  + \mathds{1}_{z_2=2}Gold_{\text{cont}} ({x_1}, {x_2}, {x_3},80, {z_3}, {z_4}, {x_5},  {w_2})  \\[12pt]
   &f_2 ({x_1}, {x_2}, {x_4}, {z_1}, {z_3}, {z_4}, {x_5}, {w_2})  = \\
       &  \mathds{1}_{z_1=0}  Gold_{\text{cont}} ({x_1}, {x_2}, 20, {x_4}, {z_3}, {z_4}, {x_5},  {w_2})  \\
      & \quad   + \mathds{1}_{z_1=1} Gold_{\text{cont}} ({x_1}, 50, {x_2}, {x_4}, {z_3}, {z_4}, {x_5},  {w_2})   \\
      & \quad  + \mathds{1}_{z_1=2}Gold_{\text{cont}} ({x_1}, {x_2}, 80, {x_4}, {z_3}, {z_4}, {x_5},  {w_2})  
\end{split}
\end{equation*}
To finish with, the function $Gold_{\text{cont}}$ is given by
\begin{equation}
\begin{split}
    & Gold_{\text{cont}} ({x_1}, {x_2}, {x_3}, {x_4}, {z_3}, {z_4}, {x_5},  {w_2}) = 
        53.3108
        + 0.184901   {x_1}  \\
        &- 5.02914   {x_1} ^3   .10    ^{-6}
        + 7.72522    {x_1} ^{z_3}    .10    ^{-8}
        - 0.0870775    {x_2}
        - 0.106959    {x_3}  \\
        &+ 7.98772    {x_3} ^{z_4}    .10    ^{-6} 
         + 0.00242482    {x_4}
        + 1.32851    {x_4} ^3    .10    ^{-6}
        - 0.00146393    {x_1}    {x_2}\\
        &- 0.00301588    {x_1}    {x_3} 
         - 0.00272291    {x_1}    {x_4}
        + 0.0017004    {x_2}    {x_3}
        + 0.0038428    {x_2}    {x_4}\\
       & - 0.000198969    {x_3}    {x_4} 
        + 1.86025    {x_1}    {x_2}    {x_3}    .10    ^{-5}
        - 1.88719    {x_1}    {x_2}    {x_4}    .10    ^{-6}\\
       & + 2.50923    {x_1}    {x_3}    {x_4}    .10    ^{-5} 
        - 5.62199    {x_2}    {x_3}    {x_4}    .10    ^{-5} 
        +  {w_2} \left( 5    \cos \left( \frac{ 2 \pi}{100}  x_5 \right) - 2\right).
\end{split}
\end{equation}

\section{Engineering test cases}
\label{app:eng_cases}
 
The following gives an in-depth description for every test case analyzed or optimized in Chapter~\ref{c5}.

\subsection{Rover path planning}
\label{app:rover}

\begin{table}[H]
    \centering
   \caption{Definition of the Rover path planning problem.}
   \small
   \resizebox{\columnwidth}{!}{%
      \begin{tabular}{lllrr}
 & Function/variable & Nature & Quantity & Range\\
\hline
\hline
Minimize & Trajectory length & cont & 1 &\\
 & \multicolumn{2}{l}{\bf Total objectives} & {\bf 1} & \\
\hline
with respect to & Variables encoding 30 control points in a 2D plane & cont  &  600 &  $[0, 1]  $ \\
 & \multicolumn{2}{r}{Total continuous variables} & { 600} & \\ \cline{2-5}
 & \multicolumn{2}{r}{\bf Total design variables} & {\bf 600} & \\
 \hline
\end{tabular}
}
   \label{tab:roverc6}
\end{table}

\subsection{Mixed cantilever beam}
\label{app:canti}

\begin{table}[H]
    \centering
   \caption{Definition of the mixed cantilever beam problem.}
   \small
   \resizebox{\columnwidth}{!}{%
      \begin{tabular}{lllrr}
 & Function/variable & Nature & Quantity & Range\\
\hline
\hline
Minimize & Displacement & cont & 1 &\\
 & \multicolumn{2}{l}{\bf Total objectives} & {\bf 1} & \\
\hline
with respect to & Length & cont  &  1 &  $[10, 20] \ (m) $ \\
&  Surface & cont  &  1 &  $[1, 2]  \ (m^2) $ \\
 & \multicolumn{2}{r}{Total continuous variables} & { 2} & \\ \cline{2-5}
&  Beam shape & cat  &  12 levels &  various beam shapes \\
& \multicolumn{2}{r}{ Total categorical variables} & { 1} & \\ \cline{2-5}
 & \multicolumn{2}{r}{\bf Total relaxed variables} & {\bf 14} & \\
 \hline
\end{tabular}
}
   \label{tab:cantic6}
\end{table}

\subsection{Neural network problem}
\label{app:nn}

\begin{table}[H]
    \centering
   \caption{Definition of the neural network problem.}
   \small
   \resizebox{\columnwidth}{!}{%

      \begin{tabular}{lllrr}
 & Function/variable & Nature & Quantity & Range\\
\hline
\hline
Minimize & Loss function & cont & 1 &\\
 & \multicolumn{2}{l}{\bf Total objectives} & {\bf 1} & \\
\hline
with respect to &  Learning rate & cont  &  1 &  $[10^{-5}, 10^{-2}]$ \\
 & \multicolumn{2}{r}{Total continuous variables} & { 1} & \\ \cline{2-5}
 & Batch size & discrete & 1 & \{8,16, $\ldots$, 256\} \\
 & Number of hidden layers & discrete & 1 & \{1,2,3\} \\
 & Number of neurons in hidden layers* & discrete & 3 & $[50,55]$ \\
 & \multicolumn{3}{l}{* only active if the corresponding hidden layer is included} \\
 & \multicolumn{2}{r}{Total discrete variables} & {5} & \\ \cline{2-5}
&  Activation function & cat  &  2 levels & \{ReLu, Sigmoid\} \\
& \multicolumn{2}{r}{ Total categorical variables} & { 1} & \\ \cline{2-5}
 & \multicolumn{2}{r}{\bf Total relaxed variables} & {\bf 8} & \\
 \hline
\end{tabular}
}
   \label{tab:NNc6}
\end{table}

\subsection{Mixed integer `\texttt{CERAS}'}
\label{app:ceras}

\begin{table}[H]
\centering
   \caption{Definition of the ``\texttt{CERAS}'' optimization problem.}
\small
\resizebox{\columnwidth}{!}{%
\small
\begin{tabular}{lllrr}
 & Function/variable & Nature & Quantity & Range\\
\hline
\hline
Minimize & Fuel mass & cont & 1 &\\
 & \multicolumn{2}{l}{\bf Total objectives} & {\bf 1} & \\
\hline
with respect to & \mbox{x position of MAC} & cont & 1 & $\left[16., 18.\right]$ ($m$)\\
 & \mbox{Wing aspect ratio}  & cont &1 & $\left[5., 11.\right]$ \\
 & \mbox{Vertical tail aspect ratio} & cont & 1 & $\left[1.5, 6.\right]$ \\
 & \mbox{Horizontal tail aspect ratio} & cont & 1 & $\left[1.5, 6.\right]$ \\
  & \mbox{Wing taper aspect ratio} & cont & 1 & $\left[0., 1.\right]$ \\
   & \mbox{Angle for swept wing} & cont & 1 & $\left[20., 30.\right]$ ($^\circ$)\\
 & \multicolumn{2}{c}{Total  continuous variables} & 6 & \\
 \cline{2-5}
& \mbox{Cruise altitude} & discrete & 1 & \{30k,32k,34k,36k\} ($ft$)\\
& \mbox{Number of engines} & discrete & 1 & \{2,3,4\} \\
 & \multicolumn{2}{c}{Total  discrete variables} & 2 & \\
 \cline{2-5}
& \mbox{Tail geometry} & cat & 2 levels & \{T-tail, no T-tail\} \\
& \mbox{Engine position}& cat & 2 levels & \{\mbox{front or rear engines}\}  \\
 & \multicolumn{2}{c}{Total  categorical variables} & 2 & \\
 \cline{2-5}
  &   \multicolumn{2}{c}{\textbf{Total relaxed variables}} & {\textbf{12}} & \\
  \hline
subject to & 0.05 \textless \ Static margin \textless \  0.1 & cont & 2 \\
 & \multicolumn{2}{c}{\textbf{Total  constraints}} & {\textbf{2}} & \\
\hline
\end{tabular}
}
\label{tab:cerasc6}
\end{table}

\subsection{Multi-objective `\texttt{CERAS}'}
\label{app:ceras_moo}

\begin{table}[H]
\centering
   \caption{Definition of the ``\texttt{CERAS}'' bi-objective optimization problem.}
\small
\resizebox{\columnwidth}{!}{%
\small
\begin{tabular}{lllrr}
 & Function/variable & Nature & Quantity & Range\\
\hline
\hline
Minimize & Fuel mass & cont &  1 &\\
&  Operating Weight Empty & cont & 1 \\
 & \multicolumn{2}{l}{\bf Total objectives} & {\bf 2} & \\
\hline
with respect to & \mbox{x position of Mean Aerodynamic Chord} & cont &  1 & $\left[16., 18.\right]$ ($m$)\\
 & \mbox{Wing aspect ratio} & cont  &1 & $\left[5., 11.\right]$ \\
 & \mbox{Horizontal tail aspect ratio} & cont & 1 & $\left[1.5, 6.\right]$ \\
  & \mbox{Wing taper aspect ratio} & cont &  1 & $\left[0., 1.\right]$ \\
   & \mbox{Angle for swept wing} & cont &  1 & $\left[20., 30.\right]$ ($^\circ$) \\
    \cline{2-5}
 & {\textbf{Total  design  variables}} &  & {\textbf{5}} & \\
 \cline{2-4}
  \hline
subject to & 0.05 \textless \ Static margin \textless \  0.1 & cont & 2 \\
& Wing span \textless \  36   & cont &1\\
 & {\textbf{Total  constraints}} &  & {\textbf{3}} & \\
\hline
\end{tabular}
}
\label{tab:cerasmoc6}
\end{table}

\subsection{`\texttt{DRAGON}' aircraft concept}
\label{app:dragon}

\begin{table}[H]
\centering
 \caption{Definition of the ``\texttt{DRAGON}'' optimization problem.}
\small

\resizebox{\columnwidth}{!}{%
\small

\begin{tabular}{lllrr}
 & Function/variable & Nature & Quantity & Range\\
\hline
\hline
Minimize & Fuel mass & cont & 1 &\\
 & \multicolumn{2}{l}{\bf Total objectives} & {\bf 1} & \\
\hline
with respect to & \mbox{Fan operating pressure ratio} & cont & 1 & $\left[1.05, 1.3\right]$ \\  
     & \mbox{Wing aspect ratio} & cont & 1 &    $\left[8, 12\right]$ \\
    & \mbox{Angle for swept wing} & cont & 1 & $\left[15, 40\right]$  ($^\circ$) \\
     & \mbox{Wing taper ratio} & cont & 1 &    $\left[0.2, 0.5\right]$ \\
     & \mbox{HT aspect ratio} & cont & 1 &    $\left[3, 6\right]$ \\
    & \mbox{Angle for swept HT} & cont & 1 & $\left[20, 40\right]$  ($^\circ$) \\
     & \mbox{HT taper ratio} & cont & 1 &    $\left[0.3, 0.5\right]$ \\
 & \mbox{TOFL for sizing}  & cont &1 & $\left[1800., 2500.\right]$ ($m$) \\
 & \mbox{Top of climb vertical speed for sizing} & cont & 1 & $\left[300., 800.\right]$($ft/min$) \\
 & \mbox{Start of climb slope angle} & cont & 1 & $\left[0.075., 0.15.\right]$($rad$) \\

 & \multicolumn{2}{l}{Total  continuous variables} & 10 & \\
 \cline{2-5}
& \mbox{Architecture} & cat & 17 levels & \{1,2,3, \ldots,15,16,17\} \\
& \mbox{Turboshaft layout} & cat & 2 levels & \{1,2\} \\

 & \multicolumn{2}{l}{Total categorical variables} & 2 & \\
 \cline{2-5}

  &   \multicolumn{2}{l}{\textbf{Total relaxed variables}} & {\textbf{29}} & \\
  \hline
  
subject to & Wing span \textless  \ 36   ($m$)  & cont & 1 \\
 & TOFL \textless  \ 2200 ($m$) & cont & 1 \\
 & Wing trailing edge occupied by fans  \textless  \ 14.4 ($m$) & cont & 1 \\
 & Climb duration \textless  \ 1740 ($s $) & cont & 1 \\
 & Top of climb slope \textgreater \ 0.0108 ($rad$) & cont & 1 \\

 & \multicolumn{2}{l}{\textbf{Total  constraints}} & {\textbf{5}} & \\
\hline
\end{tabular}
}
\label{tab:dragonc6app1}
\end{table}

\subsection{Airframe upgrade rettrofitting}
\label{app:air_retro}

\begin{table}[H]
\centering
   \caption{Definition of  the airframe upgrade design problem.}
   \small
   \resizebox{\columnwidth}{!}{
      \begin{tabular}{lllrr}
 & Function/variable & Nature & Quantity & Range\\
\hline
\hline
Minimize & Maximum Take-Off Weight & cont & 1 &\\
Minimize &  Cumulative Emission Index & cont & 1 &\\
Minimize &  Cost - savings & cont & 1 &\\
Maximize &  Max Cruise Specific Air Range & cont & 1 &\\
 & \multicolumn{2}{l}{\bf Total objectives} & {\bf 4} & \\
\hline
with respect to &  Engine Bypass Ratio & cont & 1 & $\left[9,15 \right]$ \\
 & Engine $X$ & cont & 1 & $\left[ -0.98, -0.80  \right]$ ($m$)\\
 & Engine $Z$ & cont & 1 & $\left[-0.39,-0.21\right]$  ($m$)\\
 & \multicolumn{2}{l}{ Total continuous variables} & { 3} & \\ \cline{2-5}
& OBS architecture & cat  &  4 levels & [CONV, MEA1, MEA2, AEA]\\
 & \multicolumn{2}{l}{ Total categorical variables} & { 1} & \\ \cline{2-5} 
 & \multicolumn{2}{l}{\bf Total relaxed variables} & {\bf 7} & \\
\hline
subject to & 
\multicolumn{2}{l}{ Maximum Take-Off Weight $\le 39058.5$ ($kg$) } & 1 & \\
&\multicolumn{2}{l}{ Take-Off Field Length $ \le 1500$  ($m$) }  & 1 & \\
&\multicolumn{2}{l}{ Landing Field Length $\le 1400$ ($m$) }  & 1 & \\
&\multicolumn{2}{l}{ Cumulative Noise $ \le 263$ ($dB$) } & 1 & \\
%
 & \multicolumn{2}{l}{\bf Total constraints } & {\bf 4} & \\
\hline
\end{tabular}
}
   \label{tab:pb_ac6}
\end{table}

\subsection{Family of business jets}
\label{app:jet_family}

\begin{table}[H]
    \centering
   \caption{Definition of the aircraft family design problem.}
   \small
      \begin{tabular}{lllrr}
 & Function/variable & Nature & Quantity & Range\\
\hline
\hline
Minimize & Direct Operating Costs & cont & 1 &\\
Minimize &  OEM Non-Recurring Costs& cont & 1 &\\
 & \multicolumn{2}{l}{\bf Total objectives} & {\bf 2} & \\
\hline
with respect to   &  Leading edge sweep* & cont & 3 & $\left[ 30.0, 42.0 \right]$ ($^\circ$)\\
 & Rear span location* & cont & 3 & $\left[ 0.72, 0.82 \right]$ ($\% chord$)\\
 & Wing thickness/cord ratio*  & cont & 3 & $\left[0.06, 0.11\right]$ \\
 & \multicolumn{3}{l}{* only active if the corresponding wing is \textbf{not} shared} \\
 & \multicolumn{2}{l}{Total continuous variables} & { 9} & \\ \cline{2-5}
& Engine commonality  1\&2 & cat  &  2 levels & sharing yes/no\\
& Engine commonality  2\&3 & cat  &  2 levels & sharing yes/no \\
& Wing commonality  1\&2 & cat  &  2 levels & sharing yes/no\\
& Wing commonality  2\&3 & cat  &  2 levels & sharing yes/no\\
& Landing gear commonality 1\&2 & cat  &  2 levels & sharing yes/no\\
& Landing gear commonality 2\&3 & cat  &  2 levels & sharing yes/no\\
& OBS commonality 1\&2 & cat  &  2 levels & sharing yes/no\\
& OBS commonality 2\&3 & cat  &  2 levels & sharing yes/no\\
& Empennage commonality  1\&3 & cat  &  2 levels & sharing yes/no\\
& Empennage commonality  3\&2 & cat  &  2 levels & sharing yes/no\\
& \multicolumn{2}{l}{ Total categorical variables} & { 10} & \\ \cline{2-5}
 & \multicolumn{2}{l}{\bf Total relaxed variables} & {\bf 29} & \\
\hline
subject to & 
\multicolumn{2}{l}{ Balanced  Field Length $\le 1524$ ($m$) } & 1 & \\
&\multicolumn{2}{l}{ Landing Field Length $ \le 762$ ($m$) }  & 1 & \\
 & \multicolumn{2}{l}{\bf Total constraints} & {\bf 2} & \\
\hline
\end{tabular}
   \label{tab:pb_ac7app1}
\end{table}

\subsection{Production of aircraft}
\label{app:air_prod}

\begin{table}[H]
    \centering
   \caption{Definition of the aircraft production problem.}
   \small
      \begin{tabular}{lllrr}
 & Function/variable & Nature & Quantity & Range\\
\hline
\hline
Minimize & Cost & cont & 1 &\\
Minimize &  Risk& cont & 1 &\\
Minimize & Time & cont & 1 &\\
Minimize &  Fuel Burn& cont & 1 &\\
Maximize & Quality & cont & 1 &\\
 & \multicolumn{2}{l}{\bf Total objectives} & {\bf 5} & \\
\hline
with respect to & Skin prod. location & cat  &  21 levels & geographic sites\\
 & Spar production location & cat  &  21 levels & geographic sites\\
  & Stringer production location & cat  &  21 levels & geographic sites\\
   & Rib production location & cat  &  21 levels & geographic sites\\
& Skin   Material \& Manuf. process  & cat &  6 levels & Alu-Machining,...\\
& Spar   Material \& Manuf. process  & cat  &  5 levels & Alu-Machining,...\\
& Stringer   Material \& Manuf. process  & cat  &  4 levels & Alu-Machining,...\\
& Rib  Material \& Manuf. process  & cat  &  5 levels & Alu-Machining,...\\
& \multicolumn{2}{l}{ Total categorical variables} & { 8} & \\ \cline{2-5}
 & \multicolumn{2}{l}{\bf Total relaxed variables} & {\bf 104} & \\
\hline
subject to & 
\multicolumn{2}{l}{ Material incompatibility} & 1 & \\
&\multicolumn{2}{l}{ Supply chain sites  }  & 1 & \\
 & \multicolumn{2}{l}{\bf Total constraints } & {\bf 2} & \\
\hline
\end{tabular}
   \label{tab:pb_ac2}
\end{table}
    \chapter{Chapter II appendix} 

In~\ref{apendix:EHH2CR}, we give the parameterization that allows us to obtain the continuous relaxation kernel using our proposed framework. 
In~\ref{subsec:cosine}, the cosine test case is detailed. 

\section{Continuous relaxation is a particular instance of our proposed FE Kernel.} 
\label{apendix:EHH2CR}
To show that CR is a particular instance of FE, it suffices to show that the matrix $\Phi(\Theta_i)$ is diagonal whenever $\Theta_i$ is set to a diagonal one. In fact, assume that we have, in our general model, $[\Theta_i]_{j \neq j'} = 0, \  \forall (j,j') \in \{ 1, \ldots, L_i \}.$ 
Knowing that $\cos(0)=1$ and $\sin(0)=0$, the matrix $C(\Theta_i)$ writes as  \\
$$C(\Theta_i) = \begin{bmatrix}
1 & 0 & 0  & 0  \\
1  & 0 &  \ldots & 0 \\
\vdots &\vdots & \ddots & 0 \\
1 &  0 & 0  & 0 \\
\end{bmatrix} \quad \mbox{ and } \quad  C(\Theta_i) C(\Theta_i)^\top = \begin{bmatrix}
1 & 1 & 1  & 1  \\
1  & 1 &  \ldots & 1 \\
\vdots &\vdots & \ddots & 1 \\
1 &  1 & 1  & 1 \\
\end{bmatrix} $$
Therefore, we also have 
$$[\Phi(\Theta_i)]_{j \neq j'} = \frac{\log \epsilon }{2} ([C(\Theta_i) C(\Theta_i)^\top]_{j,j'} -1) = 0, \  \forall (j,j') \in \{ 1, \ldots, L_i \} $$ that is the continuous relaxation kernel. \qed \mbox{}\\ 


    \chapter{Chapter IV appendix}

\section{A correlation kernel for decreed categorical variables}
\label{app:catdec}

Based on the recent advances done in~\cite{Mixed_Paul}, we can define a new hierarchical kernel for categorical decreed variable based on one-hot encoding and on the algebraic distance defined in~\cite{saves2023smt}. In that case:

\begin{itemize}
    \item If a categorical variable is decreed-excluded for both inputs, then all the one-hot relaxed dimensions peculiar to this variable are also decreed-excluded meaning none of them is relevant. 
    
    \item If a categorical variable is decreed-included for both inputs, then all the one-hot relaxed dimensions peculiar to this variable are also decreed-included. 
    For two given inputs in the DoE, for example, the $r^{\text{th}}$ and  $s^{\text{th}}$ points, let $c^r_{i} $ and $c^s_{i} $ be the associated categorical variables taking respectively the $\ell^i_r$ and the $\ell^i_s$ level on the categorical variable $c_i$.
    If $\ell^i_r = \ell^i_s, [R_i(\Theta_i)]_{{\ell^i_r},{\ell^i_r}} =1$. Otherwise, assuming that the chosen kernel is the exponential kernel, because of the algebraic distance related to \texttt{Alg-Kernel},  this lead to the  kernel $$ [R_i(\Theta_i)]_{{\ell^i_r},{\ell^i_s}} = \exp \left( -  \sqrt{2} [\Phi(\Theta_i)]_{{\ell^i_r},{\ell^i_r}}- \sqrt{2} [\Phi(\Theta_i)]_{{\ell^i_s},{\ell^i_s}}  \right). $$
    This kernel differs from the non-decreed one by a factor $\sqrt{2}$ over the hyperparameters but, theoretically, these hyperparameters could be whatever positive value so it should not impact the kernel.    
    
    \item If a categorical variable is decreed-excluded for an input but decreed-included for the other one, then, for all the relaxed dimensions peculiar to this variable, because of the algebraic distance related to \texttt{Alg-Kernel}, there is a induced distance $1$ between both inputs. Assuming that the chosen kernel is the exponential kernel, this lead to the  kernel $$ [R_i(\Theta_i)]_{{\ell^i_r},{\ell^i_s}} = \exp \left( -  \displaystyle{ \sum^{L_i}_{j=1} \    [\Phi(\Theta_i)]_{j,j}}  \right). $$ 
\end{itemize}
As we have seen before, continuous relaxation (or one-hot encoding) generalizes Gower distance. Therefore we can define the following adaptation for the particular case in which we are using the Gower distance kernel for a categorical variable $c_i$ with $L_i$ levels.
\begin{itemize}
    \item $d_i(x,x') = 0$ if both $x$ and $x'$ are decreed-excluded.
    \item $d_i(x,x') = \sqrt{2} \ \theta_{\text{cov}}$ if both $x$ and $x'$ are decreed-included.
    \item $d_i(x,x') = \frac{1}{2} L_i \ \theta_{\text{cov}}$ if either $x$ or $x'$ are decreed-included.
\end{itemize}

To finish with, for homoscedastic hypersphere and exponential homoscedastic hypersphere we opt for an imputation method. The imputed value is the first in the list. 

\section{Symmetric Positive Definiteness of hierachical kernels}
\label{app:SPD}

\begin{theorem}
The kernel $K^{naive}$ defined as follows is not a SPD kernel. 
\begin{equation}
 \begin{cases}
  K^{naive}(u,v) = k_{\decreed}(u_{\acting} (u_\meta), v_{\acting} (v_\meta)) k_{\neutral}(u_{\neutral}, v_{\neutral}) , & \text{if}\ u_\meta = v_\meta \\  
  K^{naive}(u,v) =  k_{\meta}(u_\meta, v_\meta) k_{\neutral}(u_{\neutral}, v_{\neutral}), & \text{if}\ u_\meta \neq v_\meta \\
\end{cases} 
\end{equation}
\end{theorem}

\begin{proof}

$K^{naive}$ is a kernel if and only if $K^{naive}$ can be written as $K^{naive} (u,v) = \kappa ( d(u,v) ) $ with $ \kappa $ being an SPD function and $d$ being a distance in a given Hilbert space.
If the two continuous decreed inputs $u_{dec}$ and $v_{dec}$ are in the same subspace ($u_\meta=v_\meta$), then  $k_\meta(u_\meta,u_\meta) =1$ and $ K^{naive}(u,v) = k_{\decreed}(u_\acting(u_\meta), v_\acting(v_\meta) )$. This works without any problem, a categorical variable being fully corellated with itself. 

On the contrary, if the two continuous inputs $u_{\decreed}$ and $v_{\decreed}$ lie in two different subspaces ($u_\meta \neq v_\meta$), then, we have $ k_{dec}( u_{\decreed}, {v}_{\decreed} ) =1 \implies  d(u_{\acting} (u_\meta), v_{\acting} (v_\meta) ) = 0 $ .

We know that a distance is defined such that $ d(A,B) = 0 \Leftrightarrow A=B  $ (Identity of indiscernibles).
Yet, we have  $ d( u_{\acting} (u_\meta), v_{\acting} (v_\meta) ) = 0 $ and  $u_{\acting} (u_\meta) \neq  v_{\acting}(v_\meta)$. $d$ is not a distance because it is always equal to $0$ for points in distinct decreed spaces. Therefore the kernel vanishes to 1 and the correlation matrix is degenerated and can not be SPD.

Based on this proof, the idea of the arc-kernel is to have a constant residual distance between distinct subspaces. In other words, there is a non-zero distance between $\mathcal{X}(u_\meta)$ and $\mathcal{X}(v_\meta)$ if $u_\meta \neq v_\meta  $.

\end{proof}

\begin{theorem}
Our \texttt{Alg-Kernel} kernel is SPD and so is $k$ defined in~\eqnref{eq:hier_ker}.\end{theorem}
\begin{proof}

Let $I_{u}$ be the subset of indices $ i \in I_\decreed$ that are decreed-included by $ u_{\meta}$ such that $ I^{inter}_{u,v}  = I_{u} \bigcap I_{v}$. Let $d_E ( [x_u,x_v],[y_u,y_v] ) = \theta_i \sqrt{ ( x_u -x_v )^2 +(y_u- y_v )^2 }  $ be an Euclidean distance in $\mathbb{R}^2  \times \mathbb{R}^2  $.
Due to~\cite[Proposition 2]{Hutter}, we only need to show that, for any two inputs $u,v \in \mathcal{X}$ , the isometry condition $d_E (f^{alg}_{i} (u) ,f^{alg}_{i} (v )) $  holds for a given function $f^{alg}_{i}$, that is equivalent to having a Hilbertian metric. In other words, $d_E$ is isomorphic to an $L^2$ norm~\cite{haasdonk2004learning}. Such kernels are well-known and referred as "substitution kernels with euclidean distance"~\cite{sow2023learning}.

$\forall i \in I_\decreed, f_i^{alg}(u_{\decreed}) $ is defined by
\begin{equation}
 \begin{cases}
      f_i^{alg}(u_\acting ) = [\frac{1-((  u_{dec})_i)^2}{1+((  u_{dec})_i)^2 },\frac{2 (  u_{dec})_i   }{1+ ((u_{dec})_i)^2 }]  \text{ if } i \in I_{u}  \\
      f_i^{alg}(u_{\decreed}) = [0,0] \text{ otherwise }  \\
\end{cases} 
\end{equation}  

Case 1: $i \in I_\decreed, i \notin I_{u}, i \notin I_{v} $. 
$$d_E ( f_i^{alg} ( \overline{u}_{dec} ) ,f_i^{alg}( \overline{v}_{dec}  )) =d_E ([0,0],[0,0]) = 0. $$
Therefore, the complementary space is not relevant, as for the arc-kernel. \medbreak

Case 2: $i \in I_\decreed, i \in I_{u}, i \notin I_{v} $.

$$d_E ( f_i^{alg}( u_{dec} ) ,f_i^{alg} (  \overline{v}_{dec})) =d_E ( [\frac{1-((  u_{dec})_i)^2}{1+((  u_{dec})_i)^2 },\frac{2 (  u_{dec})_i }{1+ ((u_{dec})_i)^2 }] , [0,0] ) $$
$$=  \theta_i . $$
This case corresponds to the kernel $K_\meta^{alg}(u_\meta,v_\meta)$ when $u_\meta \neq v_\meta $.  \medbreak

Case 3:  $i \in I_\decreed, i \in I^{inter}_{u,v}$.

\begin{equation}
    \begin{split}
         &d_E ( f_i^{alg}( u_{dec}  ) ,f_i^{alg}( v_{dec}  ))  \\
        &= d_E \left( \left[  \frac{1-((  u_{dec})_i)^2}{1+((  u_{dec})_i)^2 },\frac{2 (  u_{dec})_i }{1+ ((u_{dec})_i)^2 } \right], \left[\frac{1-((  v_{dec})_i)^2}{1+((  v_{dec})_i)^2 },\frac{2 (  v_{dec})_i }{1+ ((v_{dec})_i)^2 }\right] \right) \\
        &= 2 \theta_i \frac{ |(u_{dec})_i- (v_{dec})_i|}{ \sqrt{{( (u_{dec})_i })^2+1}\sqrt{{ ((v_{dec})_i })^2+1}}
    \end{split}
\end{equation}
This case corresponds to the kernel $K_{dec}^{alg}(u_{dec}(u_m),v_{dec}(v_m))$ when $u_m = v_m $.  \medbreak

Therefore, there exists an isometry between our algebraic distance and the euclidean distance. The distance being well-defined for any given kernel, the matrix obtained with our model is SPD. 

Although, in~\cite{Hutter}, they also do the demonstration of the metric property of their distance formally. The non-negativity and symmetry of $d^\text{alg}$ are trivially proven knowing that the hyperparameters $\theta$ are strictly positive. To prove the triangle inequality, consider $(u,v,w) \in \mathcal{X}^3.$

Case 1: $i \in I_\decreed, i \notin I_{u}, i \notin I_{v} $. 
$$d^\text{alg} (u_i,v_i) = 0 \leq d^\text{alg} (u_i,w_i) + d^\text{alg} (w_i,v_i) \text{ by non negativity.}$$   \medbreak

Case 2: $i \in I_\decreed, i \in I_{u}, i \notin I_{v} $.

\begin{itemize} [leftmargin=3cm]
    \item $   i \in I_{w}, \ d^\text{alg} (u_i,v_i) = \theta_i \leq d^\text{alg} (u_i,w_i)  + \theta_i\text{ by non negativity.}$   
 \item $  i \notin I_{w}, \ d^\text{alg} (u_i,v_i) = \theta_i \leq \theta_i  + d^\text{alg} (w_i,v_i)  \text{ by non negativity.}$   

\end{itemize}
\medbreak

Case 3: $i \in I_\decreed, i \in I^{inter}_{u,v}$.

\begin{itemize} [leftmargin=3cm]
    \item $  i \notin I_{w}, \ d^\text{alg} (u_i,v_i) = 2 \theta_i \frac{ |u_i- v_i|}{ \sqrt{{( u_i })^2+1}\sqrt{{ (v_i })^2+1}}  \leq  \theta_i  + \theta_i \ \text{ for } (u_i,v_i) \in [0,1]^2.$   
 \item $  i \in I_{w}, \text{ Knowing that}  \frac{ |a- c|}{ \sqrt{{ a }^2+1}\sqrt{{ c }^2+1}} +\frac{ |c- b|}{ \sqrt{{ c }^2+1}\sqrt{{ b}^2+1}}   \geq \frac{ |a- b|}{ \sqrt{{ a }^2+1}\sqrt{{ b }^2+1}}   $, we have $\ d^\text{alg} (u_i,v_i) =  2 \theta_i \frac{ |u_i- v_i|}{ \sqrt{{( u_i })^2+1}\sqrt{{ (v_i })^2+1}}   \leq  d^\text{alg} (u_i,w_i) +  d^\text{alg} (w_i,v_i) $.
\end{itemize}


\medbreak

\end{proof}

\section{A new algebraic distance}
\label{app:dalg}
Compared to our new algebraic distance, we acknowledge that simpler distances have been developed to tackle one-hot encoded variables in the context of computer vision and bioengineering~\cite{gardner2017definiteness, gomez2018automatic}. Still our work is more general as we are not considering binary variables but more complex hierarchical variables either included or excluded. In that context, the multidimensional formulation of our distance writes
\begin{equation}
d^{\text{alg}}(x,x') = 
\begin{cases}
1 , & \parbox{3cm}{
$\text{if}\ x^\top x' = 0 $}  \\ 
\frac{ ||x-x'||}{ \sqrt{||x||^2+1}\sqrt{||x'||^2+1}}, & \text{otherwise.}
\end{cases}
\end{equation}
\\

\end{appendices}

\addcontentsline{toc}{section}{\bf REFERENCES}
\bibliography{main.bib}

\end{document}